\documentclass[letterpaper,twoside]{article}
\usepackage[OT1,T1]{fontenc}
\usepackage[dvips]{graphicx,color}
\usepackage{citesort}
\usepackage{pst-3d}
\usepackage{pst-char}
\usepackage{pst-coil}
\usepackage{pst-eps}
\usepackage{pst-fill}
\usepackage{pst-grad}
\usepackage{pst-node}
\usepackage{pst-plot}
\usepackage{pst-text}
\usepackage{pst-tree}
\usepackage{paralist}

\usepackage{ccfonts}
\usepackage{beton}
\usepackage{fancybox}
\usepackage{amssymb}
\usepackage{amsmath,bm}
\usepackage[mathscr]{eucal}
\usepackage{amsfonts}
\usepackage{latexsym}
\usepackage{amsxtra}
\usepackage{amsbsy}
\usepackage{amsthm}
\usepackage{amscd}
\usepackage{amsopn}
\usepackage{amstext}
\usepackage{oldgerm}
\usepackage[varumlaut]{yfonts}
\newcounter{z0}
\newcounter{y0}

\newtheorem{bbbbb}{Proposition}[section]
\newtheorem{ccccc}{Lemma}[section]
\newtheorem{ddddd}{Theorem}[section]

\newtheorem{dy}{Theorem}[section]

\theoremstyle{definition}

\theoremstyle{definition}

\theoremstyle{definition}
\newtheorem{eeeee}{Remark}[section]

\theoremstyle{definition}
\newtheorem{ey}{Remark}[section]

\newtheorem*{rhp1}{RHP1}
\newtheorem*{rhp2}{RHP2}
\newcommand{\me}{\mathrm{e}}
\newcommand{\mi}{\mathrm{i}}
\newcommand{\md}{\mathrm{d}}
\renewcommand{\Im}{\mathrm{Im}}

\newcommand{\id}{\pmb{\mi \md}}
\newcommand{\norm}[1]{\lVert#1\rVert}
\numberwithin{equation}{section}
\allowdisplaybreaks[4]

\usepackage{relsize,exscale}
\usepackage{stmaryrd}
\usepackage{pxfonts}
\begin{document}
\fontsize{10}{12}\selectfont
\fontencoding{T1}\selectfont
\baselineskip=12pt
\frenchspacing
\title{Asymptotics of Recurrence Relation Coefficients, Hankel Determinant
Ratios, and Root Products Associated with Laurent Polynomials Orthogonal
with Respect to Varying Exponential Weights}
\author{K.~T.-R.~McLaughlin\thanks{\texttt{E-mail: mcl@math.arizona.edu}} \\
Department of Mathematics \\
The University of Arizona \\
617 N.~Santa Rita Ave. \\
P.~O.~Box 210089 \\
Tucson, Arizona 85721--0089 \\
U.~S.~A. \and
A.~H.~Vartanian\thanks{\texttt{E-mail: arthurv@math.ucf.edu}} \\
Department of Mathematics \\
University of Central Florida \\
P.~O.~Box 161364 \\
Orlando, Florida 32816--1364 \\
U.~S.~A. \and
X.~Zhou\thanks{\texttt{E-mail: zhou@math.duke.edu}} \\
Department of Mathematics \\
Duke University \\
Box 90320 \\
Durham, North Carolina 27708--0320 \\
U.~S.~A.}
\date{8 February 2006}
\maketitle
\begin{abstract}
\noindent
Let $\Lambda^{\mathbb{R}}$ denote the linear space over $\mathbb{R}$ spanned 
by $z^{k}$, $k \! \in \! \mathbb{Z}$. Define the real inner product (with 
varying exponential weights) $\langle \boldsymbol{\cdot},\boldsymbol{\cdot} 
\rangle_{\mathscr{L}} \colon \Lambda^{\mathbb{R}} \times \Lambda^{\mathbb{
R}} \! \to \! \mathbb{R}$, $(f,g) \! \mapsto \! \int_{\mathbb{R}}f(s)g(s) 
\exp (-\mathscr{N} \, V(s)) \, \md s$, $\mathscr{N} \! \in \! \mathbb{N}$, 
where the external field $V$ satisfies: (i) $V$ is real analytic on $\mathbb{
R} \setminus \{0\}$; (ii) $\lim_{\vert x \vert \to \infty}(V(x)/ \ln (x^{2} 
\! + \! 1)) \! = \! +\infty$; and (iii) $\lim_{\vert x \vert \to 0}(V(x)/\ln
(x^{-2} \! + \! 1)) \! = \! +\infty$. Orthogonalisation of the (ordered) base 
$\lbrace 1,z^{-1},z,z^{-2},z^{2},\dotsc,z^{-k},z^{k},\dotsc \rbrace$ with 
respect to $\langle \boldsymbol{\cdot},\boldsymbol{\cdot} \rangle_{\mathscr{
L}}$ yields the even degree and odd degree orthonormal Laurent polynomials 
(OLPs) $\lbrace \phi_{m}(z) \rbrace_{m=0}^{\infty}$: $\phi_{2n}(z) \! = \! 
\sum_{k=-n}^{n} \xi^{(2n)}_{k}z^{k}$, $\xi^{(2n)}_{n} \! > \! 0$, and 
$\phi_{2n+1}(z) \! = \! \sum_{k=-n-1}^{n} \xi^{(2n+1)}_{k}z^{k}$, 
$\xi^{(2n+1)}_{-n-1} \! > \! 0$. Associated with the even degree and odd 
degree OLPs are the following two pairs of recurrence relations: $z \phi_{2n}
(z) \! = \! c_{2n}^{\sharp} \phi_{2n-2}(z) \! + \! b_{2n}^{\sharp} \phi_{2n-1}
(z) \! + \! a_{2n}^{\sharp} \phi_{2n}(z) \! + \! b_{2n+1}^{\sharp} \phi_{2n+1}
(z) \! + \! c_{2n+2}^{\sharp} \phi_{2n+2}(z)$ and $z \phi_{2n+1}(z) \! = \! 
b_{2n+1}^{\sharp} \phi_{2n}(z) \! + \! a_{2n+1}^{\sharp} \phi_{2n+1}(z) \! + 
\! b_{2n+2}^{\sharp} \phi_{2n+2}(z)$, where $c_{0}^{\sharp} \! = \! b_{0}^{
\sharp} \! = \! 0$, and $c_{2k}^{\sharp} \! > \! 0$, $k \! \in \! \mathbb{N}$, 
and $z^{-1} \phi_{2n+1}(z) \! = \! \gamma_{2n+1}^{\sharp} \phi_{2n-1}(z) 
\! + \! \beta_{2n+1}^{\sharp} \phi_{2n}(z) \! + \! \alpha_{2n+1}^{\sharp} 
\phi_{2n+1}(z) \! + \! \beta_{2n+2}^{\sharp} \phi_{2n+2}(z) \! + \! 
\gamma_{2n+3}^{\sharp} \phi_{2n+3}(z)$ and $z^{-1} \phi_{2n}(z) \! = \! 
\beta_{2n}^{\sharp} \phi_{2n-1}(z) \! + \! \alpha_{2n}^{\sharp} \phi_{2n}
(z) \! + \! \beta_{2n+1}^{\sharp} \phi_{2n+1}(z)$, where $\beta_{0}^{\sharp} 
\! = \! \gamma_{1}^{\sharp} \! = \! 0$, $\beta_{1}^{\sharp} \! > \! 0$, and 
$\gamma_{2l+1}^{\sharp} \! > \! 0$, $l \! \in \! \mathbb{N}$. Asymptotics 
in the double-scaling limit as $\mathscr{N},n \! \to \! \infty$ such that 
$\mathscr{N}/n \! = \! 1 \! + \! o(1)$ of the coefficients of these two 
pairs of recurrence relations, Hankel determinant ratios associated with 
the real-valued, bi-infinite, strong moment sequence $\left\lbrace c_{k} 
\! = \! \int_{\mathbb{R}}s^{k} \exp (-\mathscr{N} \, V(s)) \, \md s 
\right\rbrace_{k \in \mathbb{Z}}$, and the products of the (real) roots of 
the OLPs are obtained by formulating the even degree and odd degree OLP 
problems as matrix Riemann-Hilbert problems on $\mathbb{R}$, and then 
extracting the large-$n$ behaviours by applying the non-linear 
steepest-descent method introduced in \cite{a1} and further developed in 
\cite{a2,a3}.

\vspace{0.65cm}

30C70, 30C15, 41A60, 31A99

\vspace{0.50cm}
\textbf{Abbreviated Title.} Asymptotics of Recurrence Coefficients, Hankel
Determinant Ratios, and

Root Products

\vspace{0.50cm}

ices, orthogonal Laurent polynomials, recurrence relations, Riemann-Hilbert
problems, stro-

ng moment problems, variational problems
\end{abstract}
\clearpage
\section{Introduction and Background}
The \emph{strong Stieltjes moment problem} (SSMP) considers the relation 
between a given doubly- or bi-infinite (moment) sequence $\{c_{n}\}_{n \in 
\mathbb{Z}}$ of real numbers and associated measures. There are, in fact, 
three general types of problems in this area:
\begin{compactenum}
\item[(i)] find necessary and sufficient conditions for the existence of a 
non-negative Borel measure $\mu \! := \! \mu^{\text{\tiny SS}}_{\text{\tiny 
MP}}$ on $[0,+\infty)$, and with infinite support, such that $c_{n}^{\text{
\tiny SSMP}} \! = \! \int_{0}^{+\infty} \lambda^{n} \, \md \mu^{\text{\tiny 
SS}}_{\text{\tiny MP}}(\lambda)$, $n \! \in \! \mathbb{Z}$, where the 
(improper) integral is to be understood in the Riemann-Stieltjes sense;
\item[(ii)] when there is a solution, in which case the SSMP is 
\emph{determinate}, find conditions for the uniqueness of the solution; and
\item[(iii)] when there is more than one solution, in which case the SSMP 
is \emph{indeterminate}, describe the family of all solutions.
\end{compactenum}
Similarly, when the support of the non-negative Borel measure $\mu \! := \! 
\mu^{\text{\tiny SH}}_{\text{\tiny MP}}$ is only required to be contained in 
$(-\infty,+\infty)$, the analogous problem is referred to as the \emph{strong 
Hamburger moment problem} (SHMP). (In this latter case, $c_{n}^{\text{\tiny 
SHMP}} \! = \! \int_{-\infty}^{+\infty} \lambda^{n} \, \md \mu^{\text{\tiny 
SH}}_{\text{\tiny MP}}(\lambda)$, $n \! \in \! \mathbb{Z}.)$ The SSMP (resp., 
SHMP) was introduced in $1980$ (resp., $1981)$ by Jones \emph{et al.} 
\cite{a4} (resp., Jones \emph{et al.} \cite{a5}), and studied further in 
\cite{a6,a7,a8,a9,a10} (see, also, \cite{a11}). Unlike the moment theory 
for the \emph{classical Stieltjes} (resp., \emph{classical Hamburger}) 
\emph{moment problem} (SMP) \cite{a12} (resp., (HMP) \cite{a13}), wherein 
the theory of orthogonal polynomials \cite{a14} played a prominent r\^{o}le 
(see, for example, \cite{a15}), the extension of the moment theory to 
the SSMP and the SHMP introduced a `rational generalisation' of the 
orthogonal polynomials, namely, the \emph{orthogonal Laurent} (or $L$-) 
\emph{polynomials} \cite{a4,a5,a6,a7,a8,a9,a10,a11,a16,a17,a18,a19,a20}.

Orthogonal $L$-polynomials are defined in terms of certain orthogonality 
relations involving an associated measure $\mu$. In this work, measures 
that satisfy the following conditions are considered:
\begin{compactenum}
\item[$\circledast$] $\widetilde{\mu} \! \in \! \mathcal{M}_{1}(\mathbb{R})$,
where
\begin{equation*}
\mathcal{M}_{1}(\mathbb{R}) \! := \! \left\{\mathstrut \mu; \, \int_{\mathbb{
R}} \md \mu (s) \! = \! 1, \, \int_{\mathbb{R}}s^{k} \md \mu (s) \! < \! 
\infty, \, k \! \in \! \mathbb{Z} \setminus \{0\} \right\};
\end{equation*}
\item[$\circledast$] $\md \widetilde{\mu}(z) \! = \! \widetilde{w}(z) \, \md
z$, where
\begin{equation*}
\widetilde{w}(z) \! := \! \exp \! \left(-\mathscr{N} \, V(z) \right), \quad
\mathscr{N} \! \in \! \mathbb{N}.
\end{equation*}
\end{compactenum}
Here, $\widetilde{w}(z)$ is referred to as the varying exponential weight 
function, and the function $V \colon \mathbb{R} \setminus \{0\} \! \to \! 
\mathbb{R}$, which is referred to as the \emph{external field}, satisfies:
\begin{gather}
V \, \, \text{is real analytic on $\mathbb{R} \setminus \{0\}$}; \tag{V1} \\
\lim_{\vert x \vert \to \infty} \! \left(V(x)/\ln (x^{2} \! + \! 1) \right) 
\! = \! +\infty; \tag{V2} \\
\lim_{\vert x \vert \to 0} \! \left(V(x)/\ln (x^{-2} \! + \! 1) \right) \! 
= \! +\infty. \tag{V3}
\end{gather}
(For example, a rational function of the form $V(z) \! = \! \sum_{k=-2m_{1}
}^{2m_{2}} \varrho_{k}z^{k}$, with $\varrho_{k} \! \in \! \mathbb{R}$, $k \! 
= \! -2m_{1},\dotsc,2m_{2}$, $m_{1,2} \! \in \! \mathbb{N}$, and $\varrho_{
-2m_{1}},\varrho_{2m_{2}} \! > \! 0$ would satisfy conditions~(V1)--(V3).)

Rational functions with poles at $0$ and at $\infty$ are seminal for obtaining 
the results of this work: they are now discussed. It will be important to have 
the ordered sequence of (complex) monomials $\lbrace 1,z^{-1},z,z^{-2},z^{2},
\dotsc,z^{-k},z^{k},\dotsc \rbrace$ corresponding to the 
\emph{cyclically-repeated pole sequence} $\lbrace \text{no pole},0,\infty,
\linebreak[4]
0,\infty,\dotsc,0,\infty,\dotsc \rbrace$. This is a basis for the set of 
functions $\Lambda^{\mathbb{R}} \! = \! \cup_{m \in \mathbb{Z}_{0}^{+}}
(\Lambda^{\mathbb{R}}_{2m} \cup \Lambda^{\mathbb{R}}_{2m+1})$, with $\mathbb{
Z}_{0}^{+} \! := \! \{0\} \cup \mathbb{N}$, where $\Lambda^{\mathbb{R}}_{2m} 
\! := \! \Lambda^{\mathbb{R}}_{-m,m}$, and $\Lambda^{\mathbb{R}}_{2m+1} \! 
:= \! \Lambda^{\mathbb{R}}_{-m-1,m}$, where, for any pair $(p,q) \! \in \! 
\mathbb{Z} \times \mathbb{Z}$, with $p \! \leqslant \! q$, $\Lambda^{\mathbb{
R}}_{p,q} \! := \! \left\lbrace \mathstrut f \colon \mathbb{C}^{\ast} \! \to 
\! \mathbb{C}; \, f(z) \! = \! \sum_{k=p}^{q} \widetilde{\lambda}_{k}z^{k}, 
\, \widetilde{\lambda}_{k} \! \in \! \mathbb{R}, \, k \! = \! p,\dotsc,q 
\right\rbrace$, where $\mathbb{C}^{\ast} \! := \! \mathbb{C} \setminus 
\{0\}$. For each $l \! \in \! \mathbb{Z}_{0}^{+}$ and $0 \! \not\equiv \! 
f \! \in \! \Lambda^{\mathbb{R}}_{l}$, the $L$-\emph{degree} of $f$, 
symbolically $LD(f)$, is defined as
\begin{equation*}
LD(f) \! := \! l.
\end{equation*}
The \emph{leading coefficient} of $0 \! \not\equiv \! f \! \in \! \Lambda^{
\mathbb{R}}_{l}$, symbolically $LC(f)$, is defined as
\begin{equation*}
LC(f) \! := \!
\begin{cases}
\widetilde{\lambda}_{m}, &\text{$l \! = \! 2m$,} \\
\widetilde{\lambda}_{-m-1}, &\text{$l \! = \! 2m \! + \! 1$,}
\end{cases}
\end{equation*}
and the \emph{trailing coefficient} of $0 \! \not\equiv \! f \! \in \!
\Lambda^{\mathbb{R}}_{l}$, symbolically $TC(f)$, is defined as
\begin{equation*}
TC(f) \! := \!
\begin{cases}
\widetilde{\lambda}_{-m}, &\text{$l \! = \! 2m$,} \\
\widetilde{\lambda}_{m}, &\text{$l \! = \! 2m \! + \! 1$.}
\end{cases}
\end{equation*}
Thus, for $0 \! \not\equiv \! f \! \in \! \Lambda^{\mathbb{R}}_{2m}$, one 
writes $f(z) \! = \! TC(f)z^{-m} \! + \! \cdots \! + \! LC(f)z^{m}$, and, for 
$0 \! \not\equiv \! f \! \in \! \Lambda^{\mathbb{R}}_{2m+1}$, one writes 
$f(z) \! = \! LC(f)z^{-m-1} \! + \! \cdots \! + \! TC(f)z^{m}$. For each $l 
\! \in \! \mathbb{Z}_{0}^{+}$, $0 \! \not\equiv \! f \! \in \! \Lambda^{
\mathbb{R}}_{l}$ is called \emph{monic} if $LC(f) \! = \! 1$.

Define (uniquely) the \emph{strong moment linear functional} $\mathscr{L}$ 
by its action on the basis elements of $\Lambda^{\mathbb{R}}$ as follows: 
$\mathscr{L} \colon \Lambda^{\mathbb{R}} \! \to \! \Lambda^{\mathbb{R}}$, 
$f \! = \! \sum_{k \in \mathbb{Z}} \widetilde{\lambda}_{k}z^{k} \! \mapsto 
\! \mathscr{L}(f) \! := \! \sum_{k \in \mathbb{Z}} \widetilde{\lambda}_{k}
c_{k}$, where $c_{k} \! = \! \mathscr{L}(z^{k}) \! = \! \int_{\mathbb{R}}
s^{k} \me^{-\mathscr{N} \, V(s)} \, \md s$, $(k,\mathscr{N}) \! \in \! 
\mathbb{Z} \times \mathbb{N}$. (Note that $\{c_{k}\}_{k \in \mathbb{Z}}$ is 
a bi-infinite, real-valued, \emph{strong moment sequence}: $c_{k}$ is called 
the \emph{$k$th strong moment of $\mathscr{L}$}.) Associated with the 
above-defined bi-infinite, real-valued, strong moment sequence $\{c_{k}\}_{k 
\in \mathbb{Z}}$ are the \emph{Hankel determinants} $H^{(m)}_{k}$, $(m,k) \! 
\in \! \mathbb{Z} \times \mathbb{N}$ (see, for example, \cite{a11}):
\begin{equation}
H^{(m)}_{0} \! := \! 1 \qquad \quad \text{and} \qquad \quad 
H^{(m)}_{k} \! := \!
\begin{vmatrix}
c_{m} & c_{m+1} & \cdots & c_{m+k-2} & c_{m+k-1} \\
c_{m+1} & c_{m+2} & \cdots & c_{m+k-1} & c_{m+k} \\
c_{m+2} & c_{m+3} & \cdots & c_{m+k} & c_{m+k+1} \\
\vdots & \vdots & \ddots & \vdots & \vdots \\
c_{m+k-1} & c_{m+k} & \cdots & c_{m+2k-3} & c_{m+2k-2}
\end{vmatrix}.
\end{equation}

Define the bilinear form (with varying exponential weight) $\langle 
\boldsymbol{\cdot},\boldsymbol{\cdot} \rangle_{\mathscr{L}}$ as follows: 
$\langle \boldsymbol{\cdot},\boldsymbol{\cdot} \rangle_{\mathscr{L}} \colon 
\Lambda^{\mathbb{R}} \times \Lambda^{\mathbb{R}} \! \to \! \mathbb{R}$, $(f,
g) \! \mapsto \! \langle f,g \rangle_{\mathscr{L}} \! := \! \mathscr{L}(f(z)
g(z)) \! = \! \int_{\mathbb{R}}f(s)g(s) \exp (-\mathscr{N} \, V(s)) \, \md 
s$, $\mathscr{N} \! \in \! \mathbb{N}$. It is a fact \cite{a6,a7,a11,a17} 
(see, also, the proofs of Lemmas~2.2.1 and~2.2.2 in \cite{a21} and \cite{a22}, 
respectively) that the bilinear form $\langle \boldsymbol{\cdot},\boldsymbol{
\cdot} \rangle_{\mathscr{L}}$ thus defined is an inner product if and only 
if $H^{(-2m)}_{2m} \! > \! 0$ and $H^{(-2m)}_{2m+1} \! > \! 0$ for each 
$m \! \in \! \mathbb{Z}_{0}^{+}$ (see Equations~(1.8) below).
\begin{eeeee}
These latter two (Hankel determinant) inequalities also appear when the 
question of the solvability of the SHMP is posed (in this case, the $c_{k}$, 
$k \! \in \! \mathbb{Z}$, which appear in Equations~(1.1) should be replaced 
by $c_{k}^{\text{\tiny SHMP}}$, $k \! \in \! \mathbb{Z})$; indeed, if these 
two inequalities are true $\forall \, \, m \! \in \! \mathbb{Z}_{0}^{+}$, then 
there is a non-negative measure $\mu^{\text{\tiny SH}}_{\text{\tiny MP}}$ (on 
$\mathbb{R})$ with the given (real) moments. For the case of the SSMP, there 
are four (Hankel determinant) inequalities (in this latter case, the $c_{k}$, 
$k \! \in \! \mathbb{Z}$, which appear in Equations~(1.1) should be replaced 
by $c_{k}^{\text{\tiny SSMP}}$, $k \! \in \! \mathbb{Z})$ which guarantee the 
existence of a non-negative measure $\mu^{\text{\tiny SS}}_{\text{\tiny MP}}$ 
(on $[0,+\infty))$ with the given moments, namely \cite{a4} (see, also, 
\cite{a6,a7}): for each $m \! \in \! \mathbb{Z}_{0}^{+}$, $H^{(-2m)}_{2m} 
\! > \! 0$, $H^{(-2m)}_{2m+1} \! > \! 0$, $H^{(-2m+1)}_{2m} \! > \! 0$, and 
$H^{(-2m-1)}_{2m+1} \! < \! 0$. It is interesting to note that the former 
solvability conditions do not automatically imply that the positive (real) 
moments $\{c_{k}^{\text{\tiny SHMP}}\}_{k \in \mathbb{Z}_{0}^{+}}$ determine 
a measure via the HMP: a similar statement holds true for the SMP (see the 
latter four solvability conditions). \hfill $\blacksquare$
\end{eeeee}

The \emph{norm of $f \! \in \! \Lambda^{\mathbb{R}}$ with respect to 
$\mathscr{L}$} is now defined:
\begin{equation*}
\norm{f(\pmb{\cdot})}_{\mathscr{L}} \! := \! (\langle f,f \rangle_{\mathscr{
L}})^{1/2}.
\end{equation*}
Given an inner product $\langle \boldsymbol{\cdot},\boldsymbol{\cdot} \rangle_{
\mathscr{L}}$, one may define the (real) orthonormal Laurent polynomial
sequence $\lbrace \phi_{n}^{\natural}(z) \rbrace_{n \in \mathbb{Z}_{0}^{+}}$
with respect to $\mathscr{L}$ as follows; $\forall \, \, m,n \! \in \!
\mathbb{Z}_{0}^{+}$:
\begin{compactenum}
\item[(i)] $\phi_{n}^{\natural}(z) \! \in \! \Lambda^{\mathbb{R}}_{n}$, 
that is, $LD(\phi_{n}^{\natural}) \! = \! n$;
\item[(ii)] $\langle \phi_{m}^{\natural},\phi_{m^{\prime}}^{\natural} 
\rangle_{\mathscr{L}} \! = \! 0$, $m \! \not= \! m^{\prime}$, or, 
alternatively, $\langle f,\phi_{m}^{\natural} \rangle_{\mathscr{L}} \! 
= \! 0 \, \, \forall \, \, f \! \in \! \Lambda^{\mathbb{R}}_{m-1}$;
\item[(iii)] $\langle \phi_{m}^{\natural},\phi_{m}^{\natural} \rangle_{
\mathscr{L}} \! =: \! \norm{\phi_{m}^{\natural}(\cdot)}^{2}_{\mathscr{L}} 
\! = \! 1$.
\end{compactenum}
Hereafter, these shall be referred to as orthonormal $L$-polynomials. 
Orthonormal $L$-polynomials may be obtained by Gram-Schmidt 
orthogonalisation. To do this, however, one must have an ordered sequence 
of complex monomials, $z^{j}$, $j \! \in \! \mathbb{Z}$, and the orthonormal 
$L$-polynomials, themselves, will depend on how this sequence of complex 
monomials is ordered. In this work, orthonormalisation with respect to 
$\langle \boldsymbol{\cdot},\boldsymbol{\cdot} \rangle_{\mathscr{L}}$ of the 
base sequence of $\Lambda^{\mathbb{R}}$ is considered, that is, $\lbrace 1,
z^{-1},z,z^{-2},z^{2},\dotsc,z^{-k},z^{k},\dotsc \rbrace$. The result of this 
procedure is the orthonormal $L$-polynomials $\lbrace \phi_{m}(z) \rbrace_{m 
\in \mathbb{Z}_{0}^{+}}$:
\begin{itemize}
\item for $m \! = \! 2n$,
\begin{equation}
\phi_{2n}(z) \! = \! \sum_{k=-n}^{n} \xi^{(2n)}_{k}z^{k}, \qquad \quad
\xi^{(2n)}_{n} \! > \! 0;
\end{equation}
\item for $m \! = \! 2n \! + \! 1$,
\begin{equation}
\phi_{2n+1}(z) \! = \! \sum_{k=-n-1}^{n} \xi^{(2n+1)}_{k}z^{k}, \qquad 
\quad \xi^{(2n+1)}_{-n-1} \! > \! 0.
\end{equation}
\end{itemize}
The $\phi_{n}$'s are normalised so that they all have real coefficients; in
particular, the leading coefficients, $LC(\phi_{2n}) \! = \! \xi^{(2n)}_{n}$
and $LC(\phi_{2n+1}) \! = \! \xi^{(2n+1)}_{-n-1}$, $n \! \in \! \mathbb{Z}_{
0}^{+}$, are both positive, $\xi^{(0)}_{0} \! = \! 1$, and $\phi_{0}(z) \!
\equiv \! 1$. Even though the leading coefficients, $\xi^{(2n)}_{n}$ and
$\xi^{(2n+1)}_{-n-1}$, $n \! \in \! \mathbb{Z}_{0}^{+}$, are non-zero (in
particular, they are positive), no such restriction applies to the trailing
coefficients, $TC(\phi_{2n}) \! = \! \xi^{(2n)}_{-n}$ and $TC(\phi_{2n+1})
\! = \! \xi^{(2n+1)}_{n}$, $n \! \in \! \mathbb{Z}_{0}^{+}$ (see below).
Furthermore, note that, by construction:
\begin{compactenum}
\item[(i)] $\langle \phi_{2n},z^{j} \rangle_{\mathscr{L}} \! = \! 0$, $j \! 
= \! -n,\dotsc,n \! - \! 1$;
\item[(ii)] $\langle \phi_{2n+1},z^{j} \rangle_{\mathscr{L}} \! = \! 0$, $j 
\! = \! -n,\dotsc,n$;
\item[(iii)] $\langle \phi_{j},\phi_{k} \rangle_{\mathscr{L}} \! = \! \delta_{
jk}$, $j,k \! \in \! \mathbb{Z}_{0}^{+}$, where $\delta_{jk}$ is the Kronecker
delta.
\end{compactenum}

It is convenient to introduce the monic orthogonal Laurent (or $L$-)
polynomials, $\boldsymbol{\pi}_{j}(z)$, $j \! \in \! \mathbb{Z}_{0}^{+}$: (i)
for $j \! = \! 2n$, with $\boldsymbol{\pi}_{0}(z) \! \equiv \! 1$,
\begin{equation}
\boldsymbol{\pi}_{2n}(z) \! := \! (\xi^{(2n)}_{n})^{-1} \phi_{2n}(z) \! =
\nu^{(2n)}_{-n}z^{-n} \! + \! \dotsb \! + \! z^{n}, \qquad \quad \nu^{(2n)}_{-
n} \! := \! \xi^{(2n)}_{-n}/\xi^{(2n)}_{n};
\end{equation}
and (ii) for $j \! = \! 2n \! + \! 1$,
\begin{equation}
\boldsymbol{\pi}_{2n+1}(z) \! := \! (\xi^{(2n+1)}_{-n-1})^{-1} \phi_{2n+1}(z)
\! = \! z^{-n-1} \! + \! \dotsb \! + \nu^{(2n+1)}_{n}z^{n}, \qquad \quad
\nu^{(2n+1)}_{n} \! := \! \xi^{(2n+1)}_{n}/\xi^{(2n+1)}_{-n-1}.
\end{equation}
The monic orthogonal $L$-polynomials, $\boldsymbol{\pi}_{j}(z)$, $j \! \in \!
\mathbb{Z}_{0}^{+}$, satisfy the following orthogonality and normalisation
conditions:
\begin{compactenum}
\item[(1)] $\langle \boldsymbol{\pi}_{2n},z^{j} \rangle_{\mathscr{L}} \! = 
\! 0$, $j \! = \! -n,\dotsc,n \! - \! 1$;
\item[(2)] $\langle \boldsymbol{\pi}_{2n+1},z^{j} \rangle_{\mathscr{L}} \! 
= \! 0$, $j \! = \! -n,\dotsc,n$;
\item[(3)] $\langle \boldsymbol{\pi}_{2n},\boldsymbol{\pi}_{2n} \rangle_{
\mathscr{L}} \! =: \! \norm{\boldsymbol{\pi}_{2n}(\cdot)}^{2}_{\mathscr{L}} 
\! = \! (\xi^{(2n)}_{n})^{-2}$;
\item[(4)] $\langle \boldsymbol{\pi}_{2n+1},\boldsymbol{\pi}_{2n+1} \rangle_{
\mathscr{L}} \! =: \! \norm{\boldsymbol{\pi}_{2n+1}(\cdot)}^{2}_{\mathscr{L}} 
\! = \! (\xi^{(2n+1)}_{-n-1})^{-2}$.
\end{compactenum}
Conditions~(3) and~(4) above imply, respectively, that $\xi^{(2n)}_{n} \! = 
\! 1/\norm{\boldsymbol{\pi}_{2n}(\pmb{\cdot})}_{\mathscr{L}}$ $(> \! 0)$ and 
$\xi^{(2n+1)}_{-n-1} \! = \! 1/\norm{\boldsymbol{\pi}_{2n+1}(\pmb{\cdot})}_{
\mathscr{L}}$ $(> \! 0)$.

In terms of the Hankel determinants, $H^{(m)}_{k}$, $(m,k) \! \in \! \mathbb{
Z} \times \mathbb{N}$, associated with the real-valued, bi-infinite, strong
moment sequence $\left\lbrace c_{k} \! = \! \int_{\mathbb{R}}s^{k} \exp 
(-\mathscr{N} \, V(s)) \, \md s, \, \mathscr{N} \! \in \! \mathbb{N} 
\right\rbrace_{k \in \mathbb{Z}}$, the monic orthogonal $L$-polynomials, 
$\boldsymbol{\pi}_{j}(z)$, $j \! \in \! \mathbb{Z}_{0}^{+}$, are represented 
via the following determinantal formulae \cite{a6,a7,a11,a17}: for $m \! 
\in \! \mathbb{Z}_{0}^{+}$,
\begin{gather}
\boldsymbol{\pi}_{2m}(z) \! = \! \dfrac{1}{H^{(-2m)}_{2m}}
\begin{vmatrix}
c_{-2m} & c_{-2m+1} & \cdots & c_{-1} & z^{-m} \\
c_{-2m+1} & c_{-2m+2} & \cdots & c_{0} & z^{-m+1} \\
\vdots & \vdots & \ddots & \vdots & \vdots \\
c_{-1} & c_{0} & \cdots & c_{2m-2} & z^{m-1} \\
c_{0} & c_{1} & \cdots & c_{2m-1} & z^{m}
\end{vmatrix}, \\
\intertext{and}
\boldsymbol{\pi}_{2m+1}(z) \! = \! -\dfrac{1}{H^{(-2m)}_{2m+1}}
\begin{vmatrix}
c_{-2m-1} & c_{-2m} & \cdots & c_{-1} & z^{-m-1} \\
c_{-2m} & c_{-2m+1} & \cdots & c_{0} & z^{-m} \\
\vdots & \vdots & \ddots & \vdots & \vdots \\
c_{-1} & c_{0} & \cdots & c_{2m-1} & z^{m-1} \\
c_{0} & c_{1} & \cdots & c_{2m} & z^{m}
\end{vmatrix};
\end{gather}
moreover, it can be shown that (see, for example, \cite{a11,a17}), for $n \!
\in \! \mathbb{Z}_{0}^{+}$,
\begin{gather}
\xi^{(2n)}_{n} \! \left(= \! \dfrac{1}{\norm{\boldsymbol{\pi}_{2n}(\cdot)}_{
\mathscr{L}}} \right) \! = \sqrt{\dfrac{H^{(-2n)}_{2n}}{H^{(-2n)}_{2n+1}}},
\qquad \qquad \xi^{(2n+1)}_{-n-1} \! \left(= \! \dfrac{1}{\norm{\boldsymbol{
\pi}_{2n+1}(\cdot)}_{\mathscr{L}}} \right) \! = \sqrt{\dfrac{H^{(-2n)}_{2n+
1}}{H^{(-2n-2)}_{2n+2}}}, \\
\nu^{(2n)}_{-n} \! \left(:= \! \dfrac{\xi^{(2n)}_{-n}}{\xi^{(2n)}_{n}} \right)
\! = \dfrac{H^{(-2n+1)}_{2n}}{H^{(-2n)}_{2n}}, \qquad \qquad \nu^{(2n+1)}_{n}
\! \left(:= \! \dfrac{\xi^{(2n+1)}_{n}}{\xi^{(2n+1)}_{-n-1}} \right) \! =
-\dfrac{H^{(-2n-1)}_{2n+1}}{H^{(-2n)}_{2n+1}}.
\end{gather}
\begin{eeeee}
It is important to note \cite{a10} that there is a more general theory, 
where moments corresponding to an arbitrary, countable sequence of (fixed) 
points are involved, and where \emph{orthogonal rational functions} \cite{a23} 
(see, also, \cite{a24}) play the r\^{o}le of orthogonal $L$-polynomials. 
Furthermore, since orthogonal $L$-polynomials are rational functions with 
(fixed) poles at the origin and at the point at infinity, the step towards 
a more general theory where poles are at arbitrary, but fixed, 
positions/locations in $\mathbb{C} \cup \{\infty\}$ is natural. \hfill 
$\blacksquare$
\end{eeeee}

A basic fact about orthogonal $L$-polynomials is \cite{a6,a7,a17}:
\begin{enumerate}
\item[$\circledast$] for each $m \! \in \! \mathbb{Z}_{0}^{+}$, all the roots
of $\boldsymbol{\pi}_{2m}(z)$ and $\boldsymbol{\pi}_{2m+1}(z)$ are real,
simple, and non-zero.
\end{enumerate}
In the literature (see, for example, \cite{a6,a7,a17}), there are two 
distinguished cases, depending on whether or not the trailing coefficients 
vanish. Recall that, for each $m \! \in \! \mathbb{Z}_{0}^{+}$, the trailing 
coefficient, $TC(\boldsymbol{\pi}_{m})$, is defined as follows:
$TC(\boldsymbol{\pi}_{m}) \! := \!
\begin{cases}
\nu^{(2n)}_{-n}, &\text{$m \! = \! 2n$,} \\
\nu^{(2n+1)}_{n}, &\text{$m \! = \! 2n \! + \! 1$.}
\end{cases}$ If, for each $m \! \in \! \mathbb{Z}_{0}^{+}$, $TC(\boldsymbol{
\pi}_{m}) \! \neq \! 0$, then $\boldsymbol{\pi}_{m}(z)$ and the index $m$ are
called \emph{non-singular}; otherwise, if $TC(\boldsymbol{\pi}_{m}) \! = \!
0$, then $\boldsymbol{\pi}_{m}(z)$ and the index $m$ are called
\emph{singular}. {}From Equations~(1.9), it can be seen that, for each $m 
\! \in \! \mathbb{Z}_{0}^{+}$:
\begin{compactenum}
\item[(i)] $\boldsymbol{\pi}_{2m}(z)$ is non-singular (resp., singular) if
$H^{(-2m+1)}_{2m} \! \not= \! 0$ (resp., $H^{(-2m+1)}_{2m} \! = \! 0)$;
\item[(ii)] $\boldsymbol{\pi}_{2m+1}(z)$ is non-singular (resp., singular) 
if $H^{(-2m-1)}_{2m+1} \! \not= \! 0$ (resp., $H^{(-2m-1)}_{2m+1} \! = \! 
0)$.
\end{compactenum}
It is an established fact \cite{a6,a7,a17} that, for $m \! \in \! \mathbb{
Z}_{0}^{+}$:
\begin{compactenum}
\item[(1)] if $\boldsymbol{\pi}_{2m}(z)$ is non-singular, then it possesses 
$2m$ roots and $\boldsymbol{\pi}_{2m}(0) \! \neq \! 0$, and if $\boldsymbol{
\pi}_{2m}(z)$ is singular, then it possesses $2m \! - \! 1$ roots;
\item[(2)] if $\boldsymbol{\pi}_{2m+1}(z)$ is non-singular, then it possesses 
$2m \! + \! 1$ roots and $\boldsymbol{\pi}_{2m+1}(0) \! \neq \! 0$, and if 
$\boldsymbol{\pi}_{2m+1}(z)$ is singular, then it possesses $2m$ roots.
\end{compactenum}

For each $n \! \in \! \mathbb{Z}_{0}^{+}$, it can be shown that, via a 
straightforward factorisation argument and using Equations~(1.6) and~(1.7):
\begin{compactenum}
\item[(i)] if $\boldsymbol{\pi}_{2n}(z)$ is non-singular, setting $\left\{
\mathstrut \alpha^{(2n)}_{k}, \, k \! = \! 1,\dotsc,2n \right\} \! := \!
\lbrace \mathstrut z; \, \boldsymbol{\pi}_{2n}(z) \! = \! 0 \rbrace$,
\begin{equation}
\prod_{k=1}^{2n} \alpha^{(2n)}_{k} \! = \nu^{(2n)}_{-n} \qquad (\in \!
\mathbb{R});
\end{equation}
\item[(ii)] if $\boldsymbol{\pi}_{2n}(z)$ is singular, setting (with abuse 
of notation) $\left\{\mathstrut \alpha^{(2n)}_{k}, \, k \! = \! 1,\dotsc,2n 
\! - \! 1 \right\} \! := \! \lbrace \mathstrut z; \, \boldsymbol{\pi}_{2n}
(z) \! = \! 0 \rbrace$,
\begin{equation}
\quad \prod_{k=1}^{2n-1} \alpha^{(2n)}_{k} \! = \! -\nu^{(2n)}_{-(n-1)} \!
:= \! -\xi^{(2n)}_{-(n-1)}(\xi^{(2n)}_{n})^{-1} \qquad (\in \! \mathbb{R});
\end{equation}
\item[(iii)] if $\boldsymbol{\pi}_{2n+1}(z)$ is non-singular, setting $\left\{
\mathstrut \alpha^{(2n+1)}_{k}, \, k \! = \! 1,\dotsc,2n \! + \! 1 \right\}
\! := \! \lbrace \mathstrut z; \, \boldsymbol{\pi}_{2n+1}(z) \! = \! 0
\rbrace$,
\begin{equation}
\left(\prod_{k=1}^{2n+1} \alpha^{(2n+1)}_{k} \right)^{-1} \! = \! -\nu^{(2n+1)
}_{n} \qquad (\in \! \mathbb{R});
\end{equation}
\item[(iv)] if $\boldsymbol{\pi}_{2n+1}(z)$ is singular, setting (with abuse
of notation) $\left\{\mathstrut \alpha^{(2n+1)}_{k}, \, k \! = \! 1,\dotsc,
2n \right\} \! := \! \lbrace \mathstrut z; \, \boldsymbol{\pi}_{2n+1}(z) \! = 
\! 0 \rbrace$,
\begin{equation}
\left(\prod_{k=1}^{2n} \alpha^{(2n+1)}_{k} \right)^{-1} \! =\nu^{(2n+1)}_{n-1}
\! := \! \xi^{(2n+1)}_{n-1}(\xi^{(2n+1)}_{-n-1})^{-1} \qquad (\in \! \mathbb{
R}).
\end{equation}
\end{compactenum}
Via Equations~(1.8) and~(1.9), and using the fact that $\xi^{(0)}_{0} \!
= \! H^{(0)}_{0} \! = \! H^{(0)}_{1} \! = \! 1$ (cf. Equations~(1.1)), a
straightforward calculation shows that, for $n \! \in \! \mathbb{Z}_{0}^{+}$,
\begin{gather}
H^{(-2n)}_{2n} \! = \! \dfrac{(\xi^{(2n)}_{n})^{2}(\xi^{(2n+1)}_{-n-1})^{2}}{
\prod_{j=0}^{n}(\xi^{(2j)}_{j})^{2} \prod_{k=0}^{n}(\xi^{(2k+1)}_{-k-1})^{2}}
\qquad (> \! 0), \\
H^{(-2n)}_{2n+1} \! = \! \dfrac{(\xi^{(2n+1)}_{-n-1})^{2}}{\prod_{j=0}^{n}
(\xi^{(2j)}_{j})^{2} \prod_{k=0}^{n}(\xi^{(2k+1)}_{-k-1})^{2}} \qquad (> \!
0), \\
H^{(-2n+1)}_{2n} \! = \! \dfrac{\nu^{(2n)}_{-n}(\xi^{(2n)}_{n})^{2}(\xi^{(2n+
1)}_{-n-1})^{2}}{\prod_{j=0}^{n}(\xi^{(2j)}_{j})^{2} \prod_{k=0}^{n}(\xi^{(2k
+1)}_{-k-1})^{2}} \qquad (\in \! \mathbb{R}), \\
H^{(-2n-1)}_{2n+1} \! = \! -\dfrac{\nu^{(2n+1)}_{n}(\xi^{(2n+1)}_{-n-1})^{2}}{
\prod_{j=0}^{n}(\xi^{(2j)}_{j})^{2} \prod_{k=0}^{n}(\xi^{(2k+1)}_{-k-1})^{2}}
\qquad (\in \! \mathbb{R}), \\
\dfrac{H^{(-2n-2)}_{2n+2}}{H^{(-2n)}_{2n}} \! = \! \dfrac{1}{(\xi^{(2n)}_{n}
)^{2}(\xi^{(2n+1)}_{-n-1})^{2}} \qquad (> \! 0), \\
\dfrac{H^{(-2n-2)}_{2n+3}}{H^{(-2n)}_{2n+1}} \! = \! \dfrac{1}{(\xi^{(2n+2)}_{
n+1})^{2}(\xi^{(2n+1)}_{-n-1})^{2}}  \qquad (> \! 0),
\end{gather}
with
\begin{equation*}
\xi^{(2n+2)}_{n+1} \! := \! h^{+} \! \left[\xi^{(2n)}_{n} \right],
\end{equation*}
where the homomorphisms $h^{\pm}$ (shift-up, $+$, and shift-down, $-$, 
operators) are defined as follows:
\begin{equation*}
h^{\pm} \colon \mathbb{Z} \! \to \! \mathbb{Z}, \, \, n \! \mapsto \! h^{\pm}
[n] \! := \! n \! \pm \! 1 \qquad \, \, \text{and} \qquad \, \, f(n) \! 
\mapsto \! h^{\pm}[f(n)] \! := \! f(n \! \pm \! 1);
\end{equation*}
furthermore, denoting by $\id_{\mathbb{Z}}$ the identity operator on $\mathbb{
Z}$, that is, $\id_{\mathbb{Z}} \colon \mathbb{Z} \! \to \! \mathbb{Z}$, $n 
\! \mapsto \! \id_{\mathbb{Z}}[n] \! := \! n$ and $f(n) \! \mapsto \! \id_{
\mathbb{Z}}[f(n)] \! := \! f(n)$, one arrives at (the Abelian relation) 
$h^{+}h^{-} \! = \! h^{-}h^{+} \! = \! \id_{\mathbb{Z}}$.
\begin{eeeee}
It turns out that, for the presentation of the asymptotic results of this 
work, it is convenient to re-write $\md \widetilde{\mu}(z) \! = \! \exp (-
\mathscr{N} \, V(z)) \, \md z \! = \! \exp (-n \widetilde{V}(z)) \, \md z 
\! =: \! \md \mu (z)$, $n \! \in \! \mathbb{N}$, where
\begin{equation*}
\widetilde{V}(z) \! = \! z_{o}V(z),
\end{equation*}
with
\begin{equation*}
z_{o} \colon \mathbb{N} \times \mathbb{N} \! \to \! (0,+\infty), \, \,
(\mathscr{N},n) \! \mapsto \! z_{o} \! := \mathscr{N}/n,
\end{equation*}
and where the `scaled' external field $\widetilde{V} \colon \mathbb{R}
\setminus \{0\} \! \to \! \mathbb{R}$ satisfies:
\begin{gather}
\widetilde{V} \, \, \text{is real analytic on} \, \, \mathbb{R} \setminus
\{0\}; \\
\lim_{\vert x \vert \to \infty} \! \left(\widetilde{V}(x)/\ln (x^{2} \! + \!
1) \right) \! = \! +\infty; \\
\lim_{\vert x \vert \to 0} \! \left(\widetilde{V}(x)/\ln (x^{-2} \! + \! 1)
\right) \! = \! +\infty.
\end{gather}
(For example, a rational function of the form $\widetilde{V}(z) \! = \! \sum_{
k=-2m_{1}}^{2m_{2}} \widetilde{\varrho}_{k}z^{k}$, with $\widetilde{\varrho}_{
k} \! \in \! \mathbb{R}$, $k \! = \! -2m_{1},\dotsc,2m_{2}$, $m_{1,2} \! \in
\! \mathbb{N}$, and $\widetilde{\varrho}_{-2m_{1}},\widetilde{\varrho}_{2
m_{2}} \! > \! 0$ would satisfy conditions~(1.20)--(1.22).) Note that, under
the homomorphisms $h^{\pm} \colon \mathbb{Z} \! \to \! \mathbb{Z}$,
\begin{equation*}
z_{o} \! \mapsto \! h^{\pm}[z_{o}] \! = \! z_{o} \! \left(1 \! \pm \! 
\dfrac{1}{n} \right)^{-1}.
\end{equation*}
Hereafter, the double-scaling limit as $\mathscr{N},n \! \to \! \infty$ such 
that $z_{o} \! = \! 1 \! + \! o(1)$ is studied: the simplified `notation' 
$n \! \to \! \infty$ will be adopted. \hfill $\blacksquare$
\end{eeeee}

\begin{eeeee}
The reader is reminded that, in the more classical setting of standard 
Hankel determinants, they are directly related to the theory of random 
matrices; indeed, as is well known, the partition function in Random Matrix 
Theory is actually a Hankel determinant. The asymptotic (as $n \! \to \! 
\infty)$ behaviour of Hankel determinants with respect to varying exponential 
weights has only been proved for the `single-interval/one-cut' case 
\cite{a25}. In the present, more general, setting, the asymptotic (as 
$n \! \to \! \infty)$ behaviour of the generalised Hankel 
determinants~(1.14)--(1.17) with respect to the varying exponential weight 
$w(z) \! = \! \exp (-n \widetilde{V}(z))$, $n \! \in \! \mathbb{N}$, where 
$\widetilde{V} \colon \mathbb{R} \setminus \{0\} \! \to \! \mathbb{R}$ 
satisfies conditions~(1.20)--(1.22), is a very challenging undertaking; 
however, with the asymptotic (as $n \! \to \! \infty)$ results of this work, 
one can easily obtain (real-valued) asymptotics (as $n \! \to \! \infty)$ 
for $\mathfrak{P}_{\xi}H^{(-2n)}_{2n}$, $\mathfrak{P}_{\xi}H^{(-2n)}_{2n+1}$, 
$\mathfrak{P}_{\xi}H^{(-2n+1)}_{2n}$, and $\mathfrak{P}_{\xi}H^{(-2n-1)}_{2n
+1}$, where $\mathfrak{P}_{\xi} \! := \! \prod_{j=0}^{n}(\xi^{(2j)}_{j})^{2} 
\prod_{k=0}^{n}(\xi^{(2k+1)}_{-k-1})^{2}$. \hfill $\blacksquare$
\end{eeeee}

Unlike orthogonal polynomials, which satisfy a system of three-term recurrence 
relations, orthonormal/orthogonal $L$-polynomials satisfy more complicated 
recurrence relations: there are several different versions of these recurrence 
relations \cite{a11,a16,a17,a26,a27,a28}. The general form of these (system 
of) recurrence relations is a pair of three- and five-term recurrence 
relations (see, for example, \cite{a26}): for $n \! \in \! \mathbb{Z}_{
0}^{+}$,
\begin{gather}
z \phi_{2n+1}(z) \! = \! b_{2n+1}^{\sharp} \phi_{2n}(z) \! + \! a_{2n+1}^{
\sharp} \phi_{2n+1}(z) \! + \! b_{2n+2}^{\sharp} \phi_{2n+2}(z), \\
z \phi_{2n}(z) \! = \! c_{2n}^{\sharp} \phi_{2n-2}(z) \! + \! b_{2n}^{\sharp}
\phi_{2n-1}(z) \! + \! a_{2n}^{\sharp} \phi_{2n}(z) \! + \! b_{2n+1}^{\sharp}
\phi_{2n+1}(z) \! + \! c_{2n+2}^{\sharp} \phi_{2n+2}(z),
\end{gather}
where $c_{0}^{\sharp} \! = \! b_{0}^{\sharp} \! = \! 0$, and $c_{2k}^{\sharp}
\! > \! 0$, $k \! \in \! \mathbb{N}$, and
\begin{gather}
z^{-1} \phi_{2n}(z) \! = \! \beta_{2n}^{\sharp} \phi_{2n-1}(z) \! + \!
\alpha_{2n}^{\sharp} \phi_{2n}(z) \! + \! \beta_{2n+1}^{\sharp} \phi_{2n+1}
(z), \\
z^{-1} \phi_{2n+1}(z) \! = \! \gamma_{2n+1}^{\sharp} \phi_{2n-1}(z) \! + \!
\beta_{2n+1}^{\sharp} \phi_{2n}(z) \! + \! \alpha_{2n+1}^{\sharp} \phi_{2n+1}
(z) \! + \! \beta_{2n+2}^{\sharp} \phi_{2n+2}(z) \! + \! \gamma_{2n+3}^{
\sharp} \phi_{2n+3}(z),
\end{gather}
where $\beta_{0}^{\sharp} \! = \! \gamma_{1}^{\sharp} \! = \! 0$, $\beta_{
1}^{\sharp} \! > \! 0$, and $\gamma_{2k+1}^{\sharp} \! > \! 0$, $k \! \in 
\! \mathbb{N}$, leading, respectively, to the \emph{real-symmetric}, 
\emph{tri-penta-diagonal-type Laurent-Jacobi matrices}, $\mathcal{F}$ and 
$\mathcal{G}$, for the mappings $\mathscr{F} \colon \Lambda^{\mathbb{R}} \! 
\to \! \Lambda^{\mathbb{R}}$, $\mathfrak{f}(z) \! \mapsto \! z \mathfrak{f}
(z)$, and $\mathscr{G} \colon \Lambda^{\mathbb{R}} \! \to \! \Lambda^{
\mathbb{R}}$, $\mathfrak{g}(z) \! \mapsto \! z^{-1} \mathfrak{g}(z)$,
\begin{align}
\setcounter{MaxMatrixCols}{18}
\mathcal{F} & = \!
\begin{pmatrix}
a_{0}^{\sharp} & b_{1}^{\sharp} & c_{2}^{\sharp} & & & & & & & & & & & & & \\
b_{1}^{\sharp} & a_{1}^{\sharp} & b_{2}^{\sharp} & & & & & & & & & & & & & \\
c_{2}^{\sharp} & b_{2}^{\sharp} & a_{2}^{\sharp} & b_{3}^{\sharp} & c_{4}^{
\sharp} & & & & & & & & & & & \\
 & & b_{3}^{\sharp} & a_{3}^{\sharp} & b_{4}^{\sharp} & & & & & & & & & & & \\
 & & c_{4}^{\sharp} & b_{4}^{\sharp} & a_{4}^{\sharp} & b_{5}^{\sharp} & c_{
6}^{\sharp} & & & & & & & & & \\
 & & & & b_{5}^{\sharp} & a_{5}^{\sharp} & b_{6}^{\sharp} & & & & & & & & & \\
 & & & & c_{6}^{\sharp} & b_{6}^{\sharp} & a_{6}^{\sharp} & b_{7}^{\sharp} &
c_{8}^{\sharp} & & & & & & & \\
 & & & & & & b_{7}^{\sharp} & a_{7}^{\sharp} & b_{8}^{\sharp} & & & & & & & \\
 & & & & & & c_{8}^{\sharp} & b_{8}^{\sharp} & a_{8}^{\sharp} & b_{9}^{\sharp}
& c_{10}^{\sharp} & & & & & \\
 & & & & & & & & \ddots & \ddots & \ddots & & & & & \\
 & & & & & & & & b_{2k+1}^{\sharp} & a_{2k+1}^{\sharp} & b_{2k+2}^{\sharp} & &
& & & \\
 & & & & & & & & c_{2k+2}^{\sharp} & b_{2k+2}^{\sharp} & a_{2k+2}^{\sharp} &
b_{2k+3}^{\sharp} & c_{2k+4}^{\sharp} & & & \\
 & & & & & & & & & & \ddots & \ddots & \ddots & & &
\end{pmatrix},
\end{align}
and
\begin{align}
\setcounter{MaxMatrixCols}{18}
\mathcal{G} &= \!
\begin{pmatrix}
\alpha_{0}^{\sharp} & \beta_{1}^{\sharp} & & & & & & & & & & & & & & \\
\beta_{1}^{\sharp} & \alpha_{1}^{\sharp} & \beta_{2}^{\sharp} & \gamma_{3}^{
\sharp}  & & & & & & & & & & & & \\
 & \beta_{2}^{\sharp} & \alpha_{2}^{\sharp} & \beta_{3}^{\sharp} & & & & & & &
& & & & & \\
 & \gamma_{3}^{\sharp} & \beta_{3}^{\sharp} & \alpha_{3}^{\sharp} & \beta_{4}^{
\sharp}  & \gamma_{5}^{\sharp} & & & & & & & & & & \\
 & & & \beta_{4}^{\sharp} & \alpha_{4}^{\sharp} & \beta_{5}^{\sharp} & & & & &
& & & & & \\
 & & & \gamma_{5}^{\sharp} & \beta_{5}^{\sharp} & \alpha_{5}^{\sharp} & \beta_{
6}^{\sharp} & \gamma_{7}^{\sharp} & & & & & & & & \\
 & & & & & \beta_{6}^{\sharp} & \alpha_{6}^{\sharp} & \beta_{7}^{\sharp} & & &
& & & & & \\
 & & & & & \gamma_{7}^{\sharp} & \beta_{7}^{\sharp} & \alpha_{7}^{\sharp} &
\beta_{8}^{\sharp} & \gamma_{9}^{\sharp} & & & & & & \\
 & & & & & & \beta_{8}^{\sharp} & \alpha_{8}^{\sharp} & \beta_{9}^{\sharp} & &
& & & & & \\
 & & & & & & \ddots & \ddots & \ddots & & & & & & & \\
 & & & & & & & & \gamma_{2k+1}^{\sharp} & \beta_{2k+1}^{\sharp} & \alpha_{2k+
1}^{\sharp} & \beta_{2k+2}^{\sharp} & \gamma_{2k+3}^{\sharp} & & & \\
 & & & & & & & & & & \beta_{2k+2}^{\sharp} & \alpha_{2k+2}^{\sharp} & \beta_{2
k+3}^{\sharp} & & & \\
 & & & & & & & & & & \ddots & \ddots & \ddots & & &
\end{pmatrix},
\end{align}
with zeros outside the indicated diagonals; furthermore, as shown in 
\cite{a26}, $\mathcal{F}$ and $\mathcal{G}$ are formal inverses, that 
is, $\mathcal{F} \mathcal{G} \! = \! \mathcal{G} \mathcal{F} \! = \! 
\operatorname{diag}(1,1,\dotsc,1,\dotsc)$.
\begin{eeeee}
There are multitudinous applications of orthogonal $L$-polynomials and their 
associated Laurent-Jacobi matrices. To list a few: (i) numerical analysis 
(quadrature formulae) and trigonometric moment problems \cite{a29}; (ii) the 
spectral theory of one-to-one self-adjoint operators in infinite-dimensional 
(necessarily separable) Hilbert spaces \cite{a26,a30}; (iii) complex 
approximation theory (two-point Pad\'{e} approximants) \cite{a31}; (iv) the 
direct/inverse scattering theory for the finite relativistic Toda lattice 
\cite{a32,a33} (see, also, the recent work \cite{a34}, related to an extension 
of the Toda lattice); and (v) the classical Pick-Nevanlinna problem 
\cite{a35}. Some further comments on these connections can be found in 
Section~1 of \cite{a21} (see, also, \cite{a23}, and the many references 
therein). It turns out that $n \! \to \! \infty$ asymptotics of orthogonal 
$L$-polynomials are an essential calculational ingredient in analyses related 
to the above-mentioned, seemingly disparate, topics. \hfill $\blacksquare$
\end{eeeee}

It is, by now, a well-known, if not established, mathematical fact that 
variational conditions for minimisation problems in logarithmic potential 
theory, via the \emph{equilibrium measure} \cite{a36,a37,a38,a39}, play a 
crucial r\^{o}le in the asymptotic analysis of (matrix) Riemann-Hilbert 
problems (RHPs) associated with (continuous and discrete) orthogonal 
polynomials, their roots, and corresponding recurrence relation coefficients 
(see, for example, \cite{a40,a41,a42,a43,a44,a45,a46,a47}). The situation 
with respect to the large-$n$ asymptotic analysis for the monic orthogonal 
$L$-polynomials, $\boldsymbol{\pi}_{n}(z)$, is analogous; however, unlike 
the asymptotic analysis for the orthogonal polynomials case, the asymptotic 
analysis for $\boldsymbol{\pi}_{n}(z)$ requires the consideration of two 
different families of RHPs, one for even degree (see \textbf{RHP1} below) 
and one for odd degree (see \textbf{RHP2} below). Thus, one must consider 
two sets of variational conditions for two (suitably posed) minimisation 
problems. Recall that the orthonormal $L$-polynomials, $\lbrace \phi_{m}(z) 
\rbrace_{m \in \mathbb{Z}_{0}^{+}}$ (cf. Equations~(1.2) and~(1.3)), are 
related to the monic orthogonal $L$-polynomials, $\lbrace \boldsymbol{\pi}_{m}
(z) \rbrace_{m \in \mathbb{Z}_{0}^{+}}$, via Equations~(1.4) and~(1.5). The 
even degree and odd degree monic orthogonal $L$-polynomial problems for 
$\boldsymbol{\pi}_{2n}(z)$ $(m \! = \! 2n)$ and $\boldsymbol{\pi}_{2n+1}(z)$ 
$(m \! = \! 2n \! + \! 1)$, respectively, are formulated, \emph{\`{a} la} 
Fokas-Its-Kitaev \cite{a48,a49}, as the following matrix RHPs on $\mathbb{R}$.

Before stating them, however, the notations/nomenclaturae used throughout 
the text are now summarised.
\begin{center}
\Ovalbox{\textsc{Notational Conventions}}
\end{center}
\begin{enumerate}
\item[(1)] $\mathrm{I} \! = \!
\left(
\begin{smallmatrix}
1 & 0 \\
0 & 1
\end{smallmatrix}
\right)$ is the $2 \times 2$ identity matrix, $\sigma_{1} \! = \!
\left(
\begin{smallmatrix}
0 & 1 \\
1 & 0
\end{smallmatrix}
\right)$, $\sigma_{2} \! = \!
\left(
\begin{smallmatrix}
0 & -\mi \\
\mi & 0
\end{smallmatrix}
\right)$, and $\sigma_{3} \! = \!
\left(
\begin{smallmatrix}
1 & 0 \\
0 & -1
\end{smallmatrix}
\right)$ are the Pauli matrices, $\sigma_{\pm} \! := \! \tfrac{1}{2}(\sigma_{
1} \! \pm \! \mi \sigma_{2})$, $\mathbb{R}_{\pm} \! := \! \left\lbrace 
\mathstrut x \! \in \! \mathbb{R}; \, \pm x \! > \! 0 \right\rbrace$, and 
$\mathbb{C}_{\pm} \! := \! \left\lbrace \mathstrut z \! \in \! \mathbb{C}; \, 
\pm \Im (z) \! > \! 0 \right\rbrace$;
\item[(2)] a contour $\mathcal{D}$ which is the finite union of 
piecewise-smooth, simple curves (as closed sets) is said to be 
\emph{orientable} if its complement, $\mathbb{C} \setminus \mathcal{D}$, can 
always be divided into two, possibly disconnected, disjoint open sets $\pmb{
\mho^{+}}$ and $\pmb{\mho^{-}}$, either of which has finitely many components, 
such that $\mathcal{D}$ admits an orientation so that it can either be viewed 
as a positively oriented boundary $\mathcal{D}^{+}$ for $\pmb{\mho^{+}}$ or 
as a negatively oriented boundary $\mathcal{D}^{-}$ for $\pmb{\mho^{-}}$ 
\cite{a50}, that is, the (possibly disconnected) components of $\mathbb{C} 
\setminus \mathcal{D}$ can be coloured by $+$ or $-$ in such a way that the 
$+$ regions do not share boundary with the $-$ regions, except, possibly, 
at finitely many points \cite{a51};
\item[(3)] for each segment of an oriented contour $\mathcal{D}$, according to 
the given orientation, the ``+'' side is to the left and the ``-'' side is to 
the right as one traverses the contour in the direction of orientation; so, 
for a matrix-valued function $\mathcal{A}(z)$, $\mathcal{A}_{\pm}(z)$ denote 
the non-tangential limits $\mathcal{A}_{\pm}(z) \! := \! \lim_{\genfrac{}{}
{0pt}{2}{z^{\prime} \to z}{z^{\prime} \in \pm \, \mathrm{side} \, \mathrm{of} 
\, \mathcal{D}}} \mathcal{A}(z^{\prime})$;
\item[(4)] unless stated otherwise, exponents such as $z^{\nu}$ shall be
interpreted as principal branches;
\item[(5)] for some point set $\mathcal{D} \subset \mathcal{X}$, with
$\mathcal{X} \! = \! \mathbb{C}$ or $\mathbb{R}$, $\overline{\mathcal{D}}:= \!
\mathcal{D} \cup \partial \mathcal{D}$, and $\mathcal{D}^{c} \! := \! \mathcal{
X} \setminus \overline{\mathcal{D}}$;
\item[(6)] unless stated otherwise, $\mathbb{R}$ is oriented {}from $-\infty$
to $+\infty$;
\item[(7)] for a $2 \times 2$ matrix-valued function $\mathscr{A}(z)$,
$\mathscr{A}_{ij}(z)$, $i,j \! = \! 1,2$, denotes the $(i \, j)$-element of
$\mathscr{A}(z)$;
\item[(8)] for a $2 \times 2$ matrix-valued function $\mathfrak{Y}(z)$, the
notation $\mathfrak{Y}(z) \! =_{z \to z_{0}} \! \mathcal{O}(\ast)$ means
$\mathfrak{Y}_{ij}(z) \! =_{z \to z_{0}} \! \mathcal{O}(\ast_{ij})$, $i,j \! 
= \! 1,2$ (\emph{mutatis mutandis} for $o(1))$.
\end{enumerate}

The RHP which characterises the even degree monic orthogonal $L$-polynomials 
$\boldsymbol{\pi}_{2n}(z)$, $n \! \in \! \mathbb{Z}_{0}^{+}$, is now stated.
\begin{rhp1}[{\rm \cite{a21}}]
Let $V \colon \mathbb{R} \setminus \{0\} \! \to \! \mathbb{R}$ satisfy 
conditions~{\rm (V1)--(V3)}. Find $\overset{e}{\mathrm{Y}} \colon \mathbb{C} 
\setminus \mathbb{R} \! \to \! \operatorname{SL}_{2}(\mathbb{C})$ solving: 
{\rm (i)} $\overset{e}{\mathrm{Y}}(z)$ is holomorphic for $z \! \in \! 
\mathbb{C} \setminus \mathbb{R};$ {\rm (ii)} for $z \! \in \! \mathbb{R}$, 
the boundary values $\overset{e}{\mathrm{Y}}_{\pm}(z)$  satisfy the jump 
condition
\begin{equation*}
\overset{e}{\mathrm{Y}}_{+}(z)= \overset{e}{\mathrm{Y}}_{-}(z) \! \left(
\mathrm{I} \! + \! \exp (-\mathscr{N} \, V(z)) \sigma_{+} \right);
\end{equation*}
{\rm (iii)} $\overset{e}{\mathrm{Y}}(z)z^{-n \sigma_{3}} \! =_{\underset{z \in
\mathbb{C} \setminus \mathbb{R}}{z \to \infty}} \! \mathrm{I} \! + \! \mathcal{
O}(z^{-1});$ and {\rm (iv)} $\overset{e}{\mathrm{Y}}(z)z^{n \sigma_{3}} \!
=_{\underset{z \in \mathbb{C} \setminus \mathbb{R}}{z \to 0}} \! \mathcal{O}
(1)$.
\end{rhp1}
\begin{ccccc}[{\rm \cite{a21}}]
Let $\overset{e}{\mathrm{Y}} \colon \mathbb{C} \setminus \mathbb{R} \! \to 
\! \mathrm{SL}_{2}(\mathbb{C})$ solve {\rm \pmb{RHP1}}. {\rm \pmb{RHP1}} 
possesses a unique solution given by: {\rm (i)} for $n \! = \! 0$,
\begin{equation*}
\overset{e}{\mathrm{Y}}(z) \! = \!
\begin{pmatrix}
1 & \int_{\mathbb{R}} \frac{\exp (-\mathscr{N} \, V(s))}{s-z} \, \frac{\md
s}{2 \pi \mi} \\
0 & 1
\end{pmatrix}, \quad z \! \in \! \mathbb{C} \setminus \mathbb{R},
\end{equation*}
where $\bm{\pi}_{0}(z) \! := \! \overset{e}{\mathrm{Y}}_{11}(z) \! \equiv \!
1;$ and {\rm (ii)} for $n \! \in \! \mathbb{N}$,
\begin{equation}
\overset{e}{\mathrm{Y}}(z) \! = \!
\begin{pmatrix}
\boldsymbol{\pi}_{2n}(z) & \int_{\mathbb{R}} \frac{\boldsymbol{\pi}_{2n}(s)
\exp (-\mathscr{N} \, V(s))}{s-z} \, \frac{\md s}{2 \pi \mi} \\
\overset{e}{\mathrm{Y}}_{21}(z) & \int_{\mathbb{R}} \frac{\overset{e}{\mathrm{
Y}}_{21}(s) \exp (-\mathscr{N} \, V(s))}{s-z} \, \frac{\md s}{2 \pi \mi}
\end{pmatrix}, \quad z \! \in \! \mathbb{C} \setminus \mathbb{R},
\end{equation}
where $\overset{e}{\mathrm{Y}}_{21} \colon \mathbb{C}^{\ast} \! \to \!
\mathbb{C}$, and $\boldsymbol{\pi}_{2n}(z)$ is the even degree monic {\rm OLP}
defined in Equation~{\rm (1.4)}.
\end{ccccc}

The RHP which characterises the odd degree monic orthogonal $L$-polynomials
$\boldsymbol{\pi}_{2n+1}(z)$, $n \! \in \! \mathbb{Z}_{0}^{+}$, is now stated.
\begin{rhp2}[{\rm \cite{a22}}]
Let $V \colon \mathbb{R} \setminus \{0\} \! \to \! \mathbb{R}$ satisfy
conditions~{\rm (V1)--(V3)}. Find $\overset{o}{\mathrm{Y}} \colon \mathbb{C}
\setminus \mathbb{R} \! \to \! \operatorname{SL}_{2}(\mathbb{C})$ solving:
{\rm (i)} $\overset{o}{\mathrm{Y}}(z)$ is holomorphic for $z \! \in \!
\mathbb{C} \setminus \mathbb{R};$ {\rm (ii)} for $z \! \in \! \mathbb{R}$, the
boundary values $\overset{o}{\mathrm{Y}}_{\pm}(z)$ satisfy the jump condition
\begin{equation*}
\overset{o}{\mathrm{Y}}_{+}(z)=\overset{o}{\mathrm{Y}}_{-}(z) \! \left(
\mathrm{I} \! + \! \exp (-\mathscr{N} \, V(z)) \sigma_{+} \right);
\end{equation*}
{\rm (iii)} $\overset{o}{\mathrm{Y}}(z)z^{n \sigma_{3}} \! =_{\underset{z \in
\mathbb{C} \setminus \mathbb{R}}{z \to 0}} \! \mathrm{I} \! + \! \mathcal{O}
(z);$ and {\rm (iv)} $\overset{o}{\mathrm{Y}}(z)z^{-(n+1) \sigma_{3}} \! =_{
\underset{z \in \mathbb{C} \setminus \mathbb{R}}{z \to \infty}} \! \mathcal{
O}(1)$.
\end{rhp2}
\begin{ccccc}[{\rm \cite{a22}}]
Let $\overset{o}{\mathrm{Y}} \colon \mathbb{C} \setminus \mathbb{R} \! \to 
\! \mathrm{SL}_{2}(\mathbb{C})$ solve {\rm \pmb{RHP2}}. {\rm \pmb{RHP2}} 
possesses a unique solution given by: {\rm (i)} for $n \! = \! 0$,
\begin{equation*}
\overset{o}{\mathrm{Y}}(z) \! = \!
\begin{pmatrix}
z \bm{\pi}_{1}(z) & z \int_{\mathbb{R}} \frac{(s \bm{\pi}_{1}(s)) \exp (-
\mathscr{N} \, V(s))}{s(s-z)} \, \frac{\md s}{2 \pi \mi} \\
2 \pi \mi z & 1 \! + \! z \int_{\mathbb{R}} \frac{\exp (-\mathscr{N} \, V(s))
}{s-z} \, \md s
\end{pmatrix}, \quad z \! \in \! \mathbb{C} \setminus \mathbb{R},
\end{equation*}
where $\bm{\pi}_{1}(z) \! = \! \tfrac{1}{z} \! + \! \tfrac{\xi^{(1)}_{0}}{
\xi^{(1)}_{-1}}$, with $\tfrac{\xi^{(1)}_{0}}{\xi^{(1)}_{-1}} \! = \! -\int_{
\mathbb{R}}s^{-1} \exp (-\mathscr{N} \, V(s)) \, \md s$, $\mathscr{N} \! \in
\! \mathbb{N};$ and {\rm (ii)} for $n \! \in \! \mathbb{N}$,
\begin{equation}
\overset{o}{\mathrm{Y}}(z) \! = \!
\begin{pmatrix}
z \boldsymbol{\pi}_{2n+1}(z) & z \int_{\mathbb{R}} \frac{(s \boldsymbol{\pi}_{
2n+1}(s)) \exp (-\mathscr{N} \, V(s))}{s(s-z)} \, \frac{\md s}{2 \pi \mi} \\
\overset{o}{\mathrm{Y}}_{21}(z) & z \int_{\mathbb{R}} \frac{\overset{o}{
\mathrm{Y}}_{21}(s) \exp (-\mathscr{N} \, V(s))}{s(s-z)} \, \frac{\md s}{2
\pi \mi}
\end{pmatrix}, \quad z \! \in \! \mathbb{C} \setminus \mathbb{R},
\end{equation}
where $\overset{o}{\mathrm{Y}}_{21} \colon \mathbb{C}^{\ast} \! \to \! 
\mathbb{C}$, and $\boldsymbol{\pi}_{2n+1}(z)$ is the odd degree monic 
{\rm OLP} defined in Equation~{\rm (1.5)}.
\end{ccccc}

This paper, which is not completely self-contained (see Remark~1.6 below), is
a continuation, proper, of \cite{a21,a22}, where, for example, asymptotics in 
the double-scaling limit as $\mathscr{N},n \! \to \! \infty$ such that $z_{o} 
\! = \! 1 \! + \! o(1)$ of the monic orthogonal $L$-poly\-n\-o\-m\-i\-a\-l\-s, 
$\lbrace \boldsymbol{\pi}_{n}(z) \rbrace_{n \in \mathbb{Z}_{0}^{+}}$, and 
orthonormal $L$-polynomials, $\lbrace \phi_{n}(z) \rbrace_{n \in \mathbb{
Z}_{0}^{+}}$, in the entire complex plane, were obtained. In \cite{a21,a22}, 
the $L$-polynomials were taken to be orthogonal with respect to the varying
exponential measure (cf. Remark~1.3)\footnote{Note that $LD(\boldsymbol{
\pi}_{m}) \! = \! LD(\phi_{m}) \! = \!
\begin{cases}
2n, &\text{$m \! = \! \text{even}$,} \\
2n \! + \! 1, &\text{$m \! = \! \text{odd}$,}
\end{cases}$ coincides with the parameter in the measure of orthogonality:
the large parameter, $n$, enters simultaneously into the $L$-degree of the
$L$-polynomials and the (varying exponential) weight; thus, asymptotics of
the $L$-polynomials are studied along a `diagonal strip' of a doubly-indexed
sequence.} $\md \mu (z) \! = \! \me^{-n \widetilde{V}(z)} \, \md z$, where
the external field\footnote{For real non-analytic external fields, see the
recent work \cite{a52}.} $\widetilde{V} \colon \mathbb{R} \setminus \{0\}
\! \to \! \mathbb{R}$ satisfies conditions~(1.20)--(1.22). Furthermore, 
asymptotics of the `even' and `odd' norming constants, $\xi^{(2n)}_{n}$ and 
$\xi^{(2n+1)}_{-n-1}$, respectively, as well as, subsequently, those of the 
(`even' and `odd') norms $\norm{\boldsymbol{\pi}_{2n}(\pmb{\cdot})}_{\mathscr{
L}}$, $\norm{\boldsymbol{\pi}_{2n+1}(\pmb{\cdot})}_{\mathscr{L}}$, and the 
Hankel determinant ratios $H^{(-2n)}_{2n}/H^{(-2n)}_{2n+1}$, $H^{(-2n)}_{2n
+1}/H^{(-2n-2)}_{2n+2}$, were also obtained. The aforementioned asymptotics 
were derived by extracting the large-$n$ behaviours of the respective 
solutions of \pmb{RHP1} and~\pmb{RHP2}. The paradigm for the asymptotic 
analysis of \pmb{RHP1} and~\pmb{RHP2} is a union of the Deift-Zhou (DZ) 
non-linear steepest-descent method \cite{a1,a2}, used for the asymptotic 
analysis of undulatory RHPs, and the extension of Deift-Venakides-Zhou 
\cite{a3}, incorporating into the DZ method a non-linear analogue of the 
WKB method, making the asymptotic analysis of fully non-linear problems 
tractable. It should be re-emphasized that, in this context, the 
equilibrium measure plays an absolutely crucial r\^{o}le in the analysis 
\cite{a36,a37,a38,a39}; multitudinous extensions and applications of the DZ 
method can be found in \cite{a40,a41,a42,a43,a44,a45,a46,a47,a53,a54,a55,%
a56,a57,a58,a59,a60,a61,a62,a63,a64,a65,a66,a67}.

In this paper, asymptotics for the following quantities are obtained:
\begin{compactenum}
\item[\textbullet] the (real) root products of the $L$-polynomials, for both
the non-singular and singular cases,
\begin{gather*}
\prod_{k=1}^{2n} \alpha^{(2n)}_{k}, \qquad \qquad \prod_{k=1}^{2n+1} \alpha^{
(2n+1)}_{k}, \qquad \qquad \prod_{k=1}^{2n-1} \alpha^{(2n)}_{k}, \qquad \qquad
\prod_{k=1}^{2n} \alpha^{(2n+1)}_{k};
\end{gather*}
\item[\textbullet] Hankel determinant ratios (all real valued) associated with
the bi-infinite, real-valued, strong moment sequence $\left\lbrace c_{k} \! =
\! \int_{\mathbb{R}}s^{k} \exp (-n \widetilde{V}(s)) \, \md s, \, n \! \in \!
\mathbb{N} \right\rbrace_{k \in \mathbb{Z}}$,
\begin{gather*}
\dfrac{H^{(-2n-2)}_{2n+2}}{H^{(-2n)}_{2n}}, \qquad \qquad \dfrac{H^{(-2n-2)}_{
2n+3}}{H^{(-2n)}_{2n+1}}, \qquad \qquad \dfrac{H^{(-2n+1)}_{2n}}{H^{(-2n)}_{2
n}}, \qquad \qquad \dfrac{H^{(-2n-1)}_{2n+1}}{H^{(-2n)}_{2n+1}};
\end{gather*}
\item[\textbullet] recurrence relation coefficients (all real valued),
\begin{gather*}
a_{2n}^{\sharp}, \, \, b_{2n}^{\sharp}, \, \, c_{2n}^{\sharp}, \, \, a_{2n+1}^{
\sharp}, \, \, b_{2n+1}^{\sharp}, \, \, c_{2n+2}^{\sharp}, \, \, b_{2n+2}^{
\sharp}, \, \, \alpha_{2n}^{\sharp}, \, \, \beta_{2n}^{\sharp}, \, \, \gamma_{
2n+1}^{\sharp}, \, \, \alpha_{2n+1}^{\sharp}, \, \, \beta_{2n+1}^{\sharp}, \,
\, \gamma_{2n+3}^{\sharp}, \, \, \beta_{2n+2}^{\sharp},
\end{gather*}
and subsequently those of the (real-symmetric, tri-penta-diagonal-type)
Laurent-Jacobi matrices $\mathcal{F}$ and $\mathcal{G}$.
\end{compactenum}
In order to accomplish this, the asymptotic results of \cite{a21,a22}, in 
particular, asymptotics of $\overset{e}{\operatorname{Y}}(z)$ and $\overset{
o}{\operatorname{Y}}(z)$ as $z \! \to \! 0$ and $z \! \to \! \infty$, 
where $z\! \in \! \mathbb{C} \setminus \mathbb{R}$, with $n$-dependent 
coefficients, are used.

Given the results of \cite{a21,a22}, carrying through with the above-mentioned 
programme involves lengthy algebraic calculations. A plethora of notations 
remains to be introduced and defined; in particular, basic elements associated 
with the construction of hyperelliptic and finite genus (compact) Riemann 
surfaces and the corresponding Riemann theta functions, and the construction 
of the respective `even' and `odd' equilibrium measures (see Section~2).
\begin{eeeee}
Strictly speaking, this paper is a hybrid continuation of \cite{a21,a22}, 
and is not completely self-contained: concepts, definitions, formulas, and 
theorems proven in \cite{a21,a22} are used. \hfill $\blacksquare$
\end{eeeee}
\begin{eeeee}
The superscripts ${}^{\pm}$, and sometimes the subscripts ${}_{\pm}$, used
throughout the paper should not be confused with the subscripts ${}_{\pm}$
appearing in the various RHPs. \hfill $\blacksquare$
\end{eeeee}
\begin{eeeee}
Although $\overline{\mathbb{C}}$ (or $\mathbb{C} \mathbb{P}^{1})$ $:= \! 
\mathbb{C} \cup \{\infty\}$ (resp., $\overline{\mathbb{R}} \! := \! \mathbb{
R} \cup \{-\infty\} \cup \{+\infty\})$ is the standard definition for the 
(closed) Riemann sphere (resp., closed real line), the simplified, and 
somewhat abusive, notation $\mathbb{C}$ (resp., $\mathbb{R})$ is used to 
denote both the (closed) Riemann sphere, $\overline{\mathbb{C}}$ (resp., 
closed real line, $\overline{\mathbb{R}})$, and the (open) complex field, 
$\mathbb{C}$ (resp., open real line, $\mathbb{R})$, and the context(s) should 
make clear which object(s) the notation $\mathbb{C}$ (resp., $\mathbb{R})$ 
represents. \hfill $\blacksquare$
\end{eeeee}
\section{Preliminary Considerations and Auxiliary Functions}
In this section, auxiliary information which is necessary in order to obtain 
the asymptotic (as $n \! \to \! \infty)$ results presented in Section~3 
is gathered. Towards this end, a terse presentation of all the relevant 
particulars follows (see \cite{a21,a22} for complete details and proofs).

It turns out that, for the asymptotic analysis of this work, the following 
multi-valued functions, and their associated hyperelliptic Riemman surfaces 
of \emph{genus} $N$ $(\in \! \mathbb{N})$, assumed finite, are essential; for 
further details, and proofs, see, for example, \cite{a68,a69}. There are two 
cases to consider: the case of even degree, and the case of odd degree. The 
even degree case is treated first.
\begin{enumerate}
\item[\fbox{$\pmb{\sqrt{\smash[b]{R_{e}(z)}}}$}] Given the set of real points
$-\infty \! < \! a_{0}^{e} \! < \! b_{0}^{e} \! < \! a_{1}^{e} \! < \! b_{1}^{
e} \! < \! \cdots \! < \! a_{N}^{e} \! < \! b_{N}^{e} \! < \! +\infty$, where
$a_{N+1}^{e} \! \equiv \! a_{0}^{e}$ (as points on the Riemann sphere,
$\overline{\mathbb{C}})$, define the multi-valued function $(R_{e}(z))^{1/2}$
as follows:
\begin{equation}
(R_{e}(z))^{1/2} \! := \! \left(\prod_{k=0}^{N} \! \left(z \! - \! b_{k}^{e}
\right) \! \left(z \! - \! a_{k+1}^{e} \right) \right)^{1/2}
\end{equation}
(see \cite{a21}, Figure~1). It must be noted that the set of points $\lbrace
b_{j-1}^{e},a_{j}^{e} \rbrace_{j=1}^{N+1}$, which form the end-points of the
support of the associated `even' equilibrium measure (see below), are not
arbitrary; rather, they satisfy the (real) $n$-dependent and locally solvable
system of $2(N \! + \! 1)$ \emph{moment conditions}~(2.12). The function
$R_{e}(z)$ $(\in \! \mathbb{R}[z]$: the algebra of polynomials in $z$ with
coefficients in $\mathbb{R})$ is a unital polynomial of even degree $(\deg
(R_{e}(z)) \! = \! 2(N \! + \! 1))$ whose simple roots/zeros are $\lbrace
b_{j-1}^{e},a_{j}^{e} \rbrace_{j=1}^{N+1}$.

\hspace*{0.50cm}
Let $\mathcal{Y}_{e} \! := \! \left\lbrace \mathstrut (y,z); \, y^{2} \! = \! 
R_{e}(z) \right\rbrace$ denote the two-sheeted Riemann surface of genus $N$ 
associated with $y^{2} \! = \! R_{e}(z)$: the first sheet of $\mathcal{Y}_{
e}$ is denoted by $\mathcal{Y}_{e}^{+}$, and points on the first sheet are 
represented as $z^{+} \! := \! (z,+(R_{e}(z))^{1/2})$; similarly, the second 
sheet of $\mathcal{Y}_{e}$ is denoted by $\mathcal{Y}_{e}^{-}$, and points on 
the second sheet are represented as $z^{-} \! := \! (z,-(R_{e}(z))^{1/2})$. 
As points on the plane $\mathbb{C}$, $z^{+} \! = \! z^{-} \! := \! z$. The
single-valued branch for the square root of $(R_{e}(z))^{1/2}$ is chosen so 
that $z^{-(N+1)}(R_{e}(z))^{1/2} \! \sim_{\underset{z \in \mathcal{Y}_{e}^{\pm}
}{z \to \infty}} \! \pm 1$. $\mathcal{Y}_{e}$ is realised as a (two-sheeted) 
branched/ramified covering of the Riemann sphere such that its two sheets 
are two identical copies of $\mathbb{C}$ with branch cuts (slits) along the 
intervals $(a_{0}^{e},b_{0}^{e}),(a_{1}^{e},b_{1}^{e}),\dotsc,(a_{N}^{e},
b_{N}^{e})$ and glued together `cross-wise' in the canonical manner along 
$\cup_{j=1}^{N+1}(a_{j-1}^{e},b_{j-1}^{e})$. The \emph{canonical 
$\mathbf{1}$-homology basis} for $\mathcal{Y}_{e}$ are the cycles (closed 
contours) $\boldsymbol{\alpha}_{0}^{e}$ and $\lbrace \boldsymbol{\alpha}_{
j}^{e},\boldsymbol{\beta}_{j}^{e} \rbrace$, $j \! = \! 1,\dotsc,N$, which are 
characterised by the fact that the cycles $\boldsymbol{\alpha}_{j}^{e}$, 
$j \! = \! 0,\dotsc,N$, lie on $\mathcal{Y}_{e}^{+}$, and the cycles 
$\boldsymbol{\beta}_{j}^{e}$, $j \! = \! 1,\dotsc,N$, pass {}from $\mathcal{
Y}_{e}^{+}$ (starting {}from the slit $(a_{j}^{e},b_{j}^{e}))$, through the 
slit $(a_{0}^{e},b_{0}^{e})$ to $\mathcal{Y}_{e}^{-}$, and back again to 
$\mathcal{Y}_{e}^{+}$ through the slit $(a_{j}^{e},b_{j}^{e})$ (see 
\cite{a21}, Figure~3). (Note: in this work, $\mathcal{Y}_{e}^{\pm} \! = \! 
\mathbb{C}_{\pm}.)$

\hspace*{0.50cm}
The canonical $\mathbf{1}$-homology basis $\lbrace \boldsymbol{\alpha}_{j}^{
e},\boldsymbol{\beta}_{j}^{e} \rbrace$, $j \! = \! 1,\dotsc,N$, generates, on 
$\mathcal{Y}_{e}$, the (corresponding) $\boldsymbol{\alpha}^{e}$-normalised 
basis of holomorphic Abelian differentials (one-forms) $\lbrace \omega_{1}^{
e},\omega_{2}^{e},\dotsc,\omega_{N}^{e} \rbrace$, where $\omega_{j}^{e} \! 
:= \! \sum_{k=1}^{N} \tfrac{c_{jk}^{e}z^{N-k}}{\sqrt{\smash[b]{R_{e}(z)}}} \, 
\md z$, $c_{jk}^{e} \! \in \! \mathbb{C}$, $j \! = \! 1,\dotsc,N$, and $\oint_{
\boldsymbol{\alpha}_{k}^{e}} \omega_{j}^{e} \! = \! \delta_{kj}$, $k,j \! = 
\! 1,\dotsc,N$: $\omega_{l}^{e}$, $l \! = \! 1,\dotsc,N$, is real valued on 
$\cup_{j=1}^{N+1}(a_{j-1}^{e},b_{j-1}^{e})$, and has exactly one (real) root 
in any (open) interval $(a_{j-1}^{e},b_{j-1}^{e})$, $j \! = \! 1,\dotsc,N \! 
+ \! 1$; furthermore, in the intervals $(b_{j-1}^{e},a_{j}^{e})$, $j \! = \! 
1,\dotsc,N$, $\omega_{l}^{e}$, $l \! = \! 1,\dotsc,N$, take non-zero, pure 
imaginary values. Let $\boldsymbol{\omega}^{e} \! := \! (\omega_{1}^{e},
\omega_{2}^{e},\dotsc,\omega_{N}^{e})$ denote the basis of holomorphic 
one-forms on $\mathcal{Y}_{e}$ as normalised above with the associated $N 
\! \times \! N$ \emph{Riemann matrix} of $\boldsymbol{\beta}^{e}$-periods, 
$\tau^{e} \! = \! (\tau^{e})_{i,j=1,\dotsc,N} \! := \! (\oint_{\boldsymbol{
\beta}_{j}^{e}} \omega_{i}^{e})_{i,j=1,\dotsc,N}$. The Riemann matrix, 
$\tau^{e}$, is symmetric and pure imaginary, $-\mi \tau^{e}$ is positive 
definite $(\Im (\tau^{e}_{ij}) \! > \! 0)$, and $\det (\tau^{e}) \! \not= \! 
0$. For the holomorphic Abelian differential (one-form) $\boldsymbol{\omega}^{
e}$ defined above, choose $a_{N+1}^{e}$ as the \emph{base point}, and set 
$\boldsymbol{u}^{e} \colon \mathcal{Y}_{e} \! \to \! \operatorname{Jac}
(\mathcal{Y}_{e}),z \! \mapsto \! \boldsymbol{u}^{e}(z) \! := \! \int_{a_{N
+1}^{e}}^{z} \boldsymbol{\omega}^{e}$, where $\operatorname{Jac}(\mathcal{
Y}_{e}) \! := \! \mathbb{C}^{N}/\lbrace N \! + \! \tau^{e}M \rbrace$, $(N,M) 
\! \in \! \mathbb{Z}^{N} \! \times \! \mathbb{Z}^{N}$, is the \emph{Jacobi 
variety} of $\mathcal{Y}_{e}$, and the integration {}from $a_{N+1}^{e}$ to 
$z$ $(\in \mathcal{Y}_{e})$ is taken along any path on $\mathcal{Y}_{e}^{+}$.
\begin{eeeee}
{}From the representation $\omega_{j}^{e} \! = \! \sum_{k=1}^{N} \tfrac{c_{j
k}^{e}z^{N-k}}{\sqrt{\smash[b]{R_{e}(z)}}} \, \md z$, $j \! = \! 1,\dotsc,N$, 
and the normalisation condition $\oint_{\boldsymbol{\alpha}_{k}^{e}} \omega_{
j}^{e} \! = \! \delta_{kj}$, $k,j \! = \! 1,\dotsc,N$, one shows that $c_{j
k}^{e}$, $k,j \! = \! 1,\dotsc,N$, are obtained from
\begin{equation}
\begin{pmatrix}
c_{11}^{e} & c_{12}^{e} & \dotsb & c_{1N}^{e} \\
c_{21}^{e} & c_{22}^{e} & \dotsb & c_{2N}^{e} \\
\vdots     & \vdots     & \ddots & \vdots     \\
c_{N1}^{e} & c_{N2}^{e} & \dotsb & c_{NN}^{e}
\end{pmatrix} \! = \! \widetilde{\mathfrak{S}}_{e}^{-1} \tag{E1},
\end{equation}
where
\begin{equation}
\widetilde{\mathfrak{S}}_{e} \! := \!
\begin{pmatrix}
\oint_{\boldsymbol{\alpha}_{1}^{e}} \frac{\md s_{1}}{\sqrt{\smash[b]{R_{e}
(s_{1})}}} & \oint_{\boldsymbol{\alpha}_{2}^{e}} \frac{\md s_{2}}{\sqrt{
\smash[b]{R_{e}(s_{2})}}} & \dotsb & \oint_{\boldsymbol{\alpha}_{N}^{e}}
\frac{\md s_{N}}{\sqrt{\smash[b]{R_{e}(s_{N})}}} \\
\oint_{\boldsymbol{\alpha}_{1}^{e}} \frac{s_{1} \md s_{1}}{\sqrt{\smash[b]{
R_{e}(s_{1})}}} & \oint_{\boldsymbol{\alpha}_{2}^{e}} \frac{s_{2} \md s_{2}}{
\sqrt{\smash[b]{R_{e}(s_{2})}}} & \dotsb & \oint_{\boldsymbol{\alpha}_{N}^{e}}
\frac{s_{N} \md s_{N}}{\sqrt{\smash[b]{R_{e}(s_{N})}}} \\
\vdots & \vdots & \ddots & \vdots \\
\oint_{\boldsymbol{\alpha}_{1}^{e}} \frac{s_{1}^{N-1} \md s_{1}}{\sqrt{
\smash[b]{R_{e}(s_{1})}}} & \oint_{\boldsymbol{\alpha}_{2}^{e}} \frac{s_{
2}^{N-1} \md s_{2}}{\sqrt{\smash[b]{R_{e}(s_{2})}}} & \dotsb & \oint_{
\boldsymbol{\alpha}_{N}^{e}} \frac{s_{N}^{N-1} \md s_{N}}{\sqrt{\smash[b]{
R_{e}(s_{N})}}}
\end{pmatrix} \tag{E2}.
\end{equation}
For a (representation-independent) proof of the fact that $\det (\widetilde{
\mathfrak{S}}_{e}) \! \not= \! 0$, see, for example, Chapter~10, 
Section~10--2, of \cite{a68}. \hfill $\blacksquare$
\end{eeeee}

\hspace*{0.50cm}
Set \cite{a21}
\begin{equation}
\gamma^{e}(z) \! := \!
\begin{cases}
\left(\mathlarger{\prod_{k=1}^{N+1}}(z \! - \! b_{k-1}^{e})(z \! - \! a_{k}^{
e})^{-1} \right)^{1/4}, &\text{$z \! \in \! \mathbb{C}_{+}$,} \\
-\mi \! \left(\mathlarger{\prod_{k=1}^{N+1}}(z \! - \! b_{k-1}^{e})(z \! - \!
a_{k}^{e})^{-1} \right)^{1/4}, &\text{$z \! \in \! \mathbb{C}_{-}$.}
\end{cases}
\end{equation}
It is shown in Section~4 of \cite{a21} that $\gamma^{e}(z) \! =_{\underset{z
\in \mathcal{Y}_{e}^{\pm}}{z \to \infty}} \! (-\mi)^{(1 \mp 1)/2}(1 \! + \!
\mathcal{O}(z^{-1}))$, and that
\begin{equation}
\left\{z_{j}^{e,\pm} \right\}_{j=1}^{N} \! = \! \left\{\mathstrut z^{\pm} \!
\in \! \mathcal{Y}_{e}^{\pm}; \, (\gamma^{e}(z) \! \mp \! (\gamma^{e}(z))^{-
1}) \vert_{z=z^{\pm}} \! = \! 0 \right\},
\end{equation}
with $z_{j}^{e,\pm} \! \in \! (a_{j}^{e},b_{j}^{e})^{\pm}$ $(\subset \!
\mathcal{Y}_{e}^{\pm})$, $j \! = \! 1,\dotsc,N$, where, as points on the
plane, $z_{j}^{e,+} \! = \! z_{j}^{e,-} \! := \! z_{j}^{e}$, $j \! = \! 1,
\dotsc,N$ (on the plane, $z_{j}^{e} \! \in \! (a_{j}^{e},b_{j}^{e})$, $j \!
= \! 1,\dotsc,N)$.

\hspace*{0.50cm}
Corresponding to $\mathcal{Y}_{e}$, define $\boldsymbol{d}_{e} \! := \! 
-\boldsymbol{K}_{e} \! - \! \sum_{j=1}^{N} \int_{a_{N+1}^{e}}^{z_{j}^{e,-}} 
\boldsymbol{\omega}^{e}$ $(\in \! \mathbb{C}^{N})$, where $\boldsymbol{K}_{
e}$ is the associated `even' \emph{vector of Riemann constants}, and the 
integration {}from $a_{N+1}^{e}$ to $z_{j}^{e,-}$, $j \! = \! 1,\dotsc,N$, is 
taken along a fixed path in $\mathcal{Y}_{e}^{-}$. It is shown in Chapter~VII 
of \cite{a69} that $\boldsymbol{K}_{e} \! = \! \sum_{j=1}^{N} \int_{a_{j}^{e}
}^{a_{N+1}^{e}} \boldsymbol{\omega}^{e}$; furthermore, $\boldsymbol{K}_{e}$ 
is a point of order $2$, that is, $2 \boldsymbol{K}_{e} \! = \! 0$ and $s 
\boldsymbol{K}_{e} \! \not= \! 0$ for $0 \! < \! s \! < \! 2$. Recalling the 
definition of $\boldsymbol{\omega}^{e}$ and that $z^{-(N+1)}(R_{e}(z))^{1/2} 
\! \sim_{\underset{z \in \mathbb{C}_{\pm}}{z \to \infty}} \! \pm 1$, using 
the fact that $\boldsymbol{K}_{e}$ is a point of order $2$, one arrives at
\begin{align}
\boldsymbol{d}_{e} =& \, -\boldsymbol{K}_{e} \! - \! \sum_{j=1}^{N} \int_{a_{N
+1}^{e}}^{z_{j}^{e,-}} \boldsymbol{\omega}^{e} \! = \! \boldsymbol{K}_{e} \! -
\! \sum_{j=1}^{N} \int_{a_{N+1}^{e}}^{z_{j}^{e,-}} \boldsymbol{\omega}^{e} \!
= \! -\boldsymbol{K}_{e} \! + \! \sum_{j=1}^{N} \int_{a_{N+1}^{e}}^{z_{j}^{e,
+}} \boldsymbol{\omega}^{e} \! = \! \boldsymbol{K}_{e} \! + \! \sum_{j=1}^{N}
\int_{a_{N+1}^{e}}^{z_{j}^{e,+}} \boldsymbol{\omega}^{e} \nonumber \\
=& \, -\sum_{j=1}^{N} \int_{a_{j}^{e}}^{z_{j}^{e,-}} \boldsymbol{\omega}^{e}
\! = \! \sum_{j=1}^{N} \int_{a_{j}^{e}}^{z_{j}^{e,+}} \boldsymbol{\omega}^{e}.
\end{align}

\hspace*{0.50cm}
Associated with the Riemann matrix of $\boldsymbol{\beta}^{e}$-periods,
$\tau^{e}$, is the `even' \emph{Riemann theta function}:
\begin{equation}
\boldsymbol{\theta}(z;\tau^{e}) \! =: \! \boldsymbol{\theta}^{e}(z) \! = \!
\sum_{m \in \mathbb{Z}^{N}} \me^{2 \pi \mi (m,z)+\pi \mi (m,\tau^{e}m)}, \quad
z \! \in \! \mathbb{C}^{N},
\end{equation}
where $(\pmb{\cdot},\pmb{\cdot})$ denotes the real Euclidean inner/scalar
product\footnote{For $A \! = \! (A_{1},A_{2},\dotsc,A_{N}) \! \in \! \mathbb{
E}^{N}$ and $B \! = \! (B_{1},B_{2},\dotsc,B_{N}) \! \in \! \mathbb{E}^{N}$,
$(A,B) \! := \! \sum_{k=1}^{N}A_{k}B_{k}$.}; $\boldsymbol{\theta}^{e}(z)$ has
the following evenness and (quasi-) periodicity properties,
\begin{equation*}
\boldsymbol{\theta}^{e}(-z) \! = \! \boldsymbol{\theta}^{e}(z), \qquad
\boldsymbol{\theta}^{e}(z \! + \! e_{j}) \! = \! \boldsymbol{\theta}^{e}(z),
\qquad \mathrm{and} \qquad \boldsymbol{\theta}^{e}(z \! \pm \! \tau_{j}^{e})
\! = \! \me^{\mp 2 \pi \mi z_{j}-\mi \pi \tau_{jj}^{e}} \boldsymbol{\theta}^{e}
(z),
\end{equation*}
where $e_{j}$ is the standard (basis) column vector in $\mathbb{C}^{N}$ with
$1$ in the $j$th entry and $0$ elsewhere, and $\tau_{j}^{e} \! := \! \tau^{e}
e_{j}$ $(\in \! \mathbb{C}^{N})$, $j \! = \! 1,\dotsc,N$.

The odd degree case is now treated.
\item[\fbox{$\pmb{\sqrt{\smash[b]{R_{o}(z)}}}$}] Given the set of real points
$-\infty \! < \! a_{0}^{o} \! < \! b_{0}^{o} \! < \! a_{1}^{o} \! < \! b_{
1}^{o} \! < \! \cdots \! < \! a_{N}^{o} \! < \! b_{N}^{o} \! < \! +\infty$,
where $a_{N+1}^{o} \! \equiv \! a_{0}^{o}$ (as points on the Riemann sphere,
$\overline{\mathbb{C}})$, define the multi-valued function $(R_{o}(z))^{1/2}$
as follows:
\begin{equation}
(R_{o}(z))^{1/2} \! := \! \left(\prod_{k=0}^{N} \! \left(z \! - \! b_{k}^{o}
\right) \! \left(z \! - \! a_{k+1}^{o} \right) \right)^{1/2}
\end{equation}
(see \cite{a22}, Figure~2). It must be noted that the set of points $\lbrace
b_{j-1}^{o},a_{j}^{o} \rbrace_{j=1}^{N+1}$, which form the end-points of the
support of the associated `odd' equilibrium measure (see below), are not
arbitrary; rather, they satisfy the (real) $n$-dependent and locally solvable
system of $2(N \! + \! 1)$ moment conditions~(2.18). The function $R_{o}(z)$
$(\in \! \mathbb{R}[z])$ is a unital polynomial of even degree $(\deg (R_{o}
(z)) \! = \! 2(N \! + \! 1))$ whose simple roots/zeros are $\lbrace b_{j-1}^{
o},a_{j}^{o} \rbrace_{j=1}^{N+1}$.

\hspace*{0.50cm}
Let $\mathcal{Y}_{o} \! := \! \left\lbrace \mathstrut (y,z); \, y^{2} \! = \! 
R_{o} (z) \right\rbrace$ denote the two-sheeted Riemann surface of genus $N$ 
associated with $y^{2} \! = \! R_{o}(z)$: the first sheet of $\mathcal{Y}_{o}$ 
is denoted by $\mathcal{Y}_{o}^{+}$, and points on the first sheet are 
represented as $z^{+} \! := \! (z,+(R_{o}(z))^{1/2})$; similarly, the second 
sheet of $\mathcal{Y}_{o}$ is denoted by $\mathcal{Y}_{o}^{-}$, and points on 
the second sheet are represented as $z^{-} \! := \! (z,-(R_{o}(z))^{1/2})$. 
As points on the plane $\mathbb{C}$, $z^{+} \! = \! z^{-} \! := \! z$. The 
single-valued branch for the square root of $(R_{o}(z))^{1/2}$ is chosen so 
that $z^{-(N+1)}(R_{o}(z))^{1/2} \! \sim_{\underset{z \in \mathcal{Y}_{o}^{\pm}
}{z \to \infty}} \! \pm 1$. $\mathcal{Y}_{o}$ is realised as a (two-sheeted) 
branched/ramified covering of the Riemann sphere such that its two sheets 
are two identical copies of $\mathbb{C}$ with branch cuts (slits) along the 
intervals $(a_{0}^{o},b_{0}^{o}),(a_{1}^{o},b_{1}^{o}),\dotsc,(a_{N}^{o},b_{
N}^{o})$ and glued together `cross-wise' in the canonical manner along $\cup_{
j=1}^{N+1}(a_{j-1}^{o},b_{j-1}^{o})$. The canonical $\mathbf{1}$-homology 
basis for $\mathcal{Y}_{o}$ are the cycles (closed contours) $\boldsymbol{
\alpha}_{0}^{o}$ and $\lbrace \boldsymbol{\alpha}_{j}^{o},\boldsymbol{\beta}_{
j}^{o} \rbrace$, $j \! = \! 1,\dotsc,N$, which are characterised by the fact 
that the cycles $\boldsymbol{\alpha}_{j}^{o}$, $j \! = \! 0,\dotsc,N$, lie on 
$\mathcal{Y}_{o}^{+}$, and the cycles $\boldsymbol{\beta}_{j}^{o}$, $j \! = \! 
1,\dotsc,N$, pass {}from $\mathcal{Y}_{o}^{+}$ (starting {}from the slit $(a_{
j}^{o},b_{j}^{o}))$, through the slit $(a_{0}^{o},b_{0}^{o})$ to $\mathcal{
Y}_{o}^{-}$, and back again to $\mathcal{Y}_{o}^{+}$ through the slit $(a_{
j}^{o},b_{j}^{o})$ (see \cite{a22}, Figure~4). (Note: in this work, $\mathcal{
Y}_{o}^{\pm} \! = \! \mathbb{C}_{\pm}.)$

\hspace*{0.50cm}
The canonical $\mathbf{1}$-homology basis $\lbrace \boldsymbol{\alpha}_{j}^{o},
\boldsymbol{\beta}_{j}^{o} \rbrace$, $j \! = \! 1,\dotsc,N$, generates, on
$\mathcal{Y}_{o}$, the (corresponding) $\boldsymbol{\alpha}^{o}$-normalised
basis of holomorphic Abelian differentials (one-forms) $\lbrace \omega_{1}^{
o},\omega_{2}^{o},\dotsc,\omega_{N}^{o} \rbrace$, where $\omega_{j}^{o} \! :=
\! \sum_{k=1}^{N} \tfrac{c_{jk}^{o}z^{N-k}}{\sqrt{\smash[b]{R_{o}(z)}}} \, \md
z$, $c_{jk}^{o} \! \in \! \mathbb{C}$, $j \! = \! 1,\dotsc,N$, and $\oint_{
\boldsymbol{\alpha}_{k}^{o}} \omega_{j}^{o} \! = \! \delta_{kj}$, $k,j \! =
\! 1,\dotsc,N$: $\omega_{l}^{o}$, $l \! = \! 1,\dotsc,N$, is real valued on
$\cup_{j=1}^{N+1}(a_{j-1}^{o},b_{j-1}^{o})$, and has exactly one (real) root
in any (open) interval $(a_{j-1}^{o},b_{j-1}^{o})$, $j \! = \! 1,\dotsc,N \!
+ \! 1$; furthermore, in the intervals $(b_{j-1}^{o},a_{j}^{o})$, $j \! = \!
1,\dotsc,N$, $\omega_{l}^{o}$, $l \! = \! 1,\dotsc,N$, take non-zero, pure
imaginary values. Let $\boldsymbol{\omega}^{o} \! := \! (\omega_{1}^{o},
\omega_{2}^{o},\dotsc,\omega_{N}^{o})$ denote the basis of holomorphic
one-forms on $\mathcal{Y}_{o}$ as normalised above with the associated $N \!
\times \! N$ Riemann matrix of $\boldsymbol{\beta}^{o}$-periods, $\tau^{o} \!
= \! (\tau^{o})_{i,j=1,\dotsc,N} \! := \! (\oint_{\boldsymbol{\beta}_{j}^{
o}} \omega_{i}^{o})_{i,j=1,\dotsc,N}$. The Riemann matrix, $\tau^{o}$, is
symmetric and pure imaginary, $-\mi \tau^{o}$ is positive definite $(\Im
(\tau^{o}_{ij}) \! > \! 0)$, and $\det (\tau^{o}) \! \not= \! 0$. For the
holomorphic Abelian differential (one-form) $\boldsymbol{\omega}^{o}$ defined
above, choose $a_{N+1}^{o}$ as the base point, and set $\boldsymbol{u}^{o}
\colon \mathcal{Y}_{o} \! \to \! \operatorname{Jac}(\mathcal{Y}_{o}),z \!
\mapsto \! \boldsymbol{u}^{o}(z) \! := \! \int_{a_{N+1}^{o}}^{z} \boldsymbol{
\omega}^{o}$, where $\operatorname{Jac}(\mathcal{Y}_{o}) \! := \! \mathbb{
C}^{N}/\lbrace N \! + \! \tau^{o}M \rbrace$, $(N,M) \! \in \! \mathbb{Z}^{N}
\! \times \! \mathbb{Z}^{N}$, is the Jacobi variety of $\mathcal{Y}_{o}$, and
the integration {}from $a_{N+1}^{o}$ to $z$ $(\in \mathcal{Y}_{o})$ is taken
along any path on $\mathcal{Y}_{o}^{+}$.
\begin{eeeee}
{}From the representation $\omega_{j}^{o} \! = \! \sum_{k=1}^{N} \tfrac{c_{j
k}^{o}z^{N-k}}{\sqrt{\smash[b]{R_{o}(z)}}} \, \md z$, $j \! = \! 1,\dotsc,N$,
and the normalisation condition $\oint_{\boldsymbol{\alpha}_{k}^{o}} \omega_{
j}^{o} \! = \! \delta_{kj}$, $k,j \! = \! 1,\dotsc,N$, one shows that $c_{j
k}^{o}$, $k,j \! = \! 1,\dotsc,N$, are obtained from
\begin{equation}
\begin{pmatrix}
c_{11}^{o} & c_{12}^{o} & \dotsb & c_{1N}^{o} \\
c_{21}^{o} & c_{22}^{o} & \dotsb & c_{2N}^{o} \\
\vdots     & \vdots     & \ddots & \vdots     \\
c_{N1}^{o} & c_{N2}^{o} & \dotsb & c_{NN}^{o}
\end{pmatrix} \! = \! \widetilde{\mathfrak{S}}_{o}^{-1} \tag{O1},
\end{equation}
where
\begin{equation}
\widetilde{\mathfrak{S}}_{o} \! := \!
\begin{pmatrix}
\oint_{\boldsymbol{\alpha}_{1}^{o}} \frac{\md s_{1}}{\sqrt{\smash[b]{R_{o}
(s_{1})}}} & \oint_{\boldsymbol{\alpha}_{2}^{o}} \frac{\md s_{2}}{\sqrt{
\smash[b]{R_{o}(s_{2})}}} & \dotsb & \oint_{\boldsymbol{\alpha}_{N}^{o}}
\frac{\md s_{N}}{\sqrt{\smash[b]{R_{o}(s_{N})}}} \\
\oint_{\boldsymbol{\alpha}_{1}^{o}} \frac{s_{1} \md s_{1}}{\sqrt{\smash[b]{
R_{o}(s_{1})}}} & \oint_{\boldsymbol{\alpha}_{2}^{o}} \frac{s_{2} \md s_{2}}{
\sqrt{\smash[b]{R_{o}(s_{2})}}} & \dotsb & \oint_{\boldsymbol{\alpha}_{N}^{o}}
\frac{s_{N} \md s_{N}}{\sqrt{\smash[b]{R_{o}(s_{N})}}} \\
\vdots & \vdots & \ddots & \vdots \\
\oint_{\boldsymbol{\alpha}_{1}^{o}} \frac{s_{1}^{N-1} \md s_{1}}{\sqrt{
\smash[b]{R_{o}(s_{1})}}} & \oint_{\boldsymbol{\alpha}_{2}^{o}} \frac{s_{
2}^{N-1} \md s_{2}}{\sqrt{\smash[b]{R_{o}(s_{2})}}} & \dotsb & \oint_{
\boldsymbol{\alpha}_{N}^{o}} \frac{s_{N}^{N-1} \md s_{N}}{\sqrt{\smash[b]{
R_{o}(s_{N})}}}
\end{pmatrix} \tag{O2}.
\end{equation}
For a (representation-independent) proof of the fact that $\det (\widetilde{
\mathfrak{S}}_{o}) \! \not= \! 0$, see, for example, Chapter~10,
Section~10--2, of \cite{a68}. \hfill $\blacksquare$
\end{eeeee}

\hspace*{0.50cm}
Set \cite{a22}
\begin{equation}
\gamma^{o}(z) \! := \!
\begin{cases}
\left(\mathlarger{\prod_{k=1}^{N+1}}(z \! - \! b_{k-1}^{o})(z \! - \! a_{k}^{
o})^{-1} \right)^{1/4}, &\text{$z \! \in \! \mathbb{C}_{+}$,} \\
-\mi \! \left(\mathlarger{\prod_{k=1}^{N+1}}(z \! - \! b_{k-1}^{o})(z \! - \!
a_{k}^{o})^{-1} \right)^{1/4}, &\text{$z \! \in \! \mathbb{C}_{-}$.}
\end{cases}
\end{equation}
It is shown in Section~4 of \cite{a22} that $\gamma^{o}(z) \! =_{\underset{z 
\in \mathcal{Y}_{o}^{\pm}}{z \to 0}} \! (-\mi)^{(1 \mp 1)/2} \gamma^{o}(0)(1 
\! + \! \mathcal{O}(z))$, where $\gamma^{o}(0) \! := \! (\prod_{k=1}^{N+1}
b_{k-1}^{o}(a_{k}^{o})^{-1})^{1/4}$ $(> \! 0)$, and that a set of $N$ 
upper-edge and lower-edge finite-length-gap roots/zeros are
\begin{equation}
\left\{z_{j}^{o,\pm} \right\}_{j=1}^{N} \! = \! \left\{\mathstrut z^{\pm} \!
\in \! \mathcal{Y}_{o}^{\pm}; \, ((\gamma^{o}(0))^{-1} \gamma^{o}(z) \! \mp
\! \gamma^{o}(0)(\gamma^{o}(z))^{-1}) \vert_{z=z^{\pm}} \! = \! 0 \right\},
\end{equation}
with $z_{j}^{o,\pm} \! \in \! (a_{j}^{o},b_{j}^{o})^{\pm}$ $(\subset \!
\mathcal{Y}_{o}^{\pm})$, $j \! = \! 1,\dotsc,N$, where, as points on the
plane, $z_{j}^{o,+} \! = \! z_{j}^{o,-} \! := \! z_{j}^{o}$, $j \! = \! 1,
\dotsc,N$ (on the plane, $z_{j}^{o} \! \in \! (a_{j}^{o},b_{j}^{o})$, $j \!
= \! 1,\dotsc,N)$.

\hspace*{0.50cm}
Corresponding to $\mathcal{Y}_{o}$, define $\boldsymbol{d}_{o} \! := \! - 
\boldsymbol{K}_{o} \! - \! \sum_{j=1}^{N} \int_{a_{N+1}^{o}}^{z_{j}^{o,-}} 
\boldsymbol{\omega}^{o}$ $(\in \! \mathbb{C}^{N})$, where $\boldsymbol{K}_{
o}$ is the associated `odd' vector of Riemann constants, and the integration 
{}from $a_{N+1}^{o}$ to $z_{j}^{o,-}$, $j \! = \! 1,\dotsc,N$, is taken along 
a fixed path in $\mathcal{Y}_{o}^{-}$. It is shown in Chapter~VII of 
\cite{a69} that $\boldsymbol{K}_{o} \! = \! \sum_{j=1}^{N} \int_{a_{j}^{o}}^{
a_{N+1}^{o}} \boldsymbol{\omega}^{o}$; furthermore, $\boldsymbol{K}_{o}$ 
is a point of order $2$, that is, $2 \boldsymbol{K}_{o} \! = \! 0$ and $s 
\boldsymbol{K}_{o} \! \not= \! 0$ for $0 \! < \! s \! < \! 2$. Recalling the 
definition of $\boldsymbol{\omega}^{o}$ and that $z^{-(N+1)}(R_{o}(z))^{1/2} 
\! \sim_{\underset{z \in \mathbb{C}_{\pm}}{z \to \infty}} \! \pm 1$, using 
the fact that $\boldsymbol{K}_{o}$ is a point of order $2$, one arrives at
\begin{align}
\boldsymbol{d}_{o} =& \, -\boldsymbol{K}_{o} \! - \! \sum_{j=1}^{N} \int_{a_{N
+1}^{o}}^{z_{j}^{o,-}} \boldsymbol{\omega}^{o} \! = \! \boldsymbol{K}_{o} \! -
\! \sum_{j=1}^{N} \int_{a_{N+1}^{o}}^{z_{j}^{o,-}} \boldsymbol{\omega}^{o} \!
= \! -\boldsymbol{K}_{o} \! + \! \sum_{j=1}^{N} \int_{a_{N+1}^{o}}^{z_{j}^{o,
+}} \boldsymbol{\omega}^{o} \! = \! \boldsymbol{K}_{o} \! + \! \sum_{j=1}^{N}
\int_{a_{N+1}^{o}}^{z_{j}^{o,+}} \boldsymbol{\omega}^{o} \nonumber \\
=& \, -\sum_{j=1}^{N} \int_{a_{j}^{o}}^{z_{j}^{o,-}} \boldsymbol{\omega}^{o}
\! = \! \sum_{j=1}^{N} \int_{a_{j}^{o}}^{z_{j}^{o,+}} \boldsymbol{\omega}^{o}.
\end{align}

\hspace*{0.50cm}
Associated with the Riemann matrix of $\boldsymbol{\beta}^{o}$-periods, 
$\tau^{o}$, is the `odd' Riemann theta function:
\begin{equation}
\boldsymbol{\theta}(z;\tau^{o}) \! =: \! \boldsymbol{\theta}^{o}(z) \! = \!
\sum_{m \in \mathbb{Z}^{N}} \me^{2 \pi \mi (m,z)+\pi \mi (m,\tau^{o}m)},
\quad z \! \in \! \mathbb{C}^{N};
\end{equation}
$\boldsymbol{\theta}^{o}(z)$ has the following evenness and (quasi-) 
periodicity properties,
\begin{equation*}
\boldsymbol{\theta}^{o}(-z) \! = \! \boldsymbol{\theta}^{o}(z), \qquad 
\boldsymbol{\theta}^{o}(z \! + \! e_{j}) \! = \! \boldsymbol{\theta}^{o}(z), 
\qquad \mathrm{and} \qquad \boldsymbol{\theta}^{o}(z \! \pm \! \tau_{j}^{o}) 
\! = \! \me^{\mp 2 \pi \mi z_{j}-\mi \pi \tau_{jj}^{o}} \boldsymbol{\theta}^{
o}(z),
\end{equation*}
where $\tau_{j}^{o} \! := \! \tau^{o}e_{j}$ $(\in \! \mathbb{C}^{N})$, $j \! 
= \! 1,\dotsc,N$.
\end{enumerate}

Some variational theory {}from \cite{a21} and \cite{a22} is now summarised. 
The following discussion is decomposed into two parts: one part corresponding 
to \textbf{RHP1} for $\overset{e}{\operatorname{Y}} \colon \mathbb{C} 
\setminus \mathbb{R} \! \to \! \operatorname{SL}_{2}(\mathbb{C})$, denoted 
as \fbox{$\pmb{\mathrm{P}_{1}}$}; and the other part corresponding to 
\textbf{RHP2} for $\overset{o}{\operatorname{Y}} \colon \mathbb{C} 
\setminus \mathbb{R} \! \to \! \operatorname{SL}_{2}(\mathbb{C})$, denoted as 
\fbox{$\pmb{\mathrm{P}_{2}}$}.
\begin{enumerate}
\item[\fbox{$\pmb{\mathrm{P}_{1}}$}] Let $\widetilde{V} \colon \mathbb{R}
\setminus \{0\} \! \to \! \mathbb{R}$ satisfy conditions~(1.20)--(1.22). Let
$\mathrm{I}_{V}^{e}[\mu^{e}] \colon \mathcal{M}_{1}(\mathbb{R}) \! \to \!
\mathbb{R}$ denote the functional
\begin{equation}
\mathrm{I}_{V}^{e}[\mu^{e}] \! = \! \iint_{\mathbb{R}^{2}} \ln \! \left(
\dfrac{\lvert st \rvert}{\lvert s \! - \! t \rvert^{2}} \right) \md \mu^{e}(s)
\, \md \mu^{e}(t) \! + \! 2 \int_{\mathbb{R}} \widetilde{V}(s) \, \md \mu^{e}
(s),
\end{equation}
and consider the associated minimisation problem,
\begin{equation*}
E_{V}^{e} \! = \! \inf \! \left\lbrace \mathstrut \mathrm{I}_{V}^{e}[\mu^{e}]; 
\, \mu^{e} \! \in \! \mathcal{M}_{1}(\mathbb{R}) \right\rbrace.
\end{equation*}
In \cite{a21}, it is proven that the infimum is finite, and there is a unique 
measure $\mu_{V}^{e}$ $(\in \! \mathcal{M}_{1}(\mathbb{R}))$ achieving the 
infimum. The measure $\mu_{V}^{e}$ is referred to as the `even' equilibrium 
measure. Furthermore, the following `regularity' (see below) properties of 
$\mu_{V}^{e}$ are also established in \cite{a21}:
\begin{itemize}
\item the `even' equilibrium measure has compact support which consists of the 
disjoint union of a finite number of bounded real intervals; in fact, as shown 
in \cite{a21}, $\mathrm{supp}(\mu_{V}^{e}) \! =: \! J_{e} \! = \! \cup_{j=1}^{
N+1}(b_{j-1}^{e},a_{j}^{e})$ $(\subset \! \mathbb{R} \setminus \{0\})$, where 
$\lbrace b_{j-1}^{e},a_{j}^{e} \rbrace_{j=1}^{N+1}$, as described in item 
\fbox{$\pmb{\sqrt{\smash[b]{R_{e}(z)}}}$} above, with $b_{0}^{e} \! := \! 
\min \lbrace \mathrm{supp}(\mu_{V}^{e}) \rbrace \! \notin \! \lbrace -\infty,
0 \rbrace$, $a_{N+1}^{e} \! := \! \max \lbrace \mathrm{supp}(\mu_{V}^{e}) 
\rbrace \! \notin \! \lbrace 0,+\infty \rbrace$, and $-\infty \! < \! b_{0}^{
e} \! < \! a_{1}^{e} \! < \! b_{1}^{e} \! < \! a_{2}^{e} \! < \! \dotsb \! < 
\! b_{N}^{e} \! < \! a_{N+1}^{e} \! < \! +\infty$, constitute the end-points 
of the support of the `even' equilibrium measure;
\item the end-points $\lbrace b_{j-1}^{e},a_{j}^{e} \rbrace_{j=1}^{N+1}$
satisfy the following $n$-dependent and (locally) solvable system of $2(N 
\! + \! 1)$ real moment conditions \cite{a21}:
\begin{gather}
\begin{split}
\int_{J_{e}} \dfrac{(2s^{-1} \! + \! \widetilde{V}^{\prime}(s))s^{j}}{(R_{e}
(s))^{1/2}_{+}} \, \md s \! = \! 0, \quad j \! = \! 0,\dotsc,N, \qquad \qquad
\int_{J_{e}} \dfrac{(2s^{-1} \! + \! \widetilde{V}^{\prime}(s))s^{N+1}}{(R_{e}
(s))^{1/2}_{+}} \, \md s \! = \! -4 \pi \mi, \\
\int_{a_{j}^{e}}^{b_{j}^{e}} \! \left(\dfrac{\mi (R_{e}(s))^{1/2}}{2 \pi}
\int_{J_{e}} \dfrac{(2 \xi^{-1} \! + \! \widetilde{V}^{\prime}(\xi))}{(R_{e}
(\xi))^{1/2}_{+}(\xi \! - \! s)} \, \md \xi \right) \! \md s \! = \! \ln \!
\left\vert \dfrac{a_{j}^{e}}{b_{j}^{e}} \right\vert \! + \! \dfrac{1}{2} \!
\left(\widetilde{V}(a_{j}^{e}) \! - \! \widetilde{V}(b_{j}^{e}) \right), 
\quad j \! = \! 1,\dotsc,N,
\end{split}
\end{gather}
where ${}^{\pmb{\prime}}$ denotes differentiation with respect to the
argument,
\begin{equation*}
(R_{e}(z))^{1/2} \! := \! \left(\prod_{k=1}^{N+1}(z \! - \! b_{k-1}^{e})(z \!
- \! a_{k}^{e}) \right)^{1/2},
\end{equation*}
$(R_{e}(x))^{1/2}_{\pm} \! := \! \lim_{\varepsilon \downarrow 0}(R_{e}(x \!
\pm \! \mi \varepsilon))^{1/2}$, and the branch of the square root is chosen
so that $z^{-(N+1)}(R_{e}(z))^{1/2} \! \sim_{\underset{z \in \mathbb{C}_{\pm}
}{z \to \infty}} \! \pm 1$;
\item the `even' equilibrium measure is absolutely continuous with respect to
Lebesgue measure. The \emph{density} is given by
\begin{equation}
\md \mu_{V}^{e}(x) \! := \! \psi_{V}^{e}(x) \, \md x \! = \! \dfrac{1}{2 \pi
\mi}(R_{e}(x))^{1/2}_{+}h_{V}^{e}(x) \pmb{1}_{J_{e}}(x) \, \md x,
\end{equation}
where
\begin{equation}
h_{V}^{e}(z) \! = \! \dfrac{1}{2} \oint_{C_{\mathrm{R}}^{e}} \dfrac{(\frac{
\mi}{\pi s} \! + \! \frac{\mi \widetilde{V}^{\prime}(s)}{2 \pi})}{\sqrt{
\smash[b]{R_{e}(s)}} \, (s \! - \! z)} \, \md s
\end{equation}
(real analytic for $z \! \in \! \mathbb{R} \setminus \{0\})$, $C_{\mathrm{R}
}^{e}$ $(\subset \mathbb{C}^{\ast})$ is the union of two circular contours,
one outer one of large radius $R^{\natural}$ traversed clockwise and one inner
one of small radius $r^{\natural}$ traversed counter-clockwise, with the
numbers $0 \! < \! r^{\natural} \! < \! R^{\natural} \! < \! +\infty$ chosen
such that, for (any) non-real $z$ in the domain of analyticity of $\widetilde{
V}$ (that is, $\mathbb{C}^{\ast})$, $\mathrm{int}(C_{\mathrm{R}}^{e}) \!
\supset \! J_{e} \cup \{z\}$, and $\pmb{1}_{J_{e}}(x)$ denotes the indicator
(characteristic) function of the set $J_{e}$. (Note that $\psi_{V}^{e}(x) \!
\geqslant \! 0 \, \, \forall \, \, x \! \in \! \overline{J_{e}} := \! \cup_{j
=1}^{N+1}[b_{j-1}^{e},a_{j}^{e}]$: it vanishes like a square root at the
end-points of the support of the `even' equilibrium measure, that is, $\psi_{
V}^{e}(s) \! =_{s \downarrow b_{j-1}^{e}} \! \mathcal{O}((s \! - \! b_{j-1}^{
e})^{1/2})$ and $\psi_{V}^{e}(s) \! =_{s \uparrow a_{j}^{e}} \! \mathcal{O}
((a_{j}^{e} \! - \! s)^{1/2})$, $j \! = \! 1,\dotsc,N \! + \! 1.)$;
\item the `even' equilibrium measure and its (compact) support are (uniquely)
characterised by the following Euler-Lagrange variational equations: there
exists $\ell_{e} \! \in \! \mathbb{R}$, the `even' Lagrange multiplier, and
$\mu^{e} \! \in \! \mathcal{M}_{1}(\mathbb{R})$ such that
\begin{gather*}
4 \int_{J_{e}} \ln (\vert x \! - \! s \vert) \, \md \mu^{e} (s) \! - \! 2 \ln
\vert x \vert \! - \! \widetilde{V}(x) \! - \! \ell_{e} \! = \! 0, \quad x \!
\in \! \overline{J_{e}}, \tag{$\pmb{\mathrm{P}_{1}^{(a)}}$} \\
4 \int_{J_{e}} \ln (\vert x \! - \! s \vert) \, \md \mu^{e} (s) \! - \! 2 \ln
\vert x \vert \! - \! \widetilde{V}(x) \! - \! \ell_{e} \! \leqslant \! 0,
\quad x \! \in \! \mathbb{R} \setminus \overline{J_{e}};
\tag{$\pmb{\mathrm{P}_{1}^{(b)}}$}
\end{gather*}
\item the Euler-Lagrange variational equations can be conveniently recast in
terms of the complex potential $g^{e}(z)$ of $\mu_{V}^{e}$:
\begin{equation}
g^{e}(z) \! := \! \int_{J_{e}} \! \ln \! \left((z \! - \! s)^{2}(zs)^{-1}
\right) \md \mu_{V}^{e}(s), \quad z \! \in \! \mathbb{C} \setminus (-\infty,
\max \lbrace 0,a_{N+1}^{e} \rbrace).
\end{equation}
The function $g^{e} \colon \mathbb{C} \setminus (-\infty,\max \lbrace 0,
a_{N+1}^{e} \rbrace) \! \to \! \mathbb{C}$ so defined satisfies:
\begin{enumerate}
\item[$\pmb{(\mathrm{P}_{1}^{(1)})}$] $g^{e}(z)$ is analytic for $z \! \in \!
\mathbb{C} \setminus (-\infty,\max \lbrace 0,a_{N+1}^{e} \rbrace)$;
\item[$\pmb{(\mathrm{P}_{1}^{(2)})}$] $g^{e}(z) \! =_{\underset{z \in \mathbb{
C} \setminus \mathbb{R}}{z \to \infty}} \! \ln (z) \! + \! \mathcal{O}(1)$;
\item[$\pmb{(\mathrm{P}_{1}^{(3)})}$] $g^{e}_{+}(z) \! + \! g^{e}_{-}(z) \! -
\! \widetilde{V}(z) \! - \! \ell_{e} \! + \! 2Q_{e} \! = \! 0$, $z \! \in \!
\overline{J_{e}}$, where $g^{e}_{\pm}(z) \! := \! \lim_{\varepsilon \downarrow
0}g^{e}(z \! \pm \! \mi \varepsilon)$, and
\begin{equation}
Q_{e} \! := \! \int_{J_{e}} \ln (s) \, \md \mu_{V}^{e}(s) \! = \! \int_{J_{e}}
\ln (\lvert s \rvert) \, \md \mu_{V}^{e}(s) \! + \! \mi \pi \int_{J_{e} \cap
\mathbb{R}_{-}} \md \mu_{V}^{e}(s);
\end{equation}
\item[$\pmb{(\mathrm{P}_{1}^{(4)})}$] $g^{e}_{+}(z) \! + \! g^{e}_{-}(z) \! -
\! \widetilde{V}(z) \! - \! \ell_{e} \! + \! 2Q_{e} \! \leqslant \! 0$, $z \!
\in \! \mathbb{R} \setminus \overline{J_{e}}$, where equality holds for at
most a finite number of points;
\item[$\pmb{(\mathrm{P}_{1}^{(5)})}$] $g^{e}_{+}(z) \! - \! g^{e}_{-}(z) \! =
\! \mi f_{g^{e}}^{\mathbb{R}}(z)$, $z \! \in \! \mathbb{R}$, where $f_{g^{e}}^{
\mathbb{R}} \colon \mathbb{R} \! \to \! \mathbb{R}$, and, in particular, $g^{
e}_{+}(z) \! - \! g^{e}_{-}(z) \! = \! \mi \operatorname{const.}$, $z \! \in
\! \mathbb{R} \setminus \overline{J_{e}}$, with $\operatorname{const.} \! \in
\! \mathbb{R}$;
\item[$\pmb{(\mathrm{P}_{1}^{(6)})}$] $\mi (g^{e}_{+}(z) \! - \! g^{e}_{-}
(z))^{\prime} \! \geqslant \! 0$, $z \! \in \! J_{e}$, where equality holds 
for at most a finite number of points.
\end{enumerate}
\end{itemize}
\end{enumerate}
\begin{enumerate}
\item[\fbox{$\pmb{\mathrm{P}_{2}}$}] Let $\widetilde{V} \colon \mathbb{R}
\setminus \{0\} \! \to \! \mathbb{R}$ satisfy conditions~(1.20)--(1.22). Let
$\mathrm{I}_{V}^{o}[\mu^{o}] \colon \mathcal{M}_{1}(\mathbb{R}) \! \to \!
\mathbb{R}$ denote the functional
\begin{equation}
\mathrm{I}_{V}^{o}[\mu^{o}] \! = \! \iint_{\mathbb{R}^{2}} \ln \! \left(\dfrac{
\lvert st \rvert}{\lvert s \! - \! t \rvert^{2+\frac{1}{n}}} \right) \md \mu^{
o}(s) \, \md \mu^{o}(t) \! + \! 2 \int_{\mathbb{R}} \widetilde{V}(s) \, \md
\mu^{o}(s), \quad n \! \in \! \mathbb{N},
\end{equation}
and consider the associated minimisation problem,
\begin{equation*}
E_{V}^{o} \! = \! \inf \! \left\lbrace \mathstrut \mathrm{I}_{V}^{o}[\mu^{o}]; 
\, \mu^{o} \! \in \! \mathcal{M}_{1}(\mathbb{R}) \right\rbrace.
\end{equation*}
In \cite{a22}, it is proven that the infimum is finite, and there is a unique 
measure $\mu_{V}^{o}$ $(\in \! \mathcal{M}_{1}(\mathbb{R}))$ achieving the 
infimum. The measure $\mu_{V}^{o}$ is referred to as the `odd' equilibrium 
measure. Furthermore, the following `regularity' (see below) properties of 
$\mu_{V}^{o}$ are also established in \cite{a22}:
\begin{itemize}
\item the `odd' equilibrium measure has compact support which consists of the 
disjoint union of a finite number of bounded real intervals; in fact, as shown 
in \cite{a22}, $\mathrm{supp}(\mu_{V}^{o}) \! =: \! J_{o} \! = \! \cup_{j=1}^{
N+1}(b_{j-1}^{o},a_{j}^{o})$ $(\subset \! \mathbb{R} \setminus \{0\})$, where 
$\lbrace b_{j-1}^{o},a_{j}^{o} \rbrace_{j=1}^{N+1}$, as described in item 
\fbox{$\pmb{\sqrt{\smash[b]{R_{o}(z)}}}$} above, with $b_{0}^{o} \! := \! 
\min \lbrace \mathrm{supp}(\mu_{V}^{o}) \rbrace \! \notin \! \lbrace -\infty,
0 \rbrace$, $a_{N+1}^{o} \! := \! \max \lbrace \mathrm{supp}(\mu_{V}^{o}) 
\rbrace \! \notin \! \lbrace 0,+\infty \rbrace$, and $-\infty \! < \! b_{0}^{
o} \! < \! a_{1}^{o} \! < \! b_{1}^{o} \! < \! a_{2}^{o} \! < \! \dotsb \! < 
\! b_{N}^{o} \! < \! a_{N+1}^{o} \! < \! +\infty$, constitute the end-points 
of the support of the `odd' equilibrium measure; (The number of intervals, $N 
\! + \! 1$, is the same in the `odd' case as in the `even' case, which can be 
established by a lengthy analysis similar to that contained in \cite{a38}.) 
\item the end-points $\lbrace b_{j-1}^{o},a_{j}^{o} \rbrace_{j=1}^{N+1}$
satisfy the following $n$-dependent and (locally) solvable system of $2(N 
\! + \! 1)$ real moment conditions \cite{a22}:
\begin{gather}
\begin{split}
\int_{J_{o}} \dfrac{(2s^{-1} \! + \! \widetilde{V}^{\prime}(s))s^{j}}{(R_{o}
(s))^{1/2}_{+}} \, \md s \! = \! 0, \quad j \! = \! 0,\dotsc,N, \qquad \qquad
\int_{J_{o}} \dfrac{(2s^{-1} \! + \! \widetilde{V}^{\prime}(s))s^{N+1}}{(R_{o}
(s))^{1/2}_{+}} \, \dfrac{\md s}{2 \pi \mi} \! = \! -\left(2 \! + \! \dfrac{
1}{n} \right), \\
\int_{a_{j}^{o}}^{b_{j}^{o}} \! \left(\dfrac{\mi (R_{o}(s))^{1/2}}{2 \pi}
\int_{J_{o}} \dfrac{(2 \xi^{-1} \! + \! \widetilde{V}^{\prime}(\xi))}{(R_{o}
(\xi))^{1/2}_{+}(\xi \! - \! s)} \, \md \xi \right) \! \md s \! = \! \ln \!
\left\vert \dfrac{a_{j}^{o}}{b_{j}^{o}} \right\vert \! + \! \dfrac{1}{2} \!
\left(\widetilde{V}(a_{j}^{o}) \! - \! \widetilde{V}(b_{j}^{o}) \right), \quad
j \! = \! 1,\dotsc,N,
\end{split}
\end{gather}
where
\begin{equation*}
(R_{o}(z))^{1/2} \! := \! \left(\prod_{k=1}^{N+1}(z \! - \! b_{k-1}^{o})(z \!
- \! a_{k}^{o}) \right)^{1/2},
\end{equation*}
$(R_{o}(x))^{1/2}_{\pm} \! := \! \lim_{\varepsilon \downarrow 0}(R_{o}(x \!
\pm \! \mi \varepsilon))^{1/2}$, and the branch of the square root is chosen
so that $z^{-(N+1)}(R_{o}(z))^{1/2} \! \sim_{\underset{z \in \mathbb{C}_{\pm}
}{z \to \infty}} \! \pm 1$;
\item the `odd' equilibrium measure is absolutely continuous with respect to
Lebesgue measure. The density is given by
\begin{equation}
\md \mu_{V}^{o}(x) \! := \! \psi_{V}^{o}(x) \, \md x \! = \! \dfrac{1}{2 \pi
\mi}(R_{o}(x))^{1/2}_{+}h_{V}^{o}(x) \pmb{1}_{J_{o}}(x) \, \md x,
\end{equation}
where
\begin{equation}
h_{V}^{o}(z) \! = \! \left(2 \! + \! \dfrac{1}{n} \right)^{-1} \oint_{C_{
\mathrm{R}}^{o}} \dfrac{(\frac{\mi}{\pi s} \! + \! \frac{\mi \widetilde{V}^{
\prime}(s)}{2 \pi})}{\sqrt{\smash[b]{R_{o}(s)}} \, (s \! - \! z)} \, \md s
\end{equation}
(real analytic for $z \! \in \! \mathbb{R} \setminus \{0\})$, $C_{\mathrm{R}
}^{o}$ $(\subset \mathbb{C}^{\ast})$ is the union of two circular contours,
one outer one of large radius $R^{\flat}$ traversed clockwise and one inner
one of small radius $r^{\flat}$ traversed counter-clockwise, with the numbers
$0 \! < \! r^{\flat} \! < \! R^{\flat} \! < \! +\infty$ chosen such that, for
(any) non-real $z$ in the domain of analyticity of $\widetilde{V}$ (that is,
$\mathbb{C}^{\ast})$, $\mathrm{int}(C_{\mathrm{R}}^{o}) \! \supset \! J_{o}
\cup \{z\}$, and $\pmb{1}_{J_{o}}(x)$ denotes the indicator (characteristic)
function of the set $J_{o}$. (Note that $\psi_{V}^{o}(x) \! \geqslant \! 0 \,
\, \forall \, \, x \! \in \! \overline{J_{o}} := \! \cup_{j=1}^{N+1}[b_{j-1}^{
o},a_{j}^{o}]$: it vanishes like a square root at the end-points of the
support of the `odd' equilibrium measure, that is, $\psi_{V}^{o}(s) \! =_{s
\downarrow b_{j-1}^{o}} \! \mathcal{O}((s \! - \! b_{j-1}^{o})^{1/2})$ and
$\psi_{V}^{o}(s) \! =_{s \uparrow a_{j}^{o}} \! \mathcal{O}((a_{j}^{o} \! -
\! s)^{1/2})$, $j \! = \! 1,\dotsc,N \! + \! 1.)$;
\item the `odd' equilibrium measure and its (compact) support are (uniquely)
characterised by the following Euler-Lagrange variational equations: there
exists $\ell_{o} \! \in \! \mathbb{R}$, the `odd' Lagrange multiplier, and
$\mu^{o} \! \in \! \mathcal{M}_{1}(\mathbb{R})$ such that
\begin{gather*}
2 \! \left(2 \! + \! \dfrac{1}{n} \right) \! \int_{J_{o}} \ln (\vert x \! - \!
s \vert) \, \md \mu^{o} (s) \! - \! 2 \ln \vert x \vert \! - \! \widetilde{V}
(x) \! - \! \ell_{o} \! - \! 2 \! \left(2 \! + \! \dfrac{1}{n} \right) \!
\widetilde{Q}_{o} = \! 0, \quad x \! \in \! \overline{J_{o}},
\tag{$\pmb{\mathrm{P}_{2}^{(a)}}$} \\
2 \! \left(2 \! + \! \dfrac{1}{n} \right) \! \int_{J_{o}} \ln (\vert x \! - \!
s \vert) \, \md \mu^{o} (s) \! - \! 2 \ln \vert x \vert \! - \! \widetilde{V}
(x) \! - \! \ell_{o} \! - \! 2 \! \left(2 \! + \! \dfrac{1}{n} \right) \!
\widetilde{Q}_{o} \! \leqslant \! 0, \quad x \! \in \! \mathbb{R} \setminus
\overline{J_{o}},
\tag{$\pmb{\mathrm{P}_{2}^{(b)}}$}
\end{gather*}
where $\widetilde{Q}_{o} \! := \! \int_{J_{o}} \ln (\lvert s \rvert) \, \md
\mu^{o}(s)$;
\item the Euler-Lagrange variational equations can be conveniently recast in
terms of the complex potential $g^{o}(z)$ of $\mu_{V}^{o}$:
\begin{equation}
g^{o}(z) \! := \! \int_{J_{o}} \ln \! \left((z \! - \! s)^{2+\frac{1}{n}}
(zs)^{-1} \right) \md \mu_{V}^{o}(s), \quad z \! \in \! \mathbb{C} \setminus
(-\infty,\max \lbrace 0,a_{N+1}^{o} \rbrace).
\end{equation}
The function $g^{o} \colon \mathbb{C} \setminus (-\infty,\max \lbrace 0,
a_{N+1}^{o} \rbrace) \! \to \! \mathbb{C}$ so defined satisfies:
\begin{enumerate}
\item[$\pmb{(\mathrm{P}_{2}^{(1)})}$] $g^{o}(z)$ is analytic for $z \! \in \!
\mathbb{C} \setminus (-\infty,\max \lbrace 0,a_{N+1}^{o} \rbrace)$;
\item[$\pmb{(\mathrm{P}_{2}^{(2)})}$] $g^{o}(z) \! =_{\underset{z \in \mathbb{
C} \setminus \mathbb{R}}{z \to 0}} \! -\ln (z) \! + \! \mathcal{O}(1)$;
\item[$\pmb{(\mathrm{P}_{2}^{(3)})}$] $g^{o}_{+}(z) \! + \! g^{o}_{-}(z) \! -
\! \widetilde{V}(z) \! - \! \ell_{o} \! - \! \mathfrak{Q}^{+}_{\mathscr{A}} \!
- \! \mathfrak{Q}^{-}_{\mathscr{A}} \! = \! 0$, $z \! \in \! \overline{J_{o}}
$, where $g^{o}_{\pm}(z) \! := \! \lim_{\varepsilon \downarrow 0}g^{o}(z \!
\pm \! \mi \varepsilon)$, and
\begin{equation}
\mathfrak{Q}^{\pm}_{\mathscr{A}} \! := \! \left(1 \! + \! \dfrac{1}{n} \right)
\! \int_{J_{o}} \ln (\lvert s \rvert) \, \md \mu_{V}^{o}(s) \! - \! \mi \pi
\int_{J_{o} \cap \mathbb{R}_{-}} \md \mu_{V}^{o}(s) \! \pm \! \mi \pi \! \left(
2 \! + \! \dfrac{1}{n} \right) \! \int_{J_{o} \cap \mathbb{R}_{+}} \md \mu_{
V}^{o}(s);
\end{equation}
\item[$\pmb{(\mathrm{P}_{2}^{(4)})}$] $g^{o}_{+}(z) \! + \! g^{o}_{-}(z) \! -
\! \widetilde{V}(z) \! - \! \ell_{o} \! - \! \mathfrak{Q}^{+}_{\mathscr{A}} \!
- \! \mathfrak{Q}^{-}_{\mathscr{A}} \! \leqslant \! 0$, $z \! \in \! \mathbb{
R} \setminus \overline{J_{o}}$, where equality holds for at most a finite
number of points;
\item[$\pmb{(\mathrm{P}_{2}^{(5)})}$] $g^{o}_{+}(z) \! - \! g^{o}_{-}(z) \! -
\! \mathfrak{Q}^{+}_{\mathscr{A}} \! + \! \mathfrak{Q}^{-}_{\mathscr{A}} \! =
\! \mi f_{g^{o}}^{\mathbb{R}}(z)$, $z \! \in \! \mathbb{R}$, where $f_{g^{o}}^{
\mathbb{R}} \colon \mathbb{R} \! \to \! \mathbb{R}$, and, in particular, $g^{
o}_{+}(z) \! - \! g^{o}_{-}(z) \! - \! \mathfrak{Q}^{+}_{\mathscr{A}} \! + \!
\mathfrak{Q}^{-}_{\mathscr{A}} \! = \! \mi \operatorname{const.}$, $z \! \in
\! \mathbb{R} \setminus \overline{J_{o}}$, with $\operatorname{const.} \! \in
\! \mathbb{R}$;
\item[$\pmb{(\mathrm{P}_{2}^{(6)})}$] $\mi (g^{o}_{+}(z) \! - \! g^{o}_{-}(z)
\! - \! \mathfrak{Q}^{+}_{\mathscr{A}} \! + \! \mathfrak{Q}^{-}_{\mathscr{A}}
)^{\prime} \! \geqslant \! 0$, $z \! \in \! J_{o}$, where equality holds for
at most a finite number of points.
\end{enumerate}
\end{itemize}
\end{enumerate}

In this work, as in \cite{a21,a22}, the so-called `regular case' is 
considered, namely:
\begin{itemize}
\item $\md \mu_{V}^{e}$, or $\widetilde{V} \colon \mathbb{R} \setminus \{0\}
\! \to \! \mathbb{R}$ satisfying conditions~(1.20)--(1.22), is \emph{regular}
if
\begin{gather*}
h_{V}^{e}(x) \! \not\equiv \! 0 \quad \text{on} \quad \overline{J_{e}}, \\
4 \int_{J_{e}} \ln (\vert x \! - \! s \vert) \, \md \mu_{V}^{e}(s) \! - \! 2
\ln \vert x \vert \! - \! \widetilde{V}(x) \! - \! \ell_{e} \! < \! 0, \quad 
x \! \in \! \mathbb{R} \setminus \overline{J_{e}},
\end{gather*}
and inequalities~$\pmb{(\mathrm{P}_{1}^{(4)})}$
and~$\pmb{(\mathrm{P}_{1}^{(6)})}$ in \fbox{$\pmb{\mathrm{P}_{1}}$} are
strict, that is, $\leqslant$ (resp., $\geqslant)$ is replaced by $<$ (resp.,
$>)$;
\item $\md \mu_{V}^{o}$, or $\widetilde{V} \colon \mathbb{R} \setminus \{0\}
\! \to \! \mathbb{R}$ satisfying conditions~(1.20)--(1.22), is regular if
\begin{gather*}
h_{V}^{o}(x) \! \not\equiv \! 0 \quad \text{on} \quad \overline{J_{o}}, \\
2 \! \left(2 \! + \! \dfrac{1}{n} \right) \! \int_{J_{o}} \ln (\vert x \! - \!
s \vert) \, \md \mu_{V}^{o}(s) \! - \! 2 \ln \vert x \vert \! - \! \widetilde{
V}(x) \! - \! \ell_{o} \! - \! 2 \! \left(2 \! + \! \dfrac{1}{n} \right) \!
Q_{o} \! < \! 0, \quad x \! \in \! \mathbb{R} \setminus \overline{J_{o}},
\end{gather*}
where $Q_{o} \! := \! \int_{J_{o}} \ln (\lvert s \rvert) \, \md \mu_{V}^{o}
(s)$, and inequalities~$\pmb{(\mathrm{P}_{2}^{(4)})}$ 
and~$\pmb{(\mathrm{P}_{2}^{(6)})}$ in \fbox{$\pmb{\mathrm{P}_{2}}$} are 
strict, that is, $\leqslant$ (resp., $\geqslant)$ is replaced by $<$ (resp., 
$>)$\footnote{There are three distinct situations in which these conditions
may fail: (i) for at least one $\widetilde{x}_{e} \! \in \! \mathbb{R}
\setminus \overline{J_{e}}$ (resp., $\widetilde{x}_{o} \! \in \! \mathbb{R}
\setminus \widetilde{J_{o}})$, $4 \int_{J_{e}} \ln (\lvert \widetilde{x}_{e}
\! - \! s \rvert) \, \md \mu_{V}^{e}(s) \! - \! 2 \ln \lvert \widetilde{x}_{e}
\rvert \! - \! \widetilde{V}(\widetilde{x}_{e}) \! - \! \ell_{e} \! = \! 0$
(resp., $2(2 \! + \! \tfrac{1}{n}) \int_{J_{o}} \ln (\lvert \widetilde{x}_{o}
\! - \! s \rvert) \, \md \mu_{V}^{o}(s) \! - \! 2 \ln \lvert \widetilde{x}_{o}
\rvert \! - \! \widetilde{V}(\widetilde{x}_{o}) \! - \! \ell_{o} \! - \! 2(2
\! + \! \tfrac{1}{n})Q_{o} \! = \! 0)$, that is, for $n$ even (resp., $n$
odd) equality is attained for at least one point $\widetilde{x}_{e}$ (resp.,
$\widetilde{x}_{o})$ in the complement of the closure of the support of the
`even' (resp., `odd') equilibrium measure $\mu_{V}^{e}$ (resp., $\mu_{V}^{o}
)$, which corresponds to the situation in which a `band' has just closed, or
is about to open, about $\widetilde{x}_{e}$ (resp., $\widetilde{x}_{o})$; (ii)
for at least one $\widehat{x}_{e}$ (resp., $\widehat{x}_{o})$, $h_{V}^{e}
(\widehat{x}_{e}) \! = \! 0$ (resp., $h_{V}^{o}(\widehat{x}_{o}) \! = \! 0)$,
that is, for $n$ even (resp., $n$ odd) the function $h_{V}^{e}$ (resp., $h_{
V}^{o})$ vanishes for at least one point $\widehat{x}_{e}$ (resp., $\widehat{
x}_{o})$ within the support of the `even' (resp., `odd') equilibrium measure
$\mu_{V}^{e}$ (resp., $\mu_{V}^{o})$, which corresponds to the situation in
which a `gap' is about to open, or close, about $\widehat{x}_{e}$ (resp.,
$\widehat{x}_{o})$; and (iii) there exists at least one $j \! \in \! \lbrace
1,\dotsc,N \! + \! 1 \rbrace$, denoted $j_{e}$ (resp., $j_{o})$, such that
$h_{V}^{e}(b_{j_{e}-1}^{e}) \! = \! 0$ and/or $h_{V}^{e}(a_{j_{e}}^{e}) \!
= \! 0$ (resp., $h_{V}^{o}(b_{j_{o}-1}^{o}) \! = \! 0$ and/or $h_{V}^{o}(a_{
j_{o}}^{o}) \! = \! 0)$. Each of these three cases can occur only a finite
number of times due to the fact that $\widetilde{V} \colon \mathbb{R}
\setminus \lbrace 0 \rbrace \! \to \! \mathbb{R}$ satisfies
conditions~(1.20)--(1.22) \cite{a38,a40}.}.
\end{itemize}
The (density of the) `even' and `odd' equilibrium measures, $\md \mu_{V}^{e}$ 
and $\md \mu_{V}^{o}$, respectively, together with the variational problems, 
emerge naturally in the asymptotic analysis of \pmb{RHP1} and \pmb{RHP2}.

\begin{eeeee}
The following correspondences should be noted:
\begin{itemize}
\item $g^{e} \colon \mathbb{C} \setminus (-\infty,\max \lbrace 0,a_{N+1}^{e}
\rbrace) \! \to \! \mathbb{C}$ solves the phase
conditions~$\pmb{(\mathrm{P}_{1}^{(1)})}$---$\pmb{(\mathrm{P}_{1}^{(6)})} \!
\Leftrightarrow \! \mu_{V}^{e} \! \in \! \mathcal{M}_{1}(\mathbb{R})$ solves
the variational conditions~($\pmb{\mathrm{P}_{1}^{(a)}}$)
and~($\pmb{\mathrm{P}_{1}^{(b)}}$);
\item $g^{o} \colon \mathbb{C} \setminus (-\infty,\max \lbrace 0,a_{N+1}^{o}
\rbrace) \! \to \! \mathbb{C}$ solves the phase
conditions~$\pmb{(\mathrm{P}_{2}^{(1)})}$---$\pmb{(\mathrm{P}_{2}^{(6)})}
\Leftrightarrow \! \mu_{V}^{o} \! \in \! \mathcal{M}_{1}(\mathbb{R})$ solves
the variational conditions~($\pmb{\mathrm{P}_{2}^{(a)}}$)
and~($\pmb{\mathrm{P}_{2}^{(b)}}$). \hfill $\blacksquare$
\end{itemize}
\end{eeeee}
\section{Asymptotics for the Root Products, Recurrence Relation Coefficients, 
and Hankel Determinant Ratios of OLPs}
In this section, asymptotics in the double-scaling limit as $\mathscr{N},n \! 
\to \! \infty$ such that $z_{o} \! = \! 1 \! + \! o(1)$, denoted, hereafter, 
via the simplified `notation' $n \! \to \! \infty$, for the following 
quantities are obtained:
\begin{itemize}
\item the products of the (real) roots of the $L$-polynomials, for both the 
non-singular and singular cases (cf. Equations~(1.10)--(1.13));
\item the Hankel determinant ratios associated with the bi-infinite, 
real-valued, strong moment sequence $\left\lbrace \mathstrut c_{k} \! = \! 
\int_{\mathbb{R}}s^{k} \exp (-n \widetilde{V}(s)) \, \md s, \, n \! \in \! 
\mathbb{N} \right\rbrace_{k \in \mathbb{Z}}$ (cf. Equations~(1.9), (1.18), 
and~(1.19));
\item the (real-valued) recurrence relation coefficients (cf. 
Equations~(1.23)--(1.26)), and subsequently those of the (real-symmetric, 
tri-penta-diagonal-type) Laurent-Jacobi matrices $\mathcal{F}$ and $\mathcal{
G}$ (cf. Equations~(1.27) and~(1.28)).
\end{itemize}

Now that the relevant auxiliary information, functions, and results have been 
given in Section~2, one turns to the large-$n$ asymptotic description of the 
solutions to \pmb{RHP1} and~\pmb{RHP2}. As mentioned in Sections~1 and~2, this 
asymptotic analysis was the subject of \cite{a21} (for \pmb{RHP1}) and 
\cite{a22} (for \pmb{RHP2}). In this manuscript, all that is necessary {}from 
the results of \cite{a21,a22} are the large-$n$ asymptotic behaviours of the 
$n$-dependent coefficients of the following asymptotic expansions for 
$\overset{e}{\operatorname{Y}}(z)$ and $\overset{o}{\operatorname{Y}}(z)$:
\begin{gather}
\overset{e}{\mathrm{Y}}(z)z^{-n \sigma_{3}} \underset{\underset{z \in \mathbb{
C} \setminus \mathbb{R}}{z \to \infty}}{=} \mathrm{I} \! + \! \dfrac{1}{z}
\mathrm{Y}^{e,\infty}_{1} \! + \! \dfrac{1}{z^{2}} \mathrm{Y}^{e,\infty}_{2}
\! + \! \mathcal{O} \! \left(\dfrac{1}{z^{3}} \right), \\
\overset{e}{\mathrm{Y}}(z)z^{n \sigma_{3}} \underset{\underset{z \in \mathbb{
C} \setminus \mathbb{R}}{z \to 0}}{=} \mathrm{Y}^{e,0}_{0} \! + \! z \mathrm{
Y}^{e,0}_{1} \! + \! z^{2} \mathrm{Y}^{e,0}_{2} \! + \! \mathcal{O}(z^{3}), \\
\overset{o}{\mathrm{Y}}(z)z^{n \sigma_{3}} \underset{\underset{z \in \mathbb{
C} \setminus \mathbb{R}}{z \to 0}}{=} \mathrm{I} \! + \! z \mathrm{Y}^{o,0}_{
1} \! + \! z^{2} \mathrm{Y}^{o,0}_{2} \! + \! \mathcal{O}(z^{3}), \\
\overset{o}{\mathrm{Y}}(z)z^{-(n+1) \sigma_{3}} \underset{\underset{z \in
\mathbb{C} \setminus \mathbb{R}}{z \to \infty}}{=} \mathrm{Y}^{o,\infty}_{0}
\! + \! \dfrac{1}{z} \mathrm{Y}^{o,\infty}_{1} \! + \! \dfrac{1}{z^{2}}
\mathrm{Y}^{o,\infty}_{2} \! + \! \mathcal{O} \! \left(\dfrac{1}{z^{3}}
\right),
\end{gather}
where, for the convenience of the presentation which follows, as well as that
of the reader's, large-$n$ asymptotics of $\mathrm{Y}^{e,\infty}_{j}$, $j
\! = \! 1,2$, and $\mathrm{Y}^{e,0}_{k}$, $k \! = \! 0,1,2$, are given in
Appendix~A (Theorems~A.1 and~A.2, respectively), and large-$n$ asymptotics of
$\mathrm{Y}^{o,0}_{j}$, $j \! = \! 1,2$, and $\mathrm{Y}^{o,\infty}_{k}$, $k
\! = \! 0,1,2$, are given in Appendix~B (Theorems~B.1 and~B.2, respectively).
Furthermore, large-$n$ asymptotics of the `even' and `odd' leading
coefficients $\xi^{(2n)}_{n}$ and $\xi^{(2n+1)}_{-n-1}$, respectively, are
given in  Appendix~A (Theorem~A.3) and Appendix~B (Theorem~B.3).

The underlying idea of the following analysis is relatively straightforward:
the quantities of interest, namely, the coefficients in the (system of) three-
and five-term recurrence relations, Hankel determinant ratios, etc., are
all expressible in terms of the $n$-dependent coefficients of asymptotic
expansions~(3.1)--(3.4), for which rather complete large-$n$ asymptotics are
available.
\begin{eeeee}
It turns out that (cf. Equations~(1.29) and~(1.30)) only the $(1 \, 1)$- and 
$(1 \, 2)$-elements of asymptotics~(3.1)--(3.4) are necessary for obtaining 
the results of this work. In Lemma~3.1 (see below), and Propositions~3.1--3.3 
(see below), a series of formulae are presented which provide explicit 
relations between the quantities of interest and the $(n$-dependent) 
coefficients in the asymptotic expansions~(3.1)--(3.4). The statements, and 
proofs, are somewhat lengthy; but the important point should not be missed, 
namely: all quantities of interest are expressed, via straightforward 
algebraic calculations, in terms of the above asymptotic expansions. 
\hfill $\blacksquare$
\end{eeeee}
\begin{ccccc}
Let the external field $\widetilde{V} \colon \mathbb{R} \setminus \{0\} \! 
\to \! \mathbb{R}$ be regular and satisfy conditions~{\rm (1.20)--(1.22)}. 
Let the orthonormal $L$-polynomials (resp., monic orthogonal $L$-polynomials), 
$\lbrace \phi_{k}(z) \rbrace_{k \in \mathbb{Z}_{0}^{+}}$ (resp., $\lbrace 
\boldsymbol{\pi}_{k}(z) \rbrace_{k \in \mathbb{Z}_{0}^{+}})$ be as defined 
in Equations~{\rm (1.2)} and~{\rm (1.3)} (resp., Equations~{\rm (1.4)} 
and~{\rm (1.5))}, and let $\lbrace \phi_{k}(z) \rbrace_{k \in \mathbb{Z}_{
0}^{+}}$ satisfy the system of recurrence relations~{\rm (1.23)--(1.26)}. 
Let $\overset{e}{\operatorname{Y}} \colon \mathbb{C} \setminus \mathbb{R} \! 
\to \! \operatorname{SL}_{2}(\mathbb{C})$ (resp., $\overset{o}{\operatorname{
Y}} \colon \mathbb{C} \setminus \mathbb{R} \! \to \! \operatorname{SL}_{2}
(\mathbb{C}))$ be the (unique) solution of {\rm \pmb{RHP1}} (resp., 
{\rm \pmb{RHP2})} with integral representation {\rm (1.29)} (resp., 
{\rm (1.30))}. Then,
\begin{align}
z^{n} \overset{e}{\operatorname{Y}}_{12}(z) \underset{\underset{z \in \mathbb{
C} \setminus \mathbb{R}}{z \to \infty}}{=}& \, \dfrac{1}{z} \! \left(-\dfrac{(
\xi^{(2n)}_{n})^{-2}}{2 \pi \mi} \right) \! + \! \dfrac{1}{z^{2}} \! \left(-
\dfrac{a_{2n}^{\sharp}(\xi^{(2n)}_{n})^{-2}}{2 \pi \mi} \! + \! \dfrac{(\xi^{
(2n)}_{n})^{-1}}{2 \pi \mi} \! \left(\dfrac{b_{2n}^{\sharp} \nu^{(2n)}_{-n}}{
\xi^{(2n-1)}_{-n}} \right. \right. \nonumber \\
&\left. \left. + \, \dfrac{c_{2n}^{\sharp}}{\xi^{(2n-2)}_{n-1}} \! \left(
\dfrac{\xi^{(2n)}_{n-1}}{\xi^{(2n)}_{n}} \! - \! \dfrac{\nu^{(2n)}_{-n} \xi^{
(2n-1)}_{n-1}}{\xi^{(2n-1)}_{-n}} \right) \right) \right) \! + \! \mathcal{O}
\! \left(\dfrac{1}{z^{3}} \right), \\
z^{-n} \overset{e}{\operatorname{Y}}_{12}(z) \underset{\underset{z \in \mathbb{
C} \setminus \mathbb{R}}{z \to 0}}{=}& \, \dfrac{1}{2 \pi \mi} \dfrac{\beta_{2
n}^{\sharp}}{\xi^{(2n)}_{n} \xi^{(2n-1)}_{-n}} \! + \! z \left(\dfrac{(\xi^{(2
n)}_{n})^{-1}}{2 \pi \mi} \! \left(\dfrac{\beta_{2n+1}^{\sharp}}{\xi^{(2n+1)}_{
-n-1}} \! - \! \dfrac{\alpha_{2n}^{\sharp} \nu^{(2n+1)}_{n}}{\xi^{(2n)}_{n}}
\right. \right. \nonumber \\
&\left. \left. + \, \dfrac{\beta_{2n}^{\sharp}}{\xi^{(2n-1)}_{-n}} \! \left(
\nu^{(2n+1)}_{n} \nu^{(2n)}_{-n} \! - \! \dfrac{\xi^{(2n+1)}_{-n}}{\xi^{(2n+1)
}_{-n-1}} \right) \right) \right) \! + \! \mathcal{O}(z^{2}), \\
\intertext{and}
z^{-n} \overset{o}{\operatorname{Y}}_{12}(z) \underset{\underset{z \in \mathbb{
C} \setminus \mathbb{R}}{z \to 0}}{=}& \, z \left(\dfrac{(\xi^{(2n+1)}_{-n-1}
)^{-2}}{2 \pi \mi} \right) \! + \! z^{2} \! \left(\dfrac{\alpha_{2n+1}^{\sharp}
(\xi^{(2n+1)}_{-n-1})^{-2}}{2 \pi \mi} \! - \! \dfrac{(\xi^{(2n+1)}_{-n-1})^{-
1}}{2 \pi \mi} \! \left(\dfrac{\beta_{2n+1}^{\sharp} \nu^{(2n+1)}_{n}}{\xi^{(
2n)}_{n}} \right. \right. \nonumber \\
&\left. \left. + \, \dfrac{\gamma_{2n+1}^{\sharp}}{\xi^{(2n-1)}_{-n}} \! \left(
\dfrac{\xi^{(2n+1)}_{-n}}{\xi^{(2n+1)}_{-n-1}} \! -\nu^{(2n+1)}_{n} \nu^{(2n)
}_{-n} \right) \right) \right) \! + \! \mathcal{O}(z^{3}), \\
z^{n+1} \overset{o}{\operatorname{Y}}_{12}(z) \underset{\underset{z \in
\mathbb{C} \setminus \mathbb{R}}{z \to \infty}}{=}& \, -\dfrac{1}{2 \pi \mi}
\dfrac{b_{2n+1}^{\sharp}}{\xi^{(2n+1)}_{-n-1} \xi^{(2n)}_{n}} \! + \! \dfrac{
1}{z} \! \left(\dfrac{(\xi^{(2n+1)}_{-n-1})^{-1}}{2 \pi \mi} \! \left(\dfrac{
a_{2n+1}^{\sharp}}{\xi^{(2n+1)}_{-n-1}} \dfrac{\xi^{(2n+2)}_{-(n+1)}}{\xi^{(2n
+2)}_{n+1}} \! - \! \dfrac{b_{2n+2}^{\sharp}}{\xi^{(2n+2)}_{n+1}} \right.
\right. \nonumber \\
&\left. \left. + \, \dfrac{b_{2n+1}^{\sharp}}{\xi^{(2n)}_{n}} \! \left(\dfrac{
\xi^{(2n+2)}_{n}}{\xi^{(2n+2)}_{n+1}} \! -\nu^{(2n+1)}_{n} \dfrac{\xi^{(2n+2)
}_{-(n+1)}}{\xi^{(2n+2)}_{n+1}} \right) \right) \right) \! + \! \mathcal{O} \!
\left(\dfrac{1}{z^{2}} \right).
\end{align}
\end{ccccc}

\emph{Proof.} Noting {}from Lemma~1.1, Equation~(1.29) that
\begin{gather*}
\boldsymbol{\pi}_{2n}(z) \! := \! \overset{e}{\operatorname{Y}}_{11}(z) 
\qquad \quad \text{and} \qquad \quad \overset{e}{\operatorname{Y}}_{12}(z) \! 
= \! \int_{\mathbb{R}} \dfrac{\boldsymbol{\pi}_{2n}(s) \exp (-n \widetilde{V}
(s))}{s \! - \! z} \, \dfrac{\md s}{2 \pi \mi}, \quad z \! \in \! \mathbb{C} 
\setminus \mathbb{R},
\end{gather*}
and recalling that, for $\widetilde{V} \colon \mathbb{R} \setminus \{0\} \! 
\to \! \mathbb{R}$ satisfying conditions~(1.20)--(1.22) and the regularity 
assumption, $\int_{\mathbb{R}}s^{k} \exp (-n \widetilde{V}(s)) \, \md s \! 
< \! \infty$, $(k,n) \! \in \! \mathbb{Z} \times \mathbb{N}$, one uses the 
expansion $\tfrac{1}{s-z} \! = \! -\sum_{j=0}^{m} \tfrac{s^{j}}{z^{j+1}} \! 
+ \! \tfrac{s^{m+1}}{z^{m+1}(s-z)}$, $m \! \in \! \mathbb{Z}_{0}^{+}$, for 
$(\mathbb{C} \setminus \mathbb{R} \! \ni)$ $z \! \to \! \infty$, and $\tfrac{
1}{z-s} \! = \! -\sum_{j=0}^{m} \tfrac{z^{j}}{s^{j+1}} \! + \! \tfrac{z^{m+1}
}{s^{m+1}(z-s)}$ for $(\mathbb{C} \setminus \mathbb{R} \! \ni)$ $z \! \to \! 
0$, and proceeds exactly as in the proof of Proposition~5.4 of \cite{a21}. The 
proof of expansions~(3.5) and~(3.6) is now completed upon using, repeatedly, 
the orthogonality relations $\langle \boldsymbol{\pi}_{2n},z^{j} \rangle_{
\mathscr{L}} \! = \! 0$, $j \! = \! -n,\dotsc,n \! - \! 1$, $\langle 
\boldsymbol{\pi}_{2n+1},z^{i} \rangle_{\mathscr{L}} \! = \! 0$, $i \! = \! 
-n,\dotsc,n$, $\langle \boldsymbol{\pi}_{2n},\boldsymbol{\pi}_{2n} \rangle_{
\mathscr{L}} \! = \! (\xi^{(2n)}_{n})^{-2}$, and $\langle \boldsymbol{\pi}_{2
n+1},\boldsymbol{\pi}_{2n+1} \rangle_{\mathscr{L}} \! = \! (\xi^{(2n+1)}_{-
n-1})^{-2}$, as well as Equations~(1.2)--(1.5), the recurrence 
relations~(1.23)--(1.26), and (the $(1 \, 2)$-elements of) asymptotic 
expansions~(3.1) and~(3.2).

\emph{Mutatis mutandis} for expansions~(3.7) and~(3.8); but, in the latter 
case, one begins with Lemma 1.2, Equation~(1.30), that is,
\begin{gather*}
z \boldsymbol{\pi}_{2n+1}(z) \! := \! \overset{o}{\operatorname{Y}}_{11}(z) 
\qquad \quad \text{and} \qquad \quad \overset{o}{\operatorname{Y}}_{12}(z) 
\! = \! z \int_{\mathbb{R}} \dfrac{(s \boldsymbol{\pi}_{2n+1}(s)) \exp (-n 
\widetilde{V}(s))}{s(s \! - \! z)} \, \dfrac{\md s}{2 \pi \mi}, \quad z \! 
\in \! \mathbb{C} \setminus \mathbb{R},
\end{gather*}
uses (the $(1 \, 2)$-elements of) asymptotic expansions~(3.3) and~(3.4), and 
proceeds exactly as in the proof of Proposition~5.4 of \cite{a22}. \hfill 
$\qed$
\begin{bbbbb}
Let the external field $\widetilde{V} \colon \mathbb{R} \setminus \{0\} \! 
\to \! \mathbb{R}$ be regular and satisfy conditions~{\rm (1.20)--(1.22)}. 
Let the orthonormal $L$-polynomials (resp., monic orthogonal $L$-polynomials), 
$\lbrace \phi_{k}(z) \rbrace_{k \in \mathbb{Z}_{0}^{+}}$ (resp., $\lbrace 
\boldsymbol{\pi}_{k}(z) \rbrace_{k \in \mathbb{Z}_{0}^{+}})$ be as defined 
in Equations~{\rm (1.2)} and~{\rm (1.3)} (resp., Equations~{\rm (1.4)} 
and~{\rm (1.5))}, and let $\lbrace \phi_{k}(z) \rbrace_{k \in \mathbb{Z}_{
0}^{+}}$ satisfy the system of recurrence relations~{\rm (1.23)--(1.26)}. Let 
$\overset{e}{\operatorname{Y}} \colon \mathbb{C} \setminus \mathbb{R} \! \to 
\! \operatorname{SL}_{2}(\mathbb{C})$ (resp., $\overset{o}{\operatorname{Y}} 
\colon \mathbb{C} \setminus \mathbb{R} \! \to \! \operatorname{SL}_{2}
(\mathbb{C}))$ be the (unique) solution of {\rm \pmb{RHP1}} (resp., 
{\rm \pmb{RHP2})} with integral representation {\rm (1.29)} (resp., 
{\rm (1.30))}. Then a set of relations between the coefficients of 
asymptotic expansions~{\rm (3.1)--(3.4)}, the `even' and `odd' norming 
constants, and the recurrence coefficients is:
\begin{gather}
\xi^{(2n)}_{n} \! = \! \dfrac{1}{(-2 \pi \mi (\operatorname{Y}^{e,\infty}_{1}
)_{12})^{1/2}} \qquad (> \! 0), \\
2 \pi \mi (\operatorname{Y}^{e,0}_{0})_{12} \! = \! \dfrac{\beta_{2n}^{\sharp}
}{\xi^{(2n)}_{n} \xi^{(2n-1)}_{-n}}, \\
\xi^{(2n+1)}_{-n-1} \! = \! \dfrac{1}{(2 \pi \mi (\operatorname{Y}^{o,0}_{1}
)_{12})^{1/2}} \qquad (> \! 0), \\
-2 \pi \mi (\operatorname{Y}^{o,\infty}_{0})_{12} \! = \! \dfrac{b_{2n+1}^{
\sharp}}{\xi^{(2n+1)}_{-n-1} \xi^{(2n)}_{n}}, \\
-2 \pi \mi (\operatorname{Y}^{e,\infty}_{2})_{12} \! = \! \dfrac{a_{2n}^{
\sharp}}{(\xi^{(2n)}_{n})^{2}} \! - \! \dfrac{1}{\xi^{(2n)}_{n}} \! \left(
\dfrac{b_{2n}^{\sharp} \nu^{(2n)}_{-n}}{\xi^{(2n-1)}_{-n}} \! + \! \dfrac{c_{2
n}^{\sharp}}{\xi^{(2n-2)}_{n-1}} \! \left(\dfrac{\xi^{(2n)}_{n-1}}{\xi^{(2n)
}_{n}} \! -\nu^{(2n)}_{-n} \dfrac{\xi^{(2n-1)}_{n-1}}{\xi^{(2n-1)}_{-n}}
\right) \right), \\
2 \pi \mi (\operatorname{Y}^{e,0}_{1})_{12} \! = \! \dfrac{1}{\xi^{(2n)}_{n}}
\! \left(\dfrac{\beta_{2n+1}^{\sharp}}{\xi^{(2n+1)}_{-n-1}} \! - \! \dfrac{
\alpha_{2n}^{\sharp} \nu^{(2n+1)}_{n}}{\xi^{(2n)}_{n}} \! + \! \dfrac{\beta_{2
n}^{\sharp}}{\xi^{(2n-1)}_{-n}} \! \left(\nu^{(2n+1)}_{n} \nu^{(2n)}_{-n} \! -
\! \dfrac{\xi^{(2n+1)}_{-n}}{\xi^{(2n+1)}_{-n-1}} \right) \right), \\
2 \pi \mi (\operatorname{Y}^{o,0}_{2})_{12} \! = \! \dfrac{\alpha_{2n+1}^{
\sharp}}{(\xi^{(2n+1)}_{-n-1})^{2}} \! - \! \dfrac{1}{\xi^{(2n+1)}_{-n-1}} \!
\left(\dfrac{\beta_{2n+1}^{\sharp} \nu^{(2n+1)}_{n}}{\xi^{(2n)}_{n}} \! + \!
\dfrac{\gamma_{2n+1}^{\sharp}}{\xi^{(2n-1)}_{-n}} \! \left(\dfrac{\xi^{(2n+1)
}_{-n}}{\xi^{(2n+1)}_{-n-1}} \! -\nu^{(2n+1)}_{n} \nu^{(2n)}_{-n} \right)
\right), \\
-2 \pi \mi (\operatorname{Y}^{o,\infty}_{1})_{12} \! = \! \dfrac{1}{\xi^{(2n+1)
}_{-n-1}} \! \left(\dfrac{b_{2n+2}^{\sharp}}{\xi^{(2n+2)}_{n+1}} \! - \!
\dfrac{a_{2n+1}^{\sharp}}{\xi^{(2n+1)}_{-n-1}} \dfrac{\xi^{(2n+2)}_{-(n+1)}}{
\xi^{(2n+2)}_{n+1}} \! + \! \dfrac{b_{2n+1}^{\sharp}}{\xi^{(2n)}_{n}} \! \left(
\nu^{(2n+1)}_{n} \dfrac{\xi^{(2n+2)}_{-(n+1)}}{\xi^{(2n+2)}_{n+1}} \! - \!
\dfrac{\xi^{(2n+2)}_{n}}{\xi^{(2n+2)}_{n+1}} \right) \right).
\end{gather}
\end{bbbbb}

\emph{Proof.} Equate coefficients of like powers of $z^{-k}$ (resp., $z^{k})$,
$k \! = \! 1,2$, and $z^{i}$ (resp., $z^{-i})$, $i \! = \! 0,1$, respectively,
in---the $(1 \, 2)$-elements of---asymptotics~(3.1) and~(3.5) (resp.,
asymptotics~(3.3) and~(3.7)) and asymptotics~(3.2) and~(3.6) (resp.,
asymptotics~(3.4) and~(3.8)). \hfill $\qed$
\begin{eeeee}
Large-$n$ asymptotics for $\xi^{(2n)}_{n}$ are given in Appendix~A, 
Theorem~A.3; similarly, asymptotics for $\xi^{(2n+1)}_{-n-1}$ are given in 
Appendix~B, Theorem~B.3. These are obtained by using the $n \! \to \! \infty$ 
asymptotics for $(\operatorname{Y}^{e,\infty}_{1})_{12}$ and $(\operatorname{
Y}^{o,0}_{1})_{12}$, and Equations~(3.9) and~(3.11), respectively (see, also, 
\cite{a21,a22}).

It turns out that the following recurrence relation coefficients are the
essential building blocks {}from which the remaining ones can be expressed:
\begin{equation*}
a_{2n+1}^{\sharp}, \qquad b_{2n+1}^{\sharp}, \qquad b_{2n+2}^{\sharp}, \qquad
c_{2n+2}^{\sharp}, \qquad \alpha_{2n}^{\sharp}, \qquad \beta_{2n+1}^{\sharp},
\qquad \beta_{2n+2}^{\sharp}, \qquad \text{and} \qquad \gamma_{2n+3}^{\sharp};
\end{equation*}
for example, formulae for $a_{2n}^{\sharp}$ and $\alpha_{2n+1}^{\sharp}$ may
be derived: they have the functional form $a_{2n}^{\sharp} \! = \! a_{2n}^{
\sharp}(b_{2n}^{\sharp},\linebreak[4]
b_{2n+1}^{\sharp},\alpha_{2n}^{\sharp},\beta_{2n}^{\sharp},\beta_{2n+1}^{
\sharp})$ and $\alpha_{2n+1}^{\sharp} \! = \! \alpha_{2n+1}^{\sharp}(a_{2n}^{
\sharp},b_{2n}^{\sharp},b_{2n+1}^{\sharp},c_{2n+2}^{\sharp},\beta_{2n+1}^{
\sharp},\beta_{2n+2}^{\sharp},\gamma_{2n+1}^{\sharp})$. Furthermore,
e\-v\-e\-n though Equations~(3.13)--(3.16) are not used explicitly hereafter,
they give rise to several very interesting identities/relations between the
`even' and `odd' Riemann theta functions, $\boldsymbol{\theta}^{e}(z)$ and
$\boldsymbol{\theta}^{o}(z)$, respectively, and the end-points of the supports
of the `even' and `odd' equilibrium measures, $\lbrace b_{j-1}^{e},a_{j}^{e}
\rbrace_{j=1}^{N+1}$ and $\lbrace b_{j-1}^{o},a_{j}^{o} \rbrace_{j=1}^{N+1}$,
respectively: the details are left to the reader. \hfill $\blacksquare$
\end{eeeee}
\begin{bbbbb}
Let the external field $\widetilde{V} \colon \mathbb{R} \setminus \{0\} \! \to
\! \mathbb{R}$ be regular and satisfy conditions~{\rm (1.20)--(1.22)}. Let
the orthonormal $L$-polynomials (resp., monic orthogonal $L$-polynomials),
$\lbrace \phi_{k}(z) \rbrace_{k \in \mathbb{Z}_{0}^{+}}$ (resp., $\lbrace
\boldsymbol{\pi}_{k}(z) \rbrace_{k \in \mathbb{Z}_{0}^{+}})$ be as defined
in Equations~{\rm (1.2)} and~{\rm (1.3)} (resp., Equations~{\rm (1.4)}
and~{\rm (1.5))}, and let $\lbrace \phi_{k}(z) \rbrace_{k \in \mathbb{Z}_{
0}^{+}}$ satisfy the system of recurrence relations~{\rm (1.23)--(1.26)}. Let
$\overset{e}{\operatorname{Y}} \colon \mathbb{C} \setminus \mathbb{R} \! \to
\! \operatorname{SL}_{2}(\mathbb{C})$ (resp., $\overset{o}{\operatorname{Y}}
\colon \mathbb{C} \setminus \mathbb{R} \! \to \! \operatorname{SL}_{2}
(\mathbb{C}))$ be the (unique) solution of {\rm \pmb{RHP1}} (resp.,
{\rm \pmb{RHP2})} with integral representation {\rm (1.29)} (resp.,
{\rm (1.30))}. Then,
\begin{gather}
\begin{split}
\dfrac{\xi^{(2n)}_{-n}}{\xi^{(2n)}_{n}} \! =: \nu^{(2n)}_{-n} \! = \!
(\operatorname{Y}^{e,0}_{0})_{11},& \qquad \, \, \dfrac{\xi^{(2n)}_{-(n-1)}}{
\xi^{(2n)}_{n}} \! =: \nu^{(2n)}_{-(n-1)} \! = \! (\operatorname{Y}^{e,0}_{1}
)_{11}, \qquad \, \, \dfrac{\xi^{(2n)}_{-(n-2)}}{\xi^{(2n)}_{n}} \! = \!
(\operatorname{Y}^{e,0}_{2})_{11}, \\
&\dfrac{\xi^{(2n)}_{n-1}}{\xi^{(2n)}_{n}} \! = \! (\operatorname{Y}^{e,
\infty}_{1})_{11}, \qquad \, \, \dfrac{\xi^{(2n)}_{n-2}}{\xi^{(2n)}_{n}} \!
= \! (\operatorname{Y}^{e,\infty}_{2})_{11},
\end{split}
\end{gather}
\begin{gather}
\begin{split}
\dfrac{\xi^{(2n+1)}_{n}}{\xi^{(2n+1)}_{-n-1}} \! =: \nu^{(2n+1)}_{n} \! = \!
(\operatorname{Y}^{o,\infty}_{0})_{11},& \qquad \, \, \dfrac{\xi^{(2n+1)}_{n
-1}}{\xi^{(2n+1)}_{-n-1}} \! =: \nu^{(2n+1)}_{n-1} \! = \! (\operatorname{
Y}^{o,\infty}_{1})_{11}, \qquad \, \, \dfrac{\xi^{(2n+1)}_{n-2}}{\xi^{(2n+1)
}_{-n-1}} \! = \! (\operatorname{Y}^{o,\infty}_{2})_{11}, \\
&\dfrac{\xi^{(2n+1)}_{-n}}{\xi^{(2n+1)}_{-n-1}} \! = \! (\operatorname{Y}^{o,
0}_{1})_{11}, \qquad \, \, \dfrac{\xi^{(2n+1)}_{-(n-1)}}{\xi^{(2n+1)}_{-n-1}}
\! = \! (\operatorname{Y}^{o,0}_{2})_{11},
\end{split}
\end{gather}
and
\begin{gather}
\xi^{(2n+1)}_{n} \! = \! b_{2n+2}^{\sharp} \xi^{(2n+2)}_{n+1}, \\
\xi^{(2n+1)}_{n-1} \! = \! b_{2n+1}^{\sharp} \xi^{(2n)}_{n} \! + \! a_{2n+1}^{
\sharp} \xi^{(2n+1)}_{n} \! + \! b_{2n+2}^{\sharp} \xi^{(2n+2)}_{n}, \\
\xi^{(2n)}_{n} \! = \! c_{2n+2}^{\sharp} \xi^{(2n+2)}_{n+1}, \\
\xi^{(2n)}_{n-1} \! = \! a_{2n}^{\sharp} \xi^{(2n)}_{n} \! + \! b_{2n+1}^{
\sharp} \xi^{(2n+1)}_{n} \! + \! c_{2n+2}^{\sharp} \xi^{(2n+2)}_{n}, \\
b_{2n}^{\sharp} \xi^{(2n-1)}_{-n} \! + \! a_{2n}^{\sharp} \xi^{(2n)}_{-n} \! +
\! b_{2n+1}^{\sharp} \xi^{(2n+1)}_{-n} \! + \! c_{2n+2}^{\sharp} \xi^{(2n+2)}_{
-n} \! = \! 0, \\
b_{2n+1}^{\sharp} \xi^{(2n+1)}_{-n-1} \! + \! c_{2n+2}^{\sharp} \xi^{(2n+2)}_{
-(n+1)} \! = \! 0, \\
\xi^{(2n+1)}_{-n-1} \! = \! b_{2n+1}^{\sharp} \xi^{(2n)}_{-n} \! + \! a_{2n+
1}^{\sharp} \xi^{(2n+1)}_{-n} \! + \! b_{2n+2}^{\sharp} \xi^{(2n+2)}_{-n}, \\
a_{2n+1}^{\sharp} \xi^{(2n+1)}_{-n-1} \! + \! b_{2n+2}^{\sharp} \xi^{(2n+2)}_{
-(n+1)} \! = \! 0,
\end{gather}
where Equations~{\rm (3.19)}, {\rm (3.21)}, {\rm (3.24)}, and {\rm (3.26)}
are used to determine $b_{2n+2}^{\sharp}$, $c_{2n+2}^{\sharp}$, $b_{2n+1}^{
\sharp}$, and $a_{2n+1}^{\sharp}$, respectively (see Proposition~{\rm 3.3}
below),
\begin{gather}
\xi^{(2n+1)}_{-n-1} \! = \! \gamma_{2n+3}^{\sharp} \xi^{(2n+3)}_{-(n+2)}, \\
\xi^{(2n+1)}_{-n} \! = \! \alpha_{2n+1}^{\sharp} \xi^{(2n+1)}_{-n-1} \! + \!
\beta_{2n+2}^{\sharp} \xi^{(2n+2)}_{-(n+1)} \! + \! \gamma_{2n+3}^{\sharp}
\xi^{(2n+3)}_{-(n+1)}, \\
\beta_{2n+2}^{\sharp} \xi^{(2n+2)}_{n+1} \! + \! \gamma_{2n+3}^{\sharp} \xi^{
(2n+3)}_{n+1} \! = \! 0, \\
\beta_{2n+1}^{\sharp} \xi^{(2n)}_{n} \! + \! \alpha_{2n+1}^{\sharp} \xi^{(2n+
1)}_{n} \! + \! \beta_{2n+2}^{\sharp} \xi^{(2n+2)}_{n} \! + \! \gamma_{2n+3}^{
\sharp} \xi^{(2n+3)}_{n} \! = \! 0, \\
\xi^{(2n)}_{-n} \! = \! \beta_{2n+1}^{\sharp} \xi^{(2n+1)}_{-n-1}, \\
\xi^{(2n)}_{-(n-1)} \! = \! \beta_{2n}^{\sharp} \xi^{(2n-1)}_{-n} \! + \!
\alpha_{2n}^{\sharp} \xi^{(2n)}_{-n} \! + \! \beta_{2n+1}^{\sharp} \xi^{(2n+
1)}_{-n}, \\
\alpha_{2n}^{\sharp} \xi^{(2n)}_{n} \! + \! \beta_{2n+1}^{\sharp} \xi^{(2n+1)
}_{n} \! = \! 0, \\
\xi^{(2n)}_{n} \! = \! \beta_{2n}^{\sharp} \xi^{(2n-1)}_{n-1} \! + \! \alpha_{
2n}^{\sharp} \xi^{(2n)}_{n-1} \! + \! \beta_{2n+1}^{\sharp} \xi^{(2n+1)}_{n-
1},
\end{gather}
where Equations~{\rm (3.27)}, {\rm (3.29)}, {\rm (3.31)}, and {\rm (3.33)} are
used to determine $\gamma_{2n+3}^{\sharp}$, $\beta_{2n+2}^{\sharp}$, $\beta_{
2n+1}^{\sharp}$, and $\alpha_{2n}^{\sharp}$, respectively (see
Proposition~{\rm 3.3} below), with
\begin{gather}
a_{2n}^{\sharp} \! = \! \dfrac{1}{\alpha_{2n}^{\sharp}} \! - \! \dfrac{1}{
\alpha_{2n}^{\sharp}} \! \left(b_{2n}^{\sharp} \beta_{2n}^{\sharp} \! + \! b_{
2n+1}^{\sharp} \beta_{2n+1}^{\sharp} \right), \\
\alpha_{2n+1}^{\sharp} \! = \! -\dfrac{1}{b_{2n+1}^{\sharp}} \! \left(b_{2n}^{
\sharp} \gamma_{2n+1}^{\sharp} \! + \! a_{2n}^{\sharp} \beta_{2n+1}^{\sharp}
\! + \! c_{2n+2}^{\sharp} \beta_{2n+2}^{\sharp} \right).
\end{gather}
\end{bbbbb}

\emph{Proof.} One shows {}from Equation~(1.4) that
\begin{gather*}
z^{-n} \boldsymbol{\pi}_{2n}(z) \underset{\underset{z \in \mathbb{C} \setminus
\mathbb{R}}{z \to \infty}}{=} 1 \! + \! \dfrac{1}{z} \dfrac{\xi^{(2n)}_{n-1}}{
\xi^{(2n)}_{n}} \! + \! \dfrac{1}{z^{2}} \dfrac{\xi^{(2n)}_{n-2}}{\xi^{(2n)}_{
n}} \! + \! \mathcal{O} \! \left(\dfrac{1}{z^{3}} \right), \\
z^{n} \boldsymbol{\pi}_{2n}(z) \underset{\underset{z \in \mathbb{C} \setminus
\mathbb{R}}{z \to 0}}{=} \nu^{(2n)}_{-n} \! + \! z \dfrac{\xi^{(2n)}_{-(n-1)}}{
\xi^{(2n)}_{n}} \! + \! z^{2} \dfrac{\xi^{(2n)}_{-(n-2)}}{\xi^{(2n)}_{n}} \! +
\! \mathcal{O}(z^{3}),
\end{gather*}
where $\nu^{(2n)}_{-n}$ is defined in Equations~(1.9), and, {}from (the $(1 
\, 1)$-elements of) asymptotics~(3.1) and~(3.2), that
\begin{gather*}
z^{-n} \overset{e}{\operatorname{Y}}_{11}(z) \underset{\underset{z \in \mathbb{
C} \setminus \mathbb{R}}{z \to \infty}}{=} 1 \! + \! \dfrac{1}{z}
(\operatorname{Y}^{e,\infty}_{1})_{11} \! + \! \dfrac{1}{z^{2}}(\operatorname{
Y}^{e,\infty}_{2})_{11} \! + \! \mathcal{O} \! \left(\dfrac{1}{z^{3}} \right),
\\
z^{n} \overset{e}{\operatorname{Y}}_{11}(z) \underset{\underset{z \in \mathbb{
C} \setminus \mathbb{R}}{z \to 0}}{=} (\operatorname{Y}^{e,0}_{0})_{11} \! +
\! z(\operatorname{Y}^{e,0}_{1})_{11} \! + \! z^{2}(\operatorname{Y}^{e,0}_{2}
)_{11} \! + \! \mathcal{O}(z^{3});
\end{gather*}
via Lemma~1.1 (Equation~(1.29)), upon equating coefficients of like powers 
of $z^{-k}$, $k \! = \! 1,2$, and $z^{i}$, $i \! = \! 0,1,2$, in the above 
asymptotic expansions, one arrives at Equations~(3.17). 
\emph{Mutatis mutandis} for Equations~(3.18); but in the latter case, via 
Lemma~1.2 (Equation~(1.30)), one equates coefficients of like powers of $z^{
k}$, $k \! = \! 1,2$, and $z^{-i}$, $i \! = \! 0,1,2$, in the asymptotics 
(cf. Equation~(1.5))
\begin{gather*}
z^{n+1} \boldsymbol{\pi}_{2n+1}(z) \underset{\underset{z \in \mathbb{C}
\setminus \mathbb{R}}{z \to 0}}{=} 1 \! + \! z \dfrac{\xi^{(2n+1)}_{-n}}{\xi^{
(2n+1)}_{-n-1}} \! + \! z^{2} \dfrac{\xi^{(2n+1)}_{-(n-1)}}{\xi^{(2n+1)}_{-n-
1}} \! + \! \mathcal{O}(z^{3}), \\
z^{-n} \boldsymbol{\pi}_{2n+1}(z) \underset{\underset{z \in \mathbb{C}
\setminus \mathbb{R}}{z \to \infty}}{=} \nu^{(2n+1)}_{n} \! + \! \dfrac{1}{z}
\dfrac{\xi^{(2n+1)}_{n-1}}{\xi^{(2n+1)}_{-n-1}} \! + \! \dfrac{1}{z^{2}}
\dfrac{\xi^{(2n+1)}_{n-2}}{\xi^{(2n+1)}_{-n-1}} \! + \! \mathcal{O} \! \left(
\dfrac{1}{z^{3}} \right),
\end{gather*}
where $\nu^{(2n+1)}_{n}$ is defined in Equations~(1.9), and (cf. the $(1 \, 
1)$-elements of asymptotic expansions~(3.3) and (3.4))
\begin{gather*}
z^{n+1} \! \left(z^{-1} \overset{o}{\operatorname{Y}}_{11}(z) \right)
\underset{\underset{z \in \mathbb{C} \setminus \mathbb{R}}{z \to 0}}{=} 1 \!
+ \! z(\operatorname{Y}^{o,0}_{1})_{11} \! + \! z^{2}(\operatorname{Y}^{o,
0}_{2})_{11} \! + \! \mathcal{O}(z^{3}), \\
z^{-n} \! \left(z^{-1} \overset{o}{\operatorname{Y}}_{11}(z) \right)
\underset{\underset{z \in \mathbb{C} \setminus \mathbb{R}}{z \to \infty}}{=}
(\operatorname{Y}^{o,\infty}_{0})_{11} \! + \! \dfrac{1}{z}(\operatorname{
Y}^{o,\infty}_{1})_{11} \! + \! \dfrac{1}{z^{2}}(\operatorname{Y}^{o,
\infty}_{2})_{11} \! + \! \mathcal{O} \! \left(\dfrac{1}{z^{3}} \right).
\end{gather*}
Equations~(3.19)--(3.34) are derived by substituting Equations~(1.2) and~(1.3)
into recurrence relations~(1.23)--(1.26) and equating coefficients of like
powers of $z^{\pm 1}$ as $(\mathbb{C} \setminus \mathbb{R} \! \ni)$ $z \! \to
\! \infty$ and $z \! \to \! 0$. Recall that (cf. Equations~(1.27) and~(1.28))
$\mathcal{F} \mathcal{G} \! = \! \mathcal{G} \mathcal{F} \! = \! \operatorname{
diag}(1,1,\dotsc,1,\dotsc)$; forming, for example, the product $\mathcal{F}
\mathcal{G}$ (the result remains unchanged upon forming the product $\mathcal{
G} \mathcal{F})$, using Equations~(3.19)--(3.34), and recalling that (due to
the chosen normalisation for the orthonormal $L$-polynomials) $\xi^{(0)}_{0} 
\! = \! 1$, one shows, after a tedious, but otherwise straightforward, 
algebraic calculation that the only non-trivial relations between the 
recurrence relation coefficients are Equations~(3.35) and~(3.36). \hfill $\qed$
\begin{bbbbb}
Let the external field $\widetilde{V} \colon \mathbb{R} \setminus \{0\} \!
\to \! \mathbb{R}$ be regular and satisfy conditions~{\rm (1.20)--(1.22)}.
Let the orthonormal $L$-polynomials (resp., monic orthogonal $L$-polynomials),
$\lbrace \phi_{k}(z) \rbrace_{k \in \mathbb{Z}_{0}^{+}}$ (resp., $\lbrace
\boldsymbol{\pi}_{k}(z) \rbrace_{k \in \mathbb{Z}_{0}^{+}})$ be as defined
in Equations~{\rm (1.2)} and~{\rm (1.3)} (resp., Equations~{\rm (1.4)}
and~{\rm (1.5))}, and let $\lbrace \phi_{k}(z) \rbrace_{k \in \mathbb{Z}_{
0}^{+}}$ satisfy the system of recurrence relations~{\rm (1.23)--(1.26)}. Let
$\overset{e}{\operatorname{Y}} \colon \mathbb{C} \setminus \mathbb{R} \! \to
\! \operatorname{SL}_{2}(\mathbb{C})$ (resp., $\overset{o}{\operatorname{Y}}
\colon \mathbb{C} \setminus \mathbb{R} \! \to \! \operatorname{SL}_{2}
(\mathbb{C}))$ be the (unique) solution of {\rm \pmb{RHP1}} (resp.,
{\rm \pmb{RHP2})} with integral representation {\rm (1.29)} (resp.,
{\rm (1.30))} and asymptotics~{\rm (3.1)} and~{\rm (3.2)} (resp.,
asymptotics~{\rm (3.3)} and~{\rm (3.4))}. Then,
\begin{gather}
a_{2n+1}^{\sharp} \! = \! -\nu^{(2n+1)}_{n} \dfrac{(\operatorname{Y}^{o,
\infty}_{0})_{12}}{(\operatorname{Y}^{o,0}_{1})_{12}} \qquad (\in \!
\mathbb{R}), \\
b_{2n+1}^{\sharp} \! = \! -2 \pi \mi \xi^{(2n)}_{n} \xi^{(2n+1)}_{-n-1}
(\operatorname{Y}^{o,\infty}_{0})_{12} \qquad (\in \! \mathbb{R}), \\
b_{2n+2}^{\sharp} \! = \nu^{(2n+1)}_{n} \sqrt{\dfrac{-2 \pi \mi h^{+}
[(\operatorname{Y}^{e,\infty}_{1})_{12}]}{2 \pi \mi (\operatorname{Y}^{o,0}_{
1})_{12}}} \qquad (\in \! \mathbb{R}), \\
c_{2n+2}^{\sharp} \! = \! \sqrt{\dfrac{-2 \pi \mi h^{+}[(\operatorname{Y}^{e,
\infty}_{1})_{12}]}{-2 \pi \mi (\operatorname{Y}^{e,\infty}_{1})_{12}}} \qquad
(\in \! \mathbb{R}_{+}), \\
\alpha^{\sharp}_{2n} \! = \! -\nu^{(2n)}_{-n} \nu^{(2n+1)}_{n} \qquad (\in \!
\mathbb{R}), \\
\beta_{2n+1}^{\sharp} \! = \nu^{(2n)}_{-n} \dfrac{\xi^{(2n)}_{n}}{\xi^{(2n+1)
}_{-n-1}} \qquad (\in \! \mathbb{R}), \\
\beta_{2n+2}^{\sharp} \! = \! -\left(h^{+} \! \left[\nu^{(2n+1)}_{n} \right]
\right) \! \sqrt{\dfrac{-2 \pi \mi h^{+}[(\operatorname{Y}^{e,\infty}_{1})_{1
2}]}{2 \pi \mi (\operatorname{Y}^{o,0}_{1})_{12}}} \qquad (\in \! \mathbb{R}),
\\
\gamma_{2n+3}^{\sharp} \! = \sqrt{\dfrac{2 \pi \mi h^{+}[(\operatorname{Y}^{o,
0}_{1})_{12}]}{2 \pi \mi (\operatorname{Y}^{o,0}_{1})_{12}}} \qquad (\in \!
\mathbb{R}_{+}),
\end{gather}
where all the square roots are positive.
\end{bbbbb}

\emph{Proof.} One commences by establishing the following identity:
\begin{equation}
h^{+} \! \left[\nu^{(2n)}_{-n} \right] \! = \! \dfrac{(\operatorname{Y}^{o,
\infty}_{0})_{12}}{(\operatorname{Y}^{o,0}_{1})_{12}}, \quad n \! \in \!
\mathbb{Z}_{0}^{+}.
\end{equation}
{}From Equations~(3.9),\, (3.11),\, (3.21), and~(3.24), one shows that
\begin{align*}
b_{2n+1}^{\sharp} &= -c_{2n+2}^{\sharp} \dfrac{\xi^{(2n+2)}_{-(n+1)}}{\xi^{(2n
+1)}_{-n-1}} \! = \! -\dfrac{\xi^{(2n)}_{n}}{\xi^{(2n+1)}_{-n-1}} \dfrac{\xi^{
(2n+2)}_{-(n+1)}}{\xi^{(2n+2)}_{n+1}} \! = \! -\dfrac{\xi^{(2n)}_{n}}{\xi^{(2n
+1)}_{-n-1}}h^{+} \! \left[\xi^{(2n)}_{-n}/\xi^{(2n)}_{n} \right] \\
&= \, -\dfrac{\xi^{(2n)}_{n}}{\xi^{(2n+1)}_{-n-1}}h^{+} \! \left[\nu^{(2n)}_{
-n} \right] \! = \! - \! \left(h^{+} \! \left[\nu^{(2n)}_{-n} \right] \right)
\! \dfrac{(2 \pi \mi (\operatorname{Y}^{o,0}_{1})_{12})^{1/2}}{(-2 \pi \mi
(\operatorname{Y}^{e,\infty}_{1})_{12})^{1/2}};
\end{align*}
but, {}from Equation~(3.12), which, in fact, leads to Equation~(3.38) 
for $b_{2n+1}^{\sharp}$, via Equations~(3.9) and (3.11), one has that
\begin{equation*}
b_{2n+1}^{\sharp} \! = \! -\dfrac{2 \pi \mi (\operatorname{Y}^{o,\infty}_{0}
)_{12}}{(-2 \pi \mi (\operatorname{Y}^{e,\infty}_{1})_{12})^{1/2}(2 \pi \mi
(\operatorname{Y}^{o,0}_{1})_{12})^{1/2}}:
\end{equation*}
equating the above two expressions for $b_{2n+1}^{\sharp}$, for the 
\textbf{same} choices of branches, one arrives at Identity (3.45). 
{}From Equations~(3.19) and~(3.26), one has that
\begin{align*}
a_{2n+1}^{\sharp} &= -b_{2n+2}^{\sharp} \dfrac{\xi^{(2n+2)}_{-(n+1)}}{\xi^{(2n
+1)}_{-n-1}} \! = \! -\dfrac{\xi^{(2n+1)}_{n}}{\xi^{(2n+2)}_{n+1}} \dfrac{\xi^{
(2n+2)}_{-(n+1)}}{\xi^{(2n+1)}_{-n-1}} \! = \! -\dfrac{\xi^{(2n+1)}_{n}}{\xi^{
(2n+1)}_{-n-1}} \dfrac{\xi^{(2n+2)}_{-(n+1)}}{\xi^{(2n+2)}_{n+1}} \\
&= \, -\nu^{(2n+1)}_{n}h^{+} \! \left[\xi^{(2n)}_{-n}/\xi^{(2n)}_{n} \right]
\! = \! -\nu^{(2n+1)}_{n}h^{+} \! \left[\nu^{(2n)}_{-n} \right],
\end{align*}
which, via Identity~(3.45), leads to Equation~(3.37). {}From Equation~(3.19),
one has that
\begin{equation*}
b_{2n+2}^{\sharp} \! = \! \dfrac{\xi^{(2n+1)}_{n}}{\xi^{(2n+2)}_{n+1}} \! =
\! \dfrac{\xi^{(2n+1)}_{n}}{\xi^{(2n+1)}_{-n-1}} \dfrac{\xi^{(2n+1)}_{-n-1}
}{\xi^{(2n+2)}_{n+1}} \! =\nu^{(2n+1)}_{n} \dfrac{\xi^{(2n+1)}_{-n-1}}{h^{+}
\! \left[\xi^{(2n)}_{n} \right]},
\end{equation*}
which, via Equations~(3.9) and~(3.11), lead to Equation~(3.39). {}From
Equation~(3.21), one has that
\begin{equation*}
c_{2n+2}^{\sharp} \! = \! \dfrac{\xi^{(2n)}_{n}}{\xi^{(2n+2)}_{n+1}} \! = \!
\dfrac{\xi^{(2n)}_{n}}{h^{+} \! \left[\xi^{(2n)}_{n} \right]},
\end{equation*}
which, via Equation~(3.9), leads to Equation~(3.40). {}From Equations~(3.31)
and~(3.33), one shows that
\begin{equation*}
\alpha_{2n}^{\sharp} \! = \! -\beta_{2n+1}^{\sharp} \dfrac{\xi^{(2n+1)}_{n}}{
\xi^{(2n)}_{n}} \! = \! -\dfrac{\xi^{(2n)}_{-n}}{\xi^{(2n+1)}_{-n-1}} \dfrac{
\xi^{(2n+1)}_{n}}{\xi^{(2n)}_{n}} \! = \! -\dfrac{\xi^{(2n)}_{-n}}{\xi^{(2n)}_{
n}} \dfrac{\xi^{(2n+1)}_{n}}{\xi^{(2n+1)}_{-n-1}} \! = \! -\nu^{(2n)}_{-n}
\nu^{(2n+1)}_{n},
\end{equation*}
which is Equation~(3.41). {}From Equation~(3.31), one has that
\begin{equation*}
\beta_{2n+1}^{\sharp} \! = \! \dfrac{\xi^{(2n)}_{-n}}{\xi^{(2n+1)}_{-n-1}} \!
= \! \dfrac{\xi^{(2n)}_{-n}}{\xi^{(2n)}_{n}} \dfrac{\xi^{(2n)}_{n}}{\xi^{(2n+
1)}_{-n-1}} \! =\nu^{(2n)}_{-n} \dfrac{\xi^{(2n)}_{n}}{\xi^{(2n+1)}_{-n-1}},
\end{equation*}
which is Equation~(3.42). {}From Equations~(3.27) and~(3.29), one has that
\begin{align*}
\beta_{2n+2}^{\sharp} &= -\gamma_{2n+3}^{\sharp} \dfrac{\xi^{(2n+3)}_{n+1}}{
\xi^{(2n+2)}_{n+1}} \! = \! -\dfrac{\xi^{(2n+1)}_{-n-1}}{\xi^{(2n+3)}_{-(n+2)}
} \dfrac{\xi^{(2n+3)}_{n+1}}{\xi^{(2n+2)}_{n+1}} \! = \! -\dfrac{\xi^{(2n+3)}_{
n+1}}{\xi^{(2n+3)}_{-(n+2)}} \dfrac{\xi^{(2n+1)}_{-n-1}}{\xi^{(2n+2)}_{n+1}} \\
&= \, - \! \left(h^{+} \! \left[\xi^{(2n+1)}_{n}/\xi^{(2n+1)}_{-n-1} \right]
\right) \! \dfrac{\xi^{(2n+1)}_{-n-1}}{h^{+} \! \left[\xi^{(2n)}_{n} \right]}
\! = \! - \! \left(h^{+} \! \left[\nu^{(2n+1)}_{n} \right] \right) \! \dfrac{
\xi^{(2n+1)}_{-n-1}}{h^{+} \! \left[\xi^{(2n)}_{n} \right]},
\end{align*}
which, via Equations~(3.9) and~(3.11), lead to Equation~(3.43). Finally,
{}from Equation~(3.27), one has that
\begin{equation*}
\gamma_{2n+3}^{\sharp} \! = \! \dfrac{\xi^{(2n+1)}_{-n-1}}{\xi^{(2n+3)}_{-(n
+2)}} \! = \! \dfrac{\xi^{(2n+1)}_{-n-1}}{h^{+} \! \left[\xi^{(2n+1)}_{-n-1}
\right]},
\end{equation*}
which, via Equation~(3.11), leads to Equation~(3.44). \hfill $\qed$
\begin{eeeee}
Now, explicit formulae are at hand relating all quantities of interest
(such as recurrence relation coefficients, real root products, and
leading coefficients) to the $n$-dependent coefficients of asymptotic
expansions~(3.1)--(3.4). Using these, formulae for certain ratios of Hankel
determinants are also obtained. \hfill $\blacksquare$
\end{eeeee}
Recalling Equations~(1.9)--(1.13), via relations~(3.17) and~(3.18), one
arrives at: (i) for the non-singular case,
\begin{equation*}
\prod_{k=1}^{2n} \alpha^{(2n)}_{k} \! = \! \dfrac{H^{(-2n+1)}_{2n}}{H^{(-2n)}_{
2n}} \! = \! \dfrac{\xi^{(2n)}_{-n}}{\xi^{(2n)}_{n}} \! =: \nu^{(2n)}_{-n} \!
= \! (\operatorname{Y}^{e,0}_{0})_{11},
\end{equation*}
where (for $\xi^{(2n)}_{-n} \! \not= \! 0)$ $\left\{\mathstrut \alpha^{
(2n)}_{k}, \, k \! = \! 1,\dotsc,2n \right\} \! = \! \left\{\mathstrut z; \,
\boldsymbol{\pi}_{2n}(z) \! = \! 0 \right\}$, with $\alpha^{(2n)}_{k}$ real
and simple (see Theorem~3.1, Equation~(3.48)), and
\begin{equation*}
\left(\prod_{k=1}^{2n+1} \alpha^{(2n+1)}_{k} \right)^{-1} \! = \! \dfrac{H^{(-
2n-1)}_{2n+1}}{H^{(-2n)}_{2n+1}} \! = \! -\dfrac{\xi^{(2n+1)}_{n}}{\xi^{(2n+1)
}_{-n-1}} \! =: -\nu^{(2n+1)}_{n} \! = \! -(\operatorname{Y}^{o,\infty}_{0})_{
11},
\end{equation*}
where (for $\xi^{(2n+1)}_{n} \! \not= \! 0)$ $\left\{\mathstrut \alpha^{(2n+
1)}_{k}, \, k \! = \! 1,\dotsc,2n \! + \! 1 \right\} \! = \! \left\{\mathstrut
z; \, \boldsymbol{\pi}_{2n+1}(z) \! = \! 0 \right\}$, with $\alpha^{(2n+1)}_{
k}$ real and simple (see Theorem~3.1, Equation~(3.49)); and (ii) for the
singular case,
\begin{equation*}
\prod_{k=1}^{2n-1} \alpha^{(2n)}_{k} \! = \! -\dfrac{\xi^{(2n)}_{-(n-1)}}{\xi^{
(2n)}_{n}} \! =: -\nu^{(2n)}_{-(n-1)} \! = \! -(\operatorname{Y}^{e,0}_{1})_{1
1},
\end{equation*}
where (for $\xi^{(2n)}_{-n} \! = \! 0)$ $\left\{\mathstrut \alpha^{(2n)}_{k},
\, k \! = \! 1,\dotsc,2n \! - \! 1 \right\} \! = \! \left\{\mathstrut z; \,
\boldsymbol{\pi}_{2n}(z) \! = \! 0 \right\}$, with $\alpha^{(2n)}_{k}$ real
and simple (see Theorem~3.1, Equation~(3.50)), and
\begin{equation*}
\left(\prod_{k=1}^{2n} \alpha^{(2n+1)}_{k} \right)^{-1} \! = \! \dfrac{
\xi^{(2n+1)}_{n-1}}{\xi^{(2n+1)}_{-n-1}} \! =: \nu^{(2n+1)}_{n-1} \! = \!
(\operatorname{Y}^{o,\infty}_{1})_{11},
\end{equation*}
where (for $\xi^{(2n+1)}_{n} \! = \! 0)$ $\left\{\mathstrut \alpha^{(2n+1)}_{
k}, \, k \! = \! 1,\dotsc,2n \right\} \! = \! \left\{\mathstrut z; \,
\boldsymbol{\pi}_{2n+1}(z) \! = \! 0 \right\}$, with $\alpha^{(2n+1)}_{k}$
real and simple (see Theorem~3.1, Equation~(3.51)).
\begin{eeeee}
In order to eschew a flood of superfluous notation, the `symbol' $c(n)$ 
appears hereafter in the various error estimations (as $n \! \to \! \infty)$ 
of the form $\mathcal{O}(c(n)n^{-2})$ and $\mathcal{O}(c(n)(n \! + \! 1/2)^{
-2})$, where $c(n) \! =_{n \to \infty} \! \mathcal{O}(1)$. Even though, 
\emph{sensus strictu}, the symbol $c(n)$ should properly be denoted as $c_{1}
(n),c_{2}(n),\dotsc,c_{k}(n),\dotsc$, where, in general, $c_{i}(n) \! \not= 
\! c_{j}(n)$, $i \! \not= \! j \! \in \! \mathbb{N}$, and $c_{j}(n) \! =_{n 
\to \infty} \! \mathcal{O}(1)$, $j \! \in \! \mathbb{N}$, with a certain abuse 
of notation, and concomitant with the above theme, the simplified `notation' 
$c(n)$ will be retained. \hfill $\blacksquare$
\end{eeeee}
\begin{eeeee}
The notations $\infty^{\pm}$ (resp., $0^{\pm})$, which appear in the theorems 
below, are to be understood as follows: let $\mathscr{P} \! := \! (y,z)$ 
denote an arbitrary point on the (hyperelliptic) Riemann surface $\mathcal{
Y}_{e}$ or $\mathcal{Y}_{o}$, where, for $q \! \in \! \{e,o\}$, $\mathcal{
Y}_{q} \! := \! \left\lbrace \mathstrut (y,z); \, y^{2} \! = \! R_{q}(z) 
\right\rbrace$, with $R_{e}(z)$ (resp., $R_{o}(z))$ defined in 
\fbox{$\pmb{\mathrm{P}_{1}}$} (resp., \fbox{$\pmb{\mathrm{P}_{2}}$}). Then: 
(i) $\mathscr{P} \! \to \! \infty^{\pm} \! \Leftrightarrow \! z \! \to \! 
\infty$, $y \! \sim \! \pm z^{N+1}$; and (ii) $\mathscr{P} \! \to \! 0^{\pm} 
\! \Leftrightarrow \! z \! \to \! 0$, $y \! \sim \! \pm (-1)^{\mathcal{N}_{
+}^{q}}(\prod_{k=1}^{N+1} \vert b_{k-1}^{q}a_{k}^{q} \vert)^{1/2}$ for $q \! 
\in \! \{e,o\}$, where $\mathcal{N}_{+}^{q} \! \in \! \lbrace 0,\dotsc,N 
\! + \! 1 \rbrace$ is the number of bands to the right of $z \! = \! 0$. 
Furthermore, since in this manuscript $\mathcal{Y}_{q}^{\pm} \! = \! \mathbb{
C}_{\pm}$ for $q \! \in \! \{e,o\}$, $\infty^{\pm} \, \text{and} \, 0^{\pm} 
\! \in \! \mathbb{C}_{\pm}$. \hfill $\blacksquare$
\end{eeeee}

At this point, one is ready to use the preceding formulae and relations 
subsumed in Propositions~3.1--3.3 to establish large-$n$ asymptotic 
expansions for the various quantities of interest.

Using, now, $n \! \to \! \infty$ asymptotics for the $(1 \, 1)$-elements 
of $\operatorname{Y}^{e,0}_{0}$ and $\operatorname{Y}^{e,0}_{1}$ (resp., 
$\operatorname{Y}^{o,\infty}_{0}$ and $\operatorname{Y}^{o,\infty}_{1})$ 
given in Appendix~A, Theorem~A.2 (resp., Appendix~B, Theorem~B.2), one 
proves the following theorem concerning the products of the (real) roots 
of the $L$-polynomials: the asymptotic results are contained in 
Equations~(3.48)--(3.51).
\begin{ddddd}
Let the external field $\widetilde{V} \colon \mathbb{R} \setminus \{0\} \! \to 
\! \mathbb{R}$ satisfy conditions~{\rm (1.20)--(1.22)}. Let the orthonormal 
$L$-polynomials (resp., monic orthogonal $L$-polynomials) $\lbrace \phi_{k}(z) 
\rbrace_{k \in \mathbb{Z}_{0}^{+}}$ (resp., $\lbrace \boldsymbol{\pi}_{k}(z) 
\rbrace_{k \in \mathbb{Z}_{0}^{+}})$ be as defined in Equations~{\rm (1.2)} 
and~{\rm (1.3)} (resp., Equations~{\rm (1.4)} and~{\rm (1.5))}. Let 
$\overset{e}{\operatorname{Y}} \colon \mathbb{C} \setminus \mathbb{R} \! \to 
\! \operatorname{SL}_{2}(\mathbb{C})$ (resp., $\overset{o}{\operatorname{Y}} 
\colon \mathbb{C} \setminus \mathbb{R} \! \to \! \operatorname{SL}_{2}
(\mathbb{C}))$ be the (unique) solution of {\rm \pmb{RHP1}} (resp., 
{\rm \pmb{RHP2})} with integral representation {\rm (1.29)} (resp., 
{\rm (1.30))}. Define the density of the `even' (resp., `odd') equilibrium 
measure, $\md \mu_{V}^{e}(x)$ (resp., $\md \mu_{V}^{o}(x))$, as in 
Equation~{\rm (2.13)} (resp., Equation~{\rm (2.19))}, and set $J_{q} \! := 
\! \operatorname{supp}(\mu_{V}^{q}) \! = \! \cup_{j=1}^{N+1}(b_{j-1}^{q},
a_{j}^{q})$, $q \! \in \! \{e,o\}$, where $\lbrace b_{j-1}^{e},a_{j}^{e} 
\rbrace_{j=1}^{N+1}$ (resp., $\lbrace b_{j-1}^{o},a_{j}^{o} \rbrace_{j=1}^{N
+1})$ satisfy the (real) $n$-dependent and locally solvable system of $2(N 
\! + \! 1)$ moment conditions~{\rm (2.12)} (resp., {\rm (2.18))}.

Suppose, furthermore, that $\widetilde{V} \colon \mathbb{R} \setminus \{0\} \! 
\to \! \mathbb{R}$ is regular, namely: {\rm (i)} $h_{V}^{q}(z) \! \not\equiv 
\! 0$ on $\overline{J_{q}} \! := \! \cup_{k=1}^{N+1}[b_{k-1}^{q},a_{k}^{q}]$, 
$q \! \in \! \{e,o\};$ and {\rm (ii)} the variational 
inequalities~$(\pmb{\mathrm{P}_{1}^{(b)}})$ and~$(\pmb{\mathrm{P}_{2}^{(b)}})$ 
are strict, and take the following form,
\begin{gather}
4 \int_{J_{e}} \ln (\vert x \! - \! s \vert) \, \md \mu_{V}^{e}(s) \! - \! 2 
\ln \vert x \vert \! - \! \widetilde{V}(x) \! - \! \ell_{e} \! < \! 0, \quad 
x \! \in \! \mathbb{R} \setminus \overline{J_{e}}, \\
2 \! \left(2 \! + \! \dfrac{1}{n} \right) \! \int_{J_{o}} \ln (\vert x \! - \! 
s \vert) \, \md \mu_{V}^{o}(s) \! - \! 2 \ln \vert x \vert \! - \! \widetilde{
V}(x) \! - \! \ell_{o} \! - \! 2 \! \left(2 \! + \! \dfrac{1}{n} \right) \! 
Q_{o} \! < \! 0, \quad x \! \in \! \mathbb{R} \setminus \overline{J_{o}} \,,
\end{gather}
where the `even' and `odd' variational constants, $\ell_{e}$ $(\in \! \mathbb{
R})$ and $\ell_{o}$ $(\in \! \mathbb{R})$, respectively, are defined by 
Equations~$(\pmb{\mathrm{P}_{1}^{(a)}})$ and~$(\pmb{\mathrm{P}_{2}^{(a)}})$, 
and $Q_{o} \! := \! \int_{J_{o}} \ln (\lvert s \rvert) \psi_{V}^{o}(s) \, 
\md s$. Then,
\begin{align}
\nu^{(2n)}_{-n} &= \, \prod_{k=1}^{2n} \alpha^{(2n)}_{k} \underset{n \to 
\infty}{=} \, \dfrac{\varpi^{e}_{+}}{2} \dfrac{\boldsymbol{\theta}^{e}
(\boldsymbol{u}^{e}_{+}(\infty) \! + \! \boldsymbol{d}_{e}) \boldsymbol{
\theta}^{e}(\boldsymbol{u}^{e}_{+}(0) \! - \! \frac{n}{2 \pi} \boldsymbol{
\Omega}^{e} \! + \! \boldsymbol{d}_{e})}{\boldsymbol{\theta}^{e}(\boldsymbol{
u}^{e}_{+}(0) \! + \! \boldsymbol{d}_{e}) \boldsymbol{\theta}^{e}(\boldsymbol{
u}^{e}_{+}(\infty) \! - \! \frac{n}{2 \pi} \boldsymbol{\Omega}^{e} \! + \! 
\boldsymbol{d}_{e})} \nonumber \\
&\times \left\{1 \! + \! \dfrac{1}{n} \! \left((\mathscr{R}^{e,0}_{0})_{11} \! 
+ \! (\mathscr{R}^{e,0}_{0})_{12} \dfrac{\varpi^{e}_{-}}{\mi \varpi^{e}_{+}} 
\dfrac{\boldsymbol{\theta}^{e}(\boldsymbol{u}^{e}_{+}(\infty) \! - \! \frac{
n}{2 \pi} \boldsymbol{\Omega}^{e} \! + \! \boldsymbol{d}_{e})}{\boldsymbol{
\theta}^{e}(-\boldsymbol{u}^{e}_{+}(\infty) \! - \! \frac{n}{2 \pi} 
\boldsymbol{\Omega}^{e} \! - \! \boldsymbol{d}_{e})} \right. \right. \nonumber 
\\
&\left. \left. \times \, \dfrac{\boldsymbol{\theta}^{e}(\boldsymbol{u}^{e}_{+}
(0) \! - \! \frac{n}{2 \pi} \boldsymbol{\Omega}^{e} \! - \! \boldsymbol{d}_{e}
) \boldsymbol{\theta}^{e}(\boldsymbol{u}^{e}_{+}(0) \! + \! \boldsymbol{d}_{e}
)}{\boldsymbol{\theta}^{e}(\boldsymbol{u}^{e}_{+}(0) \! - \! \frac{n}{2 \pi} 
\boldsymbol{\Omega}^{e} \! + \! \boldsymbol{d}_{e}) \boldsymbol{\theta}^{e}
(\boldsymbol{u}^{e}_{+}(0) \! - \! \boldsymbol{d}_{e})} \right) \! + \! 
\mathcal{O} \! \left(\dfrac{c(n)}{n^{2}} \right) \right\} \nonumber \\
&\times \exp \! \left(2n \! \left(\int_{J_{e}} \ln (\lvert s \rvert) \psi_{
V}^{e}(s) \, \md s \! + \! \mi \pi \int_{J_{e} \cap \mathbb{R}_{+}} \psi_{
V}^{e}(s) \, \md s \right) \right),
\end{align}
\begin{align}
-\nu^{(2n+1)}_{n} &= \, \left(\prod_{k=1}^{2n+1} \alpha^{(2n+1)}_{k} \right)^{-
1} \underset{n \to \infty}{=} \, -\dfrac{1}{2} \dfrac{\varpi^{o}_{+}}{\mathbb{
E}} \dfrac{\boldsymbol{\theta}^{o}(\boldsymbol{u}^{o}_{+}(0) \! + \!
\boldsymbol{d}_{o}) \boldsymbol{\theta}^{o}(\boldsymbol{u}^{o}_{+}(\infty) \!
- \! \frac{1}{2 \pi}(n \! + \! \frac{1}{2}) \boldsymbol{\Omega}^{o} \! + \!
\boldsymbol{d}_{o})}{\boldsymbol{\theta}^{o}(\boldsymbol{u}^{o}_{+}(\infty) \!
+ \! \boldsymbol{d}_{o}) \boldsymbol{\theta}^{o}(\boldsymbol{u}^{o}_{+}(0) \!
- \! \frac{1}{2 \pi}(n \! + \! \frac{1}{2}) \boldsymbol{\Omega}^{o} \! + \!
\boldsymbol{d}_{o})} \nonumber \\
&\times \, \left\{1 \! + \! \dfrac{1}{(n \! + \! \frac{1}{2})} \! \left((
\mathscr{R}^{o,\infty}_{0})_{11} \! - \! (\mathscr{R}^{o,\infty}_{0})_{12} 
\dfrac{\varpi^{o}_{-} \mathbb{E}^{2}}{\mi \varpi^{o}_{+}} \dfrac{\boldsymbol{
\theta}^{o}(\boldsymbol{u}^{o}_{+}(\infty) \! - \! \frac{1}{2 \pi}(n \! + \! 
\frac{1}{2}) \boldsymbol{\Omega}^{o} \! - \! \boldsymbol{d}_{o})}{\boldsymbol{
\theta}^{o}(\boldsymbol{u}^{o}_{+}(\infty) \! - \! \frac{1}{2 \pi}(n \! + \! 
\frac{1}{2}) \boldsymbol{\Omega}^{o} \! + \! \boldsymbol{d}_{o})} \right. 
\right. \nonumber \\
&\left. \left. \times \dfrac{\boldsymbol{\theta}^{o}(\boldsymbol{u}^{o}_{+}
(\infty) \! + \! \boldsymbol{d}_{o}) \boldsymbol{\theta}^{o}(\boldsymbol{u}^{
o}_{+}(0) \! - \! \frac{1}{2 \pi}(n \! + \! \frac{1}{2}) \boldsymbol{\Omega}^{
o} \! + \! \boldsymbol{d}_{o})}{\boldsymbol{\theta}^{o}(\boldsymbol{u}^{o}_{+}
(\infty) \! - \! \boldsymbol{d}_{o}) \boldsymbol{\theta}^{o}(-\boldsymbol{u}^{
o}_{+}(0) \! - \! \frac{1}{2 \pi}(n \! + \! \frac{1}{2}) \boldsymbol{\Omega}^{
o} \! - \! \boldsymbol{d}_{o})} \right) \! + \! \mathcal{O} \! \left(\dfrac{
c(n)}{(n \! + \! \frac{1}{2})^{2}} \right) \right\} \nonumber \\
&\times \exp \! \left(\! -2 \! \left(n \! + \! \dfrac{1}{2} \right) \! 
\int_{J_{o}} \ln (\lvert s \rvert) \psi_{V}^{o}(s) \, \md s \right),
\end{align}
\begin{align}
-\nu^{(2n)}_{-(n-1)} &= \, \prod_{k=1}^{2n-1} \alpha^{(2n)}_{k} \underset{n 
\to \infty}{=} \, n \varpi^{e}_{+} \dfrac{\boldsymbol{\theta}^{e}(\boldsymbol{
u}^{e}_{+}(\infty) \! + \! \boldsymbol{d}_{e}) \boldsymbol{\theta}^{e}
(\boldsymbol{u}^{e}_{+}(0) \! - \! \frac{n}{2 \pi} \boldsymbol{\Omega}^{e} \! 
+ \! \boldsymbol{d}_{e})}{\boldsymbol{\theta}^{e}(\boldsymbol{u}^{e}_{+}(0) 
\! + \! \boldsymbol{d}_{e}) \boldsymbol{\theta}^{e}(\boldsymbol{u}^{e}_{+}
(\infty) \! - \! \frac{n}{2 \pi} \boldsymbol{\Omega}^{e} \! + \! \boldsymbol{
d}_{e})} \nonumber \\
&\times \left\{1 \! + \! \dfrac{1}{n} \! \left((\mathscr{R}^{e,0}_{0})_{11} \! 
+ \! (\mathscr{R}^{e,0}_{0})_{12} \dfrac{\varpi^{e}_{-}}{\mi \varpi^{e}_{+}} 
\dfrac{\boldsymbol{\theta}^{e}(\boldsymbol{u}^{e}_{+}(0) \! + \! \boldsymbol{
d}_{e}) \boldsymbol{\theta}^{e}(\boldsymbol{u}^{e}_{+}(0) \! - \! \frac{n}{2 
\pi} \boldsymbol{\Omega}^{e} \! - \! \boldsymbol{d}_{e})}{\boldsymbol{\theta}^{
e}(\boldsymbol{u}^{e}_{+}(0) \! - \! \boldsymbol{d}_{e}) \boldsymbol{\theta}^{
e}(\boldsymbol{u}^{e}_{+}(0) \! - \! \frac{n}{2 \pi} \boldsymbol{\Omega}^{e} 
\! + \! \boldsymbol{d}_{e})} \right. \right. \nonumber \\
&\left. \left. \times \, \dfrac{\boldsymbol{\theta}^{e}(\boldsymbol{u}^{e}_{+}
(\infty) \! - \! \frac{n}{2 \pi} \boldsymbol{\Omega}^{e} \! + \! \boldsymbol{
d}_{e})}{\boldsymbol{\theta}^{e}(-\boldsymbol{u}^{e}_{+}(\infty) \! - \! \frac{
n}{2 \pi} \boldsymbol{\Omega}^{e} \! - \! \boldsymbol{d}_{e})} \! + \! \dfrac{
1}{2} \! \left(\int_{J_{e}}s^{-1} \psi_{V}^{e}(s) \, \md s \right)^{-1} \! 
\left[\dfrac{\widetilde{\alpha}^{e}_{0}(1,1,\pmb{0})}{\boldsymbol{\theta}^{e}
(\boldsymbol{u}^{e}_{+}(0) \! + \! \boldsymbol{d}_{e})} \right. \right. 
\right. \nonumber \\
&\left. \left. \left. - \, \dfrac{\widetilde{\alpha}^{e}_{0}(1,1,\boldsymbol{
\Omega}^{e})}{\boldsymbol{\theta}^{e}(\boldsymbol{u}^{e}_{+}(0) \! - \! \frac{
n}{2 \pi} \boldsymbol{\Omega}^{e} \! + \! \boldsymbol{d}_{e})} \! + \! \dfrac{
\varpi^{e}_{-} \eta^{e}_{-}}{4 \varpi^{e}_{+}} \right] \right) \! + \! 
\mathcal{O} \! \left(\dfrac{c(n)}{n^{2}} \right) \right\} \! \left(\int_{J_{
e}}s^{-1} \psi_{V}^{e}(s) \, \md s \right) \nonumber \\
&\times \exp \! \left(2n \! \left(\int_{J_{e}} \ln (\lvert s \rvert) \psi_{
V}^{e}(s) \, \md s \! + \! \mi \pi \int_{J_{e} \cap \mathbb{R}_{+}} \psi_{
V}^{e}(s) \, \md s \right) \right),
\end{align}
\begin{align}
\nu^{(2n+1)}_{n-1} &= \left(\prod_{k=1}^{2n} \alpha^{(2n+1)}_{k} \right)^{-
1} \underset{n \to \infty}{=} \, -\left(n \! + \! \dfrac{1}{2} \right) \! 
\dfrac{\varpi^{o}_{+}}{\mathbb{E}} \dfrac{\boldsymbol{\theta}^{o}(\boldsymbol{
u}^{o}_{+}(0) \! + \! \boldsymbol{d}_{o}) \boldsymbol{\theta}^{o}(\boldsymbol{
u}^{o}_{+}(\infty) \! - \! \frac{1}{2 \pi}(n \! + \! \frac{1}{2}) \boldsymbol{
\Omega}^{o} \! + \! \boldsymbol{d}_{o})}{\boldsymbol{\theta}^{o}(\boldsymbol{
u}^{o}_{+}(\infty) \! + \! \boldsymbol{d}_{o}) \boldsymbol{\theta}^{o}
(\boldsymbol{u}^{o}_{+}(0) \! - \! \frac{1}{2 \pi}(n \! + \! \frac{1}{2}) 
\boldsymbol{\Omega}^{o} \! + \! \boldsymbol{d}_{o})} \nonumber \\
&\times \left\{1 \! + \! \dfrac{1}{(n \! + \! \frac{1}{2})} \! \left((\mathscr{
R}^{o,\infty}_{0})_{11} \! - \! (\mathscr{R}^{o,\infty}_{0})_{12} \dfrac{
\varpi^{o}_{-} \mathbb{E}^{2}}{\mi \varpi^{o}_{+}} \dfrac{\boldsymbol{
\theta}^{o}(\boldsymbol{u}^{o}_{+}(\infty) \! - \! \frac{1}{2 \pi}(n \! + \! 
\frac{1}{2}) \boldsymbol{\Omega}^{o} \! - \! \boldsymbol{d}_{o})}{\boldsymbol{
\theta}^{o}(\boldsymbol{u}^{o}_{+}(\infty) \! - \! \frac{1}{2 \pi}(n \! + \! 
\frac{1}{2}) \boldsymbol{\Omega}^{o} \! + \! \boldsymbol{d}_{o})} \right. 
\right. \nonumber \\
&\left. \left. \times \, \dfrac{\boldsymbol{\theta}^{o}(\boldsymbol{u}^{o}_{+}
(0) \! - \! \frac{1}{2 \pi}(n \! + \! \frac{1}{2}) \boldsymbol{\Omega}^{o} 
\! + \! \boldsymbol{d}_{o}) \boldsymbol{\theta}^{o}(\boldsymbol{u}^{o}_{+}
(\infty) \! + \! \boldsymbol{d}_{o})}{\boldsymbol{\theta}^{o}(-\boldsymbol{
u}^{o}_{+}(0) \! - \! \frac{1}{2 \pi}(n \! + \! \frac{1}{2}) \boldsymbol{
\Omega}^{o} \! - \! \boldsymbol{d}_{o}) \boldsymbol{\theta}^{o}(\boldsymbol{
u}^{o}_{+}(\infty) \! - \! \boldsymbol{d}_{o})} \! + \! \dfrac{1}{2} \! 
\left(\int_{J_{o}}s \psi_{V}^{o}(s) \, \md s \right)^{-1} \right. \right. 
\nonumber \\
&\left. \left. \times \left[\dfrac{\widetilde{\alpha}^{o}_{\infty}(1,1,
\boldsymbol{\Omega}^{o})}{\boldsymbol{\theta}^{o}(\boldsymbol{u}^{o}_{+}
(\infty) \! - \! \frac{1}{2 \pi}(n \! + \! \frac{1}{2}) \boldsymbol{\Omega}^{
o} \! + \! \boldsymbol{d}_{o})} \! - \! \dfrac{\widetilde{\alpha}^{o}_{\infty}
(1,1,\pmb{0})}{\boldsymbol{\theta}^{o}(\boldsymbol{u}^{o}_{+}(\infty) \! + \!
\boldsymbol{d}_{o})} \! + \! \dfrac{\varpi^{o}_{-} \eta^{o}_{+}}{4 \varpi^{
o}_{+}} \right] \right) \! + \! \mathcal{O} \! \left(\dfrac{c(n)}{(n \! + \! 
\frac{1}{2})^{2}} \right) \right\} \nonumber \\
&\times \left(\int_{J_{o}}s \psi_{V}^{o}(s) \, \md s \right) \! \exp \left(
\! -2 \! \left(n \! + \! \dfrac{1}{2} \right) \! \int_{J_{o}} \ln (\lvert s
\rvert) \psi_{V}^{o}(s) \, \md s \right).
\end{align}
In the above asymptotic expansions~{\rm (3.48)--(3.51)}, there are a number 
of new parameters, which are now defined:
\begin{gather}
\varpi^{q}_{\pm} \! := \! \gamma^{q}_{0} \! \pm \! (\gamma^{q}_{0})^{-1}, 
\quad q \! \in \! \{e,o\}, \\
\gamma^{q}_{0} \! := \! \gamma^{q}(0) \! = \! \left(\prod_{k=1}^{N+1}b_{k-
1}^{q}(a_{k}^{q})^{-1} \right)^{1/4} \quad (> \! 0), \\
\eta^{q}_{+} \! := \! \sum_{k=1}^{N+1} \! \left(a_{k}^{q} \! - \! b_{k-1}^{q} 
\right), \qquad \quad \eta^{q}_{-} \! := \! \sum_{k=1}^{N+1} \! \left((b_{k-
1}^{q})^{-1} \! - \! (a_{k}^{q})^{-1} \right),
\end{gather}
$\gamma^{e}(z)$ (resp., $\gamma^{o}(z))$ is defined in Equation~{\rm (2.2)} 
(resp., Equation~{\rm (2.7))}, $\boldsymbol{\theta}^{e}(z)$ (resp., 
$\boldsymbol{\theta}^{o}(z))$ is defined in Equation~{\rm (2.5)} (resp., 
Equation~{\rm (2.10))}, $\boldsymbol{u}^{q}(z) \! := \! \int_{a_{N+1}^{q}}^{
z} \boldsymbol{\omega}^{q}$, $q \! \in \! \{e,o\}$, $\boldsymbol{u}^{q}_{+}
(\infty) \! := \! \int_{a_{N+1}^{q}}^{\infty^{+}} \boldsymbol{\omega}^{q}$, 
$\boldsymbol{u}^{q}_{+}(0) \! := \! \int_{a_{N+1}^{q}}^{0^{+}} \boldsymbol{
\omega}^{q}$, $\boldsymbol{d}_{e}$ (resp., $\boldsymbol{d}_{o})$ is given in 
Equations~{\rm (2.4)} (resp., Equations~{\rm (2.9))}, with $z_{j}^{e,\pm} \! 
\in \! (a_{j}^{e},b_{j}^{e})^{\pm} \! \subset \! \mathbb{C}_{\pm}$ (resp., 
$z_{j}^{o,\pm} \! \in \! (a_{j}^{o},b_{j}^{o})^{\pm} \! \subset \! \mathbb{
C}_{\pm})$, $j \! = \! 1,\dotsc,N$, given in Equation~{\rm (2.3)} (resp., 
Equation~{\rm (2.8))}, $\boldsymbol{\Omega}^{q} \! := \! (\Omega^{q}_{1},
\Omega^{q}_{2},\dotsc,\Omega^{q}_{N})^{\mathrm{T}}$ $(\in \! \mathbb{R}^{N})$, 
with $\Omega^{q}_{j} \! = \! 4 \pi \int_{b_{j}^{q}}^{a_{N+1}^{q}} \psi_{V}^{q}
(s) \, \md s$, $j \! = \! 1,\dotsc,N$,
\begin{align}
\mathscr{R}^{q,p}_{0} :=& \, \operatorname{sgn}(q) \, \mathlarger{\sum_{j=1}^{
N+1}} \! \left(\dfrac{\left(\mathscr{A}^{q}(b_{j-1}^{q}) \! \left(\widehat{
\alpha}_{1}^{q}(b_{j-1}^{q}) \! + \! (b_{j-1}^{q})^{-1} \widehat{\alpha}_{
0}^{q}(b_{j-1}^{q}) \right) \! - \! \mathscr{B}^{q}(b_{j-1}^{q}) \widehat{
\alpha}_{0}^{q}(b_{j-1}^{q}) \right)}{(\widehat{\alpha}_{0}^{q}(b_{j-1}^{q})
)^{2}b_{j-1}^{q}} \right. \nonumber \\
+&\left. \, \dfrac{\left(\mathscr{A}^{q}(a_{j}^{q}) \! \left(\widehat{
\alpha}_{1}^{q}(a_{j}^{q}) \! + \! (a_{j}^{q})^{-1} \widehat{\alpha}_{0}^{q}
(a_{j}^{q}) \right) \! - \! \mathscr{B}^{q}(a_{j}^{q}) \widehat{\alpha}_{0}^{
q}(a_{j}^{q}) \right)}{(\widehat{\alpha}_{0}^{q}(a_{j}^{q}))^{2}a_{j}^{q}} 
\right), \qquad q \! \in \! \{e,o\}, \quad p \! \in \! \{0,\infty\},
\end{align}
\begin{equation}
\operatorname{sgn}(q) \! := \!
\begin{cases}
+1, &\text{$q \! = \! e$,} \\
-1, &\text{$q \! = \! o$,}
\end{cases}
\end{equation}
where, for $j \! = \! 1,\dotsc,N \! + \! 1$,
\begin{gather}
\mathscr{A}^{e}(b_{j-1}^{e}) \! = \! -s_{1}(Q_{0}^{e}(b_{j-1}^{e}))^{-1} 
\me^{\mi n \mho_{j-1}^{e}} \! 
\begin{pmatrix}
\varkappa^{e}_{1}(b_{j-1}^{e}) \varkappa^{e}_{2}(b_{j-1}^{e}) & \mi 
(\varkappa^{e}_{1}(b_{j-1}^{e}))^{2} \\
\mi (\varkappa^{e}_{2}(b_{j-1}^{e}))^{2} & 
-\varkappa^{e}_{1}(b_{j-1}^{e}) \varkappa^{e}_{2}(b_{j-1}^{e})
\end{pmatrix}, \\
\mathscr{A}^{e}(a_{j}^{e}) \! = \! s_{1}Q_{0}^{e}(a_{j}^{e}) \me^{\mi n 
\mho_{j}^{e}} \! 
\begin{pmatrix}
-\varkappa^{e}_{1}(a_{j}^{e}) \varkappa^{e}_{2}(a_{j}^{e}) & \mi 
(\varkappa^{e}_{1}(a_{j}^{e}))^{2} \\
\mi (\varkappa^{e}_{2}(a_{j}^{e}))^{2} & \varkappa^{e}_{1}(a_{j}^{e}) 
\varkappa^{e}_{2}(a_{j}^{e})
\end{pmatrix}, \\
\mathscr{A}^{o}(b_{j-1}^{o}) \! = \! -\dfrac{s_{1}(\gamma^{o}(0))^{2} \me^{
\mi (n+\frac{1}{2}) \mho_{j-1}^{o}}}{Q_{0}^{o}(b_{j-1}^{o})} \! 
\begin{pmatrix}
\varkappa^{o}_{1}(b_{j-1}^{o}) \varkappa^{o}_{2}(b_{j-1}^{o}) & \mi 
(\varkappa^{o}_{1}(b_{j-1}^{o}))^{2} \\
\mi (\varkappa^{o}_{2}(b_{j-1}^{o}))^{2} & -\varkappa^{o}_{1}(b_{j-1}^{o}) 
\varkappa^{o}_{2}(b_{j-1}^{o})
\end{pmatrix}, \\
\mathscr{A}^{o}(a_{j}^{o}) \! = \! \dfrac{s_{1}Q_{0}^{o}(a_{j}^{o}) \me^{\mi 
(n+\frac{1}{2}) \mho_{j}^{o}}}{(\gamma^{o}(0))^{2}} \! 
\begin{pmatrix}
-\varkappa^{o}_{1}(a_{j}^{o}) \varkappa^{o}_{2}(a_{j}^{o}) & \mi (\varkappa^{
o}_{1}(a_{j}^{o}))^{2} \\
\mi (\varkappa^{o}_{2}(a_{j}^{o}))^{2} & \varkappa^{o}_{1}(a_{j}^{o}) 
\varkappa^{o}_{2}(a_{j}^{o})
\end{pmatrix}, \\
\dfrac{\mathscr{B}^{e}(b_{j-1}^{e})}{\me^{\mi n \mho_{j-1}^{e}}} \! = \! 
\begin{pmatrix}
\boxed{\begin{matrix} \varkappa_{1}^{e}(b_{j-1}^{e}) \varkappa_{2}^{e}
(b_{j-1}^{e}) \! \left(-s_{1}(Q_{0}^{e}(b_{j-1}^{e}))^{-1} \right. \\
\left. \times \left\{
\overset{e}{\daleth}^{\raise-1.0ex\hbox{$\scriptstyle 1$}}_{1}(b_{j-1}^{e}) 
\! + \!
\overset{e}{\daleth}^{\raise-1.0ex\hbox{$\scriptstyle 1$}}_{-1}(b_{j-1}^{e}) 
\! - \! Q_{1}^{e}(b_{j-1}^{e}) \right. \right. \\
\left. \left. \times \, (Q_{0}^{e}(b_{j-1}^{e}))^{-1} \right\} \! - \! t_{1} 
\! \left\{Q_{0}^{e}(b_{j-1}^{e}) \right. \right. \\
\left. \left. + \, (Q_{0}^{e}(b_{j-1}^{e}))^{-1} 
\overset{e}{\aleph}^{\raise-1.0ex\hbox{$\scriptstyle 1$}}_{1}(b_{j-1}^{e}) 
\overset{e}{\aleph}^{\raise-1.0ex\hbox{$\scriptstyle 1$}}_{-1} (b_{j-1}^{e}) 
\right\} \right. \\
\left. + \, \mi (s_{1} \! + \! t_{1}) \! \left\{
\overset{e}{\aleph}^{\raise-1.0ex\hbox{$\scriptstyle 1$}}_{-1}(b_{j-1}^{e}) 
\! - \! \overset{e}{\aleph}^{\raise-1.0ex\hbox{$\scriptstyle 1$}}_{1}(b_{j-
1}^{e}) \right\} \right)
\end{matrix}} & 
\boxed{\begin{matrix} (\varkappa_{1}^{e}(b_{j-1}^{e}))^{2} \! \left(-\mi s_{1}
(Q_{0}^{e}(b_{j-1}^{e}))^{-1} \! \left\{2
\overset{e}{\daleth}^{\raise-1.0ex\hbox{$\scriptstyle 1$}}_{1}(b_{j-1}^{e}) 
\right. \right. \\
\left. \left. -\, Q_{1}^{e}(b_{j-1}^{e})(Q_{0}^{e}(b_{j-1}^{e}))^{-1}\right\} 
\! + \! \mi t_{1} \! \left\{Q_{0}^{e}(b_{j-1}^{e}) \right. \right. \\
\left. \left. -\, (Q_{0}^{e}(b_{j-1}^{e}))^{-1}
(\overset{e}{\aleph}^{\raise-1.0ex\hbox{$\scriptstyle 1$}}_{1}(b_{j-1}^{e})
)^{2} \right\} \right. \\
\left. + \, 2(s_{1} \! - \! t_{1}) 
\overset{e}{\aleph}^{\raise-1.0ex\hbox{$\scriptstyle 1$}}_{1}(b_{j-1}^{e}) 
\right)
\end{matrix}} \\
\boxed{\begin{matrix} (\varkappa_{2}^{e}(b_{j-1}^{e}))^{2} \! \left(-\mi s_{1}
(Q_{0}^{e}(b_{j-1}^{e}))^{-1} \! \left\{2
\overset{e}{\daleth}^{\raise-1.0ex\hbox{$\scriptstyle 1$}}_{-1}(b_{j-1}^{e}) 
\right. \right. \\
\left. \left. -\, Q_{1}^{e}(b_{j-1}^{e})(Q_{0}^{e}(b_{j-1}^{e}))^{-1} 
\right\} \! + \! \mi t_{1} \! \left\{Q_{0}^{e}(b_{j-1}^{e}) \right. \right. \\
\left. \left. -\, (Q_{0}^{e}(b_{j-1}^{e}))^{-1}
(\overset{e}{\aleph}^{\raise-1.0ex\hbox{$\scriptstyle 1$}}_{-1}(b_{j-1}^{e}
))^{2} \right\} \right. \\
\left. - \, 2(s_{1} \! - \! t_{1}) 
\overset{e}{\aleph}^{\raise-1.0ex\hbox{$\scriptstyle 1$}}_{-1}(b_{j-1}^{e}) 
\right)
\end{matrix}} & 
\boxed{\begin{matrix} \varkappa_{1}^{e}(b_{j-1}^{e}) \varkappa_{2}^{e}
(b_{j-1}^{e}) \! \left(s_{1}(Q_{0}^{e}(b_{j-1}^{e}))^{-1} \right. \\
\left. \times \left\{
\overset{e}{\daleth}^{\raise-1.0ex\hbox{$\scriptstyle 1$}}_{1}(b_{j-1}^{e}) 
\! + \! 
\overset{e}{\daleth}^{\raise-1.0ex\hbox{$\scriptstyle 1$}}_{-1}(b_{j-1}^{e}) 
\! - \! Q_{1}^{e}(b_{j-1}^{e}) \right. \right. \\
\left. \left. \times \, (Q_{0}^{e}(b_{j-1}^{e}))^{-1} \right\} \! + \! t_{1} 
\! \left\{Q_{0}^{e}(b_{j-1}^{e}) \right. \right. \\
\left. \left. + \, (Q_{0}^{e}(b_{j-1}^{e}))^{-1} 
\overset{e}{\aleph}^{\raise-1.0ex\hbox{$\scriptstyle 1$}}_{1}(b_{j-1}^{e}) 
\overset{e}{\aleph}^{\raise-1.0ex\hbox{$\scriptstyle 1$}}_{-1} (b_{j-1}^{e}) 
\right\} \right. \\
\left. + \, \mi (s_{1} \! + \! t_{1}) \! \left\{
\overset{e}{\aleph}^{\raise-1.0ex\hbox{$\scriptstyle 1$}}_{1}(b_{j-1}^{e}) 
\! - \! \overset{e}{\aleph}^{\raise-1.0ex\hbox{$\scriptstyle 1$}}_{-1}
(b_{j-1}^{e}) \right\} \right)
\end{matrix}}
\end{pmatrix}, \\
\dfrac{\mathscr{B}^{e}(a_{j}^{e})}{\me^{\mi n \mho_{j}^{e}}} \! = \! 
\begin{pmatrix}
\boxed{\begin{matrix} \varkappa_{1}^{e}(a_{j}^{e}) \varkappa_{2}^{e}
(a_{j}^{e}) \! \left(-s_{1} \! \left\{Q_{1}^{e}(a_{j}^{e}) \right. \right. \\
\left. \left. + \, Q_{0}^{e}(a_{j}^{e}) \left[
\overset{e}{\daleth}^{\raise-1.0ex\hbox{$\scriptstyle 1$}}_{1}(a_{j}^{e}) \! 
+ \! \overset{e}{\daleth}^{\raise-1.0ex\hbox{$\scriptstyle 1$}}_{-1}(a_{j}^{
e}) \right] \right\} \! - \! t_{1} \right. \\
\left. \times \left\{(Q_{0}^{e}(a_{j}^{e}))^{-1} \! + \! Q_{0}^{e}(a_{j}^{e}) 
\overset{e}{\aleph}^{\raise-1.0ex\hbox{$\scriptstyle 1$}}_{1}(a_{j}^{e}) 
\overset{e}{\aleph}^{\raise-1.0ex\hbox{$\scriptstyle 1$}}_{-1}(a_{j}^{e}) 
\right\} \right. \\
\left. + \, \mi (s_{1} \! + \! t_{1}) \! \left\{
\overset{e}{\aleph}^{\raise-1.0ex\hbox{$\scriptstyle 1$}}_{-1}(a_{j}^{e}) \! 
- \! \overset{e}{\aleph}^{\raise-1.0ex\hbox{$\scriptstyle 1$}}_{1}(a_{j}^{e}) 
\right\} \right)
\end{matrix}} & 
\boxed{\begin{matrix} (\varkappa_{1}^{e}(a_{j}^{e}))^{2} \! \left(\mi s_{1} 
\! \left\{Q_{1}^{e}(a_{j}^{e}) \! + \! 2Q_{0}^{e}(a_{j}^{e}) \right. \right. 
\\
\left. \left. \times \, 
\overset{e}{\daleth}^{\raise-1.0ex\hbox{$\scriptstyle 1$}}_{1}(a_{j}^{e}) 
\right\} \! + \! \mi t_{1} \! \left\{Q_{0}^{e}(a_{j}^{e})(\overset{e}{
\aleph}^{\raise-1.0ex\hbox{$\scriptstyle 1$}}_{1}(a_{j}^{e}))^{2} \right. 
\right. \\
\left. \left. - \, (Q_{0}^{e}(a_{j}^{e}))^{-1} \right\} \! - \! 2 (s_{1} \! 
- \! t_{1}) \overset{e}{\aleph}^{\raise-1.0ex\hbox{$\scriptstyle 1$}}_{1}
(a_{j}^{e}) \right)
\end{matrix}} \\
\boxed{\begin{matrix} (\varkappa_{2}^{e}(a_{j}^{e}))^{2} \! \left(\mi s_{1} 
\! \left\{Q_{1}^{e}(a_{j}^{e}) \! + \! 2Q_{0}^{e}(a_{j}^{e}) \right. \right.
\\
\left. \left. \times \,
\overset{e}{\daleth}^{\raise-1.0ex\hbox{$\scriptstyle 1$}}_{-1}(a_{j}^{e}) 
\right\} \! + \! \mi t_{1} \! \left\{Q_{0}^{e}(a_{j}^{e})
(\overset{e}{\aleph}^{\raise-1.0ex\hbox{$\scriptstyle 1$}}_{-1}(a_{j}^{e}))^{
2}\right. \right. \\
\left. \left. - \, (Q_{0}^{e}(a_{j}^{e}))^{-1} \right\} \! + \! 2(s_{1} \! 
- \! t_{1}) \overset{e}{\aleph}^{\raise-1.0ex\hbox{$\scriptstyle 1$}}_{-1}
(a_{j}^{e}) \right)
\end{matrix}} & 
\boxed{\begin{matrix} \varkappa_{1}^{e}(a_{j}^{e})\varkappa_{2}^{e}(a_{j}^{e}) 
\! \left(s_{1} \! \left\{Q_{1}^{e}(a_{j}^{e}) \right. \right. \\
\left. \left. + \, Q_{0}^{e}(a_{j}^{e}) \left[
\overset{e}{\daleth}^{\raise-1.0ex\hbox{$\scriptstyle 1$}}_{1}(a_{j}^{e}) \! 
+ \! \overset{e}{\daleth}^{\raise-1.0ex\hbox{$\scriptstyle 1$}}_{-1}(a_{j}^{
e}) \right] \right\} \! + \! t_{1} \right. \\
\left. \times \left\{(Q_{0}^{e}(a_{j}^{e}))^{-1} \! + \! Q_{0}^{e}(a_{j}^{e}) 
\overset{e}{\aleph}^{\raise-1.0ex\hbox{$\scriptstyle 1$}}_{1}(a_{j}^{e}) 
\overset{e}{\aleph}^{\raise-1.0ex\hbox{$\scriptstyle 1$}}_{-1}(a_{j}^{e}) 
\right\} \right. \\
\left. + \, \mi (s_{1} \! + \! t_{1}) \! \left\{
\overset{e}{\aleph}^{\raise-1.0ex\hbox{$\scriptstyle 1$}}_{1}(a_{j}^{e}) \! - 
\! \overset{e}{\aleph}^{\raise-1.0ex\hbox{$\scriptstyle 1$}}_{-1}(a_{j}^{e}) 
\right\} \right)
\end{matrix}}
\end{pmatrix}, \\
\dfrac{\mathscr{B}^{o}(b_{j-1}^{o})}{\me^{\mi (n+\frac{1}{2}) \mho_{j-1}^{o}}} 
\! = \!
\begin{pmatrix}
\boxed{\begin{matrix} \varkappa_{1}^{o}(b_{j-1}^{o}) \varkappa_{2}^{o}(b_{j-
1}^{o}) \! \left(-\frac{s_{1}(\gamma^{o}(0))^{2}}{Q_{0}^{o}(b_{j-1}^{o})} 
\right. \\
\left. \times \left\{
\overset{o}{\daleth}^{\raise-1.0ex\hbox{$\scriptstyle 1$}}_{1}(b_{j-1}^{o}) \! 
+ \! \overset{o}{\daleth}^{\raise-1.0ex\hbox{$\scriptstyle 1$}}_{-1}(b_{j-1}^{
o}) \! - \! Q_{1}^{o}(b_{j-1}^{o}) \right. \right. \\
\left. \left. \times \, (Q_{0}^{o}(b_{j-1}^{o}))^{-1} \right\} \! - \! t_{1}
(\gamma^{o}(0))^{2} \! \left\{Q_{0}^{o}(b_{j-1}^{o}) \right. \right. \\
\left. \left. + \, (Q_{0}^{o}(b_{j-1}^{o}))^{-1} 
\overset{o}{\aleph}^{\raise-1.0ex\hbox{$\scriptstyle 1$}}_{1}(b_{j-1}^{o}) 
\overset{o}{\aleph}^{\raise-1.0ex\hbox{$\scriptstyle 1$}}_{-1}(b_{j-1}^{o}) 
\right\} \right. \\
\left. + \, \mi (s_{1} \! + \! t_{1}) \! \left\{
\overset{o}{\aleph}^{\raise-1.0ex\hbox{$\scriptstyle 1$}}_{-1}(b_{j-1}^{o}) 
\! - \! \overset{o}{\aleph}^{\raise-1.0ex\hbox{$\scriptstyle 1$}}_{1}(b_{j-
1}^{o}) \right\} \right)
\end{matrix}} & 
\boxed{\begin{matrix} (\varkappa_{1}^{o}(b_{j-1}^{o}))^{2} \! \left(-\frac{
\mi s_{1}(\gamma^{o}(0))^{2}}{Q_{0}^{o}(b_{j-1}^{o})} \! \left\{2 
\overset{o}{\daleth}^{\raise-1.0ex\hbox{$\scriptstyle 1$}}_{1}(b_{j-1}^{o}) 
\right. \right. \\
\left. \left. -\, Q_{1}^{o}(b_{j-1}^{o})(Q_{0}^{o}(b_{j-1}^{o}))^{-1} \right\} 
\! + \! \mi t_{1} \! \left\{Q_{0}^{o}(b_{j-1}^{o}) \right. \right. \\
\left. \left. \times \, (\gamma^{o}(0))^{-2} \! - \! (Q_{0}^{o}(b_{j-1}^{o})
)^{-1}(\gamma^{o}(0))^{2} \right. \right. \\
\left. \left. \times \, 
(\overset{o}{\aleph}^{\raise-1.0ex\hbox{$\scriptstyle 1$}}_{1}(b_{j-1}^{o})
)^{2} \right\} \! + \! 2(s_{1} \! - \! t_{1}) 
\overset{o}{\aleph}^{\raise-1.0ex\hbox{$\scriptstyle 1$}}_{1}(b_{j-1}^{o}) 
\right)
\end{matrix}} \\
\boxed{\begin{matrix} (\varkappa_{2}^{o}(b_{j-1}^{o}))^{2} \! \left(-\frac{
\mi s_{1}(\gamma^{o}(0))^{2}}{Q_{0}^{o}(b_{j-1}^{o})} \! \left\{2 
\overset{o}{\daleth}^{\raise-1.0ex\hbox{$\scriptstyle 1$}}_{-1}(b_{j-1}^{o}) 
\right. \right. \\
\left. \left. -\, Q_{1}^{o}(b_{j-1}^{o})(Q_{0}^{o}(b_{j-1}^{o}))^{-1} \right\} 
\! + \! \mi t_{1} \! \left\{Q_{0}^{o}(b_{j-1}^{o}) \right. \right. \\
\left. \left. \times \, (\gamma^{o}(0))^{-2} \! - \! (Q_{0}^{o}(b_{j-1}^{o})
)^{-1}(\gamma^{o}(0))^{2} \right. \right. \\
\left. \left. \times \, 
(\overset{o}{\aleph}^{\raise-1.0ex\hbox{$\scriptstyle 1$}}_{-1}(b_{j-1}^{o})
)^{2} \right\} \!- \! 2(s_{1} \! - \! t_{1}) 
\overset{o}{\aleph}^{\raise-1.0ex\hbox{$\scriptstyle 1$}}_{-1}(b_{j-1}^{o}) 
\right)
\end{matrix}} & 
\boxed{\begin{matrix} \varkappa_{1}^{o}(b_{j-1}^{o}) \varkappa_{2}^{o}(b_{j-
1}^{o}) \! \left(\frac{s_{1}(\gamma^{o}(0))^{2}}{Q_{0}^{o}(b_{j-1}^{o})} 
\right. \\
\left. \times \left\{
\overset{o}{\daleth}^{\raise-1.0ex\hbox{$\scriptstyle 1$}}_{1}(b_{j-1}^{o}) 
\! + \! \overset{o}{\daleth}^{\raise-1.0ex\hbox{$\scriptstyle 1$}}_{-1}(b_{j
-1}^{o}) \! - \! Q_{1}^{o}(b_{j-1}^{o}) \right. \right. \\
\left. \left. \times \, (Q_{0}^{o}(b_{j-1}^{o}))^{-1} \right\} \! + \! t_{1}
(\gamma^{o}(0))^{2} \! \left\{Q_{0}^{o}(b_{j-1}^{o}) \right. \right. \\
\left. \left. + \, (Q_{0}^{o}(b_{j-1}^{o}))^{-1}
\overset{o}{\aleph}^{\raise-1.0ex\hbox{$\scriptstyle 1$}}_{1}(b_{j-1}^{o}) 
\overset{o}{\aleph}^{\raise-1.0ex\hbox{$\scriptstyle 1$}}_{-1}(b_{j-1}^{o}) 
\right\} \right. \\
\left. + \, \mi (s_{1} \! + \! t_{1}) \! \left\{
\overset{o}{\aleph}^{\raise-1.0ex\hbox{$\scriptstyle 1$}}_{1}(b_{j-1}^{o}) 
\! - \! \overset{o}{\aleph}^{\raise-1.0ex\hbox{$\scriptstyle 1$}}_{-1}(b_{j
-1}^{o}) \right\} \right)
\end{matrix}}
\end{pmatrix}, \\
\dfrac{\mathscr{B}^{o}(a_{j}^{o})}{\me^{\mi (n+\frac{1}{2}) \mho_{j}^{o}}} \! 
= \!
\begin{pmatrix}
\boxed{\begin{matrix} \varkappa_{1}^{o}(a_{j}^{o}) \varkappa_{2}^{o}(a_{j}^{
o}) \! \left(-\frac{s_{1}}{(\gamma^{o}(0))^{2}} \! \left\{Q_{1}^{o}(a_{j}^{
o}) \right. \right. \\
\left. \left. + \, Q_{0}^{o}(a_{j}^{o}) \left[
\overset{o}{\daleth}^{\raise-1.0ex\hbox{$\scriptstyle 1$}}_{1}(a_{j}^{o}) \! 
+ \! \overset{o}{\daleth}^{\raise-1.0ex\hbox{$\scriptstyle 1$}}_{-1}(a_{j}^{
o}) \right] \right\} \! - \! t_{1} \right. \\
\left. \times \left\{(\gamma^{o}(0))^{-2}Q_{0}^{o}(a_{j}^{o}) 
\overset{o}{\aleph}^{\raise-1.0ex\hbox{$\scriptstyle 1$}}_{1}(a_{j}^{o}) 
\overset{o}{\aleph}^{\raise-1.0ex\hbox{$\scriptstyle 1$}}_{-1}(a_{j}^{o}) 
\right. \right. \\
\left. \left. + \, (\gamma^{o}(0))^{2}(Q_{0}^{o}(a_{j}^{o}))^{-1} \right\} 
\right. \\
\left. + \, \mi (s_{1} \! + \! t_{1}) \! \left\{
\overset{o}{\aleph}^{\raise-1.0ex\hbox{$\scriptstyle 1$}}_{-1}(a_{j}^{o}) \! 
- \! \overset{o}{\aleph}^{\raise-1.0ex\hbox{$\scriptstyle 1$}}_{1}(a_{j}^{o}) 
\right\} \right)
\end{matrix}} & 
\boxed{\begin{matrix} (\varkappa_{1}^{o}(a_{j}^{o}))^{2} \! \left(\frac{\mi 
s_{1}}{(\gamma^{o}(0))^{2}} \! \left\{Q_{1}^{o}(a_{j}^{o}) \! + \! 2Q_{0}^{o}
(a_{j}^{o}) \right. \right. \\
\left. \left. \times \, 
\overset{o}{\daleth}^{\raise-1.0ex\hbox{$\scriptstyle 1$}}_{1}(a_{j}^{o}) 
\right\} \! + \! \mi t_{1} \! \left\{Q_{0}^{o}(a_{j}^{o})
(\overset{o}{\aleph}^{\raise-1.0ex\hbox{$\scriptstyle 1$}}_{1}(a_{j}^{o}))^{
2} \right. \right. \\
\left. \left. \times \, (\gamma^{o}(0))^{-2} \! - \! (Q_{0}^{o}(a_{j}^{o}))^{
-1}(\gamma^{o}(0))^{2} \right\} \right. \\
\left. - \, 2(s_{1} \! - \! t_{1}) 
\overset{o}{\aleph}^{\raise-1.0ex\hbox{$\scriptstyle 1$}}_{1}(a_{j}^{o}) 
\right)
\end{matrix}} \\
\boxed{\begin{matrix} (\varkappa_{2}^{o}(a_{j}^{o}))^{2} \! \left(\frac{\mi 
s_{1}}{(\gamma^{o}(0))^{2}} \! \left\{Q_{1}^{o}(a_{j}^{o}) \! + \! 2Q_{0}^{o}
(a_{j}^{o}) \right. \right. \\
\left. \left. \times \, 
\overset{o}{\daleth}^{\raise-1.0ex\hbox{$\scriptstyle 1$}}_{-1}(a_{j}^{o}) 
\right\} \! + \! \mi t_{1} \! \left\{Q_{0}^{o}(a_{j}^{o})
(\overset{o}{\aleph}^{\raise-1.0ex\hbox{$\scriptstyle 1$}}_{-1}(a_{j}^{o}))^{
2} \right. \right. \\
\left. \left. \times \, (\gamma^{o}(0))^{-2} \! - \! (Q_{0}^{o}(a_{j}^{o}))^{
-1}(\gamma^{o}(0))^{2} \right\} \right. \\
\left. + \, 2(s_{1} \! - \! t_{1}) 
\overset{o}{\aleph}^{\raise-1.0ex\hbox{$\scriptstyle 1$}}_{-1}(a_{j}^{o}) 
\right)
\end{matrix}} &
\boxed{\begin{matrix} \varkappa_{1}^{o}(a_{j}^{o}) \varkappa_{2}^{o}(a_{j}^{
o}) \! \left(\frac{s_{1}}{(\gamma^{o}(0))^{2}} \! \left\{Q_{1}^{o}(a_{j}^{o}) 
\right. \right. \\
\left. \left. + \, Q_{0}^{o}(a_{j}^{o}) \left[
\overset{o}{\daleth}^{\raise-1.0ex\hbox{$\scriptstyle 1$}}_{1}(a_{j}^{o}) \! 
+ \! \overset{o}{\daleth}^{\raise-1.0ex\hbox{$\scriptstyle 1$}}_{-1}(a_{j}^{
o}) \right] \right\} \! + \! t_{1} \right. \\
\left. \times \left\{(\gamma^{o}(0))^{-2}Q_{0}^{o}(a_{j}^{o}) 
\overset{o}{\aleph}^{\raise-1.0ex\hbox{$\scriptstyle 1$}}_{1}(a_{j}^{o}) 
\overset{o}{\aleph}^{\raise-1.0ex\hbox{$\scriptstyle 1$}}_{-1}(a_{j}^{o}) 
\right. \right. \\
\left. \left. + \, (\gamma^{o}(0))^{2}(Q_{0}^{o}(a_{j}^{o}))^{-1} \right\} 
\right. \\
\left. + \, \mi (s_{1} \! + \! t_{1}) \! \left\{
\overset{o}{\aleph}^{\raise-1.0ex\hbox{$\scriptstyle 1$}}_{1}(a_{j}^{o}) \! - 
\! \overset{o}{\aleph}^{\raise-1.0ex\hbox{$\scriptstyle 1$}}_{-1}(a_{j}^{o}) 
\right\} \right)
\end{matrix}}
\end{pmatrix}, \\
\mathbb{E} \! := \! \exp \! \left(\mi 2 \pi \! \left(n \! + \! \dfrac{1}{2} 
\right) \! \int_{J_{o} \cap \mathbb{R}_{+}} \psi_{V}^{o}(s) \, \md s \right),
\end{gather}
and, for $q \! \in \! \{e,o\}$,
\begin{gather}
s_{1} \! = \! \dfrac{5}{72}, \qquad \qquad \quad t_{1} \! = \! -\dfrac{7}{72}, 
\qquad \qquad \quad \mho_{k}^{q} \! = \! 
\begin{cases}
\Omega_{k}^{q}, &\text{$k \! = \! 1,\dotsc,N$,} \\
0, &\text{$k \! = \! 0,N \! + \! 1$,}
\end{cases} \\
Q_{0}^{q}(b_{0}^{q}) \! = \! -\mi \! \left(\! (a_{N+1}^{q} \! - \! b_{0}^{q})^{
-1} \prod_{k=1}^{N} \! \left(\dfrac{b_{k}^{q} \! - \! b_{0}^{q}}{a_{k}^{q} \!
- \! b_{0}^{q}} \right) \right)^{1/2}, \\
Q_{1}^{q}(b_{0}^{q}) \! = \! \dfrac{1}{2}Q_{0}^{q}(b_{0}^{q}) \! \left(\sum_{
k=1}^{N} \! \left(\dfrac{1}{b_{0}^{q} \! - \! b_{k}^{q}} \! - \! \dfrac{1}{b_{
0}^{q} \! - \! a_{k}^{q}} \right) \! - \! \dfrac{1}{b_{0}^{q} \! - \! a_{N+
1}^{q}} \right), \\
Q_{0}^{q}(a_{N+1}^{q}) \! = \! \left(\! (a_{N+1}^{q} \! - \! b_{0}^{q}) \prod_{
k=1}^{N} \! \left(\dfrac{a_{N+1}^{q} \! - \! b_{k}^{q}}{a_{N+1}^{q} \! - \!
a_{k}^{q}} \right) \right)^{1/2}, \\
Q_{1}^{q}(a_{N+1}^{q}) \! = \! \dfrac{1}{2}Q_{0}^{q}(a_{N+1}^{q}) \! \left(
\sum_{k=1}^{N} \! \left(\dfrac{1}{a_{N+1}^{q} \! - \! b_{k}^{q}} \! - \!
\dfrac{1}{a_{N+1}^{q} \! - \! a_{k}^{q}} \right) \! + \! \dfrac{1}{a_{N+1}^{q}
\! - \! b_{0}^{q}} \right), \\
Q_{0}^{q}(b_{j}^{q}) \! = \! -\mi \! \left(\! \dfrac{(b_{j}^{q} \! - \! b_{0}^{
q})}{(a_{N+1}^{q} \! - \! b_{j}^{q})(b_{j}^{q} \! - \! a_{j}^{q})} \prod_{k=
1}^{j-1} \! \left(\dfrac{b_{j}^{q} \! - \! b_{k}^{q}}{b_{j}^{q} \! - \! a_{k}^{
q}} \right) \! \prod_{l=j+1}^{N} \! \left(\! \dfrac{b_{l}^{q} \! - \! b_{j}^{q}
}{a_{l}^{q} \! - \! b_{j}^{q}} \right) \right)^{1/2}, \\
Q_{1}^{q}(b_{j}^{q}) \! = \! \dfrac{1}{2}Q_{0}^{q}(b_{j}^{q}) \! \left(\sum_{
\substack{k=1\\k \not= j}}^{N} \! \left(\! \dfrac{1}{b_{j}^{q} \! - \! b_{k}^{
q}} \! - \! \dfrac{1}{b_{j}^{q} \! - \! a_{k}^{q}} \right) \! + \! \dfrac{1}{
b_{j}^{q} \! - \! b_{0}^{q}} \! - \! \dfrac{1}{b_{j}^{q} \! - \! a_{N+1}^{q}}
\! - \! \dfrac{1}{b_{j}^{q} \! - \! a_{j}^{q}} \right), \\
Q_{0}^{q}(a_{j}^{q}) \! = \! \left(\! \dfrac{(a_{j}^{q} \! - \! b_{0}^{q})(b_{
j}^{q} \! - \! a_{j}^{q})}{(a_{N+1}^{q} \! - \! a_{j}^{q})} \prod_{k=1}^{j-1}
\! \left(\dfrac{a_{j}^{q} \! - \! b_{k}^{q}}{a_{j}^{q} \! - \! a_{k}^{q}}
\right) \! \prod_{l=j+1}^{N} \! \left(\! \dfrac{b_{l}^{q} \! - \! a_{j}^{q}}{
a_{l}^{q} \! - \! a_{j}^{q}} \right) \right)^{1/2}, \\
Q_{1}^{q}(a_{j}^{q}) \! = \! \dfrac{1}{2}Q_{0}^{q}(a_{j}^{q}) \! \left(
\sum_{\substack{k=1\\k \not= j}}^{N} \! \left(\! \dfrac{1}{a_{j}^{q} \! - \!
b_{k}^{q}} \! - \! \dfrac{1}{a_{j}^{q} \! - \! a_{k}^{q}} \right) \! + \!
\dfrac{1}{a_{j}^{q} \! - \! b_{0}^{q}} \! - \! \dfrac{1}{a_{j}^{q} \! - \!
a_{N+1}^{q}} \! + \! \dfrac{1}{a_{j}^{q} \! - \! b_{j}^{q}} \right),
\end{gather}
with $\mi Q_{0}^{q}(b_{j-1}^{q}),Q_{0}^{q}(a_{j}^{q}) \! > \! 0$, $j \! = \! 
1,\dotsc,N \! + \! 1$,
\begin{align}
\varkappa_{1}^{e}(\xi) &= \, \dfrac{\bm{\theta}^{e}(\bm{u}^{e}_{+}(\infty) 
\! + \! \bm{d}_{e}) \bm{\theta}^{e}(\bm{u}^{e}_{+}(\xi) \! - \! \frac{n}{2 
\pi} \bm{\Omega}^{e} \! + \! \bm{d}_{e})}{\bm{\theta}^{e}(\bm{u}^{e}_{+}
(\infty) \! - \! \frac{n}{2 \pi} \bm{\Omega}^{e} \! + \! \bm{d}_{e}) \bm{
\theta}^{e}(\bm{u}^{e}_{+}(\xi) \! + \! \bm{d}_{e})}, \\
\varkappa_{2}^{e}(\xi) &= \, \dfrac{\bm{\theta}^{e}(-\bm{u}^{e}_{+}(\infty) 
\! - \! \bm{d}_{e}) \bm{\theta}^{e}(\bm{u}^{e}_{+}(\xi) \! - \! \frac{n}{2 
\pi} \bm{\Omega}^{e} \! - \! \bm{d}_{e})}{\bm{\theta}^{e}(-\bm{u}^{e}_{+}
(\infty) \! - \! \frac{n}{2 \pi} \bm{\Omega}^{e} \! - \! \bm{d}_{e}) \bm{
\theta}^{e}(\bm{u}^{e}_{+}(\xi) \! - \! \bm{d}_{e})}, \\
\varkappa_{1}^{o}(\xi) &= \, \dfrac{1}{\mathbb{E}} \dfrac{\bm{\theta}^{o}
(\bm{u}^{o}_{+}(0) \! + \! \bm{d}_{o}) \bm{\theta}^{o}(\bm{u}^{o}_{+}(\xi) \! 
- \! \frac{1}{2 \pi}(n \! + \! \frac{1}{2}) \bm{\Omega}^{o} \! + \! \bm{d}_{
o})}{\bm{\theta}^{o}(\bm{u}^{o}_{+}(0) \! - \! \frac{1}{2 \pi}(n \! + \! 
\frac{1}{2}) \bm{\Omega}^{o} \! + \! \bm{d}_{o}) \bm{\theta}^{o}(\bm{u}^{o}_{
+}(\xi) \! + \! \bm{d}_{o})}, \\
\varkappa_{2}^{o}(\xi) &= \, \mathbb{E} \dfrac{\bm{\theta}^{o}(-\bm{u}^{o}_{
+}(0) \! - \! \bm{d}_{o}) \bm{\theta}^{o}(\bm{u}^{o}_{+}(\xi) \! - \! \frac{
1}{2 \pi}(n \! + \! \frac{1}{2}) \bm{\Omega}^{o} \! - \! \bm{d}_{o})}{\bm{
\theta}^{o}(-\bm{u}^{o}_{+}(0) \! - \! \frac{1}{2 \pi}(n \! + \! \frac{1}{2}) 
\bm{\Omega}^{o} \! - \! \bm{d}_{o}) \bm{\theta}^{o}(\bm{u}^{o}_{+}(\xi) \! - 
\! \bm{d}_{o})}, \\
\overset{e}{\aleph}^{\raise-1.0ex\hbox{$\scriptstyle \varepsilon_{1}$}}_{
\varepsilon_{2}}(\xi) &= \, -\dfrac{\mathfrak{u}^{e}(\varepsilon_{1},
\varepsilon_{2},\bm{0};\xi)}{\bm{\theta}^{e}(\varepsilon_{1} \bm{u}^{e}_{+}
(\xi) \! + \! \varepsilon_{2} \bm{d}_{e})} \! + \! \dfrac{\mathfrak{u}^{e}
(\varepsilon_{1},\varepsilon_{2},\bm{\Omega}^{e};\xi)}{\bm{\theta}^{e}
(\varepsilon_{1} \bm{u}^{e}_{+}(\xi) \! - \! \frac{n}{2 \pi} \bm{\Omega}^{e} 
\! + \! \varepsilon_{2} \bm{d}_{e})}, \quad \varepsilon_{1},\varepsilon_{2} 
\! = \! \pm 1, \\
\overset{o}{\aleph}^{\raise-1.0ex\hbox{$\scriptstyle \varepsilon_{1}$}}_{
\varepsilon_{2}}(\xi) &= \, -\dfrac{\mathfrak{u}^{o}(\varepsilon_{1},
\varepsilon_{2},\bm{0};\xi)}{\bm{\theta}^{o}(\varepsilon_{1} \bm{u}^{o}_{+}
(\xi) \! + \! \varepsilon_{2} \bm{d}_{o})} \! + \! \dfrac{\mathfrak{u}^{o}
(\varepsilon_{1},\varepsilon_{2},\bm{\Omega}^{o};\xi)}{\bm{\theta}^{o}
(\varepsilon_{1} \bm{u}^{o}_{+}(\xi) \! - \! \frac{1}{2 \pi}(n \! + \! 
\frac{1}{2}) \bm{\Omega}^{o} \! + \! \varepsilon_{2} \bm{d}_{o})}, \\
\overset{e}{\daleth}^{\raise-1.0ex\hbox{$\scriptstyle \varepsilon_{1}$}}_{
\varepsilon_{2}}(\xi) &= \, -\dfrac{\mathfrak{v}^{e}(\varepsilon_{1},
\varepsilon_{2},\bm{0};\xi)}{\bm{\theta}^{e}(\varepsilon_{1} \bm{u}^{e}_{+}
(\xi) \! + \! \varepsilon_{2} \bm{d}_{e})} \! + \! \dfrac{\mathfrak{v}^{e}
(\varepsilon_{1},\varepsilon_{2},\bm{\Omega}^{e};\xi)}{\bm{\theta}^{e}
(\varepsilon_{1} \bm{u}^{e}_{+}(\xi) \! - \! \frac{n}{2 \pi} \bm{\Omega}^{e} 
\! + \! \varepsilon_{2} \bm{d}_{e})} \! - \! \left(\dfrac{\mathfrak{u}^{e}
(\varepsilon_{1},\varepsilon_{2},\bm{0};\xi)}{\bm{\theta}^{e}(\varepsilon_{
1} \bm{u}^{e}_{+}(\xi) \! + \! \varepsilon_{2} \bm{d}_{e})} \right)^{2} 
\nonumber \\
&+ \, \dfrac{\mathfrak{u}^{e}(\varepsilon_{1},\varepsilon_{2},\bm{0};\xi) 
\mathfrak{u}^{e}(\varepsilon_{1},\varepsilon_{2},\bm{\Omega}^{e};\xi)}{\bm{
\theta}^{e}(\varepsilon_{1} \bm{u}^{e}_{+}(\xi) \! + \! \varepsilon_{2} \bm{
d}_{e}) \bm{\theta}^{e}(\varepsilon_{ 1} \bm{u}^{e}_{+}(\xi) \! - \! \frac{
n}{2 \pi} \bm{\Omega}^{e} \! + \! \varepsilon_{2} \bm{d}_{e})}, \\
\overset{o}{\daleth}^{\raise-1.0ex\hbox{$\scriptstyle \varepsilon_{1}$}}_{
\varepsilon_{2}}(\xi) &= \, -\dfrac{\mathfrak{v}^{o}(\varepsilon_{1},
\varepsilon_{2},\bm{0};\xi)}{\bm{\theta}^{o}(\varepsilon_{1} \bm{u}^{o}_{+}
(\xi) \! + \! \varepsilon_{2} \bm{d}_{o})} \! + \! \dfrac{\mathfrak{v}^{o}
(\varepsilon_{1},\varepsilon_{2},\bm{\Omega}^{o};\xi)}{\bm{\theta}^{o}
(\varepsilon_{1} \bm{u}^{o}_{+}(\xi) \! - \! \frac{1}{2 \pi}(n \! + \! 
\frac{1}{2}) \bm{\Omega}^{o} \! + \! \varepsilon_{2} \bm{d}_{o})} \! - \! 
\left(\dfrac{\mathfrak{u}^{o}(\varepsilon_{1},\varepsilon_{2},\bm{0};\xi)}{
\bm{\theta}^{o}(\varepsilon_{1} \bm{u}^{o}_{+}(\xi) \! + \! \varepsilon_{2} 
\bm{d}_{o})} \right)^{2} \nonumber \\
&+ \, \dfrac{\mathfrak{u}^{o}(\varepsilon_{1},\varepsilon_{2},\bm{0};\xi) 
\mathfrak{u}^{o}(\varepsilon_{1},\varepsilon_{2},\bm{\Omega}^{o};\xi)}{\bm{
\theta}^{o}(\varepsilon_{1} \bm{u}^{o}_{+}(\xi) \! + \! \varepsilon_{2} 
\bm{d}_{o}) \bm{\theta}^{o}(\varepsilon_{1} \bm{u}^{o}_{+}(\xi) \! - \! 
\frac{1}{2 \pi}(n \! + \! \frac{1}{2}) \bm{\Omega}^{o} \! + \! \varepsilon_{
2} \bm{d}_{o})},
\end{align}
\begin{gather}
\mathfrak{u}^{q}(\varepsilon_{1},\varepsilon_{2},\bm{\Omega}^{q},\xi) \! := 
\! 2 \pi \Lambda^{\raise-1.0ex\hbox{$\scriptstyle 1$}}_{q}(\varepsilon_{1},
\varepsilon_{2},\bm{\Omega}^{q},\xi), \qquad \quad \mathfrak{v}^{q}
(\varepsilon_{1},\varepsilon_{2},\bm{\Omega}^{q},\xi) \! := \! -2 \pi^{2} 
\Lambda^{\raise-1.0ex\hbox{$\scriptstyle 2$}}_{q}(\varepsilon_{1},
\varepsilon_{2},\bm{\Omega}^{q},\xi), \\
\Lambda^{\raise-1.0ex\hbox{$\scriptstyle j^{\prime}$}}_{e}(\varepsilon_{1},
\varepsilon_{2},\bm{\Omega}^{e},\xi) \! = \! \sum_{m \in \mathbb{Z}^{N}}
(\mathfrak{r}_{e}(\xi))^{j^{\prime}} \me^{2 \pi \mi (m,\varepsilon_{1} \bm{
u}^{e}_{+}(\xi)-\frac{n}{2 \pi} \bm{\Omega}^{e}+\varepsilon_{2} \bm{d}_{e}) 
+ \pi \mi (m,\bm{\tau}^{e}m)}, \quad j^{\prime} \! = \! 1,2, \\
\Lambda^{\raise-1.0ex\hbox{$\scriptstyle j^{\prime}$}}_{o}(\varepsilon_{1},
\varepsilon_{2},\bm{\Omega}^{o},\xi) \! = \! \sum_{m \in \mathbb{Z}^{N}}
(\mathfrak{r}_{o}(\xi))^{j^{\prime}} \me^{2 \pi \mi (m,\varepsilon_{1} 
\bm{u}^{o}_{+}(\xi)-\frac{1}{2 \pi}(n+\frac{1}{2}) \bm{\Omega}^{o}+
\varepsilon_{2} \bm{d}_{o})+ \pi \mi(m,\bm{\tau}^{o}m)}, \\
\mathfrak{r}_{q}(\xi) \! := \! \dfrac{2(m,\vec{\moo}_{q}(\xi))}{
\leftthreetimes^{\raise+0.3ex\hbox{$\scriptstyle q$}}(\xi)}, \qquad \qquad 
\quad \vec{\moo}_{q}(\xi) \! = \! 
\left(\rightthreetimes^{\raise-0.9ex\hbox{$\scriptstyle q$}}_{1}(\xi),
\rightthreetimes^{\raise-0.9ex\hbox{$\scriptstyle q$}}_{2}(\xi),
\dotsc,\rightthreetimes^{\raise-0.9ex\hbox{$\scriptstyle q$}}_{N}(\xi) 
\right), \\
\rightthreetimes^{\raise-0.9ex\hbox{$\scriptstyle q$}}_{j^{\prime}}(\xi) \! 
:= \! \sum_{k=1}^{N} c_{j^{\prime}k}^{q} \xi^{N-k}, \quad j^{\prime} \! = \! 
1,\dotsc,N, \\
\leftthreetimes^{\raise+0.3ex\hbox{$\scriptstyle q$}}(b_{0}^{q}) \! = \! \mi 
(-1)^{N} \eta_{b_{0}^{q}}, \quad 
\leftthreetimes^{\raise+0.3ex\hbox{$\scriptstyle q$}}(a_{N+1}^{q}) \! = \! 
\eta_{a_{N+1}^{q}}, \quad 
\leftthreetimes^{\raise+0.3ex\hbox{$\scriptstyle q$}}(b_{j}^{q}) \! = \! \mi 
(-1)^{N-j} \eta_{b_{j}^{q}}, \quad 
\leftthreetimes^{\raise+0.3ex\hbox{$\scriptstyle q$}}(a_{j}^{q}) \! = \! 
(-1)^{N-j+1} \eta_{a_{j}^{q}}, \\
\eta_{b_{0}^{q}} \! := \! \left((a_{N+1}^{q} \! - \! b_{0}^{q}) \prod_{k=1}^{N}
(b_{k}^{q} \! - \! b_{0}^{q})(a_{k}^{q} \! - \! b_{0}^{q}) \right)^{1/2}, \\
\eta_{a_{N+1}^{q}} \! := \! \left((a_{N+1}^{q} \! - \! b_{0}^{q}) \prod_{k=1
}^{N}(a_{N+1}^{q} \! - \! b_{k}^{q})(a_{N+1}^{q} \! - \! a_{k}^{q}) \right)^{
1/2}, \\
\eta_{b_{j}^{q}} \! := \! \left(\! (b_{j}^{q} \! - \! a_{j}^{q})(a_{N+1}^{q} \!
- \! b_{j}^{q})(b_{j}^{q} \! - \! b_{0}^{q}) \prod_{k=1}^{j-1}(b_{j}^{q} \! -
\! b_{k}^{q})(b_{j}^{q} \! - \! a_{k}^{q}) \prod_{l=j+1}^{N}(b_{l}^{q} \! - \!
b_{j}^{q})(a_{l}^{q} \! - \! b_{j}^{q}) \! \right)^{1/2}, \\
\eta_{a_{j}^{q}} \! := \! \left(\! (b_{j}^{q} \! - \! a_{j}^{q})(a_{N+1}^{q} \!
- \! a_{j}^{q})(a_{j}^{q} \! - \! b_{0}^{q}) \prod_{k=1}^{j-1}(a_{j}^{q} \! -
\! b_{k}^{q})(a_{j}^{q} \! - \! a_{k}^{q}) \prod_{l=j+1}^{N}(b_{l}^{q} \! - \!
a_{j}^{q})(a_{l}^{q} \! - \! a_{j}^{q}) \! \right)^{1/2},
\end{gather}
where $c^{e}_{j^{\prime}k^{\prime}}$ (resp., $c^{o}_{j^{\prime}k^{\prime}})$, 
$j^{\prime},k^{\prime} \! = \! 1,\dotsc,N$, are obtained {}from 
Equations~{\rm (E1)} and~{\rm (E2)} (resp., 
Equations~{\rm (O1)} and~{\rm (O2))}, $\eta_{b_{j-1}^{q}},\eta_{a_{j}^{q}} \! 
> \! 0$, $j \! = \! 1,\dotsc,N \! + \! 1$, and
\begin{align}
\widehat{\alpha}^{q}_{0}(b_{0}^{q}) \! =& \dfrac{4}{3} \mi (-1)^{N}h_{V}^{q}
(b_{0}^{q}) \eta_{b_{0}^{q}}, \\
\widehat{\alpha}^{q}_{1}(b_{0}^{q}) \! =& \mi (-1)^{N} \! \left(\dfrac{2}{5}
h_{V}^{q}(b_{0}^{q}) \eta_{b_{0}^{q}} \! \left(\sum_{l=1}^{N} \! \left(\dfrac{
1}{b_{0}^{q} \! - \! b_{l}^{q}} \! + \! \dfrac{1}{b_{0}^{q} \! - \! a_{l}^{q}}
\right) \! + \! \dfrac{1}{b_{0}^{q} \! - \! a_{N+1}^{q}} \right) \! + \!
\dfrac{4}{5}(h_{V}^{q}(b_{0}^{q}))^{\prime} \eta_{b_{0}^{q}} \right), \\
\widehat{\alpha}^{q}_{0}(a_{N+1}^{q}) \! =& \dfrac{4}{3}h_{V}^{q}(a_{N+1}^{q})
\eta_{a_{N+1}^{q}}, \\
\widehat{\alpha}^{q}_{1}(a_{N+1}^{q}) \! =& \dfrac{2}{5}h_{V}^{q}(a_{N+1}^{q})
\eta_{a_{N+1}^{q}} \! \left(\sum_{l=1}^{N} \! \left(\dfrac{1}{a_{N+1}^{q} \! -
\! b_{l}^{q}} \! + \! \dfrac{1}{a_{N+1}^{q} \! - \! a_{l}^{q}} \right) \! + \!
\dfrac{1}{a_{N+1}^{q} \! - \! b_{0}^{q}} \right) \! + \! \dfrac{4}{5}(h_{V}^{
q}(a_{N+1}^{q}))^{\prime} \eta_{a_{N+1}^{q}}, \\
\widehat{\alpha}^{q}_{0}(b_{j}^{q}) \! =& \dfrac{4}{3} \mi (-1)^{N-j}h_{V}^{q}
(b_{j}^{q}) \eta_{b_{j}^{q}}, \\
\widehat{\alpha}^{q}_{1}(b_{j}^{q}) \! =& \mi (-1)^{N-j} \! \left(\dfrac{2}{5}
h_{V}^{q}(b_{j}^{q}) \eta_{b_{j}^{q}} \! \left(\sum_{\substack{k=1\\k \not= j}
}^{N} \! \left(\dfrac{1}{b_{j}^{q} \! - \! b_{k}^{q}} \! + \! \dfrac{1}{b_{j}^{
q} \! - \! a_{k}^{q}} \right) \! + \! \dfrac{1}{b_{j}^{q} \! - \! a_{j}^{q}}
\! + \! \dfrac{1}{b_{j}^{q} \! - \! b_{0}^{q}} \! + \! \dfrac{1}{b_{j}^{q} \!
- \! a_{N+1}^{q}} \right) \! + \! \dfrac{4}{5}(h_{V}^{q}(b_{j}^{q}))^{\prime} 
\eta_{b_{j}^{q}} \right), \\
\widehat{\alpha}^{q}_{0}(a_{j}^{q}) \! =& \dfrac{4}{3}(-1)^{N-j+1}h_{V}^{q}
(a_{j}^{q}) \eta_{a_{j}^{q}}, \\
\widehat{\alpha}^{q}_{1}(a_{j}^{q}) \! =& (-1)^{N-j+1} \! \left(\dfrac{2}{5}
h_{V}^{q}(a_{j}^{q}) \eta_{a_{j}^{q}} \! \left(\sum_{\substack{k=1\\k \not= j}
}^{N} \! \left(\dfrac{1}{a_{j}^{q} \! - \! b_{k}^{q}} \! + \! \dfrac{1}{a_{j}^{
q} \! - \! a_{k}^{q}} \right) \! + \! \dfrac{1}{a_{j}^{q} \! - \! b_{j}^{q}}
\! + \! \dfrac{1}{a_{j}^{q} \! - \! a_{N+1}^{q}} \! + \! \dfrac{1}{a_{j}^{q}
\! - \! b_{0}^{q}} \right) \! + \! \dfrac{4}{5}(h_{V}^{q}(a_{j}^{q}))^{\prime} 
\eta_{a_{j}^{q}} \right),
\end{align}
and
\begin{gather}
\widetilde{\alpha}^{e}_{0}(1,1,\pmb{\mathscr{Z}}) \! := \! 2 \pi \mi \sum_{m 
\in \mathbb{Z}^{N}}(m,\widehat{\boldsymbol{\alpha}}^{e}_{0}) \me^{2 \pi \mi 
(m,\boldsymbol{u}^{e}_{+}(0)-\frac{n}{2 \pi} \pmb{\mathscr{Z}}+\boldsymbol{
d}_{e})+ \pi \mi (m,\tau^{e}m)}, \\
\widetilde{\alpha}^{o}_{\infty}(1,1,\pmb{\mathscr{Z}}) \! := \! 2 \pi \mi 
\sum_{m \in \mathbb{Z}^{N}}(m,\widehat{\boldsymbol{\alpha}}^{o}_{\infty}) 
\me^{2 \pi \mi (m,\boldsymbol{u}^{o}_{+}(\infty)-\frac{1}{2 \pi}(n+\frac{
1}{2}) \pmb{\mathscr{Z}}+\boldsymbol{d}_{o})+ \pi \mi (m,\tau^{o}m)},
\end{gather}
where $\widehat{\boldsymbol{\alpha}}^{q}_{p} \! = \! (\widehat{\alpha}_{p,1}^{
q},\widehat{\alpha}^{q}_{p,2},\dotsc,\widehat{\alpha}^{q}_{p,N})$, $q \! \in 
\! \{e,o\}$, $p \! \in \! \{0,\infty\}$, with
\begin{equation}
\widehat{\alpha}^{e}_{0,j} \! := \! (-1)^{\mathcal{N}_{+}^{e}} \! \left(
\prod_{i=1}^{N+1} \lvert b_{i-1}^{e}a_{i}^{e} \rvert \right)^{-1/2}c_{jN}^{
e}, \qquad \qquad \widehat{\alpha}^{o}_{\infty,j} \! := \! c_{j1}^{o}, \quad 
j \! = \! 1,\dotsc,N,
\end{equation}
and $\pmb{0} \! := \! (0,0,\dotsc,0)^{\mathrm{T}}$ $(\in \! \mathbb{R}^{N})$.
\end{ddddd}

Via Proposition~3.3, the formulae
\begin{equation*}
b_{2n}^{\sharp} \! := \! h^{-} \! \left[b_{2n+2}^{\sharp} \right], \qquad 
c_{2n}^{\sharp} \! := \! h^{-} \! \left[c_{2n+2}^{\sharp} \right], \qquad 
\beta_{2n}^{\sharp} \! := \! h^{-} \! \left[\beta_{2n+2}^{\sharp} \right], 
\qquad \gamma_{2n+1}^{\sharp} \! := \! h^{-} \! \left[\gamma_{2n+3}^{\sharp} 
\right],
\end{equation*}
and Equations~(3.35) and~(3.36), one proves the following theorem for 
asymptotics (as $n \! \to \! \infty)$ of the coefficients of the recurrence 
relations~(1.23)--(1.26). Theorem~3.2 (see below) presents asymptotics 
for $a_{2n+1}^{\sharp}$ (see Equation~(3.104)), $b_{2n+1}^{\sharp}$ 
(see Equation~(3.106)), $b_{2n+2}^{\sharp}$ (see Equation~(3.109)), 
$c_{2n+2}^{\sharp}$ (see Equation~(3.114)), $\alpha_{2n}^{\sharp}$ (see 
Equation~(3.115)), $\beta_{2n+1}^{\sharp}$ (see Equation~(3.116)), 
$\beta_{2n+2}^{\sharp}$ (see Equation~(3.117)), $\gamma_{2n+3}^{\sharp}$ 
(see Equation~(3.118)), $b_{2n}^{\sharp}$ (see Equation~(3.119)), 
$c_{2n}^{\sharp}$ (see Equation~(3.120)), $\beta_{2n}^{\sharp}$ (see 
Equation~(3.121)), and $\gamma_{2n+1}^{\sharp}$ (see Equation~(3.122)). 
Furthermore, asymptotics for $\xi^{(2n)}_{n-1}/\xi^{(2n)}_{n}$ (cf. 
Equations~(3.17)) and $\xi^{(2n+1)}_{-n}/\xi^{(2n+1)}_{-n-1}$ (cf. 
Equations~(3.18)) are presented in Equations~(3.123) and~(3.124), 
respectively.
\begin{ddddd}
Let the external field $\widetilde{V} \colon \mathbb{R} \setminus \{0\} 
\! \to \! \mathbb{R}$ satisfy conditions~{\rm (1.20)--(1.22)}. Let the 
orthonormal $L$-polynomials (resp., monic orthogonal $L$-polynomials) 
$\lbrace \phi_{k}(z) \rbrace_{k \in \mathbb{Z}_{0}^{+}}$ (resp., $\lbrace 
\boldsymbol{\pi}_{k}(z) \rbrace_{k \in \mathbb{Z}_{0}^{+}})$ be as defined 
in Equations~{\rm (1.2)} and~{\rm (1.3)} (resp., Equations~{\rm (1.4)} 
and~{\rm (1.5))}, and let $\lbrace \phi_{k}(z) \rbrace_{k \in \mathbb{Z}_{
0}^{+}}$ satisfy the system of recurrence relations~{\rm (1.23)--(1.26)}. 
Let $\overset{e}{\operatorname{Y}} \colon \mathbb{C} \setminus \mathbb{R} \! 
\to \! \operatorname{SL}_{2}(\mathbb{C})$ (resp., $\overset{o}{\operatorname{
Y}} \colon \mathbb{C} \setminus \mathbb{R} \! \to \! \operatorname{SL}_{2}
(\mathbb{C}))$ be the (unique) solution of {\rm \pmb{RHP1}} (resp., 
{\rm \pmb{RHP2})} with integral representation {\rm (1.29)} (resp., 
{\rm (1.30))}. Define the density of the `even' (resp., `odd') equilibrium 
measure, $\md \mu_{V}^{e}(x)$ (resp., $\md \mu_{V}^{o}(x))$, as in 
Equation~{\rm (2.13)} (resp., Equation~{\rm (2.19))}, and set $J_{q} \! := 
\! \operatorname{supp}(\mu_{V}^{q}) \! = \! \cup_{j=1}^{N+1}(b_{j-1}^{q},
a_{j}^{q})$, $q \! \in \! \{e,o\}$, where $\lbrace b_{j-1}^{e},a_{j}^{e} 
\rbrace_{j=1}^{N+1}$ (resp., $\lbrace b_{j-1}^{o},a_{j}^{o} \rbrace_{j=1}^{
N+1})$ satisfy the (real) $n$-dependent and locally solvable system of 
$2(N \! + \! 1)$ moment conditions~{\rm (2.12)} (resp., {\rm (2.18))}.

Suppose, furthermore, that $\widetilde{V} \colon \mathbb{R} \setminus \{0\} 
\! \to \! \mathbb{R}$ is regular, namely: {\rm (i)} $h_{V}^{q}(z) \! 
\not\equiv \! 0$ on $\overline{J_{q}} \! := \! \cup_{k=1}^{N+1}[b_{k-1}^{q},
a_{k}^{q}]$, $q \! \in \! \{e,o\};$ and {\rm (ii)} the strict variational 
inequalities of Theorem~{\rm 3.1}, Equations~{\rm (3.46)} and~{\rm (3.47)}, 
are valid. Then,
\begin{align}
a_{2n+1}^{\sharp} \underset{n \to \infty}{=}& \, -\dfrac{\varpi^{o}_{+} 
\varpi^{o}_{-}}{\eta^{o}_{-}} \dfrac{\boldsymbol{\theta}^{o}(\boldsymbol{u}^{
o}_{+}(\infty) \! - \! \frac{1}{2 \pi}(n \! + \! \frac{1}{2}) \boldsymbol{
\Omega}^{o} \! + \! \boldsymbol{d}_{o}) \boldsymbol{\theta}^{o}(-\boldsymbol{
u}^{o}_{+}(\infty) \! - \! \frac{1}{2 \pi}(n \! + \! \frac{1}{2}) \boldsymbol{
\Omega}^{o} \! + \! \boldsymbol{d}_{o})}{\boldsymbol{\theta}^{o}(\boldsymbol{
u}^{o}_{+}(0) \! - \! \frac{1}{2}(n \! + \! \frac{1}{2}) \boldsymbol{\Omega}^{
o} \! + \! \boldsymbol{d}_{o}) \boldsymbol{\theta}^{o}(-\boldsymbol{u}^{o}_{+}
(0) \! - \! \frac{1}{2 \pi}(n \! + \! \frac{1}{2}) \boldsymbol{\Omega}^{o} \! 
+ \! \boldsymbol{d}_{o})} \nonumber \\
&\times \dfrac{\boldsymbol{\theta}^{o}(\boldsymbol{u}^{o}_{+}(0) \! + \! 
\boldsymbol{d}_{o}) \boldsymbol{\theta}^{o}(-\boldsymbol{u}^{o}_{+}(0) \! + \! 
\boldsymbol{d}_{o})}{\boldsymbol{\theta}^{o}(\boldsymbol{u}^{o}_{+}(\infty) 
\! + \! \boldsymbol{d}_{o}) \boldsymbol{\theta}^{o}(-\boldsymbol{u}^{o}_{+}
(\infty) \! + \! \boldsymbol{d}_{o})} \! \left\{1 \! + \! \dfrac{1}{(n \! + \! 
\frac{1}{2})} \! \left((\mathscr{R}^{o,\infty}_{0})_{11} \! - \! (\mathscr{
R}^{o,\infty}_{0})_{12} \dfrac{\varpi^{o}_{-} \mathbb{E}^{2}}{\mi \varpi^{
o}_{+}} \right. \right. \nonumber \\
&\left. \left. \times \, \dfrac{\boldsymbol{\theta}^{o}(\boldsymbol{u}^{o}_{+}
(\infty) \! + \! \boldsymbol{d}_{o}) \boldsymbol{\theta}^{o}(\boldsymbol{u}^{
o}_{+}(0) \! - \! \frac{1}{2 \pi}(n \! + \! \frac{1}{2}) \boldsymbol{\Omega}^{
o} \! + \! \boldsymbol{d}_{o}) \boldsymbol{\theta}^{o}(\boldsymbol{u}^{o}_{+}
(\infty) \! - \! \frac{1}{2 \pi}(n \! + \! \frac{1}{2}) \boldsymbol{\Omega}^{
o} \! - \! \boldsymbol{d}_{o})}{\boldsymbol{\theta}^{o}(\boldsymbol{u}^{o}_{+}
(\infty) \! - \! \boldsymbol{d}_{o}) \boldsymbol{\theta}^{o}(-\boldsymbol{u}^{
o}_{+}(0) \! - \! \frac{1}{2 \pi}(n \! + \! \frac{1}{2}) \boldsymbol{\Omega}^{
o} \! - \! \boldsymbol{d}_{o}) \boldsymbol{\theta}^{o}(\boldsymbol{u}^{o}_{+}
(\infty) \! - \! \frac{1}{2 \pi}(n \! + \! \frac{1}{2}) \boldsymbol{\Omega}^{
o} \! + \! \boldsymbol{d}_{o})} \right. \right. \nonumber \\
&\left. \left. + \, (\mathscr{R}^{o,\infty}_{0})_{11} \! - \! (\mathscr{R}^{o,
\infty}_{0})_{12} \dfrac{\varpi^{o}_{+} \mathbb{E}^{2}}{\mi \varpi^{o}_{-}} 
\dfrac{\boldsymbol{\theta}^{o}(-\boldsymbol{u}^{o}_{+}(\infty) \! + \! 
\boldsymbol{d}_{o}) \boldsymbol{\theta}^{o}(-\boldsymbol{u}^{o}_{+}(\infty) 
\! - \! \frac{1}{2 \pi}(n \! + \! \frac{1}{2}) \boldsymbol{\Omega}^{o} \! - \! 
\boldsymbol{d}_{o})}{\boldsymbol{\theta}^{o}(-\boldsymbol{u}^{o}_{+}(\infty) 
\! - \! \boldsymbol{d}_{o}) \boldsymbol{\theta}^{o}(-\boldsymbol{u}^{o}_{+}
(\infty) \! - \! \frac{1}{2 \pi}(n \! + \! \frac{1}{2}) \boldsymbol{\Omega}^{
o} \! + \! \boldsymbol{d}_{o})} \right. \right. \nonumber \\
&\left. \left. \times \, \dfrac{\boldsymbol{\theta}^{o}(\boldsymbol{u}^{o}_{+}
(0) \! - \! \frac{1}{2 \pi}(n \! + \! \frac{1}{2}) \boldsymbol{\Omega}^{o} \! 
+ \! \boldsymbol{d}_{o})}{\boldsymbol{\theta}^{o}(-\boldsymbol{u}^{o}_{+}(0) 
\! - \! \frac{1}{2 \pi}(n \! + \! \frac{1}{2}) \boldsymbol{\Omega}^{o} \! - 
\! \boldsymbol{d}_{o})} \! - \! \dfrac{4 \mi (\mathscr{R}^{o,0}_{1})_{12} 
\mathbb{E}^{2}}{\eta^{o}_{-}} \dfrac{\boldsymbol{\theta}^{o}(\boldsymbol{u}^{
o}_{+}(0) \! - \! \frac{1}{2 \pi}(n \! + \! \frac{1}{2}) \boldsymbol{\Omega}^{
o} \! + \! \boldsymbol{d}_{o})}{\boldsymbol{\theta}^{o}(-\boldsymbol{u}^{o}_{
+}(0) \! - \! \frac{1}{2 \pi}(n \! + \! \frac{1}{2}) \boldsymbol{\Omega}^{o} 
\! + \! \boldsymbol{d}_{o})} \right. \right. \nonumber \\
&\left. \left. \times \, \dfrac{\boldsymbol{\theta}^{o}(-\boldsymbol{u}^{o}_{
+}(0) \! + \! \boldsymbol{d}_{o})}{\boldsymbol{\theta}^{o}(\boldsymbol{u}^{
o}_{+}(0) \! + \! \boldsymbol{d}_{o})} \right) \! + \! \mathcal{O} \! \left(
\dfrac{c(n)}{(n \! + \! \frac{1}{2})^{2}} \right) \right\},
\end{align}
where
\begin{align}
\mathscr{R}^{o,0}_{1} :=& \, \mathlarger{\sum_{j=1}^{N+1}} \! \left(\dfrac{
\left(\mathscr{A}^{o}(b_{j-1}^{o}) \! \left(\widehat{\alpha}_{1}^{o}(b_{j-1}^{
o}) \! + \! 2(b_{j-1}^{o})^{-1} \widehat{\alpha}_{0}^{o}(b_{j-1}^{o}) \right) 
\! - \! \mathscr{B}^{o}(b_{j-1}^{o}) \widehat{\alpha}_{0}^{o}(b_{j-1}^{o}) 
\right)}{(b_{j-1}^{o})^{2}(\widehat{\alpha}_{0}^{o}(b_{j-1}^{o}))^{2}} 
\right. \nonumber \\
+&\left. \, \dfrac{\left(\mathscr{A}^{o}(a_{j}^{o}) \! \left(\widehat{
\alpha}_{1}^{o}(a_{j}^{o}) \! + \! 2(a_{j}^{o})^{-1} \widehat{\alpha}_{0}^{o}
(a_{j}^{o}) \right) \! - \! \mathscr{B}^{o}(a_{j}^{o}) \widehat{\alpha}_{0}^{
o}(a_{j}^{o}) \right)}{(a_{j}^{o})^{2}(\widehat{\alpha}_{0}^{o}(a_{j}^{o}))^{
2}} \right),
\end{align}
with all parameters defined in Theorem~{\rm 3.1}, 
Equations~{\rm (3.52)--(3.100)},
\begin{align}
b_{2n+1}^{\sharp} \underset{n \to \infty}{=}& \, -\dfrac{2 \varpi^{o}_{-}}{
\mathbb{E}} \dfrac{\boldsymbol{\theta}^{o}(\boldsymbol{u}^{o}_{+}(0) \! + \! 
\boldsymbol{d}_{o}) \boldsymbol{\theta}^{o}(-\boldsymbol{u}^{o}_{+}(\infty) 
\! - \! \frac{1}{2 \pi}(n \! + \! \frac{1}{2}) \boldsymbol{\Omega}^{o} \! + \! 
\boldsymbol{d}_{o})}{\boldsymbol{\theta}^{o}(-\boldsymbol{u}^{o}_{+}(\infty) 
\! + \! \boldsymbol{d}_{o}) \boldsymbol{\theta}^{o}(\boldsymbol{u}^{o}_{+}(0) 
\! - \! \frac{1}{2 \pi}(n \! + \! \frac{1}{2}) \boldsymbol{\Omega}^{o} \! + \! 
\boldsymbol{d}_{o})} \nonumber \\
&\times \sqrt{\dfrac{1}{\eta^{e}_{+}} \dfrac{\boldsymbol{\theta}^{e}
(\boldsymbol{u}^{e}_{+}(\infty) \! - \! \frac{n}{2 \pi} \boldsymbol{\Omega}^{
e} \! + \! \boldsymbol{d}_{e}) \boldsymbol{\theta}^{e}(-\boldsymbol{u}^{e}_{+}
(\infty) \! + \! \boldsymbol{d}_{e})}{\boldsymbol{\theta}^{e}(-\boldsymbol{
u}^{e}_{+}(\infty) \! - \! \frac{n}{2 \pi} \boldsymbol{\Omega}^{e} \! + \! 
\boldsymbol{d}_{e}) \boldsymbol{\theta}^{e}(\boldsymbol{u}^{e}_{+}(\infty) 
\! + \! \boldsymbol{d}_{e})}} \nonumber \\
&\times \sqrt{\dfrac{\mathbb{E}^{2}}{\eta^{o}_{-}} \dfrac{\boldsymbol{
\theta}^{o}(\boldsymbol{u}^{o}_{+}(0) \! - \! \frac{1}{2 \pi}(n \! + \! 
\frac{1}{2}) \boldsymbol{\Omega}^{o} \! + \! \boldsymbol{d}_{o}) \boldsymbol{
\theta}^{o}(-\boldsymbol{u}^{o}_{+}(0) \! + \! \boldsymbol{d}_{o})}{
\boldsymbol{\theta}^{o}(-\boldsymbol{u}^{o}_{+}(0) \! - \! \frac{1}{2 \pi}
(n \! + \! \frac{1}{2}) \boldsymbol{\Omega}^{o} \! + \! \boldsymbol{d}_{o}) 
\boldsymbol{\theta}^{o}(\boldsymbol{u}^{o}_{+}(0) \! + \! \boldsymbol{d}_{
o})}} \nonumber \\
&\times \left\{1 \! + \! \dfrac{1}{n} \! \left(\dfrac{2 \mathfrak{Q}^{\flat}_{
12}}{\eta^{e}_{+}} \dfrac{\boldsymbol{\theta}^{e}(\boldsymbol{u}^{e}_{+}
(\infty) \! - \! \frac{n}{2 \pi} \boldsymbol{\Omega}^{e} \! + \! \boldsymbol{
d}_{e}) \boldsymbol{\theta}^{e}(-\boldsymbol{u}^{e}_{+}(\infty) \! + \! 
\boldsymbol{d}_{e})}{\boldsymbol{\theta}^{e}(-\boldsymbol{u}^{e}_{+}(\infty) 
\! - \! \frac{n}{2 \pi} \boldsymbol{\Omega}^{e} \! + \! \boldsymbol{d}_{e}) 
\boldsymbol{\theta}^{e}(\boldsymbol{u}^{e}_{+}(\infty) \! + \! \boldsymbol{
d}_{e})} \right) \right. \nonumber \\
&\left. + \, \dfrac{1}{(n \! + \! \frac{1}{2})} \! \left(\dfrac{2 \mathbb{E}^{
2} \mathfrak{Q}^{\natural}_{12}}{\eta^{o}_{-}} \dfrac{\boldsymbol{\theta}^{o}
(\boldsymbol{u}^{o}_{+}(0) \! - \! \frac{1}{2 \pi}(n \! + \! \frac{1}{2}) 
\boldsymbol{\Omega}^{o} \! + \! \boldsymbol{d}_{o}) \boldsymbol{\theta}^{o}
(-\boldsymbol{u}^{o}_{+}(0) \! + \! \boldsymbol{d}_{o})}{\boldsymbol{\theta}^{
o}(-\boldsymbol{u}^{o}_{+}(0) \! - \! \frac{1}{2 \pi}(n \! + \! \frac{1}{2}) 
\boldsymbol{\Omega}^{o} \! + \! \boldsymbol{d}_{o}) \boldsymbol{\theta}^{o}
(\boldsymbol{u}^{o}_{+}(0) \! + \! \boldsymbol{d}_{o})} \right. \right. 
\nonumber \\
&\left. \left. + \, (\mathscr{R}^{o,\infty}_{0})_{11} \! + \! (\mathscr{R}^{
o,\infty}_{0})_{12} \dfrac{\mi \varpi^{o}_{+} \mathbb{E}^{2}}{\varpi^{o}_{-}} 
\dfrac{\boldsymbol{\theta}^{o}(-\boldsymbol{u}^{o}_{+}(\infty) \! - \! \frac{
1}{2 \pi}(n \! + \! \frac{1}{2}) \boldsymbol{\Omega}^{o} \! - \! \boldsymbol{
d}_{o})}{\boldsymbol{\theta}^{o}(-\boldsymbol{u}^{o}_{+}(\infty) \! - \! 
\frac{1}{2 \pi}(n \! + \! \frac{1}{2}) \boldsymbol{\Omega}^{o} \! + \! 
\boldsymbol{d}_{o})} \right. \right. \nonumber \\
&\left. \left. \times \, \dfrac{\boldsymbol{\theta}^{o}(-\boldsymbol{u}^{o}_{
+}(\infty) \! + \! \boldsymbol{d}_{o}) \boldsymbol{\theta}^{o}(\boldsymbol{
u}^{o}_{+}(0) \! - \! \frac{1}{2 \pi}(n \! + \! \frac{1}{2}) \boldsymbol{
\Omega}^{o} \! + \! \boldsymbol{d}_{o})}{\boldsymbol{\theta}^{o}(-\boldsymbol{
u}^{o}_{+}(\infty) \! - \! \boldsymbol{d}_{o}) \boldsymbol{\theta}^{o}
(-\boldsymbol{u}^{o}_{+}(0) \! - \! \frac{1}{2 \pi}(n \! + \! \frac{1}{2}) 
\boldsymbol{\Omega}^{o} \! - \! \boldsymbol{d}_{o})} \right) \! + \! \mathcal{
O} \! \left(\dfrac{c(n)}{n^{2}} \right) \right\} \nonumber \\
&\times \, \exp \! \left(2 \! \left(n \! + \! \dfrac{1}{2} \right) \! \int_{
J_{o}} \ln (\lvert s \rvert) \psi_{V}^{o}(s) \, \md s \! + \! \dfrac{n}{2} 
\! \left(\ell_{o} \! - \! \ell_{e} \right) \right),
\end{align}
where
\begin{align}
\mathfrak{Q}^{\flat} \! := \! \mi \mathscr{R}^{e,\infty}_{1}, \qquad \qquad 
\qquad \mathfrak{Q}^{\natural} \! := \! -\mi \mathscr{R}^{o,0}_{1},
\end{align}
with
\begin{equation}
\mathscr{R}^{e,\infty}_{1} \! := \mathlarger{\sum_{j=1}^{N+1}} \! \left(
\dfrac{\left(\mathscr{B}^{e}(a_{j}^{e}) \widehat{\alpha}_{0}^{e}(a_{j}^{e}) \! 
- \! \mathscr{A}^{e}(a_{j}^{e}) \widehat{\alpha}_{1}^{e}(a_{j}^{e}) \right)}{
(\widehat{\alpha}_{0}^{e}(a_{j}^{e}))^{2}} \! + \! \dfrac{\left(\mathscr{B}^{
e}(b_{j-1}^{e}) \widehat{\alpha}_{0}^{e}(b_{j-1}^{e}) \! - \! \mathscr{A}^{e}
(b_{j-1}^{e}) \widehat{\alpha}_{1}^{e}(b_{j-1}^{e}) \right)}{(\widehat{
\alpha}_{0}^{e}(b_{j-1}^{e}))^{2}} \right),
\end{equation}
and all parameters defined in Theorem~{\rm 3.1}, 
Equations~{\rm (3.52)--(3.100)},
\begin{align}
b_{2n+2}^{\sharp} \underset{n \to \infty}{=}& \, \nu^{(2n+1)}_{n} \, \sqrt{
\dfrac{h^{+} \! \left[\me^{n \ell_{e}} \eta^{e}_{+} \mathfrak{a}^{\uparrow}(n)
\! \left(1 \! + \! \frac{1}{n} \mathfrak{A}^{\uparrow}(n) \! + \! \mathcal{O}
\! \left(\frac{c(n)}{n^{2}} \right) \right) \right]}{\me^{n \ell_{o}} \eta^{
o}_{-} \mathfrak{b}^{\downarrow}(n) \! \left(1 \! + \! \frac{1}{(n+\frac{1}{
2})} \mathfrak{B}^{\downarrow}(n) \! + \! \mathcal{O} \! \left(\frac{c(n)}{(
n+\frac{1}{2})^{2}} \right) \right)}},
\end{align}
where asymptotics (as $n \! \to \infty)$ of $\nu^{(2n+1)}_{n}$ is given in 
Theorem~{\rm 3.1}, Equation~{\rm (3.49)}, with
\begin{gather}
\mathfrak{a}^{\uparrow}(n) \! := \! \dfrac{\boldsymbol{\theta}^{e}(\boldsymbol{
u}^{e}_{+}(\infty) \! + \! \boldsymbol{d}_{e}) \boldsymbol{\theta}^{e}(-
\boldsymbol{u}^{e}_{+}(\infty) \! - \! \frac{n}{2 \pi} \boldsymbol{\Omega}^{e}
\! + \! \boldsymbol{d}_{e})}{\boldsymbol{\theta}^{e}(-\boldsymbol{u}^{e}_{+}
(\infty) \! + \! \boldsymbol{d}_{e}) \boldsymbol{\theta}^{e}(\boldsymbol{u}^{
e}_{+}(\infty) \! - \! \frac{n}{2 \pi} \boldsymbol{\Omega}^{e} \! + \!
\boldsymbol{d}_{e})}, \\
\mathfrak{b}^{\downarrow}(n) \! := \! \dfrac{1}{\mathbb{E}^{2}} \dfrac{
\boldsymbol{\theta}^{o}(\boldsymbol{u}^{o}_{+}(0) \! + \! \boldsymbol{d}_{o})
\boldsymbol{\theta}^{o}(-\boldsymbol{u}^{o}_{+}(0) \! - \! \frac{1}{2 \pi}(
n \! + \! \frac{1}{2}) \boldsymbol{\Omega}^{o} \! + \! \boldsymbol{d}_{o})}{
\boldsymbol{\theta}^{o}(-\boldsymbol{u}^{o}_{+}(0) \! + \! \boldsymbol{d}_{o}
) \boldsymbol{\theta}^{o}(\boldsymbol{u}^{o}_{+}(0) \! - \! \frac{1}{2 \pi}(n
\! + \! \frac{1}{2}) \boldsymbol{\Omega}^{o} \! + \! \boldsymbol{d}_{o})}, \\
\mathfrak{A}^{\uparrow}(n) \! := \! -\dfrac{4 \mi (\mathscr{R}^{e,\infty}_{1}
)_{12}}{\eta^{e}_{+}} \dfrac{\boldsymbol{\theta}^{e}(-\boldsymbol{u}^{e}_{+}
(\infty) \! + \! \boldsymbol{d}_{e}) \boldsymbol{\theta}^{e}(\boldsymbol{u}^{
e}_{+}(\infty) \! - \! \frac{n}{2 \pi} \boldsymbol{\Omega}^{e} \! + \!
\boldsymbol{d}_{e})}{\boldsymbol{\theta}^{e}(\boldsymbol{u}^{e}_{+}(\infty)
\! + \! \boldsymbol{d}_{e}) \boldsymbol{\theta}^{e}(-\boldsymbol{u}^{e}_{+}
(\infty) \! - \! \frac{n}{2 \pi} \boldsymbol{\Omega}^{e} \! + \! \boldsymbol{
d}_{e})}, \\
\mathfrak{B}^{\downarrow}(n) \! := \! \dfrac{4 \mi (\mathscr{R}^{o,0}_{1})_{1
2} \mathbb{E}^{2}}{\eta^{o}_{-}} \dfrac{\boldsymbol{\theta}^{o}(-\boldsymbol{
u}^{o}_{+}(0) \! + \! \boldsymbol{d}_{o}) \boldsymbol{\theta}^{o}(\boldsymbol{
u}^{o}_{+}(0) \! - \! \frac{1}{2 \pi}(n \! + \! \frac{1}{2}) \boldsymbol{
\Omega}^{o} \! + \! \boldsymbol{d}_{o})}{\boldsymbol{\theta}^{o}(\boldsymbol{
u}^{o}_{+}(0) \! + \! \boldsymbol{d}_{o}) \boldsymbol{\theta}^{o}(-
\boldsymbol{u}^{o}_{+}(0) \! - \! \frac{1}{2 \pi}(n \! + \! \frac{1}{2})
\boldsymbol{\Omega}^{o} \! + \! \boldsymbol{d}_{o})},
\end{gather}
\begin{align}
c_{2n+2}^{\sharp} \underset{n \to \infty}{=}& \, \sqrt{\dfrac{h^{+} \! \left[
\me^{n \ell_{e}} \eta^{e}_{+} \mathfrak{a}^{\uparrow}(n) \! \left(1 \! + \!
\frac{1}{n} \mathfrak{A}^{\uparrow}(n) \! + \! \mathcal{O} \! \left(\frac{c(n)
}{n^{2}} \right) \right) \right]}{\me^{n \ell_{e}} \eta^{e}_{+} \mathfrak{a}^{
\uparrow}(n) \! \left(1 \! + \! \frac{1}{n} \mathfrak{A}^{\uparrow}(n) \! + \!
\mathcal{O} \! \left(\frac{c(n)}{n^{2}} \right) \right)}},
\end{align}
\begin{align}
\alpha_{2n}^{\sharp} \underset{n \to \infty}{=}& \, -\dfrac{\varpi^{e}_{+}
\varpi^{o}_{+}}{4 \mathbb{E}} \dfrac{\boldsymbol{\theta}^{e}(\boldsymbol{
u}^{e}_{+}(\infty) \! + \! \boldsymbol{d}_{e}) \boldsymbol{\theta}^{o}
(\boldsymbol{u}^{o}_{+}(0) \! + \! \boldsymbol{d}_{o}) \boldsymbol{\theta}^{e}
(\boldsymbol{u}^{e}_{+}(0) \! - \! \frac{n}{2 \pi} \boldsymbol{\Omega}^{e} \!
+ \! \boldsymbol{d}_{e})}{\boldsymbol{\theta}^{e}(\boldsymbol{u}^{e}_{+}(0)
\! + \! \boldsymbol{d}_{e}) \boldsymbol{\theta}^{o}(\boldsymbol{u}^{o}_{+}
(\infty) \! + \! \boldsymbol{d}_{o}) \boldsymbol{\theta}^{e}(\boldsymbol{
u}^{e}_{+}(\infty) \! - \! \frac{n}{2 \pi} \boldsymbol{\Omega}^{e} \! + \!
\boldsymbol{d}_{e})} \nonumber \\
&\times \, \dfrac{\boldsymbol{\theta}^{o}(\boldsymbol{u}^{o}_{+}(\infty) \!
- \! \frac{1}{2 \pi}(n \! + \! \frac{1}{2}) \boldsymbol{\Omega}^{o} \! + \!
\boldsymbol{d}_{o})}{\boldsymbol{\theta}^{o}(\boldsymbol{u}^{o}_{+}(0) \! -
\! \frac{1}{2 \pi}(n \! + \! \frac{1}{2}) \boldsymbol{\Omega}^{o} \! + \!
\boldsymbol{d}_{o})} \! \left\{1 \! + \! \dfrac{1}{n} \! \left((\mathscr{R}^{
e,0}_{0})_{11} \! + \! (\mathscr{R}^{e,0}_{0})_{12} \dfrac{\varpi^{e}_{-}}{
\mi \varpi^{e}_{+}} \right. \right. \nonumber \\
&\left. \left. \times \, \dfrac{\boldsymbol{\theta}^{e}(\boldsymbol{u}^{e}_{+}
(0) \! + \! \boldsymbol{d}_{e}) \boldsymbol{\theta}^{e}(\boldsymbol{u}^{e}_{+}
(0) \! - \! \frac{n}{2 \pi} \boldsymbol{\Omega}^{e} \! - \! \boldsymbol{d}_{e}
) \boldsymbol{\theta}^{e}(\boldsymbol{u}^{e}_{+}(\infty) \! - \! \frac{n}{2
\pi} \boldsymbol{\Omega}^{e} \! + \! \boldsymbol{d}_{e})}{\boldsymbol{\theta}^{
e}(\boldsymbol{u}^{e}_{+}(0) \! - \! \boldsymbol{d}_{e}) \boldsymbol{\theta}^{
e}(\boldsymbol{u}^{e}_{+}(0) \! - \! \frac{n}{2 \pi} \boldsymbol{\Omega}^{e}
\! + \! \boldsymbol{d}_{e}) \boldsymbol{\theta}^{e}(-\boldsymbol{u}^{e}_{+}
(\infty) \! - \! \frac{n}{2 \pi} \boldsymbol{\Omega}^{e} \! - \! \boldsymbol{
d}_{e})} \right) \right. \nonumber \\
&\left. + \, \dfrac{1}{(n \! + \! \frac{1}{2})} \! \left((\mathscr{R}^{o,
\infty}_{0})_{11} \! - \! (\mathscr{R}^{o,\infty}_{0})_{12} \dfrac{\varpi^{
o}_{-} \mathbb{E}^{2}}{\mi \varpi^{o}_{+}} \dfrac{\boldsymbol{\theta}^{o}
(\boldsymbol{u}^{o}_{+}(0) \! - \! \frac{1}{2 \pi}(n \! + \! \frac{1}{2}) 
\boldsymbol{\Omega}^{o} \! + \! \boldsymbol{d}_{o})}{\boldsymbol{\theta}^{o}
(-\boldsymbol{u}^{o}_{+}(0) \! - \! \frac{1}{2 \pi}(n \! + \! \frac{1}{2}) 
\boldsymbol{\Omega}^{o} \! - \! \boldsymbol{d}_{o})} \right. \right. \nonumber 
\\
&\left. \left. \times \, \dfrac{\boldsymbol{\theta}^{o}(\boldsymbol{u}^{o}_{+}
(\infty) \! + \! \boldsymbol{d}_{o}) \boldsymbol{\theta}^{o}(\boldsymbol{u}^{
o}_{+}(\infty) \! - \! \frac{1}{2 \pi}(n \! + \! \frac{1}{2}) \boldsymbol{
\Omega}^{o} \! - \! \boldsymbol{d}_{o})}{\boldsymbol{\theta}^{o}(\boldsymbol{
u}^{o}_{+}(\infty) \! - \! \boldsymbol{d}_{o}) \boldsymbol{\theta}^{o}
(\boldsymbol{u}^{o}_{+}(\infty) \! - \! \frac{1}{2 \pi}(n \! + \! \frac{1}{2}) 
\boldsymbol{\Omega}^{o} \! + \! \boldsymbol{d}_{o})} \right) \! + \! \mathcal{
O} \! \left(\dfrac{c(n)}{n^{2}} \right) \right\}\nonumber \\
&\times \, \exp \! \left(2n \! \left(\int_{J_{e}} \ln (\lvert s \rvert) \psi_{
V}^{e}(s) \, \md s \! + \! \mi \pi \int_{J_{e} \cap \mathbb{R}_{+}} \psi_{V}^{
e}(s) \, \md s \right) \! - \! 2 \! \left(n \! + \! \dfrac{1}{2} \right) \! 
\int_{J_{o}} \ln (\lvert s \rvert) \psi_{V}^{o}(s) \, \md s \right),
\end{align}
\begin{align}
\beta_{2n+1}^{\sharp} \underset{n \to \infty}{=}& \, \dfrac{\varpi^{e}_{+}}{
2} \dfrac{\boldsymbol{\theta}^{e}(\boldsymbol{u}^{e}_{+}(\infty) \! + \!
\boldsymbol{d}_{e}) \boldsymbol{\theta}^{e}(\boldsymbol{u}^{e}_{+}(0) \! -
\! \frac{n}{2 \pi} \boldsymbol{\Omega}^{e} \! + \! \boldsymbol{d}_{e})}{
\boldsymbol{\theta}^{e}(\boldsymbol{u}^{e}_{+}(0) \! + \! \boldsymbol{d}_{e})
\boldsymbol{\theta}^{e}(\boldsymbol{u}^{e}_{+}(\infty) \! - \! \frac{n}{2 \pi}
\boldsymbol{\Omega}^{e} \! + \! \boldsymbol{d}_{e})} \nonumber \\
&\times \sqrt{\dfrac{1}{\eta^{e}_{+}} \dfrac{\boldsymbol{\theta}^{e}
(\boldsymbol{u}^{e}_{+}(\infty) \! - \! \frac{n}{2 \pi} \boldsymbol{\Omega}^{
e} \! + \! \boldsymbol{d}_{e}) \boldsymbol{\theta}^{e}(-\boldsymbol{u}^{e}_{+}
(\infty) \! + \! \boldsymbol{d}_{e})}{\boldsymbol{\theta}^{e}(-\boldsymbol{
u}^{e}_{+}(\infty) \! - \! \frac{n}{2 \pi} \boldsymbol{\Omega}^{e} \! + \!
\boldsymbol{d}_{e}) \boldsymbol{\theta}^{e}(\boldsymbol{u}^{e}_{+}(\infty)
\! + \! \boldsymbol{d}_{e})}} \nonumber \\
&\times \left(\dfrac{\mathbb{E}^{2}}{\eta^{o}_{-}} \dfrac{\boldsymbol{
\theta}^{o}(\boldsymbol{u}^{o}_{+}(0) \! - \! \frac{1}{2 \pi}(n \! + \!
\frac{1}{2}) \boldsymbol{\Omega}^{o} \! + \! \boldsymbol{d}_{o}) \boldsymbol{
\theta}^{o}(-\boldsymbol{u}^{o}_{+}(0) \! + \! \boldsymbol{d}_{o})}{
\boldsymbol{\theta}^{o}(-\boldsymbol{u}^{o}_{+}(0) \! - \! \frac{1}{2 \pi}
(n \! + \! \frac{1}{2}) \boldsymbol{\Omega}^{o} \! + \! \boldsymbol{d}_{o})
\boldsymbol{\theta}^{o}(\boldsymbol{u}^{o}_{+}(0) \! + \! \boldsymbol{d}_{o})}
\right)^{-1/2} \nonumber \\
&\times \left\{1 \! + \! \dfrac{1}{(n \! + \! \frac{1}{2})} \! \left(-\dfrac{
2 \mathbb{E}^{2} \mathfrak{Q}^{\natural}_{12}}{\eta^{o}_{-}} \dfrac{
\boldsymbol{\theta}^{o}(\boldsymbol{u}^{o}_{+}(0) \! - \! \frac{1}{2 \pi}
(n \! + \! \frac{1}{2}) \boldsymbol{\Omega}^{o} \! + \! \boldsymbol{d}_{o})
\boldsymbol{\theta}^{o}(-\boldsymbol{u}^{o}_{+}(0) \! + \! \boldsymbol{d}_{
o})}{\boldsymbol{\theta}^{o}(-\boldsymbol{u}^{o}_{+}(0) \! - \! \frac{1}{2
\pi}(n \! + \! \frac{1}{2}) \boldsymbol{\Omega}^{o} \! + \! \boldsymbol{d}_{
o}) \boldsymbol{\theta}^{o}(\boldsymbol{u}^{o}_{+}(0) \! + \! \boldsymbol{
d}_{o})} \right) \right. \nonumber \\
&\left. + \, \dfrac{1}{n} \! \left(\dfrac{2 \mathfrak{Q}^{\flat}_{12}}{\eta^{
e}_{+}} \dfrac{\boldsymbol{\theta}^{e}(\boldsymbol{u}^{e}_{+}(\infty) \!
- \! \frac{n}{2 \pi} \boldsymbol{\Omega}^{e} \! + \! \boldsymbol{d}_{e})
\boldsymbol{\theta}^{e}(-\boldsymbol{u}^{e}_{+}(\infty) \! + \! \boldsymbol{
d}_{e})}{\boldsymbol{\theta}^{e}(-\boldsymbol{u}^{e}_{+}(\infty) \! - \!
\frac{n}{2 \pi} \boldsymbol{\Omega}^{e} \! + \! \boldsymbol{d}_{e})
\boldsymbol{\theta}^{e}(\boldsymbol{u}^{e}_{+}(\infty) \! + \! \boldsymbol{
d}_{e})} \! + \! (\mathscr{R}^{e,0}_{0})_{11} \! + \! (\mathscr{R}^{e,0}_{
0})_{12} \dfrac{\varpi^{e}_{-}}{\mi \varpi^{e}_{+}} \right. \right. \nonumber 
\\
&\left. \left. \times \, \dfrac{\boldsymbol{\theta}^{e}(\boldsymbol{u}^{e}_{+}
(0) \! - \! \frac{n}{2 \pi} \boldsymbol{\Omega}^{e} \! - \! \boldsymbol{d}_{
e}) \boldsymbol{\theta}^{e}(\boldsymbol{u}^{e}_{+}(0) \! + \! \boldsymbol{d}_{
e}) \boldsymbol{\theta}^{e}(\boldsymbol{u}^{e}_{+}(\infty) \! - \! \frac{n}{2 
\pi} \boldsymbol{\Omega}^{e} \! + \! \boldsymbol{d}_{e})}{\boldsymbol{\theta}^{
e}(\boldsymbol{u}^{e}_{+}(0) \! - \! \frac{n}{2 \pi} \boldsymbol{\Omega}^{e} 
\! + \! \boldsymbol{d}_{e}) \boldsymbol{\theta}^{e}(\boldsymbol{u}^{e}_{+}(0) 
\! - \! \boldsymbol{d}_{e}) \boldsymbol{\theta}^{e}(-\boldsymbol{u}^{e}_{+}
(\infty) \! - \! \frac{n}{2 \pi} \boldsymbol{\Omega}^{e} \! - \! \boldsymbol{
d}_{e})} \right) \! + \! \mathcal{O} \! \left(\dfrac{c(n)}{n^{2}} \right) 
\right\} \nonumber \\
&\times \, \exp \! \left(2n \! \left(\int_{J_{e}} \ln (\lvert s \rvert) \psi_{
V}^{e}(s) \, \md s \! + \! \mi \pi \int_{J_{e} \cap \mathbb{R}_{+}} \psi_{V}^{
e}(s) \, \md s \right) \! + \! \dfrac{n}{2} \! \left(\ell_{o} \! - \! \ell_{
e} \right) \right),
\end{align}
\begin{align}
\beta_{2n+2}^{\sharp} \underset{n \to \infty}{=}& \, - \! \left(h^{+} \! 
\left[\nu^{(2n+1)}_{n} \right] \right) \sqrt{\dfrac{h^{+} \! \left[\me^{n 
\ell_{e}} \eta^{e}_{+} \mathfrak{a}^{\uparrow}(n) \! \left(1 \! + \! \frac{
1}{n} \mathfrak{A}^{\uparrow}(n) \! + \! \mathcal{O} \! \left(\frac{c(n)}{n^{
2}} \right) \right) \right]}{\me^{n \ell_{o}} \eta^{o}_{-} \mathfrak{b}^{
\downarrow}(n) \! \left(1 \! + \! \frac{1}{(n+\frac{1}{2})} \mathfrak{B}^{
\downarrow}(n) \! + \! \mathcal{O} \! \left(\frac{c(n)}{(n+\frac{1}{2})^{2}} 
\right) \right)}},
\end{align}
\begin{align}
\gamma_{2n+3}^{\sharp} \underset{n \to \infty}{=}& \, \sqrt{\dfrac{h^{+} \! 
\left[\me^{n \ell_{o}} \eta^{o}_{-} \mathfrak{b}^{\downarrow}(n) \! \left(1 \! 
+ \! \frac{1}{(n+\frac{1}{2})} \mathfrak{B}^{\downarrow}(n) \! + \! \mathcal{
O} \! \left(\frac{c(n)}{(n+\frac{1}{2})^{2}} \right) \right) \right]}{\me^{n 
\ell_{o}} \eta^{o}_{-} \mathfrak{b}^{\downarrow}(n) \! \left(1 \! + \! \frac{
1}{(n+\frac{1}{2})} \mathfrak{B}^{\downarrow}(n) \! + \! \mathcal{O} \! \left(
\frac{c(n)}{(n+\frac{1}{2})^{2}} \right) \right)}},
\end{align}
\begin{align}
b_{2n}^{\sharp} \underset{n \to \infty}{=}& \, h^{-} \! \left[\nu^{(2n+1)}_{
n} \right] \sqrt{\dfrac{\me^{n \ell_{e}} \eta^{e}_{+} \mathfrak{a}^{\uparrow}
(n) \! \left(1 \! + \! \frac{1}{n} \mathfrak{A}^{\uparrow}(n) \! + \! 
\mathcal{O} \! \left(\frac{c(n)}{n^{2}} \right) \right)}{h^{-} \! \left[\me^{
n \ell_{o}} \eta^{o}_{-} \mathfrak{b}^{\downarrow}(n) \! \left(1 \! + \! 
\frac{1}{(n+\frac{1}{2})} \mathfrak{B}^{\downarrow}(n) \! + \! \mathcal{O} 
\! \left(\frac{c(n)}{(n+\frac{1}{2})^{2}} \right) \right) \right]}},
\end{align}
\begin{align}
c_{2n}^{\sharp} \underset{n \to \infty}{=}& \, \sqrt{\dfrac{\me^{n \ell_{e}} 
\eta^{e}_{+} \mathfrak{a}^{\uparrow}(n) \! \left(1 \! + \! \frac{1}{n} 
\mathfrak{A}^{\uparrow}(n) \! + \! \mathcal{O} \! \left(\frac{c(n)}{n^{2}} 
\right) \right)}{h^{-} \! \left[\me^{n \ell_{e}} \eta^{e}_{+} \mathfrak{a}^{
\uparrow}(n) \! \left(1 \! + \! \frac{1}{n} \mathfrak{A}^{\uparrow}(n) \! + 
\! \mathcal{O} \! \left(\frac{c(n)}{n^{2}} \right) \right) \right]}},
\end{align}
\begin{align}
\beta_{2n}^{\sharp} \underset{n \to \infty}{=}& \, -\nu^{(2n+1)}_{n} \, \sqrt{
\dfrac{\me^{n \ell_{e}} \eta^{e}_{+} \mathfrak{a}^{\uparrow}(n) \! \left(1
\! + \! \frac{1}{n} \mathfrak{A}^{\uparrow}(n) \! + \! \mathcal{O} \! \left(
\frac{c(n)}{n^{2}} \right) \right)}{h^{-} \! \left[\me^{n \ell_{o}} \eta^{
o}_{-} \mathfrak{b}^{\downarrow}(n) \! \left(1 \! + \! \frac{1}{(n+\frac{1}{
2})} \mathfrak{B}^{\downarrow}(n) \! + \! \mathcal{O} \! \left(\frac{c(n)}{
(n+\frac{1}{2})^{2}} \right) \right) \right]}},
\end{align}
\begin{align}
\gamma_{2n+1}^{\sharp} \underset{n \to \infty}{=}& \, \sqrt{\dfrac{\me^{n 
\ell_{o}} \eta^{o}_{-} \mathfrak{b}^{\downarrow}(n) \! \left(1 \! + \! \frac{
1}{(n+\frac{1}{2})} \mathfrak{B}^{\downarrow}(n) \! + \! \mathcal{O} \! \left(
\frac{c(n)}{(n+\frac{1}{2})^{2}} \right) \right)}{h^{-} \! \left[\me^{n \ell_{
o}} \eta^{o}_{-} \mathfrak{b}^{\downarrow}(n) \! \left(1 \! + \! \frac{1}{(n+
\frac{1}{2})} \mathfrak{B}^{\downarrow}(n) \! + \! \mathcal{O} \! \left(\frac{
c(n)}{(n+\frac{1}{2})^{2}} \right) \right) \right]}},
\end{align}
and $n \! \to \! \infty$ asymptotics for $a_{2n}^{\sharp}$ and 
$\alpha_{2n+1}^{\sharp}$ are obtained by substituting, formally, 
asymptotics~{\rm (3.106)}, {\rm (3.114)--(3.117)}, {\rm (3.119)}, 
{\rm (3.121)}, and~{\rm (3.122)} into Equations~{\rm (3.35)} 
and~{\rm (3.56)}. Furthermore,
\begin{align}
\dfrac{\xi^{(2n)}_{n-1}}{\xi^{(2n)}_{n}} \underset{n \to \infty}{=}& \, -2n 
\int_{J_{e}}s \psi_{V}^{e}(s) \, \md s \! + \! \dfrac{\alpha^{e}_{\infty}(1,
1,\pmb{0})}{\boldsymbol{\theta}^{e}(\boldsymbol{u}^{e}_{+}(\infty) \! + \! 
\boldsymbol{d}_{e})} \! - \! \dfrac{\alpha^{e}_{\infty}(1,1,\boldsymbol{
\Omega}^{e})}{\boldsymbol{\theta}^{e}(\boldsymbol{u}^{e}_{+}(\infty) \! - \! 
\frac{n}{2 \pi} \boldsymbol{\Omega}^{e} \! + \! \boldsymbol{d}_{e})} \! + \! 
\dfrac{1}{n}(\mathscr{R}^{e,\infty}_{1})_{11} \! + \! \mathcal{O} \! \left(
\dfrac{c(n)}{n^{2}} \right), \\
\dfrac{\xi^{(2n+1)}_{-n}}{\xi^{(2n+1)}_{-n-1}} \underset{n \to \infty}{=}& \,
-2 \! \left(n \! + \! \dfrac{1}{2} \right) \! \int_{J_{o}}s^{-1} \psi_{V}^{o}
(s) \, \md s \! - \! \dfrac{\alpha^{o}_{0}(1,1,\pmb{0})}{\boldsymbol{\theta}^{
o}(\boldsymbol{u}^{o}_{+}(0) \! + \! \boldsymbol{d}_{o})} \! + \! \dfrac{
\alpha^{o}_{0}(1,1,\boldsymbol{\Omega}^{o})}{\boldsymbol{\theta}^{o}
(\boldsymbol{u}^{o}_{+}(0) \! - \! \frac{1}{2 \pi}(n \! + \! \frac{1}{2})
\boldsymbol{\Omega}^{o} \! + \! \boldsymbol{d}_{o})} \! + \! \dfrac{1}{(n 
\! + \! \frac{1}{2})}(\mathscr{R}^{o,0}_{1})_{11} \nonumber \\
&+ \, \mathcal{O} \! \left(\dfrac{c(n)}{(n \! + \! \frac{1}{2})^{2}} \right),
\end{align}
where $\mathscr{R}^{o,0}_{1}$ (resp., $\mathscr{R}^{e,\infty}_{1})$ is defined 
in Equation~{\rm (3.105)} (resp., Equation~{\rm (3.108))},
\begin{gather}
\alpha^{e}_{\infty}(1,1,\pmb{\mathscr{Z}}) \! := \! 2 \pi \mi \sum_{m \in 
\mathbb{Z}^{N}} (m,\widehat{\boldsymbol{\alpha}}^{e}_{\infty}) \me^{2 \pi 
\mi (m,\boldsymbol{u}^{e}_{+}(\infty)-\frac{n}{2 \pi} \pmb{\mathscr{Z}}+
\boldsymbol{d}_{e})+ \pi \mi (m,\tau^{e}m)}, \\
\alpha^{o}_{0}(1,1,\pmb{\mathscr{Z}}) \! := \! 2 \pi \mi \sum_{m \in \mathbb{
Z}^{N}}(m,\widehat{\boldsymbol{\alpha}}^{o}_{0}) \me^{2 \pi \mi (m,\boldsymbol{
u}^{o}_{+}(0)-\frac{1}{2 \pi}(n+\frac{1}{2}) \pmb{\mathscr{Z}}+\boldsymbol{
d}_{o})+ \pi \mi (m,\tau^{o}m)},
\end{gather}
where $\widehat{\boldsymbol{\alpha}}^{q}_{p} \! = \! (\widehat{\alpha}_{p,
1}^{q},\widehat{\alpha}^{q}_{p,2},\dotsc,\widehat{\alpha}^{q}_{p,N})$, $q \! 
\in \! \{e,o\}$, $p \! \in \! \{0,\infty\}$, with
\begin{equation}
\widehat{\alpha}^{e}_{\infty,j} \! := \! c_{j1}^{e}, \qquad \qquad \widehat{
\alpha}^{o}_{0,j} \! := \! (-1)^{\mathcal{N}_{+}^{o}} \! \left(\prod_{i=1}^{N
+1} \lvert b_{i-1}^{o} a_{i}^{o} \rvert \right)^{-1/2}c_{jN}^{o}, \quad j \! 
= \! 1,\dotsc,N,
\end{equation}
where $c_{j1}^{e}$ (resp., $c_{jN}^{o})$, $j \! = \! 1,\dotsc,N$, are obtained 
{}from Equations~{\rm (E1)} and~{\rm (E2)} (resp., Equations~{\rm (O1)} 
and~{\rm (O2))}, and $\pmb{0} \! := \! (0,0,\dotsc,0)^{\mathrm{T}}$ $(\in \! 
\mathbb{R}^{N})$.
\end{ddddd}

Recalling Equations~(1.9), asymptotics (as $n \! \to \! \infty)$ for the 
Hankel determinant ratios
\begin{equation*}
H^{(-2n+1)}_{2n}/H^{(-2n)}_{2n} \qquad \qquad \text{and} \qquad \qquad 
H^{(-2n-1)}_{2n+1}/H^{(-2n)}_{2n+1}
\end{equation*}
have already been presented in Theorem~3.1, Equations~(3.48) and~(3.49), 
respectively. Via Equations~(1.18) and~(1.19), and asymptotics for 
$\xi^{(2n)}_{n}$ (resp., $\xi^{(2n+1)}_{-n-1})$ given in Appendix~A, 
Theorem~A.3 (resp., Appendix~B, Theorem~B.3), one proves the following 
theorem concerning asymptotics of the remaining Hankel determinant ratios 
$H^{(-2n-2)}_{2n+2}/H^{(-2n)}_{2n}$ (cf. Equation~(1.18)) and $H^{(-2n-2)
}_{2n+3}/H^{(-2n)}_{2n+1}$ (cf. Equation~(1.19)); see, in particular, 
Equations~(3.128) and~(3.129).
\begin{ddddd}
Let the external field $\widetilde{V} \colon \mathbb{R} \setminus \{0\} \! 
\to \! \mathbb{R}$ satisfy conditions~{\rm (1.20)--(1.22)}. Let $H^{(m)}_{k}$, 
$(m,k) \! \in \! \mathbb{Z} \times \mathbb{N}$, be the Hankel determinants 
associated with the bi-infinite, real-valued, strong moment sequence 
$\left\lbrace \mathstrut c_{j} \! = \! \int_{\mathbb{R}}s^{j} \exp (-n 
\widetilde{V}(s)) \, \md s, \, n \! \in \! \mathbb{N} \right\rbrace_{j 
\in \mathbb{Z}}$ defined in Equations~{\rm (1.1)}. Let the orthonormal 
$L$-polynomials (resp., monic orthogonal $L$-polynomials) $\lbrace \phi_{k}(z) 
\rbrace_{k \in \mathbb{Z}_{0}^{+}}$ (resp., $\lbrace \boldsymbol{\pi}_{k}(z) 
\rbrace_{k \in \mathbb{Z}_{0}^{+}})$ be as defined in Equations~{\rm (1.2)} 
and {\rm (1.3)} (resp., Equations~{\rm (1.4)} and~{\rm (1.5))}. Let 
$\overset{e}{\operatorname{Y}} \colon \mathbb{C} \setminus \mathbb{R} \! \to 
\! \operatorname{SL}_{2}(\mathbb{C})$ (resp., $\overset{o}{\operatorname{Y}} 
\colon \mathbb{C} \setminus \mathbb{R} \! \to \! \operatorname{SL}_{2}
(\mathbb{C}))$ be the (unique) solution of {\rm \pmb{RHP1}} (resp., 
{\rm \pmb{RHP2})} with integral representation {\rm (1.29)} (resp., 
{\rm (1.30))}. Define the density of the `even' (resp., `odd') equilibrium 
measure, $\md \mu_{V}^{e}(x)$ (resp., $\md \mu_{V}^{o}(x))$, as in 
Equation~{\rm (2.13)} (resp., Equation {\rm (2.19))}, and set $J_{q} \! := \! 
\operatorname{supp}(\mu_{V}^{q}) \! = \! \cup_{j=1}^{N+1}(b_{j-1}^{q},a_{j}^{
q})$, $q \! \in \! \{e,o\}$, where $\lbrace b_{j-1}^{e},a_{j}^{e} \rbrace_{j
=1}^{N+1}$ (resp., $\lbrace b_{j-1}^{o},a_{j}^{o} \rbrace_{j=1}^{N+1})$ 
satisfy the (real) $n$-dependent and locally solvable system of $2(N \! + 
\! 1)$ moment conditions~{\rm (2.12)} (resp., {\rm (2.18))}.

Suppose, furthermore, that $\widetilde{V} \colon \mathbb{R} \setminus \{0\} \! 
\to \! \mathbb{R}$ is regular, namely: {\rm (i)} $h_{V}^{q}(z) \! \not\equiv 
\! 0$ on $\overline{J_{q}} \! := \! \cup_{k=1}^{N+1}[b_{k-1}^{q},a_{k}^{q}]$, 
$q \! \in \! \{e,o\};$ and {\rm (ii)} the strict variational inequalities of 
Theorem~{\rm 3.1}, Equations~{\rm (3.46)} and~{\rm (3.47)}, are valid. Then,
\begin{align}
\dfrac{H^{(-2n-2)}_{2n+2}}{H^{(-2n)}_{2n}} \underset{n \to \infty}{=}& \,
\dfrac{\pi^{2} \eta^{e}_{+} \eta^{o}_{-}}{4 \mathbb{E}^{2}} \dfrac{
\boldsymbol{\theta}^{e}(-\boldsymbol{u}^{e}_{+}(\infty) \! - \! \frac{n}{2
\pi} \boldsymbol{\Omega}^{e} \! + \! \boldsymbol{d}_{e}) \boldsymbol{\theta}^{
o}(-\boldsymbol{u}^{o}_{+}(0) \! - \! \frac{1}{2 \pi}(n \! + \! \frac{1}{2})
\boldsymbol{\Omega}^{o} \! + \! \boldsymbol{d}_{o})}{\boldsymbol{\theta}^{e}
(\boldsymbol{u}^{e}_{+}(\infty) \! - \! \frac{n}{2 \pi} \boldsymbol{\Omega}^{
e} \! + \! \boldsymbol{d}_{e}) \boldsymbol{\theta}^{o}(\boldsymbol{u}^{o}_{+}
(0) \! - \! \frac{1}{2 \pi}(n \! + \! \frac{1}{2}) \boldsymbol{\Omega}^{o} \!
+ \! \boldsymbol{d}_{o})} \nonumber \\
&\times \dfrac{\boldsymbol{\theta}^{e}(\boldsymbol{u}^{e}_{+}(\infty) \! + \!
\boldsymbol{d}_{e}) \boldsymbol{\theta}^{o}(\boldsymbol{u}^{o}_{+}(0) \! + \!
\boldsymbol{d}_{o})}{\boldsymbol{\theta}^{e}(-\boldsymbol{u}^{e}_{+}(\infty)
\! + \! \boldsymbol{d}_{e}) \boldsymbol{\theta}^{o}(-\boldsymbol{u}^{o}_{+}(0)
\! + \! \boldsymbol{d}_{o})} \! \left\{1 \! + \! \dfrac{1}{n} \! \left(-
\dfrac{4 \mathfrak{Q}_{12}^{\flat}}{\eta^{e}_{+}} \dfrac{\boldsymbol{\theta}^{
e}(-\boldsymbol{u}^{e}_{+}(\infty) \! + \! \boldsymbol{d}_{e})}{\boldsymbol{
\theta}^{e}(\boldsymbol{u}^{e}_{+}(\infty) \! + \! \boldsymbol{d}_{e})}
\right. \right. \nonumber \\
&\left. \left. \times \, \dfrac{\boldsymbol{\theta}^{e}(\boldsymbol{u}^{e}_{+}
(\infty) \! - \! \frac{n}{2 \pi} \boldsymbol{\Omega}^{e} \! + \! \boldsymbol{
d}_{e})}{\boldsymbol{\theta}^{e}(-\boldsymbol{u}^{e}_{+}(\infty) \! - \! \frac{
n}{2 \pi} \boldsymbol{\Omega}^{e} \! + \! \boldsymbol{d}_{e})} \right) \! +
\! \dfrac{1}{(n \! + \! \frac{1}{2})} \! \left(-\dfrac{4 \mathbb{E}^{2}
\mathfrak{Q}_{12}^{\natural}}{\eta^{o}_{-}} \dfrac{\boldsymbol{\theta}^{o}(-
\boldsymbol{u}^{o}_{+}(0) \! + \! \boldsymbol{d}_{o})}{\boldsymbol{\theta}^{
o}(\boldsymbol{u}^{o}_{+}(0) \! + \! \boldsymbol{d}_{o})} \right. \right.
\nonumber \\
&\left. \left. \times \, \dfrac{\boldsymbol{\theta}^{o}(\boldsymbol{u}^{o}_{+}
(0) \! - \! \frac{1}{2 \pi}(n \! + \! \frac{1}{2}) \boldsymbol{\Omega}^{o} \!
+ \! \boldsymbol{d}_{o})}{\boldsymbol{\theta}^{o}(-\boldsymbol{u}^{o}_{+}(0)
\! - \! \frac{1}{2 \pi}(n \! + \! \frac{1}{2}) \boldsymbol{\Omega}^{o} \! + \!
\boldsymbol{d}_{o})} \right) \! + \! \mathcal{O} \left(\dfrac{c(n)}{n^{2}}
\right) \right\} \! \exp \! \left(n(\ell_{e} \! + \! \ell_{o}) \right), \\
\intertext{and}
\dfrac{H^{(-2n-2)}_{2n+3}}{H^{(-2n)}_{2n+1}} \underset{n \to \infty}{=}& \,
\dfrac{\pi^{2} \me^{n \ell_{o}} \eta^{o}_{-}}{4 \mathbb{E}^{2}} \dfrac{
\boldsymbol{\theta}^{o}(-\boldsymbol{u}^{o}_{+}(0) \! - \! \frac{1}{2 \pi}
(n \! + \! \frac{1}{2}) \boldsymbol{\Omega}^{o} \! + \! \boldsymbol{d}_{o})
\boldsymbol{\theta}^{o}(\boldsymbol{u}^{o}_{+}(0) \! + \! \boldsymbol{d}_{o}
)}{\boldsymbol{\theta}^{o}(\boldsymbol{u}^{o}_{+}(0) \! - \! \frac{1}{2 \pi}
(n \! + \! \frac{1}{2}) \boldsymbol{\Omega}^{o} \! + \! \boldsymbol{d}_{o})
\boldsymbol{\theta}^{o}(-\boldsymbol{u}^{o}_{+}(0) \! + \! \boldsymbol{d}_{o}
)} \nonumber \\
&\times \left\{1 \! + \! \dfrac{1}{(n \! + \! \frac{1}{2})} \! \left(-\dfrac{
4 \mathbb{E}^{2} \mathfrak{Q}^{\natural}_{12}}{\eta^{o}_{-}} \dfrac{
\boldsymbol{\theta}^{o}(\boldsymbol{u}^{o}_{+}(0) \! - \! \frac{1}{2 \pi}
(n \! + \! \frac{1}{2}) \boldsymbol{\Omega}^{o} \! + \! \boldsymbol{d}_{o})
\boldsymbol{\theta}^{o}(-\boldsymbol{u}^{o}_{+}(0) \! + \! \boldsymbol{d}_{o}
)}{\boldsymbol{\theta}^{o}(-\boldsymbol{u}^{o}_{+}(0) \! - \! \frac{1}{2 \pi}
(n \! + \! \frac{1}{2}) \boldsymbol{\Omega}^{o} \! + \! \boldsymbol{d}_{o})
\boldsymbol{\theta}^{o}(\boldsymbol{u}^{o}_{+}(0) \! + \! \boldsymbol{d}_{o})}
\right) \right. \nonumber \\
&\left. + \, \mathcal{O} \! \left(\dfrac{c(n)}{(n \! + \! \frac{1}{2})^{2}}
\right) \right\} \! \left\{h^{+} \! \left[\dfrac{\me^{-n \ell_{e}}}{\eta^{
e}_{+}} \dfrac{\boldsymbol{\theta}^{e}(\boldsymbol{u}^{e}_{+}(\infty) \! -
\! \frac{n}{2 \pi} \boldsymbol{\Omega}^{e} \! + \! \boldsymbol{d}_{e})
\boldsymbol{\theta}^{e}(-\boldsymbol{u}^{e}_{+}(\infty) \! + \! \boldsymbol{
d}_{e})}{\boldsymbol{\theta}^{e}(-\boldsymbol{u}^{e}_{+}(\infty) \! - \!
\frac{n}{2 \pi} \boldsymbol{\Omega}^{e} \! + \! \boldsymbol{d}_{e})
\boldsymbol{\theta}^{e}(\boldsymbol{u}^{e}_{+}(\infty) \! + \! \boldsymbol{
d}_{e})} \right. \right. \nonumber \\
&\left. \left. \times \left(1 \! + \! \dfrac{1}{n} \left(\dfrac{4 \mathfrak{
Q}^{\flat}_{12}}{\eta^{e}_{+}} \dfrac{\boldsymbol{\theta}^{e}(\boldsymbol{
u}^{e}_{+}(\infty) \! - \! \frac{n}{2 \pi} \boldsymbol{\Omega}^{e} \! + \!
\boldsymbol{d}_{e}) \boldsymbol{\theta}^{e}(-\boldsymbol{u}^{e}_{+}(\infty)
\! + \! \boldsymbol{d}_{e})}{\boldsymbol{\theta}^{e}(-\boldsymbol{u}^{e}_{+}
(\infty) \! - \! \frac{n}{2 \pi} \boldsymbol{\Omega}^{e} \! + \! \boldsymbol{
d}_{e}) \boldsymbol{\theta}^{e}(\boldsymbol{u}^{e}_{+}(\infty) \! + \!
\boldsymbol{d}_{e})} \right) \! + \! \mathcal{O} \! \left(\dfrac{c(n)}{n^{2}}
\right) \right) \right] \right\}^{-1},
\end{align}
with all parameters defined in Theorem~{\rm 3.1}, 
Equations~{\rm (3.52)--(3.100)}, and Theorem~{\rm 3.2}, Equations 
{\rm (3.105)}, {\rm (3.107)}, and~{\rm (3.108)}.
\end{ddddd}

\vspace*{1.70cm}

\pmb{\large Acknowledgements}

K.~T.-R.~McLaughlin was supported, in part, by National Science Foundation 
Grant Nos. DMS--0200749 and DMS--0451495; he also acknowledges the partial 
support of NATO Collaborative Linkage Grant `Orthogonal Polynomials: Theory, 
Applications and Generalizations', Ref. No. PST.CLG.979738. X.~Zhou was 
supported, in part, by National Science Foundation Grant No. DMS--0300844.
\clearpage
\section*{Appendix A. Asymptotics for $\overset{e}{\operatorname{Y}}(z)$}
\setcounter{section}{1}
\setcounter{z0}{1}
\setcounter{y0}{1}
\setcounter{equation}{0}
\renewcommand{\theequation}{\Alph{section}.\arabic{equation}}
In this appendix, large-$n$ asymptotics for the coefficients of the 
expansions~(3.1) and~(3.2) for the solution of \pmb{RHP1}, with integral 
representation~(1.29), are presented (for complete details and proofs, see 
\cite{a21}).
\begin{ey}
In \cite{a21}, the following formula was established for complex $z$ away 
from an open neighbourhood around $\overline{J_{e}} \! := \! \cup_{j=1}^{N+1}
[b_{j-1}^{e},a_{j}^{e}]$:
\begin{equation*}
\overset{e}{\mathrm{Y}}(z) \! = \! \me^{\frac{n \ell_{e}}{2} \sigma_{3}} 
\mathscr{R}^{e}(z) 
\overset{e}{m}^{\raise-1.0ex\hbox{$\scriptstyle \infty$}}(z) \me^{n(g^{e}(z)
-\frac{\ell_{e}}{2}+\int_{J_{e}} \ln (s) \psi_{V}^{e}(s) \, \md s) \sigma_{
3}},
\end{equation*}
where
$\ell_{e}$, $g^{e}(z)$, and $\psi_{V}^{e}(z)$ have been defined heretofore, 
the matrix-valued function 
$\overset{e}{m}^{\raise-1.0ex\hbox{$\scriptstyle \infty$}}(z)$ is determined 
explicitly, and the matrix-valued function $\mathscr{R}^{e}(z)$ solves a 
``small norm'' RHP for which a complete asymptotic expansion as $n \! \to \! 
\infty$ for $z \! \in \! \mathbb{C}$ is computed. The above explicit formula 
for the solution of \pmb{RHP1} yields an explicit asymptotic description 
for the even degree OLPs in question (for complete details and proofs, see 
\cite{a21}). \hfill $\blacksquare$
\end{ey}
\addtocounter{y0}{1}
\begin{ey}
In Theorems~A.1 and~A.2 below, only $n \! \to \! \infty$ asymptotics for 
$\mathscr{R}^{e,\infty}_{j}(n)$, $j \! = \! 1,2$, and $\mathscr{R}^{e,0}_{k-1}
(n)$, $k \! = \! 1,2,3$, are given, as they are all that are actually 
necessary in order to obtain the results of Theorems~3.1--3.3: finite-$n$ 
formulae for $\mathscr{R}^{e,\infty}_{j}(n)$, $j \! = \! 1,2$, and $\mathscr{
R}^{e,0}_{k-1}(n)$, $k \! = \! 1,2,3$, are given in \cite{a21}; furthermore, 
the `symbol' $c(n)$ appearing in the asymptotic (as $n \! \to \! \infty)$ 
expressions for $\mathscr{R}^{e,\infty}_{j}(n)$, $j \! = \! 1,2$, and 
$\mathscr{R}^{e,0}_{k-1}(n)$, $k \! = \! 1,2,3$, is to be understood as 
a uniformly bounded $(\mathcal{O}(1))$, 
$\operatorname{M}_{2}(\mathbb{C})$-valued, $n$-dependent function, that 
is, $\operatorname{M}_{2}(\mathbb{C}) \! \ni \! c(n) \! =_{n \to \infty} \! 
\mathcal{O}(1)$. \hfill $\blacksquare$
\end{ey}
\begin{dy}[{\rm \cite{a21}}]
Let the external field $\widetilde{V} \colon \mathbb{R} \setminus \{0\} \! \to 
\! \mathbb{R}$ satisfy conditions~{\rm (1.20)--(1.22)}. Let the orthonormal 
$L$-polynomials (resp., monic orthogonal $L$-polynomials) $\lbrace \phi_{k}(z) 
\rbrace_{k \in \mathbb{Z}_{0}^{+}}$ (resp., $\lbrace \boldsymbol{\pi}_{k}(z) 
\rbrace_{k \in \mathbb{Z}_{0}^{+}})$ be as defined in Equations~{\rm (1.2)} 
and~{\rm (1.3)} (resp., Equations~{\rm (1.4)} and~{\rm (1.5))}. Let 
$\overset{e}{\operatorname{Y}} \colon \mathbb{C} \setminus \mathbb{R} \! 
\to \! \operatorname{SL}_{2}(\mathbb{C})$ be the unique solution of 
{\rm \pmb{RHP1}} with integral representation {\rm (1.29)}. Define the 
density of the `even' equilibrium measure, $\md \mu_{V}^{e}(x)$, as in 
Equation~{\rm (2.13)}, and set $J_{e} \! := \! \operatorname{supp}(\mu_{V}^{
e}) \! = \! \cup_{j=1}^{N+1}(b_{j-1}^{e},a_{j}^{e})$, where $\lbrace b_{j-
1}^{e},a_{j}^{e} \rbrace_{j=1}^{N+1}$ satisfy the (real) $n$-dependent and 
locally solvable system of $2(N \! + \! 1)$ moment conditions~{\rm (2.12)}.

Suppose, furthermore, that $\widetilde{V} \colon \mathbb{R} \setminus \{ 0\} 
\! \to \! \mathbb{R}$ is regular, namely: {\rm (i)} $h_{V}^{e}(z) \! 
\not\equiv \! 0$ on $\overline{J_{e}} := \! \cup_{j=1}^{N+1}[b_{j-1}^{e},
a_{j}^{e}];$ {\rm (ii)}
\begin{equation*}
4 \int_{J_{e}} \ln (\vert x \! - \! s \vert) \, \md \mu_{V}^{e}(s) \! - \! 2 
\ln \vert x \vert \! - \! \widetilde{V}(x) \! - \! \ell_{e} \! = \! 0, \quad 
x \! \in \! \overline{J_{e}},
\end{equation*}
which defines the `even' variational constant, $\ell_{e}$ $(\in \! \mathbb{
R})$ (the same on each compact interval $[b_{j-1}^{e},a_{j}^{e}]$, $j \! = 
\! 1,\dotsc,N \! + \! 1)$, and
\begin{equation*}
4 \int_{J_{e}} \ln (\vert x \! - \! s \vert) \, \md \mu_{V}^{e}(s) \! - \! 2 
\ln \vert x \vert \! - \! \widetilde{V}(x) \! - \! \ell_{e} \! < \! 0, \quad 
x \! \in \! \mathbb{R} \setminus \overline{J_{e}} \,;
\end{equation*}
{\rm (iii)}
\begin{equation*}
g^{e}_{+}(z) \! + \! g^{e}_{-}(z) \! - \! \widetilde{V}(z) \! - \! \ell_{e} \! 
+ \! 2Q_{e} \! < \! 0, \quad z \! \in \! \mathbb{R} \setminus \overline{J_{e}},
\end{equation*}
where $g^{e}(z)$ is defined by Equation~{\rm (2.15)}, $g^{e}_{\pm}(z) \! := 
\! \lim_{\varepsilon \downarrow 0}g^{e}(z \! \pm \! \mi \varepsilon)$, and 
$Q_{e}$ is defined in Equation~{\rm (2.16);} and {\rm (iv)}
\begin{equation*}
\mi (g^{e}_{+}(z) \! - \! g^{e}_{-}(z))^{\prime} \! > \! 0, \quad z \! \in \! 
J_{e}.
\end{equation*}
Then,
\begin{equation*}
\overset{e}{\mathrm{Y}}(z)z^{-n \sigma_{3}} \underset{z \to \infty}{=}
\mathrm{I} \! + \! \dfrac{1}{z} \mathrm{Y}^{e,\infty}_{1} \! + \! \dfrac{1}{
z^{2}} \mathrm{Y}^{e,\infty}_{2} \! + \! \mathcal{O} \! \left(\dfrac{1}{z^{
3}} \right),
\end{equation*}
where
\begin{align*}
(\mathrm{Y}^{e,\infty}_{1})_{11} &= -2n \int_{J_{e}}s \psi_{V}^{e}(s) \, \md 
s \! + \! (\overset{e}{m}_{1}^{\raise-1.0ex\hbox{$\scriptstyle \infty$}})_{11} 
\! + \! (\mathscr{R}^{e,\infty}_{1}(n))_{11}, \\
(\mathrm{Y}^{e,\infty}_{1})_{12} &= \me^{n \ell_{e}} \! \left(
(\overset{e}{m}_{1}^{\raise-1.0ex\hbox{$\scriptstyle \infty$}})_{12} \! + \! 
(\mathscr{R}^{e,\infty}_{1}(n))_{12} \right), \\
(\mathrm{Y}^{e,\infty}_{1})_{21} &= \me^{-n \ell_{e}} \! \left(
(\overset{e}{m}_{1}^{\raise-1.0ex\hbox{$\scriptstyle \infty$}})_{21} \! + \! 
(\mathscr{R}^{e,\infty}_{1}(n))_{21} \right), \\
(\mathrm{Y}^{e,\infty}_{1})_{22} &= 2n \int_{J_{e}}s \psi_{V}^{e}(s) \, \md s 
\! + \! (\overset{e}{m}_{1}^{\raise-1.0ex\hbox{$\scriptstyle \infty$}})_{22} 
\! + \! (\mathscr{R}^{e,\infty}_{1}(n))_{22}, \\
(\mathrm{Y}^{e,\infty}_{2})_{11} &= 2n^{2} \! \left(\int_{J_{e}}s \psi_{V}^{e}
(s) \, \md s \right)^{2} \! - \! n \int_{J_{e}}s^{2} \psi_{V}^{e}(s) \, \md 
s \! - \! 2n \! \left(
(\overset{e}{m}_{1}^{\raise-1.0ex\hbox{$\scriptstyle \infty$}})_{11} \! + \! 
(\mathscr{R}^{e,\infty}_{1}(n))_{11} \right) \! \int_{J_{e}}s \psi_{V}^{e}(s) 
\, \md s \\
&+(\overset{e}{m}_{2}^{\raise-1.0ex\hbox{$\scriptstyle \infty$}})_{11} \! + \! 
(\mathscr{R}^{e,\infty}_{2}(n))_{11} \! + \! (\mathscr{R}^{e,\infty}_{1}(n))_{
11}(\overset{e}{m}_{1}^{\raise-1.0ex\hbox{$\scriptstyle \infty$}})_{11} \! + 
\! (\mathscr{R}^{e,\infty}_{1}(n))_{12}
(\overset{e}{m}_{1}^{\raise-1.0ex\hbox{$\scriptstyle \infty$}})_{21}, \\
(\mathrm{Y}^{e,\infty}_{2})_{12} &= \me^{n \ell_{e}} \! \left(2n \! \left(
(\overset{e}{m}_{1}^{\raise-1.0ex\hbox{$\scriptstyle \infty$}})_{12} \! + \! 
(\mathscr{R}^{e,\infty}_{1}(n))_{12} \right) \! \int_{J_{e}}s \psi_{V}^{e}(s) 
\, \md s \! + \! 
(\overset{e}{m}_{2}^{\raise-1.0ex\hbox{$\scriptstyle \infty$}})_{12} \! + \!
(\mathscr{R}^{e,\infty}_{2}(n))_{12} \right. \\
&\left. + \, (\mathscr{R}^{e,\infty}_{1}(n))_{11}
(\overset{e}{m}_{1}^{\raise-1.0ex\hbox{$\scriptstyle \infty$}})_{12} \! + \! 
(\mathscr{R}^{e,\infty}_{1}(n))_{12} 
(\overset{e}{m}_{1}^{\raise-1.0ex\hbox{$\scriptstyle \infty$}})_{22} \right), 
\\
(\mathrm{Y}^{e,\infty}_{2})_{21} &= \me^{-n \ell_{e}} \! \left(-2n \! \left(
(\overset{e}{m}_{1}^{\raise-1.0ex\hbox{$\scriptstyle \infty$}})_{21} \! + \!
(\mathscr{R}^{e,\infty}_{1}(n))_{21} \right) \! \int_{J_{e}}s \psi_{V}^{e}(s)
\, \md s \! + \!
(\overset{e}{m}_{2}^{\raise-1.0ex\hbox{$\scriptstyle \infty$}})_{21} \! + \!
(\mathscr{R}^{e,\infty}_{2}(n))_{21} \right. \\
&\left. + \, (\mathscr{R}^{e,\infty}_{1}(n))_{21}
(\overset{e}{m}_{1}^{\raise-1.0ex\hbox{$\scriptstyle \infty$}})_{11} \! + \!
(\mathscr{R}^{e,\infty}_{1}(n))_{22}
(\overset{e}{m}_{1}^{\raise-1.0ex\hbox{$\scriptstyle \infty$}})_{21} \right),
\\
(\mathrm{Y}^{e,\infty}_{2})_{22} &= 2n^{2} \! \left(\int_{J_{e}}s \psi_{V}^{e}
(s) \, \md s \right)^{2} \! + \! n \int_{J_{e}}s^{2} \psi_{V}^{e}(s) \, \md s
\! + \! 2n \! \left(
(\overset{e}{m}_{1}^{\raise-1.0ex\hbox{$\scriptstyle \infty$}})_{22} \! + \!
(\mathscr{R}^{e,\infty}_{1}(n))_{22} \right) \! \int_{J_{e}}s \psi_{V}^{e}(s)
\, \md s \\
&+(\overset{e}{m}_{2}^{\raise-1.0ex\hbox{$\scriptstyle \infty$}})_{22} \! + \!
(\mathscr{R}^{e,\infty}_{2}(n))_{22} \! + \! (\mathscr{R}^{e,\infty}_{1}(n))_{
21}(\overset{e}{m}_{1}^{\raise-1.0ex\hbox{$\scriptstyle \infty$}})_{12} \! +
\! (\mathscr{R}^{e,\infty}_{1}(n))_{22}
(\overset{e}{m}_{1}^{\raise-1.0ex\hbox{$\scriptstyle \infty$}})_{22},
\end{align*}
with
\begin{align*}
(\overset{e}{m}_{1}^{\raise-1.0ex\hbox{$\scriptstyle \infty$}})_{11} &=\dfrac{
\boldsymbol{\theta}^{e}(\boldsymbol{u}^{e}_{+}(\infty) \! + \! \boldsymbol{d}_{
e})}{\boldsymbol{\theta}^{e}(\boldsymbol{u}^{e}_{+}(\infty) \! - \! \frac{n}{2
\pi} \boldsymbol{\Omega}^{e} \! + \! \boldsymbol{d}_{e})} \! \left(\dfrac{
\theta^{e}_{\infty}(1,1,\boldsymbol{\Omega}^{e}) \alpha^{e}_{\infty}(1,1,\pmb{
0}) \! - \! \alpha^{e}_{\infty}(1,1,\boldsymbol{\Omega}^{e}) \theta^{e}_{
\infty}(1,1,\pmb{0})}{(\theta^{e}_{\infty}(1,1,\pmb{0}))^{2}} \right), \\
(\overset{e}{m}_{1}^{\raise-1.0ex\hbox{$\scriptstyle \infty$}})_{12} &=\dfrac{
1}{4 \mi} \! \left(\sum_{k=1}^{N+1}(b_{k-1}^{e} \! - \! a_{k}^{e}) \right) \!
\dfrac{\boldsymbol{\theta}^{e}(\boldsymbol{u}^{e}_{+}(\infty) \! + \!
\boldsymbol{d}_{e}) \theta^{e}_{\infty}(-1,1,\boldsymbol{\Omega}^{e})}{
\boldsymbol{\theta}^{e}(\boldsymbol{u}^{e}_{+}(\infty) \! - \! \frac{n}{2 \pi}
\boldsymbol{\Omega}^{e} \! + \! \boldsymbol{d}_{e}) \theta^{e}_{\infty}(-1,1,
\pmb{0})}, \\
(\overset{e}{m}_{1}^{\raise-1.0ex\hbox{$\scriptstyle \infty$}})_{21}&=-\dfrac{
1}{4 \mi} \! \left(\sum_{k=1}^{N+1}(b_{k-1}^{e} \! - \! a_{k}^{e}) \right)
\! \dfrac{\boldsymbol{\theta}^{e}(\boldsymbol{u}^{e}_{+}(\infty) \! + \!
\boldsymbol{d}_{e}) \theta^{e}_{\infty}(1,-1,\boldsymbol{\Omega}^{e})}{
\boldsymbol{\theta}^{e}(-\boldsymbol{u}^{e}_{+}(\infty) \! - \! \frac{n}{2
\pi} \boldsymbol{\Omega}^{e} \! - \! \boldsymbol{d}_{e}) \theta^{e}_{\infty}
(1,-1,\pmb{0})}, \\
(\overset{e}{m}_{1}^{\raise-1.0ex\hbox{$\scriptstyle \infty$}})_{22} &=\dfrac{
\boldsymbol{\theta}^{e}(\boldsymbol{u}^{e}_{+}(\infty) \! + \! \boldsymbol{
d}_{e})}{\boldsymbol{\theta}^{e}(-\boldsymbol{u}^{e}_{+}(\infty) \! - \! 
\frac{n}{2 \pi} \boldsymbol{\Omega}^{e} \! - \! \boldsymbol{d}_{e})} \! \left(
\dfrac{\theta^{e}_{\infty}(-1,-1,\boldsymbol{\Omega}^{e}) \alpha^{e}_{\infty}
(-1,-1,\pmb{0}) \! - \! \alpha^{e}_{\infty}(-1,-1,\boldsymbol{\Omega}^{e})
\theta^{e}_{\infty}(-1,-1,\pmb{0})}{(\theta^{e}_{\infty}(-1,-1,\pmb{0}))^{2}}
\right), \\
(\overset{e}{m}_{2}^{\raise-1.0ex\hbox{$\scriptstyle \infty$}})_{11} &=\left(
\theta^{e}_{\infty}(1,1,\boldsymbol{\Omega}^{e}) \! \left(\beta^{e}_{\infty}
(1,1,\pmb{0}) \theta^{e}_{\infty}(1,1,\pmb{0}) \! + \! (\alpha^{e}_{\infty}(1,
1,\pmb{0}))^{2} \right) \! - \! \alpha^{e}_{\infty}(1,1,\boldsymbol{\Omega}^{e}
) \alpha^{e}_{\infty}(1,1,\pmb{0}) \theta^{e}_{\infty}(1,1,\pmb{0}) \right.
\nonumber \\
&\left. - \, \beta^{e}_{\infty}(1,1,\boldsymbol{\Omega}^{e})(\theta^{e}_{
\infty}(1,1,\pmb{0}))^{2} \right) \! \dfrac{(\theta^{e}_{\infty}(1,1,\pmb{0})
)^{-3} \boldsymbol{\theta}^{e}(\boldsymbol{u}^{e}_{+}(\infty) \! + \!
\boldsymbol{d}_{e})}{\boldsymbol{\theta}^{e}(\boldsymbol{u}^{e}_{+}(\infty) \!
- \! \frac{n}{2 \pi} \boldsymbol{\Omega}^{e} \! + \! \boldsymbol{d}_{e})} \!
+ \! \dfrac{1}{32} \! \left(\sum_{k=1}^{N+1}(b_{k-1}^{e} \! - \! a_{k}^{e})
\right)^{2}, \\
(\overset{e}{m}_{2}^{\raise-1.0ex\hbox{$\scriptstyle \infty$}})_{12} &=
\dfrac{\boldsymbol{\theta}^{e}(\boldsymbol{u}^{e}_{+}(\infty) \! + \!
\boldsymbol{d}_{e})}{\boldsymbol{\theta}^{e}(\boldsymbol{u}^{e}_{+}(\infty) \!
- \! \frac{n}{2 \pi} \boldsymbol{\Omega}^{e} \! + \! \boldsymbol{d}_{e})} \!
\left(\! \left(\dfrac{\theta^{e}_{\infty}(-1,1,\boldsymbol{\Omega}^{e})
\alpha^{e}_{\infty}(-1,1,\pmb{0}) \! - \! \alpha^{e}_{\infty}(-1,1,\boldsymbol{
\Omega}^{e}) \theta^{e}_{\infty}(-1,1,\pmb{0})}{(\theta^{e}_{\infty}(-1,1,\pmb{
0}))^{2}} \right) \right. \nonumber \\
&\left. \times \, \dfrac{1}{4 \mi} \! \left(\sum_{k=1}^{N+1}(b_{k-1}^{e} \! -
\! a_{k}^{e}) \right) \! + \! \dfrac{1}{8 \mi} \! \left(\sum_{k=1}^{N+1}((b_{k
-1}^{e})^{2} \! - \! (a_{k}^{e})^{2}) \right) \! \dfrac{\theta^{e}_{\infty}(-
1,1,\boldsymbol{\Omega}^{e})}{\theta^{e}_{\infty}(-1,1,\pmb{0})} \right), \\
(\overset{e}{m}_{2}^{\raise-1.0ex\hbox{$\scriptstyle \infty$}})_{21} &=-
\dfrac{\boldsymbol{\theta}^{e}(\boldsymbol{u}^{e}_{+}(\infty) \! + \!
\boldsymbol{d}_{e})}{\boldsymbol{\theta}^{e}(-\boldsymbol{u}^{e}_{+}(\infty)
\! - \! \frac{n}{2 \pi} \boldsymbol{\Omega}^{e} \! - \! \boldsymbol{d}_{e})}
\! \left(\! \left(\dfrac{\theta^{e}_{\infty}(1,-1,\boldsymbol{\Omega}^{e})
\alpha^{e}_{\infty}(1,-1,\pmb{0}) \! - \! \alpha^{e}_{\infty}(1,-1,\boldsymbol{
\Omega}^{e}) \theta^{e}_{\infty}(1,-1,\pmb{0})}{(\theta^{e}_{\infty}(1,-1,\pmb{
0}))^{2}} \right) \right. \nonumber \\
&\left. \times \, \dfrac{1}{4 \mi} \! \left(\sum_{k=1}^{N+1}(b_{k-1}^{e} \! -
\! a_{k}^{e}) \right) \! + \! \dfrac{1}{8 \mi} \! \left(\sum_{k=1}^{N+1}((b_{k
-1}^{e})^{2} \! - \! (a_{k}^{e})^{2}) \right) \! \dfrac{\theta^{e}_{\infty}(1,
-1,\boldsymbol{\Omega}^{e})}{\theta^{e}_{\infty}(1,-1,\pmb{0})} \right), \\
(\overset{e}{m}_{2}^{\raise-1.0ex\hbox{$\scriptstyle \infty$}})_{22} &=\left(
\theta^{e}_{\infty}(-1,-1,\boldsymbol{\Omega}^{e}) \! \left(\beta^{e}_{\infty}
(-1,-1,\pmb{0}) \theta^{e}_{\infty}(-1,-1,\pmb{0}) \! + \! (\alpha^{e}_{\infty}
(-1,-1,\pmb{0}))^{2} \right) \! - \! \alpha^{e}_{\infty}(-1,-1,\boldsymbol{
\Omega}^{e}) \right. \nonumber \\
&\left. \times \, \alpha^{e}_{\infty}(-1,-1,\pmb{0}) \theta^{e}_{\infty}(-1,-1,
\pmb{0}) \! - \! \beta^{e}_{\infty}(-1,-1,\boldsymbol{\Omega}^{e})(\theta^{e}_{
\infty}(-1,-1,\pmb{0}))^{2} \right) \! (\theta^{e}_{\infty}(-1,-1,\pmb{0}))^{-
3} \nonumber \\
&\times \dfrac{\boldsymbol{\theta}^{e}(\boldsymbol{u}^{e}_{+}(\infty) \! + \!
\boldsymbol{d}_{e})}{\boldsymbol{\theta}^{e}(-\boldsymbol{u}^{e}_{+}(\infty)
\! - \! \frac{n}{2 \pi} \boldsymbol{\Omega}^{e} \! - \! \boldsymbol{d}_{e})}
\! + \! \dfrac{1}{32} \! \left(\sum_{k=1}^{N+1}(b_{k-1}^{e} \! - \! a_{k}^{e})
\right)^{2},
\end{align*}
and $(\star)_{ij}$, $i,j \! = \! 1,2$, denoting the $(i \, j)$-element of 
$\star$, where, for $\varepsilon_{1},\varepsilon_{2} \! = \! \pm 1$,
\begin{gather*}
\theta^{e}_{\infty}(\varepsilon_{1},\varepsilon_{2},\pmb{\mathscr{Z}}) \! :=
\! \boldsymbol{\theta}^{e}(\varepsilon_{1} \boldsymbol{u}^{e}_{+}(\infty) \! 
- \! \tfrac{n}{2 \pi} \pmb{\mathscr{Z}} \! + \! \varepsilon_{2} \boldsymbol{
d}_{e}), \\
\alpha^{e}_{\infty}(\varepsilon_{1},\varepsilon_{2},\pmb{\mathscr{Z}}) \! 
:= \! 2 \pi \mi \varepsilon_{1} \sum_{m \in \mathbb{Z}^{N}} (m,\widehat{
\boldsymbol{\alpha}}^{e}_{\infty}) \me^{2 \pi \mi (m,\varepsilon_{1} 
\boldsymbol{u}^{e}_{+}(\infty)-\frac{n}{2 \pi} \pmb{\mathscr{Z}}+\varepsilon_{
2} \boldsymbol{d}_{e})+ \pi \mi (m,\tau^{e}m)},
\end{gather*}
where $\boldsymbol{u}^{e}_{+}(\infty) \! = \! \int_{a_{N+1}^{e}}^{\infty^{+}} 
\boldsymbol{\omega}^{e}$, $\widehat{\boldsymbol{\alpha}}^{e}_{\infty} \! = 
\! (\widehat{\alpha}_{\infty,1}^{e},\widehat{\alpha}^{e}_{\infty,2},\dotsc,
\widehat{\alpha}^{e}_{\infty,N})$, with $\widehat{\alpha}^{e}_{\infty,j} \! 
:= \! c_{j1}^{e}$, $j \! = \! 1,\dotsc,N$, and
\begin{gather*}
\beta^{e}_{\infty}(\varepsilon_{1},\varepsilon_{2},\pmb{\mathscr{Z}}) \! := \!
2 \pi \sum_{m \in \mathbb{Z}^{N}} \! \left(\pi (m,\widehat{\boldsymbol{\alpha}
}^{e}_{\infty})^{2} \! + \! \mi \varepsilon_{1}(m,\widehat{\boldsymbol{\beta}
}^{e}_{\infty}) \right) \! \me^{2 \pi \mi (m,\varepsilon_{1} \boldsymbol{u}^{
e}_{+}(\infty)-\frac{n}{2 \pi} \pmb{\mathscr{Z}}+\varepsilon_{2} \boldsymbol{
d}_{e})+ \pi \mi (m,\tau^{e}m)},
\end{gather*}
where $\widehat{\boldsymbol{\beta}}^{e}_{\infty} \! = \! (\widehat{\beta}^{e}_{
\infty,1},\widehat{\beta}^{e}_{\infty,2},\dotsc,\widehat{\beta}^{e}_{\infty,N}
)$, with $\widehat{\beta}^{e}_{\infty,j} \! := \! \tfrac{1}{2}(c^{e}_{j2} \! + 
\! \tfrac{1}{2}c_{j1}^{e} \sum_{i=1}^{N+1}(b_{i-1}^{e} \! + \! a_{i}^{e}))$, 
$j \! = \! 1,\dotsc,N$, where $c_{j1}^{e},c_{j2}^{e}$, $j \! = \! 1,\dotsc,N$, 
are obtained {}from Equations~{\rm (E1)} and~{\rm (E2)}, and
\begin{align*}
\mathscr{R}^{e,\infty}_{1}(n) \underset{n \to \infty}{=}& \, \dfrac{1}{n} 
\mathlarger{\sum_{j=1}^{N+1}} \! \left(\dfrac{\left(\mathscr{B}^{e}(a_{j}^{
e}) \widehat{\alpha}_{0}^{e}(a_{j}^{e}) \! - \! \mathscr{A}^{e}(a_{j}^{e}) 
\widehat{\alpha}_{1}^{e}(a_{j}^{e}) \right)}{(\widehat{\alpha}_{0}^{e}(a_{j}^{
e}))^{2}} \! + \! \dfrac{\left(\mathscr{B}^{e}(b_{j-1}^{e}) \widehat{\alpha}_{
0}^{e}(b_{j-1}^{e}) \! - \! \mathscr{A}^{e}(b_{j-1}^{e}) \widehat{\alpha}_{
1}^{e}(b_{j-1}^{e}) \right)}{(\widehat{\alpha}_{0}^{e}(b_{j-1}^{e}))^{2}} 
\right) \! + \! \mathcal{O} \! \left(\dfrac{c(n)}{n^{2}} \right), \\
\mathscr{R}^{e,\infty}_{2}(n) \underset{n \to \infty}{=}& \, \dfrac{1}{n} 
\mathlarger{\sum_{j=1}^{N+1}} \! \left(\dfrac{\left(\widehat{\alpha}_{0}^{
e}(b_{j-1}^{e}) \mathscr{A}^{e}(b_{j-1}^{e}) \! + \! b_{j-1}^{e} \! \left(
\mathscr{B}^{e}(b_{j-1}^{e}) \widehat{\alpha}_{0}^{e}(b_{j-1}^{e}) \! - \! 
\mathscr{A}^{e}(b_{j-1}^{e}) \widehat{\alpha}_{1}^{e}(b_{j-1}^{e}) \right) 
\right)}{(\widehat{\alpha}_{0}^{e}(b_{j-1}^{e}))^{2}} \right. \nonumber \\
+&\left. \, \dfrac{\left(\widehat{\alpha}_{0}^{e}(a_{j}^{e}) \mathscr{A}^{e}
(a_{j}^{e}) \! + \! a_{j}^{e} \! \left(\mathscr{B}^{e}(a_{j}^{e}) \widehat{
\alpha}_{0}^{e}(a_{j}^{e}) \! - \! \mathscr{A}^{e}(a_{j}^{e}) \widehat{
\alpha}_{1}^{e}(a_{j}^{e}) \right) \right)}{(\widehat{\alpha}_{0}^{e}(a_{j}^{
e}))^{2}} \right) \! + \! \mathcal{O} \! \left(\dfrac{c(n)}{n^{2}} \right),
\end{align*}
with all parameters defined in Theorem~{\rm 3.1}, 
Equations~{\rm (3.52)--(3.100)}.
\end{dy}
\addtocounter{z0}{1}
\begin{dy}[{\rm \cite{a21}}]
Let all the conditions stated in Theorem~{\rm A.1} be valid, and let all 
parameters be as defined therein. Let $\overset{e}{\operatorname{Y}} \colon 
\mathbb{C} \setminus \mathbb{R} \! \to \! \operatorname{SL}_{2}(\mathbb{C})$ 
be the unique solution of {\rm \pmb{RHP1}} with integral 
representation~{\rm (1.29)}. Then,
\begin{equation*}
\overset{e}{\mathrm{Y}}(z)z^{n \sigma_{3}} \underset{z \to 0}{=} \, \mathrm{
Y}^{e,0}_{0} \! + \! z \mathrm{Y}^{e,0}_{1} \! + \! z^{2} \mathrm{Y}^{e,0}_{
2} \! + \! \mathcal{O}(z^{3}),
\end{equation*}
where
\begin{align*}
(\mathrm{Y}^{e,0}_{0})_{11} =& \, (\widehat{Q}^{e}_{0})_{11} \me^{2n(\int_{J_{
e}} \ln (\lvert s \rvert) \psi_{V}^{e}(s) \, \md s + \mi \pi \int_{J_{e} \cap
\mathbb{R}_{+}} \psi_{V}^{e}(s) \, \md s)}, \\
(\mathrm{Y}^{e,0}_{0})_{12} =& \, (\widehat{Q}^{e}_{0})_{12} \me^{n(\ell_{e}-
2 \int_{J_{e}} \ln (\lvert s \rvert) \psi_{V}^{e}(s) \, \md s-\mi 2 \pi \int_{
J_{e} \cap \mathbb{R}_{+}} \psi_{V}^{e}(s) \, \md s)}, \\
(\mathrm{Y}^{e,0}_{0})_{21} =& \, (\widehat{Q}^{e}_{0})_{21} \me^{-n(\ell_{e}-
2 \int_{J_{e}} \ln (\lvert s \rvert) \psi_{V}^{e}(s) \, \md s-\mi 2\pi \int_{
J_{e} \cap \mathbb{R}_{+}} \psi_{V}^{e}(s) \, \md s)}, \\
(\mathrm{Y}^{e,0}_{0})_{22} =& \, (\widehat{Q}^{e}_{0})_{22} \me^{-2n(\int_{J_{
e}} \ln (\lvert s \rvert) \psi_{V}^{e}(s) \, \md s+\mi \pi \int_{J_{e} \cap
\mathbb{R}_{+}} \psi_{V}^{e}(s) \, \md s)}, \\
(\mathrm{Y}^{e,0}_{1})_{11} =& \, \left(\! (\widehat{Q}^{e}_{1})_{11} \! - \!
2n(\widehat{Q}^{e}_{0})_{11} \int_{J_{e}}s^{-1} \psi_{V}^{e}(s) \, \md s
\right) \! \me^{2n(\int_{J_{e}} \ln (\lvert s \rvert) \psi_{V}^{e}(s) \, \md
s+\mi \pi \int_{J_{e} \cap \mathbb{R}_{+}} \psi_{V}^{e}(s) \, \md s)}, \\
(\mathrm{Y}^{e,0}_{1})_{12} =& \, \left(\! (\widehat{Q}^{e}_{1})_{12} \! + \!
2n(\widehat{Q}^{e}_{0})_{12} \int_{J_{e}}s^{-1} \psi_{V}^{e}(s) \, \md s
\right) \! \me^{n(\ell_{e}-2\int_{J_{e}} \ln (\lvert s \rvert) \psi_{V}^{e}(s)
\, \md s-\mi 2 \pi \int_{J_{e} \cap \mathbb{R}_{+}} \psi_{V}^{e}(s) \, \md s)
}, \\
(\mathrm{Y}^{e,0}_{1})_{21} =& \, \left(\! (\widehat{Q}^{e}_{1})_{21} \! - \!
2n(\widehat{Q}^{e}_{0})_{21} \int_{J_{e}}s^{-1} \psi_{V}^{e}(s) \, \md s
\right) \! \me^{-n(\ell_{e}-2\int_{J_{e}} \ln (\lvert s \rvert) \psi_{V}^{e}
(s) \, \md s-\mi 2 \pi \int_{J_{e} \cap \mathbb{R}_{+}} \psi_{V}^{e}(s) \, 
\md s)}, \\
(\mathrm{Y}^{e,0}_{1})_{22} =& \, \left(\! (\widehat{Q}^{e}_{1})_{22} \! + \!
2n(\widehat{Q}^{e}_{0})_{22} \int_{J_{e}}s^{-1} \psi_{V}^{e}(s) \, \md s
\right) \! \me^{-2n(\int_{J_{e}} \ln (\lvert s \rvert) \psi_{V}^{e}(s) \, \md
s+\mi \pi \int_{J_{e} \cap \mathbb{R}_{+}} \psi_{V}^{e}(s) \, \md s)}, \\
(\mathrm{Y}^{e,0}_{2})_{11} =& \left(\! (\widehat{Q}^{e}_{2})_{11} \! - \! 2n
(\widehat{Q}^{e}_{1})_{11} \int_{J_{e}}s^{-1} \psi_{V}^{e}(s) \, \md s \! + \!
(\widehat{Q}^{e}_{0})_{11} \! \left(2n^{2} \! \left(\int_{J_{e}}s^{-1} \psi_{
V}^{e}(s) \, \md s \right)^{2} \right. \right. \\
-& \left. \left. n \int_{J_{e}}s^{-2} \psi_{V}^{e}(s) \, \md s \right)
\right) \! \me^{2n(\int_{J_{e}} \ln (\lvert s \rvert) \psi_{V}^{e}(s) \, \md
s+\mi \pi \int_{J_{e} \cap \mathbb{R}_{+}} \psi_{V}^{e}(s) \, \md s)}, \\
(\mathrm{Y}^{e,0}_{2})_{12} =& \left(\! (\widehat{Q}^{e}_{2})_{12} \! + \! 2n
(\widehat{Q}^{e}_{1})_{12} \int_{J_{e}}s^{-1} \psi_{V}^{e}(s) \, \md s \! + \!
(\widehat{Q}^{e}_{0})_{12} \! \left(2n^{2} \! \left(\int_{J_{e}}s^{-1} \psi_{
V}^{e}(s) \, \md s \right)^{2} \right. \right. \\
+& \left. \left. n \int_{J_{e}}s^{-2} \psi_{V}^{e}(s) \, \md s \right) 
\right) \! \me^{n(\ell_{e}-2\int_{J_{e}} \ln (\lvert s \rvert) \psi_{V}^{e}(s)
\, \md s-\mi 2 \pi \int_{J_{e} \cap \mathbb{R}_{+}} \psi_{V}^{e}(s) \, \md s)
}, \\
(\mathrm{Y}^{e,0}_{2})_{21} =& \left(\! (\widehat{Q}^{e}_{2})_{21} \! - \! 2n
(\widehat{Q}^{e}_{1})_{21} \int_{J_{e}}s^{-1} \psi_{V}^{e}(s) \, \md s \! + \!
(\widehat{Q}^{e}_{0})_{21} \! \left(2n^{2} \! \left(\int_{J_{e}}s^{-1} \psi_{
V}^{e}(s) \, \md s \right)^{2} \right. \right. \\
-& \left. \left. n \int_{J_{e}}s^{-2} \psi_{V}^{e}(s) \, \md s \right)
\right) \! \me^{-n(\ell_{e}-2 \int_{J_{e}} \ln (\lvert s \rvert) \psi_{V}^{e}
(s) \, \md s-\mi 2 \pi \int_{J_{e} \cap \mathbb{R}_{+}} \psi_{V}^{e}(s) \, \md
s)}, \\
(\mathrm{Y}^{e,0}_{2})_{22} =& \left(\! (\widehat{Q}^{e}_{2})_{22} \! + \! 2n
(\widehat{Q}^{e}_{1})_{22} \int_{J_{e}}s^{-1} \psi_{V}^{e}(s) \, \md s \! + \!
(\widehat{Q}^{e}_{0})_{22} \! \left(2n^{2} \! \left(\int_{J_{e}}s^{-1} \psi_{
V}^{e}(s) \, \md s \right)^{2} \right. \right. \\
+& \left. \left. n \int_{J_{e}}s^{-2} \psi_{V}^{e}(s) \, \md s \right)
\right) \! \me^{-2n(\int_{J_{e}} \ln (\lvert s \rvert) \psi_{V}^{e}(s) \, 
\md s+\mi \pi \int_{J_{e} \cap \mathbb{R}_{+}} \psi_{V}^{e}(s) \, \md s)},
\end{align*}
with
\begin{align*}
(\widehat{Q}^{e}_{0})_{11} :=& \, 
(\overset{e}{m}_{0}^{\raise-1.0ex\hbox{$\scriptstyle 0$}})_{11} \! \left(1 \!
+ \! (\mathscr{R}^{e,0}_{0}(n))_{11} \right) \! + \! (\mathscr{R}^{e,0}_{0}
(n))_{12}(\overset{e}{m}_{0}^{\raise-1.0ex\hbox{$\scriptstyle 0$}})_{21}, \\
(\widehat{Q}^{e}_{0})_{12} :=& \,
(\overset{e}{m}_{0}^{\raise-1.0ex\hbox{$\scriptstyle 0$}})_{12} \! \left(1 \!
+ \! (\mathscr{R}^{e,0}_{0}(n))_{11} \right) \! + \! (\mathscr{R}^{e,0}_{0}
(n))_{12}(\overset{e}{m}_{0}^{\raise-1.0ex\hbox{$\scriptstyle 0$}})_{22}, \\
(\widehat{Q}^{e}_{0})_{21} :=& \,
(\overset{e}{m}_{0}^{\raise-1.0ex\hbox{$\scriptstyle 0$}})_{21} \! \left(1 \!
+ \! (\mathscr{R}^{e,0}_{0}(n))_{22} \right) \! + \! (\mathscr{R}^{e,0}_{0}
(n))_{21}(\overset{e}{m}_{0}^{\raise-1.0ex\hbox{$\scriptstyle 0$}})_{11}, \\
(\widehat{Q}^{e}_{0})_{22} :=& \,
(\overset{e}{m}_{0}^{\raise-1.0ex\hbox{$\scriptstyle 0$}})_{22} \! \left(1 \!
+ \! (\mathscr{R}^{e,0}_{0}(n))_{22} \right) \! + \! (\mathscr{R}^{e,0}_{0}
(n))_{21}(\overset{e}{m}_{0}^{\raise-1.0ex\hbox{$\scriptstyle 0$}})_{12}, \\
(\widehat{Q}^{e}_{1})_{11} :=& \,
(\overset{e}{m}_{1}^{\raise-1.0ex\hbox{$\scriptstyle 0$}})_{11} \! \left(1 \!
+ \! (\mathscr{R}^{e,0}_{0}(n))_{11} \right) \! + \! (\mathscr{R}^{e,0}_{0}
(n))_{12}(\overset{e}{m}_{1}^{\raise-1.0ex\hbox{$\scriptstyle 0$}})_{21} \! +
\! (\mathscr{R}^{e,0}_{1}(n))_{11}
(\overset{e}{m}_{0}^{\raise-1.0ex\hbox{$\scriptstyle 0$}})_{11} \\
+& \, (\mathscr{R}^{e,0}_{1}(n))_{12}
(\overset{e}{m}_{0}^{\raise-1.0ex\hbox{$\scriptstyle 0$}})_{21}, \\
(\widehat{Q}^{e}_{1})_{12} :=& \,
(\overset{e}{m}_{1}^{\raise-1.0ex\hbox{$\scriptstyle 0$}})_{12} \! \left(1 \!
+ \! (\mathscr{R}^{e,0}_{0}(n))_{11} \right) \! + \! (\mathscr{R}^{e,0}_{0}
(n))_{12}(\overset{e}{m}_{1}^{\raise-1.0ex\hbox{$\scriptstyle 0$}})_{22} \! +
\! (\mathscr{R}^{e,0}_{1}(n))_{11}
(\overset{e}{m}_{0}^{\raise-1.0ex\hbox{$\scriptstyle 0$}})_{12} \\
+& \, (\mathscr{R}^{e,0}_{1}(n))_{12}
(\overset{e}{m}_{0}^{\raise-1.0ex\hbox{$\scriptstyle 0$}})_{22}, \\
(\widehat{Q}^{e}_{1})_{21} :=& \,
(\overset{e}{m}_{1}^{\raise-1.0ex\hbox{$\scriptstyle 0$}})_{21} \! \left(1 \!
+ \! (\mathscr{R}^{e,0}_{0}(n))_{22} \right) \! + \! (\mathscr{R}^{e,0}_{0}
(n))_{21}(\overset{e}{m}_{1}^{\raise-1.0ex\hbox{$\scriptstyle 0$}})_{11} \! +
\! (\mathscr{R}^{e,0}_{1}(n))_{21}
(\overset{e}{m}_{0}^{\raise-1.0ex\hbox{$\scriptstyle 0$}})_{11} \\
+& \, (\mathscr{R}^{e,0}_{1}(n))_{22}
(\overset{e}{m}_{0}^{\raise-1.0ex\hbox{$\scriptstyle 0$}})_{21}, \\
(\widehat{Q}^{e}_{1})_{22} :=& \,
(\overset{e}{m}_{1}^{\raise-1.0ex\hbox{$\scriptstyle 0$}})_{22} \! \left(1 \!
+ \! (\mathscr{R}^{e,0}_{0}(n))_{22} \right) \! + \! (\mathscr{R}^{e,0}_{0}
(n))_{21}(\overset{e}{m}_{1}^{\raise-1.0ex\hbox{$\scriptstyle 0$}})_{12} \! +
\! (\mathscr{R}^{e,0}_{1}(n))_{21}
(\overset{e}{m}_{0}^{\raise-1.0ex\hbox{$\scriptstyle 0$}})_{12} \\
+& \, (\mathscr{R}^{e,0}_{1}(n))_{22}
(\overset{e}{m}_{0}^{\raise-1.0ex\hbox{$\scriptstyle 0$}})_{22}, \\
(\widehat{Q}^{e}_{2})_{11} :=& \,
(\overset{e}{m}_{2}^{\raise-1.0ex\hbox{$\scriptstyle 0$}})_{11} \! \left(1 \!
+ \! (\mathscr{R}^{e,0}_{0}(n))_{11} \right) \! + \! (\mathscr{R}^{e,0}_{0}
(n))_{12}(\overset{e}{m}_{2}^{\raise-1.0ex\hbox{$\scriptstyle 0$}})_{21} \! +
\! (\mathscr{R}^{e,0}_{1}(n))_{11}
(\overset{e}{m}_{1}^{\raise-1.0ex\hbox{$\scriptstyle 0$}})_{11} \\
+& \, (\mathscr{R}^{e,0}_{1}(n))_{12}
(\overset{e}{m}_{1}^{\raise-1.0ex\hbox{$\scriptstyle 0$}})_{21} \! + \!
(\mathscr{R}^{e,0}_{2}(n))_{11}
(\overset{e}{m}_{0}^{\raise-1.0ex\hbox{$\scriptstyle 0$}})_{11} \! + \!
(\mathscr{R}^{e,0}_{2}(n))_{12}
(\overset{e}{m}_{0}^{\raise-1.0ex\hbox{$\scriptstyle 0$}})_{21}, \\
(\widehat{Q}^{e}_{2})_{12} :=& \,
(\overset{e}{m}_{2}^{\raise-1.0ex\hbox{$\scriptstyle 0$}})_{12} \! \left(1 \!
+ \! (\mathscr{R}^{e,0}_{0}(n))_{11} \right) \! + \! (\mathscr{R}^{e,0}_{0}
(n))_{12}(\overset{e}{m}_{2}^{\raise-1.0ex\hbox{$\scriptstyle 0$}})_{22} \! +
\! (\mathscr{R}^{e,0}_{1}(n))_{11}
(\overset{e}{m}_{1}^{\raise-1.0ex\hbox{$\scriptstyle 0$}})_{12} \\
+& \, (\mathscr{R}^{e,0}_{1}(n))_{12}
(\overset{e}{m}_{1}^{\raise-1.0ex\hbox{$\scriptstyle 0$}})_{22} \! + \!
(\mathscr{R}^{e,0}_{2}(n))_{11}
(\overset{e}{m}_{0}^{\raise-1.0ex\hbox{$\scriptstyle 0$}})_{12} \! + \!
(\mathscr{R}^{e,0}_{2}(n))_{12}
(\overset{e}{m}_{0}^{\raise-1.0ex\hbox{$\scriptstyle 0$}})_{22}, \\
(\widehat{Q}^{e}_{2})_{21} :=& \,
(\overset{e}{m}_{2}^{\raise-1.0ex\hbox{$\scriptstyle 0$}})_{21} \! \left(1 \!
+ \! (\mathscr{R}^{e,0}_{0}(n))_{22} \right) \! + \! (\mathscr{R}^{e,0}_{0}
(n))_{21}(\overset{e}{m}_{2}^{\raise-1.0ex\hbox{$\scriptstyle 0$}})_{11} \! +
\! (\mathscr{R}^{e,0}_{1}(n))_{21}
(\overset{e}{m}_{1}^{\raise-1.0ex\hbox{$\scriptstyle 0$}})_{11} \\
+& \, (\mathscr{R}^{e,0}_{1}(n))_{22}
(\overset{e}{m}_{1}^{\raise-1.0ex\hbox{$\scriptstyle 0$}})_{21} \! + \!
(\mathscr{R}^{e,0}_{2}(n))_{21}
(\overset{e}{m}_{0}^{\raise-1.0ex\hbox{$\scriptstyle 0$}})_{11} \! + \!
(\mathscr{R}^{e,0}_{2}(n))_{22}
(\overset{e}{m}_{0}^{\raise-1.0ex\hbox{$\scriptstyle 0$}})_{21}, \\
(\widehat{Q}^{e}_{2})_{22} :=& \,
(\overset{e}{m}_{2}^{\raise-1.0ex\hbox{$\scriptstyle 0$}})_{22} \! \left(1 \!
+ \! (\mathscr{R}^{e,0}_{0}(n))_{22} \right) \! + \! (\mathscr{R}^{e,0}_{0}
(n))_{21}(\overset{e}{m}_{2}^{\raise-1.0ex\hbox{$\scriptstyle 0$}})_{12} \! +
\! (\mathscr{R}^{e,0}_{1}(n))_{21}
(\overset{e}{m}_{1}^{\raise-1.0ex\hbox{$\scriptstyle 0$}})_{12} \\
+& \, (\mathscr{R}^{e,0}_{1}(n))_{22}
(\overset{e}{m}_{1}^{\raise-1.0ex\hbox{$\scriptstyle 0$}})_{22} \! + \!
(\mathscr{R}^{e,0}_{2}(n))_{21}
(\overset{e}{m}_{0}^{\raise-1.0ex\hbox{$\scriptstyle 0$}})_{12} \! + \!
(\mathscr{R}^{e,0}_{2}(n))_{22}
(\overset{e}{m}_{0}^{\raise-1.0ex\hbox{$\scriptstyle 0$}})_{22},
\end{align*}
and
\begin{align*}
(\overset{e}{m}_{0}^{\raise-1.0ex\hbox{$\scriptstyle 0$}})_{11} &= \dfrac{
\boldsymbol{\theta}^{e}(\boldsymbol{u}^{e}_{+}(\infty) \! + \! \boldsymbol{d}_{
e})}{\boldsymbol{\theta}^{e}(\boldsymbol{u}^{e}_{+}(\infty) \! - \! \frac{n}{2
\pi} \boldsymbol{\Omega}^{e} \! + \! \boldsymbol{d}_{e})} \! \left(\dfrac{
\gamma^{e}_{0} \! + \! (\gamma^{e}_{0})^{-1}}{2} \right) \! \dfrac{\theta^{
e}_{0}(1,1,\boldsymbol{\Omega}^{e})}{\theta^{e}_{0}(1,1,\pmb{0})}, \\
(\overset{e}{m}_{0}^{\raise-1.0ex\hbox{$\scriptstyle 0$}})_{12} &= -\dfrac{
\boldsymbol{\theta}^{e}(\boldsymbol{u}^{e}_{+}(\infty) \! + \! \boldsymbol{d}_{
e})}{\boldsymbol{\theta}^{e}(\boldsymbol{u}^{e}_{+}(\infty) \! - \! \frac{n}{2
\pi} \boldsymbol{\Omega}^{e} \! + \! \boldsymbol{d}_{e})} \! \left(\dfrac{
\gamma^{e}_{0} \! - \! (\gamma^{e}_{0})^{-1}}{2 \mi} \right) \! \dfrac{
\theta^{e}_{0}(-1,1,\boldsymbol{\Omega}^{e})}{\theta^{e}_{0}(-1,1,\pmb{0})}, \\
(\overset{e}{m}_{0}^{\raise-1.0ex\hbox{$\scriptstyle 0$}})_{21} &= \dfrac{
\boldsymbol{\theta}^{e}(\boldsymbol{u}^{e}_{+}(\infty) \! + \! \boldsymbol{d}_{
e})}{\boldsymbol{\theta}^{e}(-\boldsymbol{u}^{e}_{+}(\infty) \! - \! \frac{n}{
2 \pi} \boldsymbol{\Omega}^{e} \! - \! \boldsymbol{d}_{e})} \! \left(\dfrac{
\gamma^{e}_{0} \! - \! (\gamma^{e}_{0})^{-1}}{2 \mi} \right) \! \dfrac{\theta^{
e}_{0}(1,-1,\boldsymbol{\Omega}^{e})}{\theta^{e}_{0}(1,-1,\pmb{0})}, \\
(\overset{e}{m}_{0}^{\raise-1.0ex\hbox{$\scriptstyle 0$}})_{22} &= \dfrac{
\boldsymbol{\theta}^{e}(\boldsymbol{u}^{e}_{+}(\infty) \! + \! \boldsymbol{d}_{
e})}{\boldsymbol{\theta}^{e}(-\boldsymbol{u}^{e}_{+}(\infty) \! - \! \frac{n}{
2 \pi} \boldsymbol{\Omega}^{e} \! - \! \boldsymbol{d}_{e})} \! \left(\dfrac{
\gamma^{e}_{0} \! + \! (\gamma^{e}_{0})^{-1}}{2} \right) \! \dfrac{\theta^{e}_{
0}(-1,-1,\boldsymbol{\Omega}^{e})}{\theta^{e}_{0}(-1,-1,\pmb{0})}, \\
(\overset{e}{m}_{1}^{\raise-1.0ex\hbox{$\scriptstyle 0$}})_{11} &= \dfrac{
\boldsymbol{\theta}^{e}(\boldsymbol{u}^{e}_{+}(\infty) \! + \! \boldsymbol{d}_{
e})}{\boldsymbol{\theta}^{e}(\boldsymbol{u}^{e}_{+}(\infty) \! - \! \frac{n}{2
\pi} \boldsymbol{\Omega}^{e} \! + \! \boldsymbol{d}_{e})} \! \left(\! \left(\!
\dfrac{\widetilde{\alpha}_{0}^{e}(1,1,\boldsymbol{\Omega}^{e}) \theta^{e}_{0}
(1,1,\pmb{0}) \! - \! \widetilde{\alpha}_{0}^{e}(1,1,\pmb{0}) \theta^{e}_{0}
(1,1,\boldsymbol{\Omega}^{e})}{(\theta^{e}_{0}(1,1,\pmb{0}))^{2}} \right)
\right. \\
&\left. \times \left(\dfrac{\gamma^{e}_{0} \! + \! (\gamma^{e}_{0})^{-1}}{2}
\right) \! + \! \left(\dfrac{\gamma^{e}_{0} \! - \! (\gamma^{e}_{0})^{-1}}{8}
\right) \! \left(\sum_{k=1}^{N+1} \! \left(\dfrac{1}{a_{k}^{e}} \! - \! \dfrac{
1}{b_{k-1}^{e}} \right) \right) \! \dfrac{\theta^{e}_{0}(1,1,\boldsymbol{
\Omega}^{e})}{\theta^{e}_{0}(1,1,\pmb{0})} \right), \\
(\overset{e}{m}_{1}^{\raise-1.0ex\hbox{$\scriptstyle 0$}})_{12} &= -\dfrac{
\boldsymbol{\theta}^{e}(\boldsymbol{u}^{e}_{+}(\infty) \! + \! \boldsymbol{d}_{
e})}{\boldsymbol{\theta}^{e}(\boldsymbol{u}^{e}_{+}(\infty) \! - \! \frac{n}{2
\pi} \boldsymbol{\Omega}^{e} \! + \! \boldsymbol{d}_{e})} \! \left(\! \left(\!
\dfrac{\widetilde{\alpha}_{0}^{e}(-1,1,\boldsymbol{\Omega}^{e}) \theta^{e}_{0}
(-1,1,\pmb{0}) \! - \! \widetilde{\alpha}_{0}^{e}(-1,1,\pmb{0}) \theta^{e}_{0}
(-1,1,\boldsymbol{\Omega}^{e})}{(\theta^{e}_{0}(-1,1,\pmb{0}))^{2}} \right)
\right. \\
&\left. \times \left(\dfrac{\gamma^{e}_{0} \! - \! (\gamma^{e}_{0})^{-1}}{2
\mi} \right) \! + \! \left(\dfrac{\gamma^{e}_{0} \! + \! (\gamma^{e}_{0})^{-1}
}{8 \mi} \right) \! \left(\sum_{k=1}^{N+1} \! \left(\dfrac{1}{a_{k}^{e}} \! -
\! \dfrac{1}{b_{k-1}^{e}} \right) \right) \! \dfrac{\theta^{e}_{0}(-1,1,
\boldsymbol{\Omega}^{e})}{\theta^{e}_{0}(-1,1,\pmb{0})} \right), \\
(\overset{e}{m}_{1}^{\raise-1.0ex\hbox{$\scriptstyle 0$}})_{21} &= \dfrac{
\boldsymbol{\theta}^{e}(\boldsymbol{u}^{e}_{+}(\infty) \! + \! \boldsymbol{d}_{
e})}{\boldsymbol{\theta}^{e}(-\boldsymbol{u}^{e}_{+}(\infty) \! - \! \frac{n}{
2 \pi} \boldsymbol{\Omega}^{e} \! - \! \boldsymbol{d}_{e})} \! \left(\! \left(
\! \dfrac{\widetilde{\alpha}_{0}^{e}(1,-1,\boldsymbol{\Omega}^{e}) \theta^{e}_{
0}(1,-1,\pmb{0}) \! - \! \widetilde{\alpha}_{0}^{e}(1,-1,\pmb{0}) \theta^{e}_{
0}(1,-1,\boldsymbol{\Omega}^{e})}{(\theta^{e}_{0}(1,-1,\pmb{0}))^{2}} \right)
\right. \\
&\left. \times \left(\dfrac{\gamma^{e}_{0} \! - \! (\gamma^{e}_{0})^{-1}}{2
\mi} \right) \! + \! \left(\dfrac{\gamma^{e}_{0} \! + \! (\gamma^{e}_{0})^{-1}
}{8 \mi} \right) \! \left(\sum_{k=1}^{N+1} \! \left(\dfrac{1}{a_{k}^{e}} \!
- \! \dfrac{1}{b_{k-1}^{e}} \right) \right) \! \dfrac{\theta^{e}_{0}(1,-1,
\boldsymbol{\Omega}^{e})}{\theta^{e}_{0}(1,-1,\pmb{0})} \right), \\
(\overset{e}{m}_{1}^{\raise-1.0ex\hbox{$\scriptstyle 0$}})_{22} &= \dfrac{
\boldsymbol{\theta}^{e}(\boldsymbol{u}^{e}_{+}(\infty) \! + \! \boldsymbol{d}_{
e})}{\boldsymbol{\theta}^{e}(-\boldsymbol{u}^{e}_{+}(\infty) \! - \! \frac{n}{
2 \pi} \boldsymbol{\Omega}^{e} \! - \! \boldsymbol{d}_{e})} \! \left(\! \left(
\! \dfrac{\widetilde{\alpha}_{0}^{e}(-1,-1,\boldsymbol{\Omega}^{e}) \theta^{
e}_{0}(-1,-1,\pmb{0}) \! - \! \widetilde{\alpha}_{0}^{e}(-1,-1,\pmb{0})
\theta^{e}_{0}(-1,-1,\boldsymbol{\Omega}^{e})}{(\theta^{e}_{0}(-1,-1,\pmb{0})
)^{2}} \right) \right. \\
&\left. \times \left(\dfrac{\gamma^{e}_{0} \! + \! (\gamma^{e}_{0})^{-1}}{2}
\right) \! + \! \left(\dfrac{\gamma^{e}_{0} \! - \! (\gamma^{e}_{0})^{-1}}{8}
\right) \! \left(\sum_{k=1}^{N+1} \! \left(\dfrac{1}{a_{k}^{e}} \! - \! \dfrac{
1}{b_{k-1}^{e}} \right) \right) \! \dfrac{\theta^{e}_{0}(-1,-1,\boldsymbol{
\Omega}^{e})}{\theta^{e}_{0}(-1,-1,\pmb{0})} \right), \\
(\overset{e}{m}_{2}^{\raise-1.0ex\hbox{$\scriptstyle 0$}})_{11} &= \dfrac{
\boldsymbol{\theta}^{e}(\boldsymbol{u}^{e}_{+}(\infty) \! + \! \boldsymbol{d}_{
e})}{\boldsymbol{\theta}^{e}(\boldsymbol{u}^{e}_{+}(\infty) \! - \! \frac{n}{2
\pi} \boldsymbol{\Omega}^{e} \! + \! \boldsymbol{d}_{e})} \! \left(\! \left(
\theta^{e}_{0}(1,1,\boldsymbol{\Omega}^{e}) \! \left((\widetilde{\alpha}_{0}^{
e}(1,1,\pmb{0}))^{2} \! - \! \beta^{e}_{0}(1,1,\pmb{0}) \theta^{e}_{0}(1,1,
\pmb{0}) \right) \right. \right. \\
&\left. \left. - \, \widetilde{\alpha}_{0}^{e}(1,1,\boldsymbol{\Omega}^{e})
\widetilde{\alpha}_{0}^{e}(1,1,\pmb{0}) \theta^{e}_{0}(1,1,\pmb{0}) \! + \!
\beta^{e}_{0}(1,1,\boldsymbol{\Omega}^{e})(\theta^{e}_{0}(1,1,\pmb{0}))^{2}
\right) \! \tfrac{1}{(\theta^{e}_{0}(1,1,\pmb{0}))^{3}} \right. \\
&\left. \times \left(\dfrac{\gamma^{e}_{0} \! + \! (\gamma^{e}_{0})^{-1}}{2}
\right) \! + \! \left(\dfrac{\widetilde{\alpha}_{0}^{e}(1,1,\boldsymbol{
\Omega}^{e}) \theta^{e}_{0}(1,1,\pmb{0}) \! - \! \widetilde{\alpha}_{0}^{e}(1,
1,\pmb{0}) \theta^{e}_{0}(1,1,\boldsymbol{\Omega}^{e})}{(\theta^{e}_{0}(1,1,
\pmb{0}))^{2}} \right) \right. \\
&\left. \times \left(\sum_{k=1}^{N+1} \! \left(\dfrac{1}{a_{k}^{e}} \! - \!
\dfrac{1}{b_{k-1}^{e}} \right) \right) \! \left(\dfrac{\gamma^{e}_{0} \! - \!
(\gamma^{e}_{0})^{-1}}{8} \right) \! + \! \left(\! \left(\dfrac{\gamma^{e}_{0}
\! - \! (\gamma^{e}_{0})^{-1}}{16} \right) \! \sum_{k=1}^{N+1} \! \left(
\dfrac{1}{(a_{k}^{e})^{2}} \! - \! \dfrac{1}{(b_{k-1}^{e})^{2}} \right)
\right. \right. \\
&\left. \left. +\left(\dfrac{\gamma^{e}_{0} \! + \! (\gamma^{e}_{0})^{-1}}{64}
\right) \! \left(\sum_{k=1}^{N+1} \! \left(\dfrac{1}{a_{k}^{e}} \! - \!
\dfrac{1}{b_{k-1}^{e}} \right) \right)^{2} \right) \! \dfrac{\theta^{e}_{0}
(1,1,\boldsymbol{\Omega}^{e})}{\theta^{e}_{0}(1,1,\pmb{0})} \right), \\
(\overset{e}{m}_{2}^{\raise-1.0ex\hbox{$\scriptstyle 0$}})_{12} &= -\dfrac{
\boldsymbol{\theta}^{e}(\boldsymbol{u}^{e}_{+}(\infty) \! + \! \boldsymbol{d}_{
e})}{\boldsymbol{\theta}^{e}(\boldsymbol{u}^{e}_{+}(\infty) \! - \! \frac{n}{2
\pi} \boldsymbol{\Omega}^{e} \! + \! \boldsymbol{d}_{e})} \! \left(\! \left(
\theta^{e}_{0}(-1,1,\boldsymbol{\Omega}^{e}) \! \left((\widetilde{\alpha}_{0}^{
e}(-1,1,\pmb{0}))^{2} \! - \! \beta^{e}_{0}(-1,1,\pmb{0}) \theta^{e}_{0}(-1,1,
\pmb{0}) \right) \right. \right. \\
&\left. \left. - \, \widetilde{\alpha}_{0}^{e}(-1,1,\boldsymbol{\Omega}^{e})
\widetilde{\alpha}_{0}^{e}(-1,1,\pmb{0}) \theta^{e}_{0}(-1,1,\pmb{0}) \! + \!
\beta^{e}_{0}(-1,1,\boldsymbol{\Omega}^{e})(\theta^{e}_{0}(-1,1,\pmb{0}))^{2}
\right) \! \tfrac{1}{(\theta^{e}_{0}(-1,1,\pmb{0}))^{3}} \right. \\
&\left. \times \left(\dfrac{\gamma^{e}_{0} \! - \! (\gamma^{e}_{0})^{-1}}{2
\mi} \right) \! + \! \left(\dfrac{\widetilde{\alpha}_{0}^{e}(-1,1,\boldsymbol{
\Omega}^{e}) \theta^{e}_{0}(-1,1,\pmb{0}) \! - \! \widetilde{\alpha}_{0}^{e}
(-1,1,\pmb{0}) \theta^{e}_{0}(-1,1,\boldsymbol{\Omega}^{e})}{(\theta^{e}_{0}
(-1,1,\pmb{0}))^{2}} \right) \right. \\
&\left. \times \left(\sum_{k=1}^{N+1} \! \left(\dfrac{1}{a_{k}^{e}} \! - \!
\dfrac{1}{b_{k-1}^{e}} \right) \right) \! \left(\dfrac{\gamma^{e}_{0} \! + \!
(\gamma^{e}_{0})^{-1}}{8 \mi} \right) \! + \! \left(\! \left(\dfrac{\gamma^{
e}_{0} \! + \! (\gamma^{e}_{0})^{-1}}{16 \mi} \right) \! \sum_{k=1}^{N+1} \!
\left(\dfrac{1}{(a_{k}^{e})^{2}} \! - \! \dfrac{1}{(b_{k-1}^{e})^{2}} \right)
\right. \right. \\
&\left. \left. +\left(\dfrac{\gamma^{e}_{0} \! - \! (\gamma^{e}_{0})^{-1}}{64
\mi} \right) \! \left(\sum_{k=1}^{N+1} \! \left(\dfrac{1}{a_{k}^{e}} \! - \!
\dfrac{1}{b_{k-1}^{e}} \right) \right)^{2} \right) \! \dfrac{\theta^{e}_{0}
(-1,1,\boldsymbol{\Omega}^{e})}{\theta^{e}_{0}(-1,1,\pmb{0})} \right), \\
(\overset{e}{m}_{2}^{\raise-1.0ex\hbox{$\scriptstyle 0$}})_{21} &= \dfrac{
\boldsymbol{\theta}^{e}(\boldsymbol{u}^{e}_{+}(\infty) \! + \! \boldsymbol{d}_{
e})}{\boldsymbol{\theta}^{e}(-\boldsymbol{u}^{e}_{+}(\infty) \! - \! \frac{n}{
2 \pi} \boldsymbol{\Omega}^{e} \! - \! \boldsymbol{d}_{e})} \! \left(\! \left(
\theta^{e}_{0}(1,-1,\boldsymbol{\Omega}^{e}) \! \left((\widetilde{\alpha}_{0}^{
e}(1,-1,\pmb{0}))^{2} \! - \! \beta^{e}_{0}(1,-1,\pmb{0}) \theta^{e}_{0}(1,-1,
\pmb{0}) \right) \right. \right. \\
&\left. \left. - \, \widetilde{\alpha}_{0}^{e}(1,-1,\boldsymbol{\Omega}^{e})
\widetilde{\alpha}_{0}^{e}(1,-1,\pmb{0}) \theta^{e}_{0}(1,-1,\pmb{0}) \! + \!
\beta^{e}_{0}(1,-1,\boldsymbol{\Omega}^{e})(\theta^{e}_{0}(1,-1,\pmb{0}))^{2}
\right) \! \tfrac{1}{(\theta^{e}_{0}(1,-1,\pmb{0}))^{3}} \right. \nonumber \\
&\left. \times \left(\dfrac{\gamma^{e}_{0} \! - \! (\gamma^{e}_{0})^{-1}}{2
\mi} \right) \! + \! \left(\dfrac{\widetilde{\alpha}_{0}^{e}(1,-1,\boldsymbol{
\Omega}^{e}) \theta^{e}_{0}(1,-1,\pmb{0}) \! - \! \widetilde{\alpha}_{0}^{e}
(1,-1,\pmb{0}) \theta^{e}_{0}(1,-1,\boldsymbol{\Omega}^{e})}{(\theta^{e}_{0}
(1,-1,\pmb{0}))^{2}} \right) \right. \\
&\left. \times \left(\sum_{k=1}^{N+1} \! \left(\dfrac{1}{a_{k}^{e}} \! - \!
\dfrac{1}{b_{k-1}^{e}} \right) \right) \! \left(\dfrac{\gamma^{e}_{0} \! + \!
(\gamma^{e}_{0})^{-1}}{8 \mi} \right) \! + \! \left(\! \left(\dfrac{\gamma^{
e}_{0} \! + \! (\gamma^{e}_{0})^{-1}}{16 \mi} \right) \! \sum_{k=1}^{N+1} \!
\left(\dfrac{1}{(a_{k}^{e})^{2}} \! - \! \dfrac{1}{(b_{k-1}^{e})^{2}} \right)
\right. \right. \\
&\left. \left. +\left(\dfrac{\gamma^{e}_{0} \! - \! (\gamma^{e}_{0})^{-1}}{64
\mi} \right) \! \left(\sum_{k=1}^{N+1} \! \left(\dfrac{1}{a_{k}^{e}} \! - \!
\dfrac{1}{b_{k-1}^{e}} \right) \right)^{2} \right) \! \dfrac{\theta^{e}_{0}
(1,-1,\boldsymbol{\Omega}^{e})}{\theta^{e}_{0}(1,-1,\pmb{0})} \right), \\
(\overset{e}{m}_{2}^{\raise-1.0ex\hbox{$\scriptstyle 0$}})_{22} &= \dfrac{
\boldsymbol{\theta}^{e}(\boldsymbol{u}^{e}_{+}(\infty) \! + \! \boldsymbol{d}_{
e})}{\boldsymbol{\theta}^{e}(-\boldsymbol{u}^{e}_{+}(\infty) \! - \! \frac{n}{
2 \pi} \boldsymbol{\Omega}^{e} \! - \! \boldsymbol{d}_{e})} \! \left(\! \left(
\theta^{e}_{0}(-1,-1,\boldsymbol{\Omega}^{e}) \! \left((\widetilde{\alpha}_{
0}^{e}(-1,-1,\pmb{0}))^{2} \! - \! \beta^{e}_{0}(-1,-1,\pmb{0}) \theta^{e}_{0}
(-1,-1,\pmb{0}) \right) \right. \right. \\
&\left. \left. - \, \widetilde{\alpha}_{0}^{e}(-1,-1,\boldsymbol{\Omega}^{e})
\widetilde{\alpha}_{0}^{e}(-1,-1,\pmb{0}) \theta^{e}_{0}(-1,-1,\pmb{0}) \! +
\! \beta^{e}_{0}(-1,-1,\boldsymbol{\Omega}^{e})(\theta^{e}_{0}(-1,-1,\pmb{0})
)^{2} \right) \! \tfrac{1}{(\theta^{e}_{0}(-1,-1,\pmb{0}))^{3}} \right. \\
&\left. \times \left(\dfrac{\gamma^{e}_{0} \! + \! (\gamma^{e}_{0})^{-1}}{2}
\right) \! + \! \left(\dfrac{\widetilde{\alpha}_{0}^{e}(-1,-1,\boldsymbol{
\Omega}^{e}) \theta^{e}_{0}(-1,-1,\pmb{0}) \! - \! \widetilde{\alpha}_{0}^{e}
(-1,-1,\pmb{0}) \theta^{e}_{0}(-1,-1,\boldsymbol{\Omega}^{e})}{(\theta^{e}_{0}
(-1,-1,\pmb{0}))^{2}} \right) \right. \\
&\left. \times \left(\sum_{k=1}^{N+1} \! \left(\dfrac{1}{a_{k}^{e}} \! - \!
\dfrac{1}{b_{k-1}^{e}} \right) \right) \! \left(\dfrac{\gamma^{e}_{0} \! - \!
(\gamma^{e}_{0})^{-1}}{8} \right) \! + \! \left(\! \left(\dfrac{\gamma^{e}_{0}
\! - \! (\gamma^{e}_{0})^{-1}}{16} \right) \! \sum_{k=1}^{N+1} \! \left(
\dfrac{1}{(a_{k}^{e})^{2}} \! - \! \dfrac{1}{(b_{k-1}^{e})^{2}} \right)
\right. \right. \\
&\left. \left. +\left(\dfrac{\gamma^{e}_{0} \! + \! (\gamma^{e}_{0})^{-1}}{64}
\right) \! \left(\sum_{k=1}^{N+1} \! \left(\dfrac{1}{a_{k}^{e}} \! - \!
\dfrac{1}{b_{k-1}^{e}} \right) \right)^{2} \right) \! \dfrac{\theta^{e}_{0}
(-1,-1,\boldsymbol{\Omega}^{e})}{\theta^{e}_{0}(-1,-1,\pmb{0})} \right),
\end{align*}
where, for $\varepsilon_{1},\varepsilon_{2} \! = \! \pm 1$,
\begin{gather*}
\theta^{e}_{0}(\varepsilon_{1},\varepsilon_{2},\pmb{\mathscr{Z}}) \! := \!
\boldsymbol{\theta}^{e}(\varepsilon_{1} \boldsymbol{u}^{e}_{+}(0) \! - \!
\tfrac{n}{2 \pi} \pmb{\mathscr{Z}} \! + \! \varepsilon_{2} \boldsymbol{d}_{
e}), \\
\widetilde{\alpha}^{e}_{0}(\varepsilon_{1},\varepsilon_{2},\pmb{\mathscr{Z}})
\! := \! 2 \pi \mi \varepsilon_{1} \sum_{m \in \mathbb{Z}^{N}} (m,\widehat{
\boldsymbol{\alpha}}^{e}_{0}) \me^{2 \pi \mi (m,\varepsilon_{1} \boldsymbol{
u}^{e}_{+}(0)-\frac{n}{2 \pi} \pmb{\mathscr{Z}}+\varepsilon_{2} \boldsymbol{
d}_{e})+ \pi \mi (m,\tau^{e}m)},
\end{gather*}
with $\boldsymbol{u}^{e}_{+}(0) \! = \! \int_{a_{N+1}^{e}}^{0^{+}} \boldsymbol{
\omega}^{e}$, $\widehat{\boldsymbol{\alpha}}^{e}_{0} \! = \! (\widehat{
\alpha}_{0,1}^{e},\widehat{\alpha}^{e}_{0,2},\dotsc,\widehat{\alpha}^{e}_{0,
N})$, and $\widehat{\alpha}^{e}_{0,j} \! := \! (-1)^{\mathcal{N}_{+}^{e}}
(\prod_{i=1}^{N+1} \vert b_{i-1}^{e}a_{i}^{e} \vert)^{-1/2}c_{jN}^{e}$, $j \! 
= \! 1,\dotsc,N$, where $\mathcal{N}_{+}^{e} \! \in \! \lbrace 0,\dotsc,N \! 
+ \! 1 \rbrace$ is the number of bands to the right of $z \! = \! 0$,
\begin{gather*}
\beta^{e}_{0}(\varepsilon_{1},\varepsilon_{2},\pmb{\mathscr{Z}}) \! := \! 
2 \pi \sum_{m \in \mathbb{Z}^{N}} \! \left(\mi \varepsilon_{1}(m,\widehat{
\boldsymbol{\beta}}^{e}_{0}) \! - \!  \pi (m,\widehat{\boldsymbol{\alpha}}^{
e}_{0})^{2} \right) \! \me^{2 \pi \mi (m,\varepsilon_{1} \boldsymbol{u}^{e}_{
+}(0)-\frac{n}{2 \pi} \pmb{\mathscr{Z}}+\varepsilon_{2} \boldsymbol{d}_{e})+
\pi \mi (m,\tau^{e}m)},
\end{gather*}
where $\widehat{\boldsymbol{\beta}}^{e}_{0} \! = \! (\widehat{\beta}^{e}_{0,1},
\widehat{\beta}^{e}_{0,2},\dotsc,\widehat{\beta}^{e}_{0,N})$, with $\widehat{
\beta}^{e}_{0,j} \! := \! \tfrac{1}{2}(-1)^{\mathcal{N}_{+}^{e}}(\prod_{i=1}^{
N+1} \vert b_{i-1}^{e}a_{i}^{e} \vert)^{-1/2}(c^{e}_{jN-1} \! + \! \tfrac{1}{
2}c_{jN}^{e} \sum_{k=1}^{N+1}((a_{k}^{e})^{-1} \! + \! (b_{k-1}^{e})^{-1}))$, 
$j \! = \! 1,\dotsc,N$, where $c_{jN}^{e},c_{jN-1}^{e}$, $j \! = \! 1,\dotsc,
N$, are obtained {}from Equations~{\rm (E1)} and~{\rm (E2)}, and $\gamma^{e}_{
0}$ is defined in Theorem~{\rm 3.1}, Equations~{\rm (3.53)}, and, for $k \! = 
\! 1,2,3$,
\begin{align*}
\mathscr{R}^{e,0}_{k-1}(n) \underset{n \to \infty}{=}& \, \dfrac{1}{n}
\mathlarger{\sum_{j=1}^{N+1}} \! \left(\dfrac{\left(\mathscr{A}^{e}(b_{j-1}^{
e}) \! \left(\widehat{\alpha}_{1}^{e}(b_{j-1}^{e}) \! + \! k(b_{j-1}^{e})^{-
1} \widehat{\alpha}_{0}^{e}(b_{j-1}^{e}) \right) \! - \! \mathscr{B}^{e}(b_{
j-1}^{e}) \widehat{\alpha}_{0}^{e}(b_{j-1}^{e}) \right)}{(b_{j-1}^{e})^{k}
(\widehat{\alpha}_{0}^{e}(b_{j-1}^{e}))^{2}} \right. \nonumber \\
+&\left. \, \dfrac{\left(\mathscr{A}^{e}(a_{j}^{e}) \! \left(\widehat{
\alpha}_{1}^{e}(a_{j}^{e}) \! + \! k(a_{j}^{e})^{-1} \widehat{\alpha}_{0}^{e}
(a_{j}^{e}) \right) \! - \! \mathscr{B}^{e}(a_{j}^{e}) \widehat{\alpha}_{0}^{
e}(a_{j}^{e}) \right)}{(a_{j}^{e})^{k}(\widehat{\alpha}_{0}^{e}(a_{j}^{e}))^{
2}} \right) \! + \! \mathcal{O} \! \left(\dfrac{c(n)}{n^{2}} \right),
\end{align*}
with all parameters defined in Theorem~{\rm A.1}.
\end{dy}
\addtocounter{z0}{1}
\begin{dy}[{\rm \cite{a21}}]
Let all the conditions stated in Theorem~{\rm A.1} be valid, and let all 
parameters be as defined therein. Let $\overset{e}{\operatorname{Y}} \colon 
\mathbb{C} \setminus \mathbb{R} \! \to \! \operatorname{SL}_{2}(\mathbb{C})$ 
be the unique solution of \pmb{{\rm RHP1}} with integral 
representation~{\rm (1.29)}. Let $\boldsymbol{\pi}_{2n}(z)$ be the even 
degree monic orthogonal $L$-polynomial defined in Equations~{\rm (1.4)} with 
$n \! \to \! \infty$ asymptotics (in the entire complex plane) given by 
Theorem~{\rm 2.3.1} of {\rm \cite{a21}}. Then,
\begin{equation}
\left(\xi^{(2n)}_{n} \right)^{2} \! = \! \dfrac{1}{\norm{\boldsymbol{\pi}_{2n}
(\pmb{\cdot})}_{\mathscr{L}}^{2}} \! = \! \dfrac{H^{(-2n)}_{2n}}{H^{(-2n)}_{2
n+1} } \underset{n \to \infty}{=} \dfrac{\me^{-n \ell_{e}}}{\pi} \Xi^{\flat} 
\! \left(1 \! + \! \dfrac{2}{n} \Xi^{\flat}(\mathfrak{Q}^{\flat})_{12} \! + 
\! \mathcal{O} \! \left(\dfrac{c^{\flat}(n)}{n^{2}} \right) \right),
\end{equation}
where
\begin{gather*}
\Xi^{\flat} \! := \! 2 \! \left(\sum_{k=1}^{N+1} \! \left(a_{k}^{e} \! - \!
b_{k-1}^{e} \right) \right)^{-1} \dfrac{\boldsymbol{\theta}^{e}(\boldsymbol{
u}^{e}_{+}(\infty) \! - \! \frac{n}{2 \pi} \boldsymbol{\Omega}^{e} \! + \!
\boldsymbol{d}_{e}) \boldsymbol{\theta}^{e}(-\boldsymbol{u}^{e}_{+}(\infty)
\! + \! \boldsymbol{d}_{e})}{\boldsymbol{\theta}^{e}(-\boldsymbol{u}^{e}_{+}
(\infty) \! - \! \frac{n}{2 \pi} \boldsymbol{\Omega}^{e} \! + \! \boldsymbol{
d}_{e}) \boldsymbol{\theta}^{e}(\boldsymbol{u}^{e}_{+}(\infty) \! + \!
\boldsymbol{d}_{e})} \quad (> \! 0),
\end{gather*}
$\mathfrak{Q}^{\flat}$ is defined in Theorem~{\rm 3.2}, 
Equations~{\rm (3.107)}, $(\mathfrak{Q}^{\flat})_{12}$ denotes the
$(1 \, 2)$-element of $\mathfrak{Q}^{\flat}$, $c^{\flat}(n) \! =_{n \to 
\infty} \! \mathcal{O}(1)$, and all the relevant parameters are defined in 
Theorem~{\rm A.1}. Asymptotics for $\xi^{(2n)}_{n}$ are obtained by taking 
the positive square root of both sides of Equation~{\rm (A.1)}.
\end{dy}
\clearpage
\section*{Appendix B. Asymptotics of $\overset{o}{\operatorname{Y}}(z)$}
\setcounter{section}{2}
\setcounter{z0}{1}
\setcounter{y0}{1}
\setcounter{equation}{0}
\renewcommand{\theequation}{\Alph{section}.\arabic{equation}}
In this appendix, large-$n$ asymptotics for the coefficients of the 
expansions~(3.3) and~(3.4) for the solution of \pmb{RHP2}, with integral 
representation~(1.30), are presented (for complete details and proofs, see 
\cite{a22}).
\begin{ey}
In \cite{a22}, the following formula was established for complex $z$ away 
from an open neighbourhood around $\overline{J_{o}} \! := \! \cup_{j=1}^{N+1}
[b_{j-1}^{o},a_{j}^{o}]$:
\begin{equation*}
\overset{o}{\mathrm{Y}}(z) \! = \!
\begin{cases}
\me^{\frac{n \ell_{o}}{2} \sigma_{3}} \mathscr{R}^{o}(z)
\overset{o}{m}^{\raise-1.0ex\hbox{$\scriptstyle \infty$}}(z) \mathbb{E}^{
\sigma_{3}} \me^{n(g^{o}(z)-\frac{\ell_{o}}{2}-\mathfrak{Q}^{+}_{\mathscr{A}})
\sigma_{3}}, &\text{$z \! \in \! \mathbb{C}_{+}$,} \\
\me^{\frac{n \ell_{o}}{2} \sigma_{3}} \mathscr{R}^{o}(z)
\overset{o}{m}^{\raise-1.0ex\hbox{$\scriptstyle \infty$}}(z) \mathbb{E}^{-
\sigma_{3}} \me^{n(g^{o}(z)-\frac{\ell_{o}}{2}-\mathfrak{Q}^{-}_{\mathscr{A}})
\sigma_{3}}, &\text{$z \! \in \! \mathbb{C}_{-}$,}
\end{cases}
\end{equation*}
where $\ell_{o}$, $g^{o}(z)$, $\mathfrak{Q}^{\pm}_{\mathscr{A}}$, and 
$\mathbb{E}$ have been defined heretofore, the matrix-valued function 
$\overset{o}{m}^{\raise-1.0ex\hbox{$\scriptstyle \infty$}}(z)$ is determined 
explicitly, and the matrix-valued function $\mathscr{R}^{o}(z)$ solves a 
``small norm'' RHP for which a complete asymptotic expansion as $n \! \to \! 
\infty$ for $z \! \in \! \mathbb{C}$ is computed. The above explicit formula 
for the solution of \pmb{RHP2} yields an explicit asymptotic description for 
the odd degree OLPs in question (for complete details and proofs, see 
\cite{a22}). \hfill $\blacksquare$
\end{ey}
\addtocounter{y0}{1}
\begin{ey}
In Theorems~B.1 and~B.2 below, only $n \! \to \! \infty$ asymptotics for 
$\mathscr{R}^{o,0}_{k-1}(n)$, $k \! = \! 2,3$, and $\mathscr{R}^{o,\infty}_{
j}(n)$, $j \! = \! 0,1,2$, are given, as they are all that are actually 
necessary in order to obtain the results of Theorems~3.1--3.3: finite-$n$ 
formulae for $\mathscr{R}^{o,0}_{k-1}(n)$, $k \! = \! 2,3$, and $\mathscr{R}^{
o,\infty}_{j}(n)$, $j \! = \! 0,1,2$, are given in \cite{a22}; furthermore, 
the `symbol' $c(n)$ appearing in the asymptotic (as $n \! \to \! \infty)$ 
expressions for $\mathscr{R}^{o,0}_{k-1}(n)$, $k \! = \! 2,3$, and $\mathscr{
R}^{o,\infty}_{j}(n)$, $j \! = \! 0,1,2$, is to be understood as a uniformly 
bounded $(\mathcal{O}(1))$, $\operatorname{M}_{2}(\mathbb{C})$-valued, 
$n$-dependent function, that is, $\operatorname{M}_{2}(\mathbb{C}) \! \ni 
\! c(n) \! =_{n \to \infty} \! \mathcal{O}(1)$. \hfill $\blacksquare$
\end{ey}
\begin{dy}[{\rm \cite{a22}}]
Let the external field $\widetilde{V} \colon \mathbb{R} \setminus \{0\} 
\! \to \! \mathbb{R}$ satisfy conditions~{\rm (1.20)--(1.22)}. Let the 
orthonormal $L$-polynomials (resp., monic orthogonal $L$-polynomials) 
$\lbrace \phi_{k}(z) \rbrace_{k \in \mathbb{Z}_{0}^{+}}$ (resp., $\lbrace 
\boldsymbol{\pi}_{k}(z) \rbrace_{k \in \mathbb{Z}_{0}^{+}})$ be as defined 
in Equations~{\rm (1.2)} and~{\rm (1.3)} (resp., Equations~{\rm (1.4)} 
and~{\rm (1.5))}. Let $\overset{o}{\operatorname{Y}} \colon \mathbb{C} 
\setminus \mathbb{R} \! \to \! \operatorname{SL}_{2}(\mathbb{C})$ be 
the unique solution of {\rm \pmb{RHP2}} with integral representation 
{\rm (1.30)}. Define the density of the `odd' equilibrium measure, $\md 
\mu_{V}^{o}(x)$, as in Equation~{\rm (2.19)}, and set $J_{o} \! := \! 
\operatorname{supp}(\mu_{V}^{o}) \! = \! \cup_{j=1}^{N+1}(b_{j-1}^{o},
a_{j}^{o})$, where $\lbrace b_{j-1}^{o},a_{j}^{o} \rbrace_{j=1}^{N+1}$ 
satisfy the (real) $n$-dependent and locally solvable system of $2(N \! + 
\! 1)$ moment conditions~{\rm (2.18)}.

Suppose, furthermore, that $\widetilde{V} \colon \mathbb{R} \setminus \{0\} \! 
\to \! \mathbb{R}$ is regular, namely: {\rm (i)} $h_{V}^{o}(z) \! \not\equiv 
\! 0$ on $\overline{J_{o}} \! := \! \cup_{j=1}^{N+1}[b_{j-1}^{o},a_{j}^{o}];$ 
{\rm (ii)}
\begin{equation*}
2 \! \left(2 \! + \! \dfrac{1}{n} \right) \! \int_{J_{o}} \ln (\vert x \! - \!
s \vert) \, \md \mu_{V}^{o}(s) \! - \! 2 \ln \vert x \vert \! - \! \widetilde{
V}(x) \! - \! \ell_{o} \! - \! 2 \! \left(2 \! + \! \dfrac{1}{n} \right) \!
Q_{o} \! = \! 0, \quad x \! \in \! \overline{J_{o}},
\end{equation*}
which defines the `odd' variational constant, $l_{o}$ $(\in \! \mathbb{R})$ 
(the same on each compact interval $[b_{j-1}^{o},a_{j}^{o}]$, $j \! = \! 1,
\dotsc,N \! + \! 1)$, where $Q_{o} \! := \! \int_{J_{o}} \ln (\lvert s \rvert) 
\psi_{V}^{o}(s) \, \md s$, and
\begin{equation*}
2 \! \left(2 \! + \! \dfrac{1}{n} \right) \! \int_{J_{o}} \ln (\vert x \! - \!
s \vert) \, \md \mu_{V}^{o}(s) \! - \! 2 \ln \vert x \vert \! - \! \widetilde{
V}(x) \! - \! \ell_{o} \! - \! 2 \! \left(2 \! + \! \dfrac{1}{n} \right) \!
Q_{o} \! < \! 0, \quad x \! \in \! \mathbb{R} \setminus \overline{J_{o}} \,;
\end{equation*}
{\rm (iii)}
\begin{equation*}
g^{o}_{+}(z) \! + \! g^{o}_{-}(z) \! - \! \widetilde{V}(z) \! - \! \ell_{o} \!
- \! \mathfrak{Q}^{+}_{\mathscr{A}} \! - \! \mathfrak{Q}^{-}_{\mathscr{A}} \!
< \! 0, \quad z \! \in \! \mathbb{R} \setminus \overline{J_{o}},
\end{equation*}
where $g^{o}(z)$ is defined by Equation~{\rm (2.21)}, $g^{o}_{\pm}(z) \! :=
\! \lim_{\varepsilon \downarrow 0}g^{o}(z \! \pm \! \mi \varepsilon)$, and 
$\mathfrak{Q}^{\pm}_{\mathscr{A}}$ are defined in Equations~{\rm (2.22);} 
and {\rm (iv)}
\begin{equation*}
\mi (g^{o}_{+}(z) \! - \! g^{o}_{-}(z) \! - \! \mathfrak{Q}^{+}_{\mathscr{A}}
\! + \! \mathfrak{Q}^{-}_{\mathscr{A}})^{\prime} \! > \! 0, \quad z \! \in \!
J_{o}.
\end{equation*}
Then,
\begin{equation*}
\overset{o}{\mathrm{Y}}(z)z^{n \sigma_{3}} \underset{z \to 0}{=} \mathrm{I}
\! + \! z \mathrm{Y}^{o,0}_{1} \! + \! z^{2} \mathrm{Y}^{o,0}_{2} \! + \!
\mathcal{O}(z^{3}),
\end{equation*}
where
\begin{align*}
(\mathrm{Y}^{o,0}_{1})_{11} &= -(2n \! + \! 1) \int_{J_{o}}s^{-1} \psi_{V}^{o}
(s) \, \md s \! + \!
(\overset{o}{m}_{1}^{\raise-1.0ex\hbox{$\scriptstyle 0$}})_{11} \mathbb{E} \!
+ \! (\mathscr{R}^{o,0}_{1}(n))_{11}, \\
(\mathrm{Y}^{o,0}_{1})_{12} &= \me^{n \ell_{o}} \! \left(
(\overset{o}{m}_{1}^{\raise-1.0ex\hbox{$\scriptstyle 0$}})_{12} \mathbb{E}^{-
1} \! + \!
(\mathscr{R}^{o,0}_{1}(n))_{12} \right), \\
(\mathrm{Y}^{o,0}_{1})_{21} &= \me^{-n \ell_{o}} \! \left(
(\overset{o}{m}_{1}^{\raise-1.0ex\hbox{$\scriptstyle 0$}})_{21} \mathbb{E} \!
+ \! (\mathscr{R}^{o,0}_{1}(n))_{21} \right), \\
(\mathrm{Y}^{o,0}_{1})_{22} &= (2n \! + \! 1) \int_{J_{o}}s^{-1} \psi_{V}^{o}
(s) \, \md s \! + \!
(\overset{o}{m}_{1}^{\raise-1.0ex\hbox{$\scriptstyle 0$}})_{22} \mathbb{E}^{-
1} \! + \!
(\mathscr{R}^{o,0}_{1}(n))_{22}, \\
(\mathrm{Y}^{o,0}_{2})_{11} =& \, \tfrac{1}{2}(2n \! + \! 1)^{2} \! \left(
\int_{J_{o}}s^{-1} \psi_{V}^{o}(s) \, \md s \right)^{2} \! - \! \tfrac{1}{2}
(2n \! + \! 1) \int_{J_{o}}s^{-2} \psi_{V}^{o}(s) \, \md s \! - \! (2n \! + \!
1) \left((\overset{o}{m}_{1}^{\raise-1.0ex\hbox{$\scriptstyle 0$}})_{11}
\mathbb{E} \! + \! (\mathscr{R}^{o,0}_{1}(n))_{11} \right) \\
\times& \, \int_{J_{o}}s^{-1} \psi_{V}^{o}(s) \, \md s \! + \! (\overset{o}{
m}_{2}^{\raise-1.0ex\hbox{$\scriptstyle 0$}})_{11} \mathbb{E} \! + \!
(\mathscr{R}^{o,0}_{2}(n))_{11} \! + \! \left((\mathscr{R}^{o,0}_{1}(n))_{11}
(\overset{o}{m}_{1}^{\raise-1.0ex\hbox{$\scriptstyle 0$}})_{11} \! + \!
(\mathscr{R}^{o,0}_{1}(n))_{12}
(\overset{o}{m}_{1}^{\raise-1.0ex\hbox{$\scriptstyle 0$}})_{21} \right) \!
\mathbb{E}, \\
(\mathrm{Y}^{o,0}_{2})_{12} &= \me^{n \ell_{o}} \! \left((2n \! + \! 1) \!
\left((\overset{o}{m}_{1}^{\raise-1.0ex\hbox{$\scriptstyle 0$}})_{12} \mathbb{
E}^{-1} \! + \! (\mathscr{R}^{o,0}_{1}(n))_{12} \right) \! \int_{J_{o}}s^{-1}
\psi_{V}^{o}(s) \, \md s \! + \!
(\overset{o}{m}_{2}^{\raise-1.0ex\hbox{$\scriptstyle 0$}})_{12} \mathbb{E}^{-
1} \! + \! (\mathscr{R}^{o,0}_{2}(n))_{12} \right. \\
&\left. + \, \left((\mathscr{R}^{o,0}_{1}(n))_{11}
(\overset{o}{m}_{1}^{\raise-1.0ex\hbox{$\scriptstyle 0$}})_{12} \! + \!
(\mathscr{R}^{o,0}_{1}(n))_{12}
(\overset{o}{m}_{1}^{\raise-1.0ex\hbox{$\scriptstyle 0$}})_{22} \right)
\mathbb{E}^{-1} \right), \\
(\mathrm{Y}^{o,0}_{2})_{21} &= \me^{-n \ell_{o}} \! \left(-(2n \! + \! 1) \!
\left((\overset{o}{m}_{1}^{\raise-1.0ex\hbox{$\scriptstyle 0$}})_{21}
\mathbb{E} \! + \! (\mathscr{R}^{o,0}_{1}(n))_{21} \right) \! \int_{J_{o}}s^{-
1} \psi_{V}^{o}(s) \, \md s \! + \!
(\overset{o}{m}_{2}^{\raise-1.0ex\hbox{$\scriptstyle 0$}})_{21} \mathbb{E} \!
+ \! (\mathscr{R}^{o,0}_{2}(n))_{21} \right. \\
&\left. + \, \left((\mathscr{R}^{o,0}_{1}(n))_{21}
(\overset{o}{m}_{1}^{\raise-1.0ex\hbox{$\scriptstyle 0$}})_{11} \! + \!
(\mathscr{R}^{o,0}_{1}(n))_{22}
(\overset{o}{m}_{1}^{\raise-1.0ex\hbox{$\scriptstyle 0$}})_{21} \right) \!
\mathbb{E} \right), \\
(\mathrm{Y}^{o,0}_{2})_{22} =& \, \tfrac{1}{2}(2n \! + \! 1)^{2} \! \left(
\int_{J_{o}}s^{-1} \psi_{V}^{o}(s) \, \md s \right)^{2} \! + \! \tfrac{1}{2}
(2n \! + \! 1) \int_{J_{o}}s^{-2} \psi_{V}^{o}(s) \, \md s \! + \! (2n \! + \!
1) \! \left((\overset{o}{m}_{1}^{\raise-1.0ex\hbox{$\scriptstyle 0$}})_{22}
\mathbb{E}^{-1} \! + \! (\mathscr{R}^{o,0}_{1}(n))_{22} \right) \\
\times& \, \int_{J_{o}}s^{-1} \psi_{V}^{o}(s) \, \md s \! + \!
(\overset{o}{m}_{2}^{\raise-1.0ex\hbox{$\scriptstyle 0$}})_{22} \mathbb{E}^{-
1} \! + \! (\mathscr{R}^{o,0}_{2}(n))_{22} \! + \! \left((\mathscr{R}^{o,0}_{1}
(n))_{21}(\overset{o}{m}_{1}^{\raise-1.0ex\hbox{$\scriptstyle 0$}})_{12} \! +
\! (\mathscr{R}^{o,0}_{1}(n))_{22}
(\overset{o}{m}_{1}^{\raise-1.0ex\hbox{$\scriptstyle 0$}})_{22} \right) \!
\mathbb{E}^{-1},
\end{align*}
with
\begin{align*}
(\overset{o}{m}_{1}^{\raise-1.0ex\hbox{$\scriptstyle 0$}})_{11} &= -\dfrac{
\boldsymbol{\theta}^{o}(\boldsymbol{u}^{o}_{+}(0) \! + \! \boldsymbol{d}_{o})
\mathbb{E}^{-1}}{\boldsymbol{\theta}^{o}(\boldsymbol{u}^{o}_{+}(0) \! - \!
\frac{1}{2 \pi}(n \! + \! \frac{1}{2}) \boldsymbol{\Omega}^{o} \! + \!
\boldsymbol{d}_{o})} \! \left(\dfrac{\theta^{o}_{0}(1,1,\boldsymbol{\Omega}^{o}
) \alpha^{o}_{0}(1,1,\pmb{0}) \! - \! \alpha^{o}_{0}(1,1,\boldsymbol{\Omega}^{
o}) \theta^{o}_{0}(1,1,\pmb{0})}{(\theta^{o}_{0}(1,1,\pmb{0}))^{2}} \right), \\
(\overset{o}{m}_{1}^{\raise-1.0ex\hbox{$\scriptstyle 0$}})_{12} &= \dfrac{
1}{4 \mi} \! \left(\sum_{k=1}^{N+1} \! \left(\dfrac{1}{b_{k-1}^{o}} \! - \!
\dfrac{1}{a_{k}^{o}} \right) \right) \! \dfrac{\boldsymbol{\theta}^{o}
(\boldsymbol{u}^{o}_{+}(0) \! + \! \boldsymbol{d}_{o}) \theta^{o}_{0}(-1,1,
\boldsymbol{\Omega}^{o}) \mathbb{E}^{-1}}{\boldsymbol{\theta}^{o}(\boldsymbol{
u}^{o}_{+}(0) \! - \! \frac{1}{2 \pi}(n \! + \! \frac{1}{2}) \boldsymbol{
\Omega}^{o} \! + \! \boldsymbol{d}_{o}) \theta^{o}_{0}(-1,1,\pmb{0})}, \\
(\overset{o}{m}_{1}^{\raise-1.0ex\hbox{$\scriptstyle 0$}})_{21}&= -\dfrac{1}{4
\mi} \! \left(\sum_{k=1}^{N+1} \! \left(\dfrac{1}{b_{k-1}^{o}} \! - \! \dfrac{
1}{a_{k}^{o}} \right) \right) \! \dfrac{\boldsymbol{\theta}^{o}(\boldsymbol{
u}^{o}_{+}(0) \! + \! \boldsymbol{d}_{o}) \theta^{o}_{0}(1,-1,\boldsymbol{
\Omega}^{o}) \mathbb{E}}{\boldsymbol{\theta}^{o}(-\boldsymbol{u}^{o}_{+}(0) \!
- \! \frac{1}{2 \pi}(n \! + \! \frac{1}{2}) \boldsymbol{\Omega}^{o} \! - \!
\boldsymbol{d}_{o}) \theta^{o}_{0}(1,-1,\pmb{0})}, \\
(\overset{o}{m}_{1}^{\raise-1.0ex\hbox{$\scriptstyle 0$}})_{22} &= -\dfrac{
\boldsymbol{\theta}^{o}(\boldsymbol{u}^{o}_{+}(0) \! + \! \boldsymbol{d}_{o}) 
\mathbb{E}}{\boldsymbol{\theta}^{o}(-\boldsymbol{u}^{o}_{+}(0) \! - \! \frac{
1}{2 \pi}(n \! + \! \frac{1}{2}) \boldsymbol{\Omega}^{o} \! - \! \boldsymbol{
d}_{o})} \! \left(\dfrac{\theta^{o}_{0}(-1,-1,\boldsymbol{\Omega}^{o}) 
\alpha^{o}_{0}(-1,-1,\pmb{0}) \! - \! \alpha^{o}_{0}(-1,-1,\boldsymbol{
\Omega}^{o}) \theta^{o}_{0}(-1,-1,\pmb{0})}{(\theta^{o}_{0}(-1,-1,\pmb{0}))^{
2}} \right), \\
(\overset{o}{m}_{2}^{\raise-1.0ex\hbox{$\scriptstyle 0$}})_{11} &= \left(
\theta^{o}_{0}(1,1,\boldsymbol{\Omega}^{o}) \! \left(-\beta^{0}_{0}(1,1,\pmb{
0}) \theta^{o}_{0}(1,1,\pmb{0}) \! + \! (\alpha^{o}_{0}(1,1,\pmb{0}))^{2}
\right) \! - \! \alpha^{o}_{0}(1,1,\boldsymbol{\Omega}^{o}) \alpha^{o}_{0}
(1,1,\pmb{0}) \theta^{o}_{0}(1,1,\pmb{0}) \right. \\
&\left. + \, \beta^{o}_{0}(1,1,\boldsymbol{\Omega}^{o})(\theta^{o}_{0}(1,1,
\pmb{0}))^{2} \right) \! \dfrac{(\theta^{o}_{0}(1,1,\pmb{0}))^{-3} \boldsymbol{
\theta}^{o}(\boldsymbol{u}^{o}_{+}(0) \! + \! \boldsymbol{d}_{o}) \mathbb{E}^{-
1}}{\boldsymbol{\theta}^{o}(\boldsymbol{u}^{o}_{+}(0) \! - \! \frac{1}{2 \pi}
(n \! + \! \frac{1}{2}) \boldsymbol{\Omega}^{o} \! + \! \boldsymbol{d}_{o})}
\! + \! \dfrac{1}{32} \! \left(\sum_{k=1}^{N+1} \! \left(\dfrac{1}{b_{k-1}^{o}
} \! - \! \dfrac{1}{a_{k}^{o}} \right) \right)^{2} \! \mathbb{E}^{-1}, \\
(\overset{o}{m}_{2}^{\raise-1.0ex\hbox{$\scriptstyle 0$}})_{12} &= -\dfrac{
\boldsymbol{\theta}^{o}(\boldsymbol{u}^{o}_{+}(0) \! + \! \boldsymbol{d}_{o})
\mathbb{E}^{-1}}{\boldsymbol{\theta}^{o}(\boldsymbol{u}^{o}_{+}(0) \! - \!
\frac{1}{2 \pi}(n \! + \! \frac{1}{2}) \boldsymbol{\Omega}^{o} \! + \!
\boldsymbol{d}_{o})} \! \left(\! \left(\dfrac{\theta^{o}_{0}(-1,1,\boldsymbol{
\Omega}^{o}) \alpha^{o}_{0}(-1,1,\pmb{0}) \! - \! \alpha^{o}_{0}(-1,1,
\boldsymbol{\Omega}^{o}) \theta^{o}_{0}(-1,1,\pmb{0})}{(\theta^{o}_{0}(-1,1,
\pmb{0}))^{2}} \right) \right. \\
&\left. \times \dfrac{1}{4 \mi} \! \left(\sum_{k=1}^{N+1} \! \left(\dfrac{1}{
b_{k-1}^{o}} \! - \! \dfrac{1}{a_{k}^{o}} \right) \right) \! - \! \dfrac{1}{8
\mi} \! \left(\sum_{k=1}^{N+1} \! \left(\dfrac{1}{(b_{k-1}^{o})^{2}} \! - \!
\dfrac{1}{(a_{k}^{o})^{2}} \right) \right) \! \dfrac{\theta^{o}_{0}(-1,1,
\boldsymbol{\Omega}^{o})}{\theta^{o}_{0}(-1,1,\pmb{0})} \right), \\
(\overset{o}{m}_{2}^{\raise-1.0ex\hbox{$\scriptstyle 0$}})_{21} &= \dfrac{
\boldsymbol{\theta}^{o}(\boldsymbol{u}^{o}_{+}(0) \! + \! \boldsymbol{d}_{o})
\mathbb{E}}{\boldsymbol{\theta}^{o}(-\boldsymbol{u}^{o}_{+}(0) \! - \! \frac{
1}{2 \pi}(n \! + \! \frac{1}{2}) \boldsymbol{\Omega}^{o} \! - \! \boldsymbol{
d}_{o})} \! \left(\! \left(\dfrac{\theta^{o}_{0}(1,-1,\boldsymbol{\Omega}^{o})
\alpha^{o}_{0}(1,-1,\pmb{0}) \! - \! \alpha^{o}_{0}(1,-1,\boldsymbol{\Omega}^{
o}) \theta^{o}_{0}(1,-1,\pmb{0})}{(\theta^{o}_{0}(1,-1,\pmb{0}))^{2}} \right)
\right. \\
&\left. \times \dfrac{1}{4 \mi} \! \left(\sum_{k=1}^{N+1} \! \left(\dfrac{1}{
b_{k-1}^{o}} \! - \! \dfrac{1}{a_{k}^{o}} \right) \right) \! - \! \dfrac{1}{8
\mi} \! \left(\sum_{k=1}^{N+1} \! \left(\dfrac{1}{(b_{k-1}^{o})^{2}} \! - \!
\dfrac{1}{(a_{k}^{o})^{2}} \right) \right) \! \dfrac{\theta^{o}_{0}(1,-1,
\boldsymbol{\Omega}^{o})}{\theta^{o}_{0}(1,-1,\pmb{0})} \right), \\
(\overset{o}{m}_{2}^{\raise-1.0ex\hbox{$\scriptstyle 0$}})_{22} &= \left(
\theta^{o}_{0}(-1,-1,\boldsymbol{\Omega}^{o}) \! \left(-\beta^{o}_{0}(-1,-1,
\pmb{0}) \theta^{o}_{0}(-1,-1,\pmb{0}) \! + \! (\alpha^{o}_{0}(-1,-1,\pmb{0})
)^{2} \right) \! - \! \alpha^{o}_{0}(-1,-1,\boldsymbol{\Omega}^{o}) \right. \\
&\left. \times \, \alpha^{o}_{0}(-1,-1,\pmb{0}) \theta^{o}_{0}(-1,-1,\pmb{0})
\! + \! \beta^{o}_{0}(-1,-1,\boldsymbol{\Omega}^{o})(\theta^{o}_{0}(-1,-1,
\pmb{0}))^{2} \right) \! (\theta^{o}_{0}(-1,-1,\pmb{0}))^{-3} \\
&\times \dfrac{\boldsymbol{\theta}^{o}(\boldsymbol{u}^{o}_{+}(0) \! + \!
\boldsymbol{d}_{o}) \mathbb{E}}{\boldsymbol{\theta}^{o}(-\boldsymbol{u}^{o}_{+}
(0) \! - \! \frac{1}{2 \pi}(n \! + \! \frac{1}{2}) \boldsymbol{\Omega}^{o} \!
- \! \boldsymbol{d}_{o})} \! + \! \dfrac{1}{32} \! \left(\sum_{k=1}^{N+1} \!
\left(\dfrac{1}{b_{k-1}^{o}} \! - \! \dfrac{1}{a_{k}^{o}} \right) \right)^{2}
\! \mathbb{E},
\end{align*}
and $(\star)_{ij}$, $i,j \! = \! 1,2$, denoting the $(i \, j)$-element of 
$\star$, where, for $\varepsilon_{1},\varepsilon_{2} \! = \! \pm 1$,
\begin{gather*}
\theta^{o}_{0}(\varepsilon_{1},\varepsilon_{2},\pmb{\mathscr{Z}}) \! := \!
\boldsymbol{\theta}^{o}(\varepsilon_{1} \boldsymbol{u}^{o}_{+}(0) \! - \!
\tfrac{1}{2 \pi}(n \! + \! \tfrac{1}{2}) \pmb{\mathscr{Z}} \! + \!
\varepsilon_{2} \boldsymbol{d}_{o}), \\
\alpha^{o}_{0}(\varepsilon_{1},\varepsilon_{2},\pmb{\mathscr{Z}}) \! := \! 2
\pi \mi \varepsilon_{1} \sum_{m \in \mathbb{Z}^{N}} (m,\widehat{\boldsymbol{
\alpha}}^{o}_{0}) \me^{2 \pi \mi (m,\varepsilon_{1} \boldsymbol{u}^{o}_{+}(0)
-\frac{1}{2 \pi}(n+\frac{1}{2}) \pmb{\mathscr{Z}}+\varepsilon_{2} \boldsymbol{
d}_{o})+\pi \mi (m,\tau^{o}m)},
\end{gather*}
where $\boldsymbol{u}^{o}_{+}(0) \! = \! \int_{a_{N+1}^{o}}^{0^{+}} 
\boldsymbol{\omega}^{o}$, $\widehat{\boldsymbol{\alpha}}^{o}_{0} \! = \! 
(\widehat{\alpha}_{0,1}^{o},\widehat{\alpha}^{o}_{0,2},\dotsc,\widehat{
\alpha}^{o}_{0,N})$, with $\widehat{\alpha}^{o}_{0,j} \! := \! (-1)^{\mathcal{
N}_{+}^{o}}(\prod_{i=1}^{N+1} \vert b_{i-1}^{o}a_{i}^{o} \vert)^{-1/2}c_{jN}^{
o}$, $j \! = \! 1,\dotsc,N$, where $\mathcal{N}_{+}^{o} \! \in \! \lbrace 0,
\dotsc,N \! + \! 1 \rbrace$ is the number of bands to the right of $z \! = 
\! 0$, and
\begin{gather*}
\beta^{o}_{0}(\varepsilon_{1},\varepsilon_{2},\pmb{\mathscr{Z}}) \! := \! 2 
\pi \sum_{m \in \mathbb{Z}^{N}} \! \left(\mi \varepsilon_{1}(m,\widehat{
\boldsymbol{\beta}}^{o}_{0}) \! - \!  \pi (m,\widehat{\boldsymbol{\alpha}}^{
o}_{0})^{2} \right) \! \me^{2 \pi \mi (m,\varepsilon_{1} \boldsymbol{u}^{o}_{
+}(0)-\frac{1}{2 \pi}(n+\frac{1}{2}) \pmb{\mathscr{Z}}+\varepsilon_{2}
\boldsymbol{d}_{o})+ \pi \mi (m,\tau^{o}m)},
\end{gather*}
where $\widehat{\boldsymbol{\beta}}^{o}_{0} \! = \! (\widehat{\beta}^{o}_{0,1},
\widehat{\beta}^{o}_{0,2},\dotsc,\widehat{\beta}^{o}_{0,N})$, with $\widehat{
\beta}^{o}_{0,j} \! := \! \tfrac{1}{2}(-1)^{\mathcal{N}_{+}^{o}}(\prod_{i=
1}^{N+1} \vert b_{i-1}^{o}a_{i}^{o} \vert)^{-1/2}(c^{o}_{jN-1} \! + \! 
\tfrac{1}{2}c_{jN}^{o} \sum_{k=1}^{N+1}((a_{k}^{o})^{-1} \! + \! (b_{k-1}^{
o})^{-1}))$, $j \! = \! 1,\dotsc,N$, where $c_{jN}^{o},c_{jN-1}^{o}$, $j \! 
= \! 1,\dotsc,N$, are obtained {}from Equations~{\rm (O1)} and~{\rm (O2)}, 
and, for $k \! = \! 2,3$,
\begin{align*}
\mathscr{R}^{o,0}_{k-1}(n) \underset{n \to \infty}{=}& \, \dfrac{1}{(n \! + \!
\frac{1}{2})} \mathlarger{\sum_{j=1}^{N+1}} \! \left(\dfrac{\left(\mathscr{
A}^{o}(b_{j-1}^{o}) \! \left(\widehat{\alpha}_{1}^{o}(b_{j-1}^{o}) \! + \! k
(b_{j-1}^{o})^{-1} \widehat{\alpha}_{0}^{o}(b_{j-1}^{o}) \right) \! - \!
\mathscr{B}^{o}(b_{j-1}^{o}) \widehat{\alpha}_{0}^{o}(b_{j-1}^{o}) \right)}{
(b_{j-1}^{o})^{k}(\widehat{\alpha}_{0}^{o}(b_{j-1}^{o}))^{2}} \right. \\
+&\left. \, \dfrac{\left(\mathscr{A}^{o}(a_{j}^{o}) \! \left(\widehat{
\alpha}_{1}^{o}(a_{j}^{o}) \! + \! k(a_{j}^{o})^{-1} \widehat{\alpha}_{0}^{o}
(a_{j}^{o}) \right) \! - \! \mathscr{B}^{o}(a_{j}^{o}) \widehat{\alpha}_{0}^{
o}(a_{j}^{o}) \right)}{(a_{j}^{o})^{k}(\widehat{\alpha}_{0}^{o}(a_{j}^{o}))^{
2}} \right) \! + \! \mathcal{O} \! \left(\dfrac{c(n)}{(n \! + \! \frac{1}{2}
)^{2}} \right),
\end{align*}
with all parameters defined in Theorem~{\rm 3.1}, 
Equations~{\rm (3.52)--(3.100)}.
\end{dy}
\addtocounter{z0}{1}
\begin{dy}[{\rm \cite{a22}}]
Let all the conditions stated in Theorem~{\rm B.1} be valid, and let all 
parameters be as defined therein. Let $\overset{o}{\operatorname{Y}} \colon 
\mathbb{C} \setminus \mathbb{R} \! \to \! \operatorname{SL}_{2}(\mathbb{C})$ 
be the unique solution of \pmb{{\rm RHP2}} with integral 
representation~{\rm (1.30)}. Then,
\begin{equation*}
\overset{o}{\mathrm{Y}}(z)z^{-(n+1) \sigma_{3}} \underset{z \to \infty}{=}
\mathrm{Y}^{o,\infty}_{0} \! + \! \dfrac{1}{z} \mathrm{Y}^{o,\infty}_{1} \! +
\! \dfrac{1}{z^{2}} \mathrm{Y}^{o,\infty}_{2} \! + \! \mathcal{O} \! \left(
\dfrac{1}{z^{3}} \right),
\end{equation*}
where
\begin{align*}
(\mathrm{Y}^{o,\infty}_{0})_{11} =& \, (\widehat{Q}^{o}_{0})_{11} \me^{-2(n+
\frac{1}{2}) \int_{J_{o}} \ln (\lvert s \rvert) \psi_{V}^{o}(s) \, \md s}, \\
(\mathrm{Y}^{o,\infty}_{0})_{12} =& \, (\widehat{Q}^{o}_{0})_{12} \me^{n(
\ell_{o}+2(1+\frac{1}{2n}) \int_{J_{o}} \ln (\lvert s \rvert) \psi_{V}^{o}(s)
\, \md s)}, \\
(\mathrm{Y}^{o,\infty}_{0})_{21} =& \, (\widehat{Q}^{o}_{0})_{21} \me^{-n(
\ell_{o}+2(1+\frac{1}{2n}) \int_{J_{o}} \ln (\lvert s \rvert) \psi_{V}^{o}(s)
\, \md s)}, \\
(\mathrm{Y}^{o,\infty}_{0})_{22} =& \, (\widehat{Q}^{o}_{0})_{22} \me^{2(n+
\frac{1}{2}) \int_{J_{o}} \ln (\lvert s \rvert) \psi_{V}^{o}(s) \, \md s}, \\
(\mathrm{Y}^{o,\infty}_{1})_{11} =& \, \left(\! (\widehat{Q}^{o}_{1})_{11} \!
- \! (2n \! + \! 1)(\widehat{Q}^{o}_{0})_{11} \int_{J_{o}}s \psi_{V}^{o}(s) \,
\md s \right) \! \me^{-2(n+\frac{1}{2}) \int_{J_{o}} \ln (\lvert s \rvert)
\psi_{V}^{o}(s) \, \md s}, \\
(\mathrm{Y}^{o,\infty}_{1})_{12} =& \, \left(\! (\widehat{Q}^{o}_{1})_{12} \!
+ \! (2n \! + \! 1)(\widehat{Q}^{o}_{0})_{12} \int_{J_{o}}s \psi_{V}^{o}(s) \,
\md s \right) \! \me^{n(\ell_{o}+2(1+\frac{1}{2n}) \int_{J_{o}} \ln (\lvert s
\rvert) \psi_{V}^{o}(s) \, \md s)}, \\
(\mathrm{Y}^{o,\infty}_{1})_{21} =& \, \left(\! (\widehat{Q}^{o}_{1})_{21} \!
- \! (2n \! + \! 1)(\widehat{Q}^{o}_{0})_{21} \int_{J_{o}}s \psi_{V}^{o}(s) \,
\md s \right) \! \me^{-n(\ell_{o}+2(1+\frac{1}{2n}) \int_{J_{o}} \ln (\lvert s
\rvert) \psi_{V}^{o}(s) \, \md s)}, \\
(\mathrm{Y}^{o,\infty}_{1})_{22} =& \, \left(\! (\widehat{Q}^{o}_{1})_{22} \!
+ \! (2n \! + \! 1)(\widehat{Q}^{o}_{0})_{22} \int_{J_{o}}s \psi_{V}^{o}(s) \,
\md s \right) \! \me^{2(n+\frac{1}{2}) \int_{J_{o}} \ln (\lvert s \rvert)
\psi_{V}^{o}(s) \, \md s}, \\
(\mathrm{Y}^{o,\infty}_{2})_{11} =& \left(\! (\widehat{Q}^{o}_{2})_{11} \! -
\! (2n \! + \! 1)(\widehat{Q}^{o}_{1})_{11} \int_{J_{o}}s \psi_{V}^{o}(s) \,
\md s \! + \! (\widehat{Q}^{o}_{0})_{11} \! \left(\tfrac{1}{2}(2n \! + \! 
1)^{2} \! \left(\int_{J_{o}}s \psi_{V}^{o}(s) \, \md s \right)^{2} \right. 
\right. \\
-& \left. \left. \tfrac{1}{2}(2n \! + \! 1) \int_{J_{o}}s^{2} \psi_{V}^{o}(s)
\, \md s \right) \right) \! \me^{-2(n+\frac{1}{2}) \int_{J_{o}} \ln (\lvert s
\rvert) \psi_{V}^{o}(s) \, \md s}, \\
(\mathrm{Y}^{o,\infty}_{2})_{12} =& \left(\! (\widehat{Q}^{o}_{2})_{12} \! +
\! (2n \! + \! 1)(\widehat{Q}^{o}_{1})_{12} \int_{J_{o}}s \psi_{V}^{o}(s) \,
\md s \! + \! (\widehat{Q}^{o}_{0})_{12} \! \left(\tfrac{1}{2}(2n \! + \! 
1)^{2} \! \left(\int_{J_{o}}s \psi_{V}^{o}(s) \, \md s \right)^{2} \right. 
\right. \\
+& \left. \left. \tfrac{1}{2}(2n \! + \! 1) \int_{J_{o}}s^{2} \psi_{V}^{o}(s)
\, \md s \right) \right) \! \me^{n(\ell_{o}+2(1+\frac{1}{2n}) \int_{J_{o}} \ln
(\lvert s \rvert) \psi_{V}^{o}(s) \, \md s)}, \\
(\mathrm{Y}^{o,\infty}_{2})_{21} =& \left(\! (\widehat{Q}^{o}_{2})_{21} \! -
\! (2n \! + \! 1)(\widehat{Q}^{o}_{1})_{21} \int_{J_{o}}s \psi_{V}^{o}(s) \,
\md s \! + \! (\widehat{Q}^{o}_{0})_{21} \! \left(\tfrac{1}{2}(2n \! + \! 
1)^{2} \! \left(\int_{J_{o}}s \psi_{V}^{o}(s) \, \md s \right)^{2} \right. 
\right. \\
-& \left. \left. \tfrac{1}{2}(2n \! + \! 1) \int_{J_{o}}s^{2} \psi_{V}^{o}(s)
\, \md s \right) \right) \! \me^{-n(\ell_{o}+2(1+\frac{1}{2n}) \int_{J_{o}}
\ln (\lvert s \rvert) \psi_{V}^{o}(s) \, \md s)}, \\
(\mathrm{Y}^{o,\infty}_{2})_{22} =& \left(\! (\widehat{Q}^{o}_{2})_{22} \! +
\! (2n \! + \! 1)(\widehat{Q}^{o}_{1})_{22} \int_{J_{o}}s \psi_{V}^{o}(s) \,
\md s \! + \! (\widehat{Q}^{o}_{0})_{22} \! \left(\tfrac{1}{2}(2n \! + \! 
1)^{2} \! \left(\int_{J_{o}}s \psi_{V}^{o}(s) \, \md s \right)^{2} \right. 
\right. \\
+& \left. \left. \tfrac{1}{2}(2n \! + \! 1) \int_{J_{o}}s^{2} \psi_{V}^{o}(s)
\, \md s \right) \right) \! \me^{2(n+\frac{1}{2}) \int_{J_{o}} \ln (\lvert s
\rvert) \psi_{V}^{o}(s) \, \md s},
\end{align*}
with
\begin{align*}
(\widehat{Q}^{o}_{0})_{11} :=& \,
(\overset{o}{m}_{0}^{\raise-1.0ex\hbox{$\scriptstyle \infty$}})_{11} \! \left(
1 \! + \! (\mathscr{R}^{o,\infty}_{0}(n))_{11} \right) \! + \! (\mathscr{R}^{o,
\infty}_{0}(n))_{12}
(\overset{o}{m}_{0}^{\raise-1.0ex\hbox{$\scriptstyle \infty$}})_{21}, \\
(\widehat{Q}^{o}_{0})_{12} :=& \,
(\overset{o}{m}_{0}^{\raise-1.0ex\hbox{$\scriptstyle \infty$}})_{12} \! \left(
1 \! + \! (\mathscr{R}^{o,\infty}_{0}(n))_{11} \right) \! + \! (\mathscr{R}^{o,
\infty}_{0}(n))_{12}
(\overset{o}{m}_{0}^{\raise-1.0ex\hbox{$\scriptstyle \infty$}})_{22}, \\
(\widehat{Q}^{o}_{0})_{21} :=& \,
(\overset{o}{m}_{0}^{\raise-1.0ex\hbox{$\scriptstyle \infty$}})_{21} \! \left(
1 \! + \! (\mathscr{R}^{o,\infty}_{0}(n))_{22} \right) \! + \! (\mathscr{R}^{o,
\infty}_{0}(n))_{21}
(\overset{o}{m}_{0}^{\raise-1.0ex\hbox{$\scriptstyle \infty$}})_{11}, \\
(\widehat{Q}^{o}_{0})_{22} :=& \,
(\overset{o}{m}_{0}^{\raise-1.0ex\hbox{$\scriptstyle \infty$}})_{22} \! \left(
1 \! + \! (\mathscr{R}^{o,\infty}_{0}(n))_{22} \right) \! + \! (\mathscr{R}^{o,
\infty}_{0}(n))_{21}
(\overset{o}{m}_{0}^{\raise-1.0ex\hbox{$\scriptstyle \infty$}})_{12}, \\
(\widehat{Q}^{o}_{1})_{11} :=& \,
(\overset{o}{m}_{1}^{\raise-1.0ex\hbox{$\scriptstyle \infty$}})_{11} \! \left(
1 \! + \! (\mathscr{R}^{o,\infty}_{0}(n))_{11} \right) \! + \! (\mathscr{R}^{o,
\infty}_{0}(n))_{12}
(\overset{o}{m}_{1}^{\raise-1.0ex\hbox{$\scriptstyle \infty$}})_{21} \! + \!
(\mathscr{R}^{o,\infty}_{1}(n))_{11}
(\overset{o}{m}_{0}^{\raise-1.0ex\hbox{$\scriptstyle \infty$}})_{11} \\
+& \, (\mathscr{R}^{o,\infty}_{1}(n))_{12}
(\overset{o}{m}_{0}^{\raise-1.0ex\hbox{$\scriptstyle \infty$}})_{21}, \\
(\widehat{Q}^{o}_{1})_{12} :=& \,
(\overset{o}{m}_{1}^{\raise-1.0ex\hbox{$\scriptstyle \infty$}})_{12} \! \left(
1 \! + \! (\mathscr{R}^{o,\infty}_{0}(n))_{11} \right) \! + \! (\mathscr{R}^{o,
\infty}_{0}(n))_{12}
(\overset{o}{m}_{1}^{\raise-1.0ex\hbox{$\scriptstyle \infty$}})_{22} \! + \!
(\mathscr{R}^{o,\infty}_{1}(n))_{11}
(\overset{o}{m}_{0}^{\raise-1.0ex\hbox{$\scriptstyle \infty$}})_{12} \\
+& \, (\mathscr{R}^{o,\infty}_{1}(n))_{12}
(\overset{o}{m}_{0}^{\raise-1.0ex\hbox{$\scriptstyle \infty$}})_{22}, \\
(\widehat{Q}^{o}_{1})_{21} :=& \,
(\overset{o}{m}_{1}^{\raise-1.0ex\hbox{$\scriptstyle \infty$}})_{21} \! \left(
1 \! + \! (\mathscr{R}^{o,\infty}_{0}(n))_{22} \right) \! + \! (\mathscr{R}^{o,
\infty}_{0}(n))_{21}
(\overset{o}{m}_{1}^{\raise-1.0ex\hbox{$\scriptstyle \infty$}})_{11} \! + \!
(\mathscr{R}^{o,\infty}_{1}(n))_{21}
(\overset{o}{m}_{0}^{\raise-1.0ex\hbox{$\scriptstyle \infty$}})_{11} \\
+& \, (\mathscr{R}^{o,\infty}_{1}(n))_{22}
(\overset{o}{m}_{0}^{\raise-1.0ex\hbox{$\scriptstyle \infty$}})_{21}, \\
(\widehat{Q}^{o}_{1})_{22} :=& \,
(\overset{o}{m}_{1}^{\raise-1.0ex\hbox{$\scriptstyle \infty$}})_{22} \! \left(
1 \! + \! (\mathscr{R}^{o,\infty}_{0}(n))_{22} \right) \! + \! (\mathscr{R}^{o,
\infty}_{0}(n))_{21}
(\overset{o}{m}_{1}^{\raise-1.0ex\hbox{$\scriptstyle \infty$}})_{12} \! + \!
(\mathscr{R}^{o,\infty}_{1}(n))_{21}
(\overset{o}{m}_{0}^{\raise-1.0ex\hbox{$\scriptstyle \infty$}})_{12} \\
+& \, (\mathscr{R}^{o,\infty}_{1}(n))_{22}
(\overset{o}{m}_{0}^{\raise-1.0ex\hbox{$\scriptstyle \infty$}})_{22}, \\
(\widehat{Q}^{o}_{2})_{11} :=& \,
(\overset{o}{m}_{2}^{\raise-1.0ex\hbox{$\scriptstyle \infty$}})_{11} \! \left(
1 \! + \! (\mathscr{R}^{o,\infty}_{0}(n))_{11} \right) \! + \! (\mathscr{R}^{o,
\infty}_{0}(n))_{12}
(\overset{o}{m}_{2}^{\raise-1.0ex\hbox{$\scriptstyle \infty$}})_{21} \! + \!
(\mathscr{R}^{o,\infty}_{1}(n))_{11}
(\overset{o}{m}_{1}^{\raise-1.0ex\hbox{$\scriptstyle \infty$}})_{11} \\
+& \, (\mathscr{R}^{o,\infty}_{1}(n))_{12}
(\overset{o}{m}_{1}^{\raise-1.0ex\hbox{$\scriptstyle \infty$}})_{21} \! + \!
(\mathscr{R}^{o,\infty}_{2}(n))_{11}
(\overset{o}{m}_{0}^{\raise-1.0ex\hbox{$\scriptstyle \infty$}})_{11} \! + \!
(\mathscr{R}^{o,\infty}_{2}(n))_{12}
(\overset{o}{m}_{0}^{\raise-1.0ex\hbox{$\scriptstyle \infty$}})_{21}, \\
(\widehat{Q}^{o}_{2})_{12} :=& \,
(\overset{o}{m}_{2}^{\raise-1.0ex\hbox{$\scriptstyle \infty$}})_{12} \! \left(
1 \! + \! (\mathscr{R}^{o,\infty}_{0}(n))_{11} \right) \! + \! (\mathscr{R}^{o,
\infty}_{0}(n))_{12}
(\overset{o}{m}_{2}^{\raise-1.0ex\hbox{$\scriptstyle \infty$}})_{22} \! + \!
(\mathscr{R}^{o,\infty}_{1}(n))_{11}
(\overset{o}{m}_{1}^{\raise-1.0ex\hbox{$\scriptstyle \infty$}})_{12} \\
+& \, (\mathscr{R}^{o,\infty}_{1}(n))_{12}
(\overset{o}{m}_{1}^{\raise-1.0ex\hbox{$\scriptstyle \infty$}})_{22} \! + \!
(\mathscr{R}^{o,\infty}_{2}(n))_{11}
(\overset{o}{m}_{0}^{\raise-1.0ex\hbox{$\scriptstyle \infty$}})_{12} \! + \!
(\mathscr{R}^{o,\infty}_{2}(n))_{12}
(\overset{o}{m}_{0}^{\raise-1.0ex\hbox{$\scriptstyle \infty$}})_{22}, \\
(\widehat{Q}^{o}_{2})_{21} :=& \,
(\overset{o}{m}_{2}^{\raise-1.0ex\hbox{$\scriptstyle \infty$}})_{21} \! \left(
1 \! + \! (\mathscr{R}^{o,\infty}_{0}(n))_{22} \right) \! + \! (\mathscr{R}^{o,
\infty}_{0}(n))_{21}
(\overset{o}{m}_{2}^{\raise-1.0ex\hbox{$\scriptstyle \infty$}})_{11} \! + \!
(\mathscr{R}^{o,\infty}_{1}(n))_{21}
(\overset{o}{m}_{1}^{\raise-1.0ex\hbox{$\scriptstyle \infty$}})_{11} \\
+& \, (\mathscr{R}^{o,\infty}_{1}(n))_{22}
(\overset{o}{m}_{1}^{\raise-1.0ex\hbox{$\scriptstyle \infty$}})_{21} \! + \!
(\mathscr{R}^{o,\infty}_{2}(n))_{21}
(\overset{o}{m}_{0}^{\raise-1.0ex\hbox{$\scriptstyle \infty$}})_{11} \! + \!
(\mathscr{R}^{o,\infty}_{2}(n))_{22}
(\overset{o}{m}_{0}^{\raise-1.0ex\hbox{$\scriptstyle \infty$}})_{21}, \\
(\widehat{Q}^{o}_{2})_{22} :=& \,
(\overset{o}{m}_{2}^{\raise-1.0ex\hbox{$\scriptstyle \infty$}})_{22} \! \left(
1 \! + \! (\mathscr{R}^{o,\infty}_{0}(n))_{22} \right) \! + \! (\mathscr{R}^{o,
\infty}_{0}(n))_{21}
(\overset{o}{m}_{2}^{\raise-1.0ex\hbox{$\scriptstyle \infty$}})_{12} \! + \!
(\mathscr{R}^{o,\infty}_{1}(n))_{21}
(\overset{o}{m}_{1}^{\raise-1.0ex\hbox{$\scriptstyle \infty$}})_{12} \\
+& \, (\mathscr{R}^{o,\infty}_{1}(n))_{22}
(\overset{o}{m}_{1}^{\raise-1.0ex\hbox{$\scriptstyle \infty$}})_{22} \! + \!
(\mathscr{R}^{o,\infty}_{2}(n))_{21}
(\overset{o}{m}_{0}^{\raise-1.0ex\hbox{$\scriptstyle \infty$}})_{12} \! + \!
(\mathscr{R}^{o,\infty}_{2}(n))_{22}
(\overset{o}{m}_{0}^{\raise-1.0ex\hbox{$\scriptstyle \infty$}})_{22},
\end{align*}
and
\begin{align*}
(\overset{o}{m}_{0}^{\raise-1.0ex\hbox{$\scriptstyle \infty$}})_{11} &= \dfrac{
\boldsymbol{\theta}^{o}(\boldsymbol{u}^{o}_{+}(0) \! + \! \boldsymbol{d}_{o})
\mathbb{E}^{-1}}{\boldsymbol{\theta}^{o}(\boldsymbol{u}^{o}_{+}(0) \! - \!
\frac{1}{2 \pi}(n \! + \! \frac{1}{2}) \boldsymbol{\Omega}^{o} \! + \!
\boldsymbol{d}_{o})} \! \left(\dfrac{\gamma^{o}_{0} \! + \! (\gamma^{o}_{0})^{-
1}}{2} \right) \! \dfrac{\theta^{o}_{\infty}(1,1,\boldsymbol{\Omega}^{o})}{
\theta^{o}_{\infty}(1,1,\pmb{0})}, \\
(\overset{o}{m}_{0}^{\raise-1.0ex\hbox{$\scriptstyle \infty$}})_{12} &= \dfrac{
\boldsymbol{\theta}^{o}(\boldsymbol{u}^{o}_{+}(0) \! + \! \boldsymbol{d}_{o})
\mathbb{E}^{-1}}{\boldsymbol{\theta}^{o}(\boldsymbol{u}^{o}_{+}(0) \! - \!
\frac{1}{2 \pi}(n \! + \! \frac{1}{2}) \boldsymbol{\Omega}^{o} \! + \!
\boldsymbol{d}_{o})} \! \left(\dfrac{\gamma^{o}_{0} \! - \! (\gamma^{o}_{0})^{-
1}}{2 \mi} \right) \! \dfrac{\theta^{o}_{\infty}(-1,1,\boldsymbol{\Omega}^{o})
}{\theta^{o}_{\infty}(-1,1,\pmb{0})}, \\
(\overset{o}{m}_{0}^{\raise-1.0ex\hbox{$\scriptstyle \infty$}})_{21} &=
-\dfrac{\boldsymbol{\theta}^{o}(\boldsymbol{u}^{o}_{+}(0) \! + \! \boldsymbol{
d}_{o}) \mathbb{E}}{\boldsymbol{\theta}^{o}(-\boldsymbol{u}^{o}_{+}(0) \! -
\! \frac{1}{2 \pi}(n \! + \! \frac{1}{2}) \boldsymbol{\Omega}^{o} \! - \!
\boldsymbol{d}_{o})} \! \left(\dfrac{\gamma^{o}_{0} \! - \! (\gamma^{o}_{0})^{-
1}}{2 \mi} \right) \! \dfrac{\theta^{o}_{\infty}(1,-1,\boldsymbol{\Omega}^{o})
}{\theta^{o}_{\infty}(1,-1,\pmb{0})}, \\
(\overset{o}{m}_{0}^{\raise-1.0ex\hbox{$\scriptstyle \infty$}})_{22} &= \dfrac{
\boldsymbol{\theta}^{o}(\boldsymbol{u}^{o}_{+}(0) \! + \! \boldsymbol{d}_{o})
\mathbb{E}}{\boldsymbol{\theta}^{o}(-\boldsymbol{u}^{o}_{+}(0) \! - \! \frac{
1}{2 \pi}(n \! + \! \frac{1}{2}) \boldsymbol{\Omega}^{o} \! - \! \boldsymbol{
d}_{o})} \! \left(\dfrac{\gamma^{o}_{0} \! + \! (\gamma^{o}_{0})^{-1}}{2}
\right) \! \dfrac{\theta^{o}_{\infty}(-1,-1,\boldsymbol{\Omega}^{o})}{\theta^{
o}_{\infty}(-1,-1,\pmb{0})}, \\
(\overset{o}{m}_{1}^{\raise-1.0ex\hbox{$\scriptstyle \infty$}})_{11} &= \dfrac{
\boldsymbol{\theta}^{o}(\boldsymbol{u}^{o}_{+}(0) \! + \! \boldsymbol{d}_{o})
\mathbb{E}^{-1}}{\boldsymbol{\theta}^{o}(\boldsymbol{u}^{o}_{+}(0) \! - \!
\frac{1}{2 \pi}(n \! + \! \frac{1}{2}) \boldsymbol{\Omega}^{o} \! + \!
\boldsymbol{d}_{o})} \! \left(\! \left(\dfrac{\widetilde{\alpha}_{\infty}^{o}
(1,1,\pmb{0}) \theta^{o}_{\infty}(1,1,\boldsymbol{\Omega}^{o}) \! - \!
\widetilde{\alpha}_{\infty}^{o}(1,1,\boldsymbol{\Omega}^{o}) \theta^{o}_{
\infty}(1,1,\pmb{0})}{(\theta^{o}_{\infty}(1,1,\pmb{0}))^{2}} \right) \right. 
\\
&\left. \times \left(\dfrac{\gamma^{o}_{0} \! + \! (\gamma^{o}_{0})^{-1}}{2}
\right) \! - \! \left(\dfrac{\gamma^{o}_{0} \! - \! (\gamma^{o}_{0})^{-1}}{8}
\right) \! \left(\sum_{k=1}^{N+1} \! \left(a_{k}^{o} \! - \! b_{k-1}^{o}
\right) \right) \! \dfrac{\theta^{o}_{\infty}(1,1,\boldsymbol{\Omega}^{o})}{
\theta^{o}_{\infty}(1,1,\pmb{0})} \right), \\
(\overset{o}{m}_{1}^{\raise-1.0ex\hbox{$\scriptstyle \infty$}})_{12} &= \dfrac{
\boldsymbol{\theta}^{o}(\boldsymbol{u}^{o}_{+}(0) \! + \! \boldsymbol{d}_{o})
\mathbb{E}^{-1}}{\boldsymbol{\theta}^{o}(\boldsymbol{u}^{o}_{+}(0) \! - \!
\frac{1}{2 \pi}(n \! + \! \frac{1}{2}) \boldsymbol{\Omega}^{o} \! + \!
\boldsymbol{d}_{o})} \! \left(\! \left(\dfrac{\widetilde{\alpha}_{\infty}^{o}
(-1,1,\pmb{0}) \theta^{o}_{\infty}(-1,1,\boldsymbol{\Omega}^{o}) \! - \!
\widetilde{\alpha}_{\infty}^{o}(-1,1,\boldsymbol{\Omega}^{o}) \theta^{o}_{
\infty}(-1,1,\pmb{0})}{(\theta^{o}_{\infty}(-1,1,\pmb{0}))^{2}} \right)
\right. \\
&\left. \times \left(\dfrac{\gamma^{o}_{0} \! - \! (\gamma^{o}_{0})^{-1}}{2
\mi} \right) \! - \! \left(\dfrac{\gamma^{o}_{0} \! + \! (\gamma^{o}_{0})^{-1}
}{8 \mi} \right) \! \left(\sum_{k=1}^{N+1} \! \left(a_{k}^{o} \! - \! b_{k-1}^{
o} \right) \right) \! \dfrac{\theta^{o}_{\infty}(-1,1,\boldsymbol{\Omega}^{o})
}{\theta^{o}_{\infty}(-1,1,\pmb{0})} \right), \\
(\overset{o}{m}_{1}^{\raise-1.0ex\hbox{$\scriptstyle \infty$}})_{21} &=
-\dfrac{\boldsymbol{\theta}^{o}(\boldsymbol{u}^{o}_{+}(0) \! + \! \boldsymbol{
d}_{o}) \mathbb{E}}{\boldsymbol{\theta}^{o}(-\boldsymbol{u}^{o}_{+}(0) \! - \!
\frac{1}{2 \pi}(n \! + \! \frac{1}{2}) \boldsymbol{\Omega}^{o} \! - \!
\boldsymbol{d}_{o})} \! \left(\! \left(\dfrac{\widetilde{\alpha}_{\infty}^{o}
(1,-1,\pmb{0}) \theta^{o}_{\infty}(1,-1,\boldsymbol{\Omega}^{o}) \! - \!
\widetilde{\alpha}_{\infty}^{o}(1,-1,\boldsymbol{\Omega}^{o}) \theta^{o}_{
\infty}(1,-1,\pmb{0})}{(\theta^{o}_{\infty}(1,-1,\pmb{0}))^{2}} \right)
\right. \\
&\left. \times \left(\dfrac{\gamma^{o}_{0} \! - \! (\gamma^{o}_{0})^{-1}}{2
\mi} \right) \! - \! \left(\dfrac{\gamma^{o}_{0} \! + \! (\gamma^{o}_{0})^{-1}
}{8 \mi} \right) \! \left(\sum_{k=1}^{N+1} \! \left(a_{k}^{o} \! - \! b_{k-1}^{
o} \right) \right) \! \dfrac{\theta^{o}_{\infty}(1,-1,\boldsymbol{\Omega}^{o}
)}{\theta^{o}_{\infty}(1,-1,\pmb{0})} \right), \\
(\overset{o}{m}_{1}^{\raise-1.0ex\hbox{$\scriptstyle \infty$}})_{22} &= \dfrac{
\boldsymbol{\theta}^{o}(\boldsymbol{u}^{o}_{+}(0) \! + \! \boldsymbol{d}_{o})
\mathbb{E}}{\boldsymbol{\theta}^{o}(-\boldsymbol{u}^{o}_{+}(0) \! - \! \frac{
1}{2 \pi}(n \! + \! \frac{1}{2}) \boldsymbol{\Omega}^{o} \! - \! \boldsymbol{
d}_{o})} \! \left(\! \left(\dfrac{\widetilde{\alpha}_{\infty}^{o}(-1,-1,\pmb{
0}) \theta^{o}_{\infty}(-1,-1,\boldsymbol{\Omega}^{o}) \! - \! \widetilde{
\alpha}_{\infty}^{o}(-1,-1,\boldsymbol{\Omega}^{o}) \theta^{o}_{\infty}(-1,-1,
\pmb{0})}{(\theta^{o}_{\infty}(-1,-1,\pmb{0}))^{2}} \right) \right. \\
&\left. \times \left(\dfrac{\gamma^{o}_{0} \! + \! (\gamma^{o}_{0})^{-1}}{2}
\right) \! - \! \left(\dfrac{\gamma^{o}_{0} \! - \! (\gamma^{o}_{0})^{-1}}{8}
\right) \! \left(\sum_{k=1}^{N+1} \! \left(a_{k}^{o} \! - \! b_{k-1}^{o}
\right) \right) \! \dfrac{\theta^{o}_{\infty}(-1,-1,\boldsymbol{\Omega}^{o})}{
\theta^{o}_{\infty}(-1,-1,\pmb{0})} \right), \\
(\overset{o}{m}_{2}^{\raise-1.0ex\hbox{$\scriptstyle \infty$}})_{11} &= \dfrac{
\boldsymbol{\theta}^{o}(\boldsymbol{u}^{o}_{+}(0) \! + \! \boldsymbol{d}_{o})
\mathbb{E}^{-1}}{\boldsymbol{\theta}^{o}(\boldsymbol{u}^{o}_{+}(0) \! - \!
\frac{1}{2 \pi}(n \! + \! \frac{1}{2}) \boldsymbol{\Omega}^{o} \! + \!
\boldsymbol{d}_{o})} \! \left(\! \left(\theta^{o}_{\infty}(1,1,\boldsymbol{
\Omega}^{o}) \! \left((\widetilde{\alpha}_{\infty}^{o}(1,1,\pmb{0}))^{2} \! +
\! \beta^{o}_{\infty}(1,1,\pmb{0}) \theta^{o}_{\infty}(1,1,\pmb{0}) \right)
\right. \right. \\
&\left. \left. - \, \widetilde{\alpha}_{\infty}^{o}(1,1,\boldsymbol{\Omega}^{
o}) \widetilde{\alpha}_{\infty}^{o}(1,1,\pmb{0}) \theta^{o}_{\infty}(1,1,\pmb{
0}) \! - \! \beta^{o}_{\infty}(1,1,\boldsymbol{\Omega}^{o})(\theta^{o}_{\infty}
(1,1,\pmb{0}))^{2} \right) \! \tfrac{1}{(\theta^{o}_{\infty}(1,1,\pmb{0}))^{3}
} \right. \\
&\left. \times \left(\dfrac{\gamma^{o}_{0} \! + \! (\gamma^{o}_{0})^{-1}}{2}
\right) \! + \! \left(\dfrac{\widetilde{\alpha}_{\infty}^{o}(1,1,\pmb{0})
\theta^{o}_{\infty}(1,1,\boldsymbol{\Omega}^{o}) \! - \! \widetilde{\alpha}_{
\infty}^{o}(1,1,\boldsymbol{\Omega}^{o}) \theta^{o}_{\infty}(1,1,\pmb{0})}{(
\theta^{o}_{\infty}(1,1,\pmb{0}))^{2}} \right) \right. \\
&\left. \times \left(\sum_{k=1}^{N+1} \! \left(b_{k-1}^{o} \! - \! a_{k}^{o}
\right) \right) \! \left(\dfrac{\gamma^{o}_{0} \! - \! (\gamma^{o}_{0})^{-1}}{
8} \right) \! + \! \left(\! \left(\dfrac{\gamma^{o}_{0} \! - \! (\gamma^{o}_{
0})^{-1}}{16} \right) \! \sum_{k=1}^{N+1} \! \left((b_{k-1}^{o})^{2} \! - \!
(a_{k}^{o})^{2} \right) \right. \right. \\
&\left. \left. +\left(\dfrac{\gamma^{o}_{0} \! + \! (\gamma^{o}_{0})^{-1}}{64}
\right) \! \left(\sum_{k=1}^{N+1} \! \left(a_{k}^{o} \! - \! b_{k-1}^{o}
\right) \right)^{2} \right) \! \dfrac{\theta^{o}_{\infty}(1,1,\boldsymbol{
\Omega}^{o})}{\theta^{o}_{\infty}(1,1,\pmb{0})} \right), \\
(\overset{o}{m}_{2}^{\raise-1.0ex\hbox{$\scriptstyle \infty$}})_{12} &= \dfrac{
\boldsymbol{\theta}^{o}(\boldsymbol{u}^{o}_{+}(0) \! + \! \boldsymbol{d}_{o})
\mathbb{E}^{-1}}{\boldsymbol{\theta}^{o}(\boldsymbol{u}^{o}_{+}(0) \! - \!
\frac{1}{2 \pi}(n \! + \! \frac{1}{2}) \boldsymbol{\Omega}^{o} \! + \!
\boldsymbol{d}_{o})} \! \left(\! \left(\theta^{o}_{\infty}(-1,1,\boldsymbol{
\Omega}^{o}) \! \left((\widetilde{\alpha}_{\infty}^{o}(-1,1,\pmb{0}))^{2} \! +
\! \beta^{o}_{\infty}(-1,1,\pmb{0}) \theta^{o}_{\infty}(-1,1,\pmb{0}) \right)
\right. \right. \\
&\left. \left. - \, \widetilde{\alpha}_{\infty}^{o}(-1,1,\boldsymbol{\Omega}^{
o}) \widetilde{\alpha}_{\infty}^{o}(-1,1,\pmb{0}) \theta^{o}_{\infty}(-1,1,
\pmb{0}) \! - \! \beta^{o}_{\infty}(-1,1,\boldsymbol{\Omega}^{o})(\theta^{o}_{
\infty}(-1,1,\pmb{0}))^{2} \right) \! \tfrac{1}{(\theta^{o}_{\infty}(-1,1,\pmb{
0}))^{3}} \right. \\
&\left. \times \left(\dfrac{\gamma^{o}_{0} \! - \! (\gamma^{o}_{0})^{-1}}{2
\mi} \right) \! - \! \left(\dfrac{\widetilde{\alpha}_{\infty}^{o}(-1,1,\pmb{0}
) \theta^{o}_{\infty}(-1,1,\boldsymbol{\Omega}^{o}) \! - \! \widetilde{
\alpha}_{\infty}^{o}(-1,1,\boldsymbol{\Omega}^{o}) \theta^{o}_{\infty}(-1,1,
\pmb{0})}{(\theta^{o}_{\infty}(-1,1,\pmb{0}))^{2}} \right) \right. \\
&\left. \times \left(\sum_{k=1}^{N+1} \! \left(a_{k}^{o} \! - \! b_{k-1}^{o}
\right) \right) \! \left(\dfrac{\gamma^{o}_{0} \! + \! (\gamma^{o}_{0})^{-1}}{
8 \mi} \right) \! - \! \left(\! \left(\dfrac{\gamma^{o}_{0} \! + \! (\gamma^{
o}_{0})^{-1}}{16 \mi} \right) \! \sum_{k=1}^{N+1} \! \left((a_{k}^{o})^{2} \!
- \! (b_{k-1}^{o})^{2} \right) \right. \right. \\
&\left. \left. -\left(\dfrac{\gamma^{o}_{0} \! - \! (\gamma^{o}_{0})^{-1}}{64
\mi} \right) \! \left(\sum_{k=1}^{N+1} \! \left(a_{k}^{o} \! - \! b_{k-1}^{o}
\right) \right)^{2} \right) \! \dfrac{\theta^{o}_{\infty}(-1,1,\boldsymbol{
\Omega}^{o})}{\theta^{o}_{\infty}(-1,1,\pmb{0})} \right), \\
(\overset{o}{m}_{2}^{\raise-1.0ex\hbox{$\scriptstyle \infty$}})_{21} &=
-\dfrac{\boldsymbol{\theta}^{o}(\boldsymbol{u}^{o}_{+}(0) \! + \! \boldsymbol{
d}_{o}) \mathbb{E}}{\boldsymbol{\theta}^{o}(-\boldsymbol{u}^{o}_{+}(0) \! - \!
\frac{1}{2 \pi}(n \! + \! \frac{1}{2}) \boldsymbol{\Omega}^{o} \! - \!
\boldsymbol{d}_{o})} \! \left(\! \left(\theta^{o}_{\infty}(1,-1,\boldsymbol{
\Omega}^{o}) \! \left((\widetilde{\alpha}_{\infty}^{o}(1,-1,\pmb{0}))^{2} \! +
\! \beta^{o}_{\infty}(1,-1,\pmb{0}) \theta^{o}_{\infty}(1,-1,\pmb{0}) \right)
\right. \right. \\
&\left. \left. - \, \widetilde{\alpha}_{\infty}^{o}(1,-1,\boldsymbol{\Omega}^{
o}) \widetilde{\alpha}_{\infty}^{o}(1,-1,\pmb{0}) \theta^{o}_{\infty}(1,-1,
\pmb{0}) \! - \! \beta^{o}_{\infty}(1,-1,\boldsymbol{\Omega}^{o})(\theta^{o}_{
\infty}(1,-1,\pmb{0}))^{2} \right) \! \tfrac{1}{(\theta^{o}_{\infty}(1,-1,\pmb{
0}))^{3}} \right. \\
&\left. \times \left(\dfrac{\gamma^{o}_{0} \! - \! (\gamma^{o}_{0})^{-1}}{2
\mi} \right) \! - \! \left(\dfrac{\widetilde{\alpha}_{\infty}^{o}(1,-1,\pmb{0}
) \theta^{o}_{\infty}(1,-1,\boldsymbol{\Omega}^{o}) \! - \! \widetilde{
\alpha}_{\infty}^{o}(1,-1,\boldsymbol{\Omega}^{o}) \theta^{o}_{\infty}(1,-1,
\pmb{0})}{(\theta^{o}_{\infty}(1,-1,\pmb{0}))^{2}} \right) \right. \\
&\left. \times \left(\sum_{k=1}^{N+1} \! \left(a_{k}^{o} \! - \! b_{k-1}^{o}
\right) \right) \! \left(\dfrac{\gamma^{o}_{0} \! + \! (\gamma^{o}_{0})^{-1}}{
8 \mi} \right) \! - \! \left(\! \left(\dfrac{\gamma^{o}_{0} \! + \! (\gamma^{
o}_{0})^{-1}}{16 \mi} \right) \! \sum_{k=1}^{N+1} \! \left((a_{k}^{o})^{2} \!
- \! (b_{k-1}^{o})^{2} \right) \right. \right. \\
&\left. \left. -\left(\dfrac{\gamma^{o}_{0} \! - \! (\gamma^{o}_{0})^{-1}}{64
\mi} \right) \! \left(\sum_{k=1}^{N+1} \! \left(a_{k}^{o} \! - \! b_{k-1}^{o}
\right) \right)^{2} \right) \! \dfrac{\theta^{o}_{\infty}(1,-1,\boldsymbol{
\Omega}^{o})}{\theta^{o}_{\infty}(1,-1,\pmb{0})} \right), \\
(\overset{o}{m}_{2}^{\raise-1.0ex\hbox{$\scriptstyle \infty$}})_{22} &= \dfrac{
\boldsymbol{\theta}^{o}(\boldsymbol{u}^{o}_{+}(0) \! + \! \boldsymbol{d}_{o})
\mathbb{E}}{\boldsymbol{\theta}^{o}(-\boldsymbol{u}^{o}_{+}(0) \! - \! \frac{
1}{2 \pi}(n \! + \! \frac{1}{2}) \boldsymbol{\Omega}^{o} \! - \! \boldsymbol{
d}_{o})} \! \left(\! \left(\theta^{o}_{\infty}(-1,-1,\boldsymbol{\Omega}^{o})
\! \left((\widetilde{\alpha}_{\infty}^{o}(-1,-1,\pmb{0}))^{2} \! + \! \beta^{
o}_{\infty}(-1,-1,\pmb{0}) \theta^{o}_{\infty}(-1,-1,\pmb{0}) \right) \right.
\right. \\
&\left. \left. - \, \widetilde{\alpha}_{\infty}^{o}(-1,-1,\boldsymbol{\Omega}^{
o}) \widetilde{\alpha}_{\infty}^{o}(-1,-1,\pmb{0}) \theta^{o}_{\infty}(-1,-1,
\pmb{0}) \! - \! \beta^{o}_{\infty}(-1,-1,\boldsymbol{\Omega}^{o})(\theta^{o}_{
\infty}(-1,-1,\pmb{0}))^{2} \right) \! \tfrac{1}{(\theta^{o}_{\infty}(-1,-1,
\pmb{0}))^{3}} \right. \\
&\left. \times \left(\dfrac{\gamma^{o}_{0} \! + \! (\gamma^{o}_{0})^{-1}}{2}
\right) \! + \! \left(\dfrac{\widetilde{\alpha}_{\infty}^{o}(-1,-1,\pmb{0})
\theta^{o}_{\infty}(-1,-1,\boldsymbol{\Omega}^{o}) \! - \! \widetilde{\alpha}_{
\infty}^{o}(-1,-1,\boldsymbol{\Omega}^{o}) \theta^{o}_{\infty}(-1,-1,\pmb{0})}{
(\theta^{o}_{\infty}(-1,-1,\pmb{0}))^{2}} \right) \right. \\
&\left. \times \left(\sum_{k=1}^{N+1} \! \left(b_{k-1}^{o} \! - \! a_{k}^{o}
\right) \right) \! \left(\dfrac{\gamma^{o}_{0} \! - \! (\gamma^{o}_{0})^{-1}}{
8} \right) \! + \! \left(\! \left(\dfrac{\gamma^{o}_{0} \! - \! (\gamma^{o}_{
0})^{-1}}{16} \right) \! \sum_{k=1}^{N+1} \! \left((b_{k-1}^{o})^{2} \! - \!
(a_{k}^{o})^{2} \right) \right. \right. \\
&\left. \left. +\left(\dfrac{\gamma^{o}_{0} \! + \! (\gamma^{o}_{0})^{-1}}{64}
\right) \! \left(\sum_{k=1}^{N+1} \! \left(a_{k}^{o} \! - \! b_{k-1}^{o}
\right) \right)^{2} \right) \! \dfrac{\theta^{o}_{\infty}(-1,-1,\boldsymbol{
\Omega}^{o})}{\theta^{o}_{\infty}(-1,-1,\pmb{0})} \right),
\end{align*}
where, for $\varepsilon_{1},\varepsilon_{2} \! = \! \pm 1$,
\begin{gather*}
\theta^{o}_{\infty}(\varepsilon_{1},\varepsilon_{2},\pmb{\mathscr{Z}}) \! :=
\! \boldsymbol{\theta}^{o}(\varepsilon_{1} \boldsymbol{u}^{o}_{+}(\infty) \!
- \! \tfrac{1}{2 \pi}(n \! + \! \tfrac{1}{2}) \pmb{\mathscr{Z}} \! + \!
\varepsilon_{2} \boldsymbol{d}_{o}), \\
\widetilde{\alpha}^{o}_{\infty}(\varepsilon_{1},\varepsilon_{2},\pmb{\mathscr{
Z}}) \! := \! 2 \pi \mi \varepsilon_{1} \sum_{m \in \mathbb{Z}^{N}}(m,\widehat{
\boldsymbol{\alpha}}^{o}_{\infty}) \me^{2 \pi \mi (m,\varepsilon_{1}
\boldsymbol{u}^{o}_{+}(\infty)-\frac{1}{2 \pi}(n+\frac{1}{2}) \pmb{\mathscr{Z}}
+\varepsilon_{2} \boldsymbol{d}_{o})+ \pi \mi (m,\tau^{o}m)},
\end{gather*}
with $\boldsymbol{u}^{o}_{+}(\infty) \! = \! \int_{a_{N+1}^{o}}^{\infty^{+}} 
\boldsymbol{\omega}^{o}$, $\widehat{\boldsymbol{\alpha}}^{o}_{\infty} \! = \! 
(\widehat{\alpha}_{\infty,1}^{o},\widehat{\alpha}^{o}_{\infty,2},\dotsc,
\widehat{\alpha}^{o}_{\infty,N})$, and $\widehat{\alpha}^{o}_{\infty,j} \! := 
\! c_{j1}^{o}$, $j \! = \! 1,\dotsc,N$,
\begin{gather*}
\beta^{o}_{\infty}(\varepsilon_{1},\varepsilon_{2},\pmb{\mathscr{Z}}) \! := 
\! 2 \pi \sum_{m \in \mathbb{Z}^{N}} \! \left(\mi \varepsilon_{1}(m,\widehat{
\boldsymbol{\beta}}^{o}_{\infty}) \! + \!  \pi (m,\widehat{\boldsymbol{\alpha}
}^{o}_{\infty})^{2} \right) \! \me^{2 \pi \mi (m,\varepsilon_{1} \boldsymbol{
u}^{o}_{+}(\infty)-\frac{1}{2 \pi}(n+\frac{1}{2}) \pmb{\mathscr{Z}}+
\varepsilon_{2} \boldsymbol{d}_{o})+ \pi \mi (m,\tau^{o}m)},
\end{gather*}
where $\widehat{\boldsymbol{\beta}}^{o}_{\infty} \! = \! (\widehat{\beta}^{
o}_{\infty,1},\widehat{\beta}^{o}_{\infty,2},\dotsc,\widehat{\beta}^{o}_{
\infty,N})$, with $\widehat{\beta}^{o}_{\infty,j} \! := \! \tfrac{1}{2}
(c^{o}_{j2} \! + \! \tfrac{1}{2}c_{j1}^{o} \sum_{i=1}^{N+1}(b_{i-1}^{o} \! + 
\! a_{i}^{o}))$, $j \! = \! 1,\dotsc,N$, where $c_{j1}^{o},c_{j2}^{o}$, $j \! 
= \! 1,\dotsc,N$, are obtained {}from Equations~{\rm (O1)} and~{\rm (O2)}, and 
$\gamma^{o}_{0}$ is defined in Theorem~{\rm 3.1}, Equations~{\rm (3.53)}, and
\begin{align*}
\mathscr{R}^{o,\infty}_{0}(n) \underset{n \to \infty}{=}& \, \dfrac{1}{(n \! +
\! \frac{1}{2})} \mathlarger{\sum_{j=1}^{N+1}} \! \left(\dfrac{\left(\mathscr{
B}^{o}(b_{j-1}^{o}) \widehat{\alpha}_{0}^{o}(b_{j-1}^{o}) \! - \! \mathscr{
A}^{o}(b_{j-1}^{o}) \! \left(\widehat{\alpha}_{1}^{o}(b_{j-1}^{o}) \! + \!
(b_{j-1}^{o})^{-1} \widehat{\alpha}_{0}^{o}(b_{j-1}^{o}) \right) \right)}{b_{
j-1}^{o}(\widehat{\alpha}_{0}^{o}(b_{j-1}^{o}))^{2}} \right. \nonumber \\
+&\left. \, \dfrac{\left(\mathscr{B}^{o}(a_{j}^{o}) \widehat{\alpha}_{0}^{o}
(a_{j}^{o}) \! - \! \mathscr{A}^{o}(a_{j}^{o}) \! \left(\widehat{\alpha}_{1}^{
o}(a_{j}^{o}) \! + \! (a_{j}^{o})^{-1} \widehat{\alpha}_{0}^{o}(a_{j}^{o})
\right) \right)}{a_{j}^{o}(\widehat{\alpha}_{0}^{o}(a_{j}^{o}))^{2}} \right)
\! + \! \mathcal{O} \! \left(\dfrac{c(n)}{(n \! + \! \frac{1}{2})^{2}}
\right), \\
\mathscr{R}^{o,\infty}_{1}(n) \underset{n \to \infty}{=}& \, \dfrac{1}{(n \! +
\! \frac{1}{2})} \mathlarger{\sum_{j=1}^{N+1}} \! \left(\dfrac{\left(\mathscr{
B}^{o}(b_{j-1}^{o}) \widehat{\alpha}_{0}^{o}(b_{j-1}^{o}) \! - \! \mathscr{
A}^{o}(b_{j-1}^{o}) \widehat{\alpha}_{1}^{o}(b_{j-1}^{o}) \right)}{(\widehat{
\alpha}_{0}^{o}(b_{j-1}^{o}))^{2}} \right. \nonumber \\
+&\left. \, \dfrac{\left(\mathscr{B}^{o}(a_{j}^{o}) \widehat{\alpha}_{0}^{o}
(a_{j}^{o}) \! - \! \mathscr{A}^{o}(a_{j}^{o}) \widehat{\alpha}_{1}^{o}(a_{
j}^{o}) \right)}{(\widehat{\alpha}_{0}^{o}(a_{j}^{o}))^{2}} \right) \! + \!
\mathcal{O} \! \left(\dfrac{c(n)}{(n \! + \! \frac{1}{2})^{2}} \right), \\
\mathscr{R}^{o,\infty}_{2}(n) \underset{n \to \infty}{=}& \, \dfrac{1}{(n \! +
\! \frac{1}{2})} \mathlarger{\sum_{j=1}^{N+1}} \! \left(\dfrac{\left(\mathscr{
B}^{o}(b_{j-1}^{o}) \widehat{\alpha}_{0}^{o}(b_{j-1}^{o})b_{j-1}^{o} \! - \!
\mathscr{A}^{o}(b_{j-1}^{o}) \! \left(b_{j-1}^{o} \widehat{\alpha}_{1}^{o}
(b_{j-1}^{o}) \! - \! \widehat{\alpha}_{0}^{o}(b_{j-1}^{o}) \right) \right)}{
(\widehat{\alpha}_{0}^{o}(b_{j-1}^{o}))^{2}} \right. \nonumber \\
+&\left. \, \dfrac{\left(\mathscr{B}^{o}(a_{j}^{o}) \widehat{\alpha}_{0}^{o}
(a_{j}^{o})a_{j}^{o} \! - \! \mathscr{A}^{o}(a_{j}^{o}) \! \left(a_{j}^{o}
\widehat{\alpha}_{1}^{o}(a_{j}^{o}) \! - \! \widehat{\alpha}_{0}^{o}(a_{j}^{
o}) \right) \right)}{(\widehat{\alpha}_{0}^{o}(a_{j}^{o}))^{2}} \right) \! +
\! \mathcal{O} \! \left(\dfrac{c(n)}{(n \! + \! \frac{1}{2})^{2}} \right),
\end{align*}
with all parameters defined in Theorem~{\rm B.1}.
\end{dy}
\addtocounter{z0}{1}
\begin{dy}[{\rm \cite{a22}}]
Let all the conditions stated in Theorem~{\rm B.1} be valid, and let all 
parameters be as defined therein. Let $\overset{o}{\operatorname{Y}} \colon 
\mathbb{C} \setminus \mathbb{R} \! \to \! \operatorname{SL}_{2}(\mathbb{C})$ 
be the unique solution of \pmb{{\rm RHP2}} with integral 
representation~{\rm (1.30)}. Let $\boldsymbol{\pi}_{2n+1}(z)$ be the odd 
degree monic orthogonal $L$-polynomial defined in Equations~{\rm (1.5)} 
with $n \! \to \! \infty$ asymptotics (in the entire complex plane) given 
by Theorem~{\rm 2.3.1} of {\rm \cite{a22}}. Then,
\begin{equation}
\left(\xi^{(2n+1)}_{-n-1} \right)^{2} \! = \! \dfrac{1}{\norm{\boldsymbol{
\pi}_{2n+1}(\pmb{\cdot})}_{\mathscr{L}}^{2}} \! = \! \dfrac{H^{(-2n)}_{2n+1}
}{H^{(-2n-2)}_{2n+1}} \underset{n \to \infty}{=} \dfrac{\me^{-n \ell_{o}}}{
\pi} \Xi^{\natural} \! \left(1 \! + \! \frac{2}{n \! + \! \frac{1}{2}} \Xi^{
\natural}(\mathfrak{Q}^{\natural})_{12} \! + \! \mathcal{O} \! \left(\dfrac{
c^{\natural}(n)}{(n \! + \! \frac{1}{2})^{2}} \right) \right),
\end{equation}
where
\begin{align*}
\Xi^{\natural} :=& \, 2 \mathbb{E}^{2} \! \left(\sum_{k=1}^{N+1} \! \left(
(b_{k-1}^{o})^{-1} \! - \! (a_{k}^{o})^{-1} \right) \right)^{-1} \dfrac{
\boldsymbol{\theta}^{o}(\boldsymbol{u}^{o}_{+}(0) \! - \! \frac{1}{2 \pi}
(n \! + \! \frac{1}{2}) \boldsymbol{\Omega}^{o} \! + \! \boldsymbol{d}_{o}) 
\boldsymbol{\theta}^{o}(-\boldsymbol{u}^{o}_{+}(0) \! + \! \boldsymbol{d}_{o})
}{\boldsymbol{\theta}^{o}(-\boldsymbol{u}^{o}_{+}(0) \! - \! \frac{1}{2 \pi}
(n \! + \! \frac{1}{2}) \boldsymbol{\Omega}^{o} \! + \! \boldsymbol{d}_{o}) 
\boldsymbol{\theta}^{o}(\boldsymbol{u}^{o}_{+}(0) \! + \! \boldsymbol{d}_{o}
)} \quad (> \! 0),
\end{align*}
$\mathfrak{Q}^{\natural}$ is defined in Theorem~{\rm 3.2}, 
Equations~{\rm (3.107)}, $(\mathfrak{Q}^{\natural})_{12}$ denotes the $(1 
\, 2)$-element of $\mathfrak{Q}^{\natural}$, $c^{\natural}(n) \! =_{n \to 
\infty} \! \mathcal{O}(1)$, and all relevant parameters are defined in 
Theorem~{\rm B.1}. Asymptotics for $\xi^{(2n+1)}_{-n-1}$ are obtained by 
taking the positive square root of both sides of Equation~{\rm (B.1)}.
\end{dy}

\end{document}